\documentclass[11pt,reqno]{amsart}

\makeatletter

\usepackage[margin=1.2in]{geometry}
\usepackage{amssymb}
\usepackage{mathtools}
\usepackage{graphicx}
\usepackage{xcolor}
\usepackage[skins,breakable]{tcolorbox}
\usepackage{tikz}
\usepackage{pgfornament}       
\usetikzlibrary{calc,arrows.meta,intersections}
\usepackage{enumitem}
\usepackage{booktabs}
\usepackage[mathcal]{euscript}
\usepackage{hyperref}

\makeatletter

\usepackage{amssymb}
\usepackage{mathrsfs}

\newcommand{\CC}{\mathbb{C}}

\newcommand{\NN}{\mathbb{N}}

\newcommand{\QQ}{\mathbb{Q}}
\newcommand{\RR}{\mathbb{R}}

\newcommand{\ZZ}{\mathbb{Z}}

\newcommand{\cA}{\mathcal{A}}
\newcommand{\cB}{\mathcal{B}}
\newcommand{\cC}{\mathcal{C}}

\newcommand{\cF}{\mathcal{F}}
\newcommand{\cG}{\mathcal{G}}
\newcommand{\cH}{\mathcal{H}}

\newcommand{\cL}{\mathcal{L}}

\newcommand{\cP}{\mathcal{P}}

\newcommand{\cW}{\mathcal{W}}

\newcommand{\fM}{\mathfrak{M}}
\newcommand{\fN}{\mathfrak{N}}

\makeatletter

\usepackage{amsmath}
\usepackage{amssymb}
\usepackage{mathtools}  

\DeclareMathOperator{\Tr}{Tr}      

\providecommand{\cH}{\mathcal{H}}

\makeatletter

\usepackage[nosort]{cleveref}   

\newcommand{\di@crefnames}[3]{\crefname{#1}{#2}{#3}\Crefname{#1}{#2}{#3}}
\di@crefnames{section}{Section}{Sections}
\di@crefnames{subsection}{Section}{Sections}
\di@crefnames{subsubsection}{Section}{Sections}
\di@crefnames{chapter}{Chapter}{Chapters}
\di@crefnames{part}{Part}{Parts}
\di@crefnames{appendix}{Appendix}{Appendices}
\di@crefnames{table}{Table}{Tables}
\di@crefnames{figure}{Figure}{Figures}
\di@crefnames{theorem}{Theorem}{Theorems}
\di@crefnames{keytheorem}{Theorem}{Theorems}
\di@crefnames{proposition}{Proposition}{Propositions}
\di@crefnames{corollary}{Corollary}{Corollaries}
\di@crefnames{lemma}{Lemma}{Lemmas}
\di@crefnames{fact}{Fact}{Facts}
\di@crefnames{conjecture}{Conjecture}{Conjectures}
\di@crefnames{definition}{Definition}{Definitions}
\di@crefnames{example}{Example}{Examples}
\di@crefnames{remark}{Remark}{Remarks}
\di@crefnames{question}{Question}{Questions}
\di@crefnames{axiom}{Axiom}{Axioms}
\di@crefnames{notation}{Notation}{Notations}
\di@crefnames{exercise}{Exercise}{Exercises}

\AtBeginDocument{%
}
\makeatletter

\usepackage{amsthm}
\usepackage{aliascnt}   

\def\di@rudin@parent{section}

\newtheoremstyle{rudin}%
  {\topsep}{\topsep}{\itshape}{}{\bfseries}{.}{ }%
  {\thmname{#1}\thmnumber{ #2}\thmnote{ (#3)}}
\newtheoremstyle{rudindef}%
  {\topsep}{\topsep}{\upshape}{}{\bfseries}{.}{ }%
  {\thmname{#1}\thmnumber{ #2}\thmnote{ (#3)}}

\theoremstyle{rudin}
\ifx\di@rudin@parent\@empty
  \newtheorem{theorem}{Theorem}
\else
  \newtheorem{theorem}{Theorem}[\di@rudin@parent]
\fi
\newcommand{\DiSharedTheorem}[2]{%
  \newaliascnt{#1}{theorem}%
  \newtheorem{#1}[#1]{#2}%
  \aliascntresetthe{#1}%
}
\DiSharedTheorem{proposition}{Proposition}
\DiSharedTheorem{corollary}{Corollary}
\DiSharedTheorem{lemma}{Lemma}
\DiSharedTheorem{fact}{Fact}
\theoremstyle{rudindef}
\DiSharedTheorem{definition}{Definition}
\DiSharedTheorem{example}{Example}
\DiSharedTheorem{remark}{Remark}
\DiSharedTheorem{question}{Question}
\makeatletter

\newif\ifdi@hi@timeline

\usepackage{amsthm}
\usepackage{aliascnt}
\usepackage{graphicx}
\usepackage{xcolor}
\usepackage{manfnt}        
\usepackage{tikzsymbols}   

\definecolor{plum}{RGB}{94,42,110}        
\definecolor{orchid}{RGB}{142,68,173}     
\definecolor{rose}{RGB}{196,96,168}       
\definecolor{lavender}{RGB}{200,180,230}  
\definecolor{blush}{RGB}{248,215,228}     
\colorlet{keyframe}{plum!75}              
\colorlet{keyfill}{lavender!22!white}
\definecolor{iconcolor}{RGB}{218,105,160}  
\definecolor{qedcolor}{RGB}{186,140,210}   

\newcommand{\bendicon}{\llap{\scalebox{0.6}{\textcolor{iconcolor}{\dbend}}\hspace{0.9em}}}
\newtheoremstyle{rudinbend}%
  {\topsep}{\topsep}%
  {\upshape}{}%
  {\bfseries}{.}%
  { }%
  {\rlap{\normalfont\bendicon}\thmname{#1}\thmnumber{ #2}\thmnote{ (#3)}}
\@ifundefined{remark}{%
  \theoremstyle{rudinbend}%
  \newtheorem{remark}[theorem]{Remark}%
  \newtheorem{remarks}[theorem]{Remarks}%
}{%
  \let\di@hi@theremark\theremark
  \let\remark\@undefined \let\endremark\@undefined
  \theoremstyle{rudinbend}%
  \newtheorem{remark}[remark]{Remark}%
  \let\theremark\di@hi@theremark
  \newaliascnt{remarks}{theorem}%
  \aliascntresetthe{remarks}%
}

\@ifpackageloaded{cleveref}{%
  \crefname{remarks}{Remarks}{Remarks}\Crefname{remarks}{Remarks}{Remarks}%
}{}

\DeclareRobustCommand{\marginicon}[1]{%
  \llap{\raisebox{-0.2ex}{\textcolor{iconcolor}{\resizebox{!}{1.75ex}{#1}}}\hspace{0.9em}}}
\newcommand{\historypara}[1]{\paragraph{\marginicon{\faHourglassHalf}#1}}

\newcommand{\@twelvehour}[2]{%
  \ifcase\numexpr#1\relax 12\or 1\or 2\or 3\or 4\or 5\or 6\or 7\or 8\or 9\or 10\or 11%
  \or 12\or 1\or 2\or 3\or 4\or 5\or 6\or 7\or 8\or 9\or 10\or 11\fi
  :\two@digits{#2}\ \ifnum\numexpr#1\relax<12 am\else pm\fi}
\newcommand{\timestamp}{%
  \edef\@tshour{\the\numexpr(\time+30)/60-1\relax}%
  \edef\@tsmin{\the\numexpr\time-60*\@tshour\relax}%
  Last update: \today, \@twelvehour{\@tshour}{\@tsmin}~CT\,/\,%
  \@twelvehour{\ifnum\@tshour=23 0\else\numexpr\@tshour+1\relax\fi}{\@tsmin}~ET}

\ifdi@hi@timeline
\usepackage{tikz}
\usepackage{tcolorbox}
\tcbuselibrary{skins,breakable}
\colorlet{tmlaccent}{rose}
\colorlet{tmlfill}{blush!60!white}
\colorlet{tmlspine}{lavender}
\newlength\tmlspine 
\newlength\tmltext  
\newlength\tmlsep   
\newlength\tmldot                                
\makeatletter
\newcommand\tml@at[1]{([xshift=\tmlspine,yshift=#1]frame.north west)}
\newcommand\tml@dot{\filldraw[fill=tmlaccent,draw=white,line width=0.8pt]
  \tml@at{-\tmldot} circle (2.9pt);}
\newcommand\tml@bar{%
  \draw[white,line width=7.4pt,line cap=round]\tml@at{-\tmldot}--\tml@at{-\tmldot-0.9em};
  \draw[tmlaccent,line width=5.8pt,line cap=round]\tml@at{-\tmldot}--\tml@at{-\tmldot-0.9em};}
\let\tml@marker\tml@dot
\newcommand\tml@down{\draw[tmlspine,line width=0.5pt]\tml@at{-\tmldot}--([xshift=\tmlspine]frame.south west);}
\newcommand\tml@up{\draw[tmlspine,line width=0.5pt]\tml@at{0pt}--\tml@at{-\tmldot};}
\newcommand\tml@full{\draw[tmlspine,line width=0.5pt]\tml@at{0pt}--([xshift=\tmlspine]frame.south west);}
\tcbset{
  tml/base/.style={enhanced, breakable, blanker, left=\tmltext, right=0pt,
    top=0pt, bottom=\tmlsep, before skip=0pt, after skip=0pt, parbox=false},
  tml/only/.style={overlay unbroken and first={\tml@marker}, bottom=0pt},
  tml/first/.style={overlay unbroken and first={\tml@down\tml@marker},
    overlay middle and last={\tml@full}},
  tml/middle/.style={overlay unbroken and first={\tml@full\tml@marker},
    overlay middle and last={\tml@full}},
  tml/last/.style={overlay unbroken and first={\tml@up\tml@marker}, bottom=0pt},
}
\tikzset{tmlbadge/.style={fill=tmlfill, draw=tmlaccent, line width=0.4pt}}
\newcommand\tml@badge[1]{\tikz[baseline=(b.base)]\node(b)[tmlbadge,
  rounded corners=2pt, inner xsep=3pt, inner ysep=1.6pt, font=\small\bfseries]{#1};}
\newcommand\tml@entry[5]{%
  \expandafter\let\expandafter\tml@marker\csname tml@#2\endcsname
  \setlength\tmldot{\dimexpr0.72\baselineskip-0.5ex\relax}%
  \begin{tcolorbox}[tml/base, tml/#1]
  \noindent\llap{\tml@badge{#3}\hspace{\dimexpr\tmltext-\tmlspine+1em\relax}}%
  \if\relax\detokenize{#4}\relax\else\textbf{#4\@addpunct{.}}\ \fi#5
  \end{tcolorbox}}

\makeatother
\ExplSyntaxOn
\seq_new:N \g__tml_seq
\cs_new:Npn \__tml_first:nnnnn #1#2#3#4#5 {#1}
\cs_new:Npn \__tml_second:nnnnn #1#2#3#4#5 {#2}
\tl_new:N \l__tml_pos_tl
\tl_new:N \l__tml_start_tl
\tl_new:N \l__tml_end_tl
\seq_new:N \l__tml_match_seq
\bool_new:N \l__tml_strip_bool
\bool_new:N \g__tml_overview_bool
\bool_gset_true:N \g__tml_overview_bool
\NewDocumentCommand \timelineoverview { m }
  {
    \str_if_eq:nnTF {#1} { false }
      { \bool_gset_false:N \g__tml_overview_bool }
      { \bool_gset_true:N \g__tml_overview_bool }
  }
\NewDocumentCommand \event { o m o +m }
  {
    \IfNoValueTF {#1}
      { \tl_set:Nx \l__tml_start_tl { \str_range:nnn {#2} {1} {4} } }
      { \tl_set:Nn \l__tml_start_tl {#1} }
    \regex_extract_once:nnNTF { --([0-9]{4}) } {#2} \l__tml_match_seq
      {
        \tl_set:Nx \l__tml_end_tl
          { \int_max:nn { \l__tml_start_tl } { \seq_item:Nn \l__tml_match_seq {2} } }
      }
      { \tl_set_eq:NN \l__tml_end_tl \l__tml_start_tl }
    \seq_gput_right:Nx \g__tml_seq
      {
        { \l__tml_start_tl } { \l__tml_end_tl }
        \exp_not:n { {#2} }
        { \IfNoValueTF {#3} { } { \exp_not:n {#3} } }
        \exp_not:n { {#4} }
      }
  }
\NewDocumentEnvironment { timeline } { O{} }
  {
    \par \seq_gclear:N \g__tml_seq
    \bool_set_eq:NN \l__tml_strip_bool \g__tml_overview_bool
    \str_if_eq:nnT {#1} { nooverview } { \bool_set_false:N \l__tml_strip_bool }
  }
  {
    \seq_gsort:Nn \g__tml_seq
      {
        \int_compare:nNnTF { \__tml_first:nnnnn ##1 } > { \__tml_first:nnnnn ##2 }
          { \sort_return_swapped: } { \sort_return_same: }
      }
    \par\medskip
    \bool_if:NT \l__tml_strip_bool { \__tml_overview: }
    \seq_map_indexed_inline:Nn \g__tml_seq
      {
        \tl_set:Nx \l__tml_pos_tl
          {
            \int_compare:nNnTF { \seq_count:N \g__tml_seq } = { 1 } { only }
              {
                \int_compare:nNnTF {##1} = { 1 } { first }
                  {
                    \int_compare:nNnTF {##1} = { \seq_count:N \g__tml_seq }
                      { last } { middle }
                  }
              }
          }
        \__tml_print:Vnnnnn \l__tml_pos_tl ##2
      }
    \par\medskip
  }
\cs_new_protected:Npn \__tml_print:nnnnnn #1#2#3#4#5#6
  {
    \use:c { tml@entry } {#1}
      { \int_compare:nNnTF {#3} > {#2} { bar } { dot } } {#4} {#5} {#6}
  }
\cs_generate_variant:Nn \__tml_print:nnnnnn { V }
\int_new:N \l__tml_ya_int
\int_new:N \l__tml_yb_int
\int_new:N \l__tml_lane_int
\int_new:N \l__tml_nlanes_int
\dim_new:N \l__tml_unit_dim
\dim_new:N \l__tml_xs_dim
\dim_new:N \l__tml_xe_dim
\dim_new:N \l__tml_ext_dim
\dim_new:N \l__tml_xl_dim
\dim_new:N \l__tml_axis_dim
\bool_new:N \l__tml_flip_bool
\prop_new:N \l__tml_lanes_prop
\tl_new:N \l__tml_draw_tl
\box_new:N \l__tml_box
\dim_const:Nn \c__tml_lanesep_dim { 13pt }
\cs_new:Npn \__tml_pt:n #1 { \dim_to_decimal:n {#1} pt }
\cs_generate_variant:Nn \prop_put:Nnn { Nxx }
\cs_new_protected:Npn \__tml_place:nnnnn #1#2#3#4#5
  {
    \dim_set:Nn \l__tml_xs_dim { \l__tml_unit_dim * ( #1 - \l__tml_ya_int ) }
    \dim_set:Nn \l__tml_xe_dim { \l__tml_unit_dim * ( #2 - \l__tml_ya_int ) }
    \hbox_set:Nn \l__tml_box { \scriptsize \tl_if_blank:nTF {#4} {#3} {#4} }
    \bool_set_false:N \l__tml_flip_bool
    \dim_compare:nNnT
      { \l__tml_xs_dim - 2pt + \box_wd:N \l__tml_box } > { \l__tml_axis_dim + 1em }
      { \bool_set_true:N \l__tml_flip_bool }
    \bool_if:NTF \l__tml_flip_bool
      {
        \dim_set:Nn \l__tml_xl_dim { \l__tml_xe_dim + 2pt - \box_wd:N \l__tml_box }
        \dim_set:Nn \l__tml_ext_dim { \l__tml_xe_dim + 8pt }
      }
      {
        \dim_set:Nn \l__tml_xl_dim { \l__tml_xs_dim }
        \dim_set:Nn \l__tml_ext_dim
          { \dim_max:nn { \l__tml_xe_dim } { \l__tml_xs_dim + \box_wd:N \l__tml_box } + 6pt }
      }
    \int_zero:N \l__tml_lane_int
    \int_step_inline:nn { \l__tml_nlanes_int }
      {
        \int_compare:nNnT { \l__tml_lane_int } = { 0 }
          {
            \dim_compare:nNnT { \prop_item:Nn \l__tml_lanes_prop {##1} }
              < { \dim_min:nn { \l__tml_xs_dim } { \l__tml_xl_dim } - 2pt }
              { \int_set:Nn \l__tml_lane_int {##1} }
          }
      }
    \int_compare:nNnT { \l__tml_lane_int } = { 0 }
      {
        \int_incr:N \l__tml_nlanes_int
        \int_set_eq:NN \l__tml_lane_int \l__tml_nlanes_int
      }
    \prop_put:Nxx \l__tml_lanes_prop
      { \int_use:N \l__tml_lane_int } { \dim_use:N \l__tml_ext_dim }
    \tl_put_right:Nx \l__tml_draw_tl
      {
        \int_compare:nNnTF {#2} > {#1}
          {
            \exp_not:N \draw [ fill=tmlfill, draw=tmlaccent, line~width=0.5pt,
                               rounded~corners=1.8pt ]
              ( \__tml_pt:n { \l__tml_xs_dim } ,
                \__tml_pt:n { \l__tml_lane_int \c__tml_lanesep_dim - 2pt } )
              rectangle
              ( \__tml_pt:n { \l__tml_xe_dim } ,
                \__tml_pt:n { \l__tml_lane_int \c__tml_lanesep_dim + 2pt } ) ;
          }
          {
            \exp_not:N \filldraw [ fill=tmlaccent, draw=white, line~width=0.6pt ]
              ( \__tml_pt:n { \l__tml_xs_dim } ,
                \__tml_pt:n { \l__tml_lane_int \c__tml_lanesep_dim } )
              circle (2.3pt) ;
          }
        \exp_not:N \node [ anchor=south~west, inner~sep=0pt,
                           font=\exp_not:N \scriptsize ]
          at ( \__tml_pt:n { \bool_if:NTF \l__tml_flip_bool { \l__tml_xl_dim } { \l__tml_xs_dim - 2pt } } ,
               \__tml_pt:n { \l__tml_lane_int \c__tml_lanesep_dim + 3.5pt } )
          { \tl_if_blank:nTF {#4} { \exp_not:n {#3} } { \exp_not:n {#4} } } ;
      }
  }
\cs_new_protected:Npn \__tml_overview:
  {
    \int_set:Nn \l__tml_ya_int { 99999 }
    \int_set:Nn \l__tml_yb_int { 0 }
    \seq_map_inline:Nn \g__tml_seq
      {
        \int_set:Nn \l__tml_ya_int { \int_min:nn { \l__tml_ya_int } { \__tml_first:nnnnn ##1 } }
        \int_set:Nn \l__tml_yb_int { \int_max:nn { \l__tml_yb_int } { \__tml_second:nnnnn ##1 } }
      }
    \int_set:Nn \l__tml_ya_int { \int_div_truncate:nn { \l__tml_ya_int } { 10 } * 10 }
    \int_set:Nn \l__tml_yb_int { \int_div_truncate:nn { \l__tml_yb_int + 9 } { 10 } * 10 }
    \int_compare:nNnT { \l__tml_yb_int } = { \l__tml_ya_int }
      { \int_add:Nn \l__tml_yb_int { 10 } }
    \dim_set:Nn \l__tml_unit_dim
      { ( \linewidth - 4em ) / ( \l__tml_yb_int - \l__tml_ya_int ) }
    \dim_set:Nn \l__tml_axis_dim { \l__tml_unit_dim * ( \l__tml_yb_int - \l__tml_ya_int ) }
    \prop_clear:N \l__tml_lanes_prop
    \int_zero:N \l__tml_nlanes_int
    \tl_clear:N \l__tml_draw_tl
    \seq_map_inline:Nn \g__tml_seq { \__tml_place:nnnnn ##1 }
    \par\noindent\hspace*{1em}
    \begin{tikzpicture}
      \draw [ tmlspine, line~width=0.5pt ] ( 0pt , 0pt ) --
        ( \__tml_pt:n { \l__tml_unit_dim * ( \l__tml_yb_int - \l__tml_ya_int ) } , 0pt ) ;
      \int_step_inline:nnnn { \l__tml_ya_int } { 10 } { \l__tml_yb_int }
        {
          \draw [ tmlspine, line~width=0.5pt ]
            ( \__tml_pt:n { \l__tml_unit_dim * ( ##1 - \l__tml_ya_int ) } , 0pt )
            -- ++ ( 0pt , -3pt )
            node [ below, inner~sep=2pt, font=\scriptsize, text=black!55 ] { ##1 } ;
        }
      \tl_use:N \l__tml_draw_tl
    \end{tikzpicture}
    \par\bigskip
  }
\ExplSyntaxOff

\fi
\makeatletter

\hypersetup{colorlinks=true, linkcolor=plum, citecolor=orchid, urlcolor=rose,
  pdftitle={Finite refutations in countably additive quantum probability},
  pdfauthor={O. Aubone and J. Iovino}}

\newtheoremstyle{rudinletter}%
  {\topsep}{\topsep}{\itshape}{}{\bfseries}{.}{ }%
  {\thmname{#1}\thmnumber{ #2}\thmnote{ (#3)}}
\theoremstyle{rudinletter}
\newtheorem{mainthm}{Theorem}

\crefname{mainthm}{Theorem}{Theorems}
\Crefname{mainthm}{Theorem}{Theorems}
\tcolorboxenvironment{mainthm}{%
  colback=keyfill, colframe=keyframe, boxrule=0.4pt, arc=6pt,
  before skip=12pt, after skip=12pt,
  left=6pt, right=6pt, top=4pt, bottom=2pt}

\theoremstyle{rudin}
\DiSharedTheorem{keytheorem}{Theorem}
\crefalias{keytheorem}{theorem}
\tcolorboxenvironment{keytheorem}{%
  colback=keyfill, colframe=keyframe, boxrule=0.4pt, arc=6pt,
  left=6pt, right=6pt, top=2pt, bottom=0pt,
  before skip=10pt, after skip=10pt}

\colorlet{axiomframe}{rose!80!black}
\colorlet{axiomfill}{blush!35!white}
\newtcolorbox{axiomgroup}[1]{colback=axiomfill, colframe=axiomframe, boxrule=0.4pt,
  arc=6pt, left=6pt, right=6pt, top=4pt, bottom=4pt, breakable,
  before skip=8pt, after skip=8pt,
  before upper={\textbf{#1.}\par\smallskip}}

\makeatletter
\def\paragraph{\@startsection{paragraph}{4}%
  \z@\z@{-\fontdimen2\font}%
  {\normalfont\bfseries}}
\renewenvironment{proof}[1][\proofname]{\par
  \pushQED{\qed}%
  \normalfont \topsep6\p@\@plus6\p@\relax
  \trivlist
  \item[\hskip\labelsep]%
    {\itshape #1\@addpunct{.}}\hskip\labelsep\ignorespaces
}{%
  \popQED\endtrivlist\@endpefalse
}
\makeatother

\newcommand{\HH}{\mathbb{H}}                 
\newcommand{\ip}[2]{\langle #1\mid #2\rangle}  
\newcommand{\DM}{\mathrm{DM}}                
\newcommand{\sH}{\mathbf{H}}                 
\newcommand{\sL}{\mathbf{L}}                 
\newcommand{\sK}{\mathbf{K}}                 
\newcommand{\at}{\mathrm{at}}
\newcommand{\rep}{\mathrm{rep}}
\newcommand{\spa}{\mathrm{span}}             
\newcommand{\nonatom}{\mathrm{nonatom}}       
\newcommand{\inv}{^{*}}
\newcommand{\QI}{\QQ(i)}
\newcommand{\doteqq}{\mathrel{\dot{=}}}      

\usepackage{pifont}
\usepackage{fontawesome5}
\setlist[itemize]{label=\textbullet}

\AtBeginDocument{}

\newcommand{\finitesubtree}[3]{\fill[lavender!55] (#1,#2) -- (#1-#3,#2-1.1) -- (#1+#3,#2-1.1) -- cycle;
  \draw[plum!60] (#1,#2) -- (#1-#3,#2-1.1) -- (#1+#3,#2-1.1) -- cycle;
  \foreach \t in {-1,0,1}{\node[text=red!75!black,font=\small\bfseries] at (#1+0.7*\t*#3,#2-1.25) {$\times$};}}

\begin{document}

\title[Finite refutations in quantum probability]{Finite refutations in countably additive quantum probability}

\author[Aubone]{Olivia Aubone}
\address{Department of Mathematics, UT San Antonio, One UTSA Circle, San Antonio, TX 78249, USA}
\email{olivia.aubone@utsa.edu}

\author[Iovino]{José Iovino}
\address{Department of Mathematics, UT San Antonio, One UTSA Circle, San Antonio, TX 78249, USA}
\email{jose.iovino@utsa.edu}
\thanks{This version: \today. The most recent version of this paper is available at \url{https://utsa-quantum-math.github.io/resources/gleason-kochen-specker-game/}; it may supersede the version deposited on arXiv.}


\date{}  

\begin{abstract}
\setlength{\emergencystretch}{2em}
We axiomatize countably additive quantum probability and prove a completeness theorem for finite refutation beyond compactness: inconsistency is witnessed by refutations each of whose branches ends in finitely many jointly unsatisfiable commitments.
Adjoining the hypothesis that every event has a predetermined, context-independent outcome makes the theory inconsistent, although every finite subset of the resulting axioms has a model.

For a countable theory of closed bounds on sentence values in our real-valued infinitary logic, inconsistency is equivalent to a winning strategy for a challenger in a game against a model builder.
The challenger chooses each test in response to the builder's commitments, and the plays that follow a winning strategy form such a refutation.
For the quantum theory with predetermined outcomes, every winning strategy must branch infinitely at some node.
In the strategy we construct, countable additivity locates the four-dimensional subspace in which a finite Kochen--Specker parity argument applies.

In their standard models, our axioms describe the closed subspaces and the density-operator probabilities of a separable infinite-dimensional complex Hilbert space.
We prove in ZFC that the sets of probability values definable uniformly across standard models are exactly the $\sigma$-projective subsets of the unit interval.
\end{abstract}
\maketitle

\begin{center}{\small\scshape Overview}\end{center}
\begin{quote}\small
\setlength{\emergencystretch}{2em}
A central question in quantum foundations is whether a measurement reveals a value that was already determined and that does not depend on which other compatible measurements are made with it.
The Bell--Kochen--Specker theorem rules out such assignments to quantum events.
A finite configuration proves this by exhibiting finitely many compatible groups of events whose assignments cannot agree on the overlaps.
In finite dimension, logical compactness connects this finite conclusion to the impossibility of an assignment on the whole space, and Hrushovski and Pitowsky used first-order compactness to obtain effective finite approximations to the predictions of quantum states.

We axiomatize the countably additive probability theory of a separable infinite-dimensional complex Hilbert space.
Its events are the closed subspaces, and its standard models give exactly the density-operator measures of Gleason's theorem.
Countable additivity concerns orthogonal families of subspaces and their closed spans.
To describe those spans we use a language with variables for vectors and sequences of vectors; even classical infinitary logic in the lattice vocabulary alone cannot express lattice completeness in the relevant class.

In this setting a finite list of axioms can fail to witness an inconsistency.
Adjoining the hypothesis that all probabilities are $0$ or $1$ gives a theory with no models, standard or not, although each finite subset has a model.
The models of the finite subsets rest on a simple valuation: on the finite-dimensional subspaces and their orthogonal complements, it gives the value $0$ to the former and $1$ to the latter.
On any finite list of these events, the valuation agrees with the Born probabilities of a single pure state, given by a unit vector orthogonal to the finite-dimensional subspaces on the list and to the complements of the others.
No density operator gives these values on all the events, since every ray of an orthonormal basis would have probability $0$, although the probabilities over the basis sum to $1$.

We prove that a different form of finite refutation is complete for real-valued infinitary logic.
For a countable theory of closed conditions on sentence values, inconsistency is equivalent to a winning strategy in a game.
The challenger chooses the next test in light of the builder's previous answers, and every play that follows a winning strategy ends with finitely many incompatible commitments.
The tree of possible plays is well-founded, and in the quantum example it must branch infinitely at some node.

In the strategy we construct, countable additivity forces the builder to name a finite initial segment of an orthonormal basis whose rays together have probability close to $1$, and his later commitments single out one ray in that segment.
The index of that ray locates the four-dimensional subspace for a finite Kochen--Specker parity argument.
The proof identifies the projections shared across measurement contexts and records how equality of projections forces equality of their probabilities.

The same axioms determine the expressive power of the language.
In a suitably chosen standard model, a single formula can define an arbitrary subset of $(0,1)$, so the substantive question is which sets of probability values a single formula recognizes in every standard model.
These are exactly the $\sigma$-projective sets, a class strictly larger than the projective sets, and the classification is proved in ZFC.
\end{quote}

\tableofcontents

\section*{Introduction}

\begingroup
\setlength{\emergencystretch}{2em}
We give an axiomatization of countably additive quantum probability and prove that every hidden-variable hypothesis that is expressed by countably many closed conditions and is incompatible with the axioms can be refuted by a strategy that ends in a finite contradiction on every play.
In a compact logic any inconsistent theory has a finite inconsistent fragment.
Our quantum example admits no such fragment: when the hypothesis of predetermined, context-independent outcomes is adjoined to our axioms, every finite subset of the resulting axioms has a model, yet the theory as a whole has none.
We call an incompatible hypothesis an \emph{obstruction}; the criterion concerns all models of the axioms, including those in which the sorts of vectors do not have their full Hilbert-space interpretations.

In place of a finite inconsistent subset of the axioms, our certificate is a well-founded strategy in a game, each of whose terminal positions carries finitely many commitments that no structure satisfies together.
A model builder commits to bounds on sentence values, and a challenger chooses the next test in light of these commitments.
For countable theories of closed conditions, the challenger has a winning strategy exactly when the theory is inconsistent.
For our axioms, as presented, together with the zero--one condition, every winning strategy must branch infinitely at some node.
In the strategy we construct, the builder's answer to the challenge about countable additivity forces a commitment to a basis ray of probability close to $1$, and that ray specifies the four-dimensional subspace in which a Kochen--Specker parity contradiction occurs.

The mathematical setting is that of Gleason's theorem.
For a Hilbert space $H$, the closed subspaces form a lattice $L(H)$, and each closed subspace is an event, identified with the orthogonal projection onto it.
A density operator on $H$ is a positive operator of trace $1$ and represents a quantum state.
When $H$ is separable and of dimension at least three, Gleason's theorem~\cite{Gleason:1957} identifies the countably additive probability measures on $L(H)$ with the maps $P\mapsto\Tr(\rho P)$ given by density operators $\rho$, none of which takes only the values $0$ and $1$.
An outcome-deterministic noncontextual hidden-variable model assigns to every event, at each hidden state, a value $0$ or $1$; in every complete set of mutually orthogonal events exactly one event receives the value $1$, and the value of an event does not depend on which compatible events are measured together with it~\cite{Budroni-Cabello-Guhne-Kleinmann-Larsson:2022}.
Specker announced the impossibility of such assignments in dimension three~\cite{Specker:1960}.
J.~S.~Bell observed that such an assignment is a $\{0,1\}$-valued measure on $L(H)$, which Gleason's theorem excludes~\cite{Bell:1966}.
Kochen and Specker then exhibited a finite configuration of rays that no such assignment respects~\cite{Kochen-Specker:1967}.

The earlier logical formulations of contextuality obtain their finite results from compactness.
J.~L.~Bell gave the argument of Kochen and Specker its logical form~\cite[Theorem~1]{Bell-JL:1996}.
Events become propositional variables, the Boolean operations within a context become connectives, and a noncontextual assignment becomes a satisfying truth assignment.
Propositional logic is compact and the constraints on an assignment are closed; Gleason's theorem therefore yields, in every finite dimension at least three, a finite inconsistent set of constraints, and the existence of configurations of Kochen and Specker becomes a consequence of compactness.
Propositional logic records only the values $0$ and $1$, and the probabilities of quantum theory lie outside its language.

Hrushovski and Pitowsky used classical two-valued first-order logic to derive quantitative constraints on quantum probabilities and finite approximations to the predictions of a quantum state~\cite{Hrushovski-Pitowsky:2004}.
Their language has separate sorts for directions and scalars, with probabilities taking values in the scalar sort.
The quantitative constraints include Pitowsky's logical indeterminacy principle~\cite{Pitowsky:1998}.
Their formulation rests on the first-order theories of the real and complex fields: compactness yields finite sets of rays on which the probabilities specifying a state constrain further probabilities to any prescribed accuracy, and the decidability of real closed and algebraically closed fields makes the construction effective.
The formulation remains within finite dimension, where orthonormal bases are finite.
Their argument uses compactness to extract effective finite constraints; our example asks what form of finite refutation survives when there is no finite inconsistent fragment to extract.

In the infinite-dimensional setting of Gleason's theorem, the event lattice is complete and the probability measure is countably additive, and our axiomatization expresses these requirements.

Fixed finite configurations reach the finitely additive part of the subject.
Alda repeated a finite three-dimensional obstruction in every summand of an orthogonal direct sum; the resulting projections, of infinite rank and infinite corank, refute every finitely additive $\{0,1\}$-valued measure on $L(H)$~\cite{Alda:1981}.
Döring proved the Kochen--Specker theorem for von Neumann algebras without summands of types $\mathrm I_1$ and $\mathrm I_2$~\cite{Doring:2005}.
The models of finite fragments described below live on the smaller lattice of finite-dimensional subspaces and their orthogonal complements, which contains none of these projections, and give every finite-dimensional subspace probability $0$; the older finite proof and the necessity of infinite branching therefore hold together.

We work in infinitary real-valued logic~\cite{Aubone-Duenez-Iovino:2026a}, which we denote $L_{\omega_1\omega}(\RR)$: the infinitary extension of Keisler's model theory of general structures~\cite[Section~2.1]{Keisler:2023}.
In this logic the truth values of formulas are real numbers, and the probability of an event is the value of a predicate.
Countable suprema and infima express countable additivity and, with a sort for vector sequences, allow the axioms to describe closed spans.
The completeness of the event lattice is a consequence of this construction in standard models.
Over the lattice alone, completeness is not expressible: among the pre-Solèr lattices, which satisfy every condition that characterizes the lattices $L(H)$ except possibly completeness, no sentence of classical two-valued $L_{\omega_1\omega}$ in the language of ortholattices expanded by constants for countably many atoms expresses it (\cref{prop:geometric-nonexpressibility}).
The reason is elementary: a complete ortholattice with a named infinite orthogonal family has a different join for each subset of the family, so a countable elementary substructure, which the downward Löwenheim--Skolem theorem provides, cannot remain complete.
The sequence sort supplies what the lattice cannot: in standard models, completeness becomes the requirement that the closed span of every sequence of vectors be an event, which axioms can state.
Countable additivity is not implied by finite additivity: on the finite-dimensional subspaces of an infinite-dimensional space and their orthogonal complements, the assignment of $0$ to every finite-dimensional subspace and of $1$ to every complement is a finitely additive probability measure (\cref{prop:fc-measure}), which only countable additivity excludes.
In general structures, a vocabulary assigns to each predicate symbol an a priori bound, and a structure interprets function and predicate symbols by arbitrary maps that respect these bounds.
There is no equality symbol, no metric is prescribed, the a priori bounds need not have a common upper bound, and no modulus of uniform continuity is required, in contrast with the metric structures of continuous logic~\cite{BenYaacov-Berenstein-Henson-Usvyatsov:2008}; the two-valued counterpart is the logic without equality studied by Keisler and Miller~\cite[Section~4]{Keisler-Miller:2001}.
General structures admit the lattice operations with their natural interpretation, although these operations are not continuous for the operator-norm distance between projections: the meet of two distinct rays is the zero subspace however close the rays are (\cref{q:finitary}).

Our theory $T$ gives the full description.
The theory has sorts for vectors, for events and for square-summable sequences of vectors, and a predicate $\mu$ for the probability of an event (\cref{sec:vocabulary,sec:axioms}).
Countable suprema and infima express the projection onto the meet of two events, the closed span of a sequence, the spanning of the space by a named orthonormal sequence, and the countable additivity of $\mu$.
The standard models of $T$\footnote{Precisely, a model of $T$ is standard when its vector and sequence sorts are full closed unit balls of complex Hilbert spaces. Nonstandard models exist, among them the countable models given by the downward Löwenheim--Skolem theorem. The vector sort of every model maps into the closed unit ball of a complex Hilbert space by a map that preserves the named affine operations, the inner product and the norm (\cref{lem:envelope}); in a nonstandard model the map may omit vectors of the ball or identify distinct elements, and the sequence sort may fail to be a full Hilbert ball.} represent exactly the density-operator measures on the lattice of closed subspaces of a separable infinite-dimensional complex Hilbert space (\cref{thm:a}).

We return to the example of the opening, which shows that the expressive power afforded by infinitary sentences costs compactness.
The dispersion-free hypothesis, the single closed condition that every probability is $0$ or $1$, yields together with $T$ a theory $T^{01}$ with no models, standard or not, by a proof that uses neither Gleason's theorem nor standardness (\cref{thm:c}).
Yet every finite subset of the axioms of $T^{01}$ has a model, in which the sequence sort contains only finitely supported sequences and the measure vanishes on every finite-dimensional subspace (\cref{prop:bks-game}).
The full theory excludes such models through the instances of axiom (K4), which prescribe the entries of a named sequence whose span contains every vector; countable additivity along this sequence forces a dispersion-free measure to give some ray of the named basis the probability $1$, and no finite subset of the axioms suffices for this step.

Every play of a winning strategy nevertheless ends in a finite contradiction.
A countable theory of closed conditions, each of which asserts that the value of a sentence is at most or at least a given rational number, is inconsistent exactly when a well-founded tree refutes it; the leaves of the tree are finite sets of commitments, placing the values of sentences in closed intervals with rational endpoints, that no structure satisfies together (\cref{thm:d}).
Completeness theorems are known for deductive calculi with infinitary rules, in two-valued $L_{\omega_1\omega}$~\cite{Karp:1964,Lopez-Escobar:1965} and in Keisler's probability logic $L_{AP}$~\cite{Hoover:1978,Keisler:1985}.
For two-valued logic, Väänänen's model existence game characterizes consistency, and its winning strategies describe the classical completeness proofs~\cite{Vaananen:2011}.
We carry the game to real-valued logic in general structures, where neither equality nor a metric is available: the refuting trees arise from a two-player game of model existence (\cref{sec:game}).
The players are Olivia, the challenger, and Josh, the model builder;
Olivia maintains that the theory has no model, while Josh attempts to build one, answering Olivia's challenges about sentence values, witnesses, function values, and countable suprema and infima.
Models are built from Henkin types~\cite[Lemma~3]{Aubone-Duenez-Iovino:2026a} along fair plays, in which every challenge that becomes available is eventually made: a fair play that Josh survives yields a model.
Hence, when the theory is inconsistent, any strategy of Olivia that keeps every play fair is winning, and when it has a model, Josh wins by answering truthfully; the plays that follow Olivia's strategy form the refuting tree.
The game is therefore determined, without appeal to the theorem of Gale and Stewart on open games~\cite{Gale-Stewart:1953}.
Only challenges about countable suprema and infima give Josh infinitely many answers, and a well-founded tree that branches only finitely is finite by K\H{o}nig's lemma, in which case some finite set of axioms is already inconsistent.
For $T^{01}$, every refuting tree therefore has an infinitely branching node (\cref{cor:bks-branching}); for the axioms of $T$ as presented, the finiteness resides in the terminal contradictions, whose commitments concern sentences that may themselves be infinitary.
\Cref{app:henkin} reads Olivia's bookkeeping as a Baire category argument in the space of plays and relates the game to the \v{C}ech-complete space of Henkin types of Aubone, Dueñez and Iovino~\cite{Aubone-Duenez-Iovino:2026a} and to Keisler's model-theoretic forcing~\cite{Keisler:1971,Keisler:1973}.

The certificates are explicit.
Whereas a configuration of Kochen and Specker is chosen in advance and serves against every assignment, in the strategy we construct (\cref{prop:parity-strategy}) the configuration is located by Josh's own commitments in answer to the challenge about countable additivity, and these commitments need not come from any measure.
Once the four-dimensional block is fixed, a parity argument on the eighteen rays of Cabello, Estebaranz and Garc\'ia-Alcaine~\cite{Cabello-Estebaranz-GarciaAlcaine:1996} closes each play that reaches that stage.
The occurrences of a ray in different contexts are distinct terms.
Axioms that do not mention $\mu$ prove that the terms denote the same projection, and only then does the substitutivity axiom (M3) of \cref{sec:axioms} carry the probability across the identity; \cref{sec:obstructions} explains how these two steps enter the finite parity certificate.
The certificate thus records separately the geometric identities and each use of noncontextuality; the eighteen cross-context equalities are indispensable in the Boolean encoding of the configuration, and, with the rays and these equalities exact, the parity contradiction survives explicit numerical tolerances on the hypothetical valuation (\cref{thm:quantitative-parity}).
Whether the least rank of a refuting strategy for $T^{01}$ is finite remains open (\cref{q:bounded}).

We also determine exactly which sets of probabilities the language defines.
The measure takes every value in $[0,1]$ in every standard model, and the sets of probabilities uniformly definable across standard models are exactly the $\sigma$-projective subsets of $[0,1]$ (\cref{thm:definable}).
The uniformity is essential: in a suitably chosen standard model, a single formula defines an arbitrary subset of $(0,1)$ (\cref{prop:nonuniform}).
The full range of the measure is what lets quantification over events realize projection along the interval.
The classification is a theorem of ZFC; whether every set of probabilities definable over $T$ is Lebesgue measurable is a separate set-theoretic question, which we discuss after \cref{thm:e}.
\par\endgroup

\subsection*{Main results of the paper}
We first state the refutation theorem, which concerns every countable theory of closed conditions, and then the theorems on the theory $T$ of quantum probability.
The refutation theorem applies to every theory whose axioms are \emph{closed conditions}, that is, assertions $\varphi\le r$ or $\varphi\ge r$ about sentences $\varphi$, with $r$ rational; an axiom with free variables takes this form once its variables are bound by suprema or infima, and the axioms of the theory $T$ described below are all of this kind.

\begin{mainthm}[Finite obstructions]\label{thm:d}
Let $T'$ be a nonempty countable theory of real-valued $L_{\omega_1\omega}$, in a countable vocabulary, consisting of closed conditions.
Then $T'$ has no models if and only if, in the model existence game of \cref{sec:game}, the player who claims that $T'$ has no models has a winning strategy.
In that case, for every such strategy, the tree of plays in which the same player follows it is well-founded; at each leaf, the commitments made along the play, each of which says that the value of a sentence lies in a closed interval with rational endpoints, form a finite set that is satisfied in no structure.
The tree branches infinitely only at the moves that concern countable suprema and infima; for each fixed finite sequence of the other player's answers to them, the plays have bounded length.
If the tree branches only finitely, the tree is finite and some finite subset of $T'$ has no models; this is the case when no axiom of $T'$ contains a countable supremum or infimum.
\end{mainthm}

Each play that follows a winning strategy ends after finitely many rounds, at a finite contradiction, although the theorem supplies no uniform bound on the lengths of the plays.
The trees are countable and well-founded, so every winning strategy has a countable ordinal rank, which is finite exactly when the plays that follow the strategy have bounded length; the least rank of a winning strategy measures the complexity of refuting the given theory in the game (\cref{sec:game}).

The theory $T$ speaks of three kinds of objects, which the logic keeps apart as sorts.
The elements of the \emph{Hilbert sort} are vectors.
The \emph{Solèr sort} carries the operations of an ortholattice; in the intended interpretation its elements are closed subspaces, each identified, as above, with the orthogonal projection onto it.
The elements of the \emph{sequence sort} are sequences of vectors, through which the completeness of the lattice becomes expressible.
A function symbol names the orthogonal projection of a vector onto an element of the Solèr sort, and a vector belongs to such an element when the projection leaves it fixed.
The predicate $\mu$ assigns a probability to every element of the Solèr sort.

We call a model \emph{standard} when its Hilbert and sequence sorts are full closed unit balls of complex Hilbert spaces with their own operations, inner products and norms (see \cref{def:standard} below).
The theory $T$ has standard models, and, like every countable theory that has a model, it also has models whose sorts are countable, which are not standard.
The Hilbert sort of every model of $T$ is mapped into the closed unit ball of a complex Hilbert space by a map that preserves the named affine operations, the inner product and the norm (see \cref{lem:envelope}).
In a nonstandard model, the image of this map may omit vectors of the ball, distinct elements may have the same image, and the sequence sort may fail to be a full Hilbert ball.
The term ``nonstandard'' is here unrelated to nonstandard analysis.
Our theorems distinguish the conclusions that hold in every model from those that require standardness.

To state how $T$ captures the setting of Gleason's theorem, consider a standard model, whose Hilbert sort is the closed unit ball of a complex Hilbert space $H$.
For an element $a$ of its Solèr sort, let $S_a$ be the closed subspace of $H$ spanned by the vectors of the Hilbert sort that belong to $a$, and let $P_a$ be the orthogonal projection onto $S_a$; for a closed subspace $a$ of $H$, we likewise write $P_a$ for the orthogonal projection onto $a$.
Equality on the Solèr sort is the relation $[a=b]=1$ of \cref{sec:vocabulary}, which holds exactly when $S_a=S_b$ (see \cref{lem:lattice-truth}).

\begin{mainthm}[Gleason's theorem in real-valued $L_{\omega_1\omega}$]\label{thm:a}
There is a countable theory $T$ of real-valued $L_{\omega_1\omega}$ with the following two properties.
\begin{enumerate}[label=(\arabic*)]
\item In every standard model of $T$, the Hilbert sort is the closed unit ball of a separable infinite-dimensional complex Hilbert space $H$, and the map $a\mapsto S_a$ induces an isomorphism of ortholattices from the quotient of the Solèr sort by its equality onto the lattice $L(H)$ of closed subspaces of $H$; moreover, there is a unique density operator $\rho$ on $H$ such that $\mu(a)=\Tr(\rho P_a)$ for every element $a$ of the Solèr sort.
\item Conversely, for every separable infinite-dimensional complex Hilbert space $H$ and every density operator $\rho$ on $H$, there is a standard model of $T$ whose Hilbert sort is the closed unit ball of $H$, whose Solèr sort is $L(H)$, and in which $\mu(a)=\Tr(\rho P_a)$ for every $a\in L(H)$.
\end{enumerate}
\end{mainthm}

In a separable infinite-dimensional Hilbert space, the probabilities that a quantum state assigns to the closed subspaces fill the interval $[0,1]$, even when the state is pure, that is, when its density operator is the orthogonal projection onto a ray (see \cref{lem:range} below).
Thus the range of $\mu$ in every standard model of $T$ is $[0,1]$ (\cref{thm:b}).
If a Borel set $B\subseteq\RR$ meets $[0,1]$, adjoining an axiom that excludes $B$ from the range of $\mu$ therefore gives a theory with no standard models.

The hypothesis that every probability is $0$ or $1$ can be written as a single axiom of the logic,
\[
\mu(a)\wedge(\mathbf 1-\mu(a))\le0.
\]
In this axiom, the symbol $\mathbf 1$ denotes the sentence whose value is $1$ in every structure, and for formulas $\varphi$ and $\psi$ the expression $\varphi\wedge\psi$ denotes the minimum of their values.
Since the axioms of $T$ confine $\mu$ to $[0,1]$, the axiom says that $\mu$ takes only the values $0$ and $1$.
We write $T^{01}$ for the theory $T$ together with this axiom.
Our next theorem shows that the axiom fails not only in the standard models but in every model of $T$.

\begin{mainthm}[The Bell--Kochen--Specker theorem in real-valued $L_{\omega_1\omega}$]\label{thm:c}
The theory $T^{01}$ has no models, standard or not.
\end{mainthm}

The proof of \cref{thm:c}, given in \cref{sec:bks}, uses neither Gleason's theorem nor the standardness of a model.
Every finite set of axioms of $T^{01}$ has a model (\cref{prop:bks-game}), so every winning strategy for $T^{01}$ has an infinitely branching node (\cref{cor:bks-branching}), and no finite set of axioms is a certificate.
The models of finite sets of axioms are not standard, since their sequence sorts contain only finitely supported sequences; \cref{thm:c}, which excludes nonstandard models as well, is therefore more than a restatement of the Hilbert-space theorem.
Whether the least rank of a winning strategy for $T^{01}$ is finite is asked in \cref{q:bounded}.

The general theorem applies to every winning strategy; the configuration of rays and the selection of a four-dimensional subspace described above come from the strategy of \cref{sec:bks-game}.
In the plays of this strategy that reach the parity argument, the rays are those of the eighteen-ray configuration of Cabello, Estebaranz and Garc\'ia-Alcaine~\cite{Cabello-Estebaranz-GarciaAlcaine:1996}, each of which occurs in exactly two of its nine bases.
In the Boolean system that has one Boolean variable for each occurrence of a ray and requires exactly one variable of each basis to be true, the parity contradiction needs all eighteen identities, one between the two occurrences of each ray.
No finite set of axioms of $T^{01}$ forces this localization: in the models of finite sets of axioms constructed in the proof of \cref{prop:bks-game}, the measure vanishes on every finite-dimensional subspace, as in \cref{prop:fc-measure}.

In \cref{sec:obstructions}, after the mathematical development, we give the foundational interpretation.
There we distinguish quantum probabilities from the responses of a hidden-variable model, explain why the geometric identities of projection events and the probabilistic requirements on a dispersion-free valuation enter the game through separate axioms, and relate the contextual presentations to earlier approaches.

The full range of the measure also reveals the expressive power of the logic.
A subset of $[0,1]$ is \emph{definable over $T$} when some formula of the logic, with one free variable $a$ of the Solèr sort, singles it out as follows (see \cref{def:definable-set} below).
In every standard model of $T$ and at every element $a$ of its Solèr sort, the formula takes the value $1$ if $\mu(a)$ lies in the set and the value $-1$ otherwise.
Because formulas may quantify over the Solèr sort, and because the measure takes every value in $[0,1]$, the definable sets form a much larger class than the Borel sets.
The $\sigma$-projective sets, studied by Aguilera, M\"uller and Schlicht~\cite{Aguilera-Muller-Schlicht:2021}, are the members of the least class of subsets of Polish spaces that contains the Borel sets and is closed under complements, countable unions and projections along $\NN^\NN$.
This class contains every projective set, and the example following \cref{thm:definable} shows that the containment is strict.

\begin{mainthm}[Definable sets of probabilities]\label{thm:e}
A subset of $[0,1]$ is definable over $T$ if and only if it is $\sigma$-projective (see \cref{thm:definable}).
\end{mainthm}

Its proof uses both the full range of the measure (\cref{thm:b}) and a standard model whose symbols are interpreted by Borel functions on Polish spaces.
In that model, a formula of finite quantifier depth $k$ defines a $\boldsymbol\Delta^1_{k+1}$ set; the nonprojective sets therefore require countable joins of formulas of unbounded quantifier depth (\cref{prop:qd}).
The uniformity over standard models is essential: a single formula defines an arbitrary subset of $(0,1)$ in a suitably chosen standard model (\cref{prop:nonuniform}).
By \cref{thm:e}, the hypothesis that the probabilities avoid a given nonempty $\sigma$-projective subset of $[0,1]$ can be written as an axiom of the logic, and, by \cref{thm:b}, the theory $T$ together with this axiom has no standard models.
The exact $\sigma$-projective classification, including the strict extension beyond the projective sets, is proved in ZFC.

The classification has a consequence beyond the reach of ZFC, and it connects the formalism with a classical chapter of descriptive set theory.
The $\sigma$-projective sets form the least class that contains the Borel sets and is closed under complements, countable unions and projections, and the class properly contains the projective sets of Lusin and Sierpi\'nski~\cite{Lusin:1925,Sierpinski:1925}.
Lusin proved Suslin's analytic sets Lebesgue measurable and observed that the methods did not reach the next projective level~\cite{Lusin:1925}.
Gödel announced in 1938 that in the constructible universe a projective set of the second level is not Lebesgue measurable~\cite{Godel:1938}; indeed, in that universe some $\boldsymbol\Delta^1_2$ subset of $[0,1]$ is not Lebesgue measurable~\cite{Moschovakis:2009}, and this set is definable over $T$.
Every $\sigma$-projective set belongs to the model $L(\RR)$, so under large cardinal hypotheses, for instance infinitely many Woodin cardinals with a measurable cardinal above them, every subset of $[0,1]$ definable over $T$ is Lebesgue measurable, has the property of Baire and, if uncountable, contains a perfect set~\cite{Larson:2004,Kechris:1995}; measurability comes from the determinacy of the associated games~\cite{Martin-Steel:1989,Aguilera-Muller-Schlicht:2021} and the theorem of Mycielski and \'Swierczkowski~\cite{Mycielski-Swierczkowski:1964}.
Whether every set of probabilities that the logic of quantum probability can define is Lebesgue measurable is therefore not decided by ZFC, and large cardinals decide it.
The countable joins and the quantifiers over the Solèr sort produce this expressive power together; whether formulas without countable joins can define any set other than $\emptyset$ and $[0,1]$ we do not know (see \cref{q:finitary}).

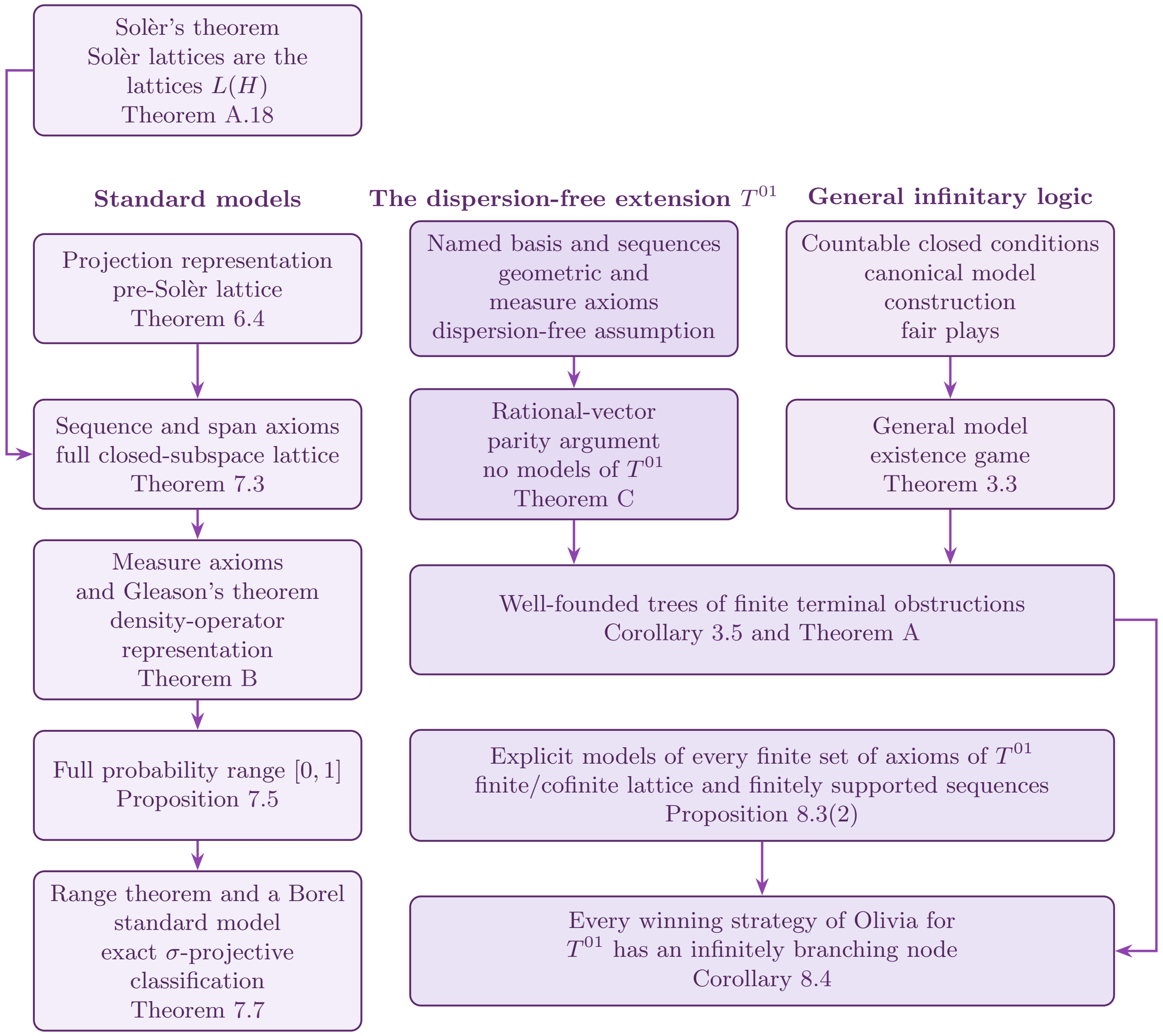
\begin{figure}[htbp]
\centering
\begin{tikzpicture}[
  font=\footnotesize,
  dep/.style={rounded corners=4pt,draw=plum,line width=0.7pt,
    fill=lavender!30,text=plum!85!black,text width=4.0cm,
    minimum height=1.45cm,align=center,inner sep=5pt},
  std/.style={dep,fill=lavender!22},
  allmod/.style={dep,fill=lavender!48},
  general/.style={dep,fill=orchid!12},
  wide/.style={dep,text width=9.0cm,fill=lavender!38},
  arr/.style={-{Stealth[length=2.3mm]},draw=orchid,line width=0.9pt},
  head/.style={text=plum,font=\footnotesize\bfseries,align=center}]
\node[head] at (0,0.55) {Standard models};
\node[head] at (5,0.55) {The dispersion-free extension $T^{01}$};
\node[head,text width=4.3cm] at (10,0.55) {General infinitary logic};

\node[std] (projection) at (0,-0.65)
  {Projection representation\\pre-Solèr lattice\\\cref{thm:t0}};
\node[std] (lattice) at (0,-2.85)
  {Sequence and span axioms\\full closed-subspace lattice\\\cref{thm:t1-iso}};
\node[std] (gleason) at (0,-5.05)
  {Measure axioms\\and Gleason's theorem\\density-operator representation\\\cref{thm:a}};
\node[std] (range) at (0,-7.25)
  {Full probability range $[0,1]$\\\cref{thm:b}};
\node[std] (definability) at (0,-9.45)
  {Range theorem and a Borel\\standard model\\exact $\sigma$-projective classification\\\cref{thm:definable}};
\node[std] (background) at (0,2.25)
  {Solèr's theorem\\Solèr lattices are the\\lattices $L(H)$\\\cref{thm:soler-lattices}};
\draw[arr] (background.west) -- ++(-0.35,0) |- (lattice.west);
\draw[arr] (projection) -- (lattice);
\draw[arr] (lattice) -- (gleason);
\draw[arr] (gleason) -- (range);
\draw[arr] (range) -- (definability);

\node[allmod] (quantum) at (5,-0.65)
  {Named basis and sequences\\geometric and measure axioms\\dispersion-free assumption};
\node[allmod] (inconsistency) at (5,-2.85)
  {Rational-vector parity argument\\no models of $T^{01}$\\\cref{thm:c}};
\node[general] (logical) at (10,-0.65)
  {Countable closed conditions\\canonical model construction\\fair plays};
\node[general] (game) at (10,-2.85)
  {General model existence game\\\cref{thm:game-existence}};
\node[wide] (obstructions) at (7.5,-5.05)
  {Well-founded trees of finite terminal obstructions\\\cref{cor:finite-obstruction,thm:d}};
\node[wide] (fragments) at (7.5,-7.25)
  {Explicit models of every finite set of axioms of $T^{01}$\\finite/cofinite lattice and finitely supported sequences\\\cref{prop:bks-game}(2)};
\node[wide] (branching) at (7.5,-9.45)
  {Every winning strategy of Olivia for $T^{01}$ has an infinitely branching node\\\cref{cor:bks-branching}};
\draw[arr] (quantum) -- (inconsistency);
\draw[arr] (logical) -- (game);
\draw[arr] (inconsistency.south) -- ([xshift=-2.5cm]obstructions.north);
\draw[arr] (game.south) -- ([xshift=2.5cm]obstructions.north);
\draw[arr] (fragments) -- (branching);
\draw[arr] (obstructions.east) -- ++(0.55,0) |- (branching.east);
\end{tikzpicture}
\caption{Direct dependencies among the results, with arrows from a result to the results that use it.}
\label{fig:proof-dependencies}
\end{figure}

The principal dependencies among the results are displayed in \cref{fig:proof-dependencies}.
The representation and definability results concern standard models, whereas the dispersion-free inconsistency holds in every model.
The general model existence theorem supplies the passage from this inconsistency to well-founded trees of finite obstructions.

\subsection*{Organization of the paper}
\begingroup\emergencystretch=1.5em
Readers from quantum foundations can pass from the Overview to the explicit refutation, \cref{prop:parity-strategy}, and to \cref{sec:obstructions}, which relates it to the contextuality of Kochen and Specker.
The refutation theorem, \cref{thm:d}, and the classification of definable sets, \cref{thm:e}, can also be read as results in mathematical logic: the first says what finite refutation means when compactness fails, and the second determines which sets of probabilities the language defines uniformly.
In \cref{sec:setting} we describe the events represented by the theory $T$ and explain, through the obstruction to expressing completeness over the lattice alone (\cref{prop:not-expressible}), why its vocabulary includes a sort of sequences of vectors.
The detailed lattice-theoretic definitions, geometric representation results and their proofs are collected in \cref{app:background}; their principal sources are MacNeille completions~\cite{MacNeille:1937}, the coordinatization of projective geometries~\cite{Baer:1952,Maeda-Maeda:1970,Holland:1995}, and the theorems on orthomodular spaces due to Amemiya and Araki~\cite{Amemiya-Araki:1966} and to Solèr~\cite{Soler:1995}.
The logic is presented in \cref{sec:logic}.
In \cref{sec:game} we organize the model existence theorem as a game, whose determinacy follows from the model construction along fair plays, and obtain the well-founded trees of finite obstructions for every countable theory of closed conditions without models, proving \cref{thm:d}; that section uses only \cref{sec:logic}.
The vocabulary and the axioms of $T$ follow in \cref{sec:vocabulary,sec:axioms}.
The standard models are studied in \cref{sec:t0,sec:t1}, which contain the proofs of \cref{thm:a,thm:b} and of the classification of definable sets of probabilities, \cref{thm:e}.
In \cref{sec:bks} we prove \cref{thm:c}, that no model of $T$ satisfies the dispersion-free axiom, and follow the proof in the game, where every winning strategy branches infinitely (\cref{cor:bks-branching}).
In \cref{sec:obstructions} we compare the finite obstructions with the configurations of Kochen and Specker: what a vertex is, which identifications are made and from which axioms they follow, the tolerances of the intervals, the location of the configuration, the role of $\sigma$-additivity for configurations of rays in infinite dimension, and the topological property of the logic behind each kind of certificate.
The remarks of \cref{sec:remarks} may be omitted on a first reading.
In \cref{app:henkin} we relate the game to Keisler's model-theoretic forcing~\cite{Keisler:1973} and compare it with the space of Henkin types~\cite{Aubone-Duenez-Iovino:2026a}.
\par\endgroup

\clearpage
\section{The geometric setting}\label{sec:setting}

The events of a quantum system whose Hilbert space is $H$ are the closed subspaces of $H$, each determined by the orthogonal projection onto it.
Ordered by inclusion, they form the lattice $L(H)$, in which the meet of two subspaces is their intersection, their join is the closed linear span of their union, and the orthocomplement of a subspace is its orthogonal complement.
Countable additivity, the hypothesis of Gleason's theorem, concerns countable orthogonal families of events, and a formalization of quantum probability must therefore describe not only the subspaces but also the closed spans of countable families of vectors, on which the hypothesis depends.

An \emph{ortholattice} is a bounded lattice with an order-reversing involution $a\mapsto a^\perp$ such that $a\wedge a^\perp=0$ and $a\vee a^\perp=1$.
Two elements $a$ and $b$ are \emph{orthogonal} when $a\le b^\perp$.
An \emph{atom} is a nonzero element with no element strictly between it and $0$; in $L(H)$ the atoms are the one-dimensional subspaces, or \emph{rays}.
A lattice is \emph{complete} when every subset has a join.
The lattice $L(H)$ is complete, but a subset of $L(H)$ that is closed under the ortholattice operations can contain every ray without containing every closed subspace.
For example, in an infinite-dimensional Hilbert space the finite-dimensional subspaces and their orthogonal complements form such a subset (see \cref{ex:finite-cofinite}).

The theory $T$ keeps vectors, events and sequences of vectors in three sorts.
A function symbol describes the projection of a vector onto an event; a vector belongs to an event when this projection leaves it fixed.
Since the logic prescribes no equality, two elements of the event sort may represent the same subspace.
The geometric axioms define an equality relation on that sort, and the representation statements concern its quotient by this relation.
In a standard model of the axioms on vectors and events alone, this quotient is a family of closed subspaces of the given complex Hilbert space, closed under the ortholattice operations and containing every ray, but possibly omitting other closed subspaces (\cref{thm:t0}).
The sequence sort and its axioms supply the missing closed spans: in a standard model of the full theory, the quotient is the whole lattice $L(H)$ (\cref{thm:t1-iso}).
Standardness specifies the interpretation of the vector sorts and is not a consequence of the axioms.

The event quotients just described are pre-Solèr lattices, whose definition is given in \cref{app:background}.
That appendix gives the coordinatization theorem and the theorem of Piron, Amemiya and Araki, and Solèr that characterizes the complete ones as lattices of closed subspaces of infinite-dimensional real, complex or quaternionic Hilbert spaces (\cref{thm:soler-lattices}).
It also describes their MacNeille completions and the obstruction to orthomodularity of those completions.
The distinction between an incomplete event lattice and the full lattice is needed here for a logical reason: completeness cannot be expressed over the lattice alone.

\begin{proposition}\label{prop:not-expressible}
In the vocabulary of ortholattices expanded by constants $p_n$, for $n\in\NN$, no sentence of classical two-valued $L_{\omega_1\omega}$ is equivalent to completeness among ortholattices in which the $p_n$ are mutually orthogonal nonzero elements.
\end{proposition}

\begin{proof}
Suppose that a sentence $\varphi$ expressed completeness in this class.
The lattice $L(\ell^2)$, expanded by the rays $p_n=\CC e_n$ of its standard orthonormal basis, satisfies $\varphi$.
Apply the downward Löwenheim--Skolem theorem to a countable fragment containing $\varphi$ and sentences expressing the ortholattice axioms and the stated conditions on the $p_n$.
It gives a countable ortholattice $\Lambda$ in the same class that satisfies $\varphi$, and hence is complete.
For each subset $A\subseteq\NN$, completeness supplies $q_A=\bigvee_{n\in A}p_n$.
If $n\notin A$, every $p_m$ with $m\in A$ lies below $p_n^\perp$, so $q_A\le p_n^\perp$; since $p_n\ne0$, it follows that $p_n\not\le q_A$.
If $n\in A$, then $p_n\le q_A$.
Thus distinct subsets give distinct joins, contradicting the countability of $\Lambda$.
\end{proof}

The same obstruction applies in real-valued $L_{\omega_1\omega}$ when the lattice equality and order are represented by two-valued predicates, with their congruence axioms.
For a countable fragment that contains these axioms and a sentence proposed to express completeness, the expanded lattice $L(\ell^2)$ has a countable substructure in which every sentence of the fragment keeps its value.
The quotient of this substructure by equality is a countable ortholattice with the named orthogonal nonzero elements.
Congruence makes every formula invariant under passage to the quotient, so the proposed sentence retains its value there, whereas the quotient cannot be complete by the argument above.
A stronger nonexpressibility statement for pre-Solèr lattices, including the condition that their MacNeille completion be orthomodular, is proved in \cref{prop:geometric-nonexpressibility}.
The vocabulary of $T$ therefore adds to the lattice a sort of vectors and a sort of sequences of vectors.
In standard models, completeness then becomes the requirement that the closed span of every sequence of vectors be an event, which axioms can state; standardness, which no sentence expresses, fixes the vectors.

\section{The logic}\label{sec:logic}

Probabilities and norms are real quantities, and so are the real and imaginary parts of inner products.
For a finitely additive measure, countable additivity amounts to continuity: the measure of the join of a sequence of events is the supremum of the measures of its finite joins.
Real-valued $L_{\omega_1\omega}$ expresses these quantities and this continuity, since its formulas take real values and countable suprema are among its connectives.
We recall this logic and the Henkin types from which its models are built; the model existence game of \cref{sec:game} uses nothing beyond this section.

\subsection{Real-valued \texorpdfstring{$L_{\omega_1\omega}$}{L-omega1-omega}}
We use real-valued infinitary logic~\cite{Aubone-Duenez-Iovino:2026a}, the infinitary extension of Keisler's model theory of general structures~\cite[Section~2.1]{Keisler:2023}.
We recall its definitions as far as we use them.

\emph{Vocabularies and structures.}
A vocabulary consists of sorts, of function symbols, each with the sorts of its arguments and of its values, and of predicate symbols.
Each predicate symbol $R$ comes with the sorts of its arguments and with a nonnegative real number $r_R$, its \emph{a priori bound}.
In this paper the a priori bound of every predicate symbol is a positive integer.
The vocabulary contains a nullary predicate symbol $\mathbf 1$, the \emph{logical unity constant}, with a priori bound $1$.
A structure for the vocabulary interprets each sort as a nonempty set, each function symbol as a function from the product of its argument sorts to its value sort, and each predicate symbol $R$ as a function from the product of its argument sorts to the interval $[-r_R,r_R]$; the logical unity constant is interpreted as $1$.
No topology is prescribed with respect to which the functions or the predicates must be continuous, and the vocabulary has no equality symbol.

\emph{Formulas.}
The atomic formulas are the expressions $R(\bar t)$, where $R$ is a predicate symbol and $\bar t$ is a tuple of terms of the appropriate sorts.
Formulas are built from the atomic formulas by the connectives $s\cdot\varphi$, for $s\in\QQ$, $\varphi+\psi$ and $\varphi\vee\psi$, by the connectives $\sup_x$, one for each variable $x$, and by countable joins $\bigvee_n\varphi_n$.
A countable join is formed from a family $(\varphi_n)_{n\in\NN}$ of formulas whose free variables lie in one finite set and whose a priori bounds, defined next, have a finite supremum.
The a priori bound of $R(\bar t)$ is $r_R$; the a priori bounds of $s\cdot\varphi$, $\varphi+\psi$, $\varphi\vee\psi$, $\sup_x\varphi$ and $\bigvee_n\varphi_n$ are $|s|r_\varphi$, $r_\varphi+r_\psi$, $\max(r_\varphi,r_\psi)$, $r_\varphi$ and $\sup_nr_{\varphi_n}$, respectively.
In a structure, the connectives $s\cdot\varphi$, $\varphi+\psi$ and $\varphi\vee\psi$ are interpreted as multiplication by $s$, as the sum and as the maximum, the connective $\sup_x$ as the supremum over the sort of $x$, and a countable join as the supremum of the values of its members.
By induction on formulas, every formula $\varphi$ takes its values in $[-r_\varphi,r_\varphi]$ in every structure, and hence in the interval $I_\varphi=[-\bar r_\varphi,\bar r_\varphi]$, where the number $\bar r_\varphi$ is the least integer not smaller than $r_\varphi$.
The number $r_\varphi$ does not change when terms are substituted for the free variables of $\varphi$.

\emph{Abbreviations.}
For a rational number $q$ we write $q$ for the sentence $q\cdot\mathbf 1$, whose value is $q$ in every structure; in particular, we write $0$ for $0\cdot\mathbf 1$.
We write $-\varphi$ for $(-1)\cdot\varphi$ and $\varphi-\psi$ for $\varphi+(-\psi)$.
The minimum, the infimum over the sort of $x$ and the countable infimum are expressed by $\varphi\wedge\psi$, which abbreviates $-(-\varphi\vee-\psi)$, by $\inf_x\varphi$, which abbreviates $-\sup_x-\varphi$, and by $\bigwedge_n\varphi_n$, which abbreviates $-\bigvee_n-\varphi_n$.
Finally, the absolute value $|\varphi|$ abbreviates $\varphi\vee-\varphi$.
The connectives do not include multiplication, but the square of a formula can be written as a countable join of affine expressions in it, and the product of two formulas as a rational combination of two squares.
For a formula $\varphi$, let $\varphi^2$ denote the countable join of the formulas $(2q)\cdot\varphi-q^2$, for the rational numbers $q$ with $|q|\le r_\varphi$.

\begin{lemma}\label{lem:squares}
For every formula $\varphi$, the join $\varphi^2$ is a formula with a priori bound at most $3r_\varphi^2$, and its value in every structure, at every tuple, is the square of the value of $\varphi$.
\end{lemma}

\begin{proof}
The member of index $q$ has a priori bound $2|q|r_\varphi+q^2\le3r_\varphi^2$, and its free variables are those of $\varphi$; the join is therefore a formula.
The values of $\varphi$ lie in $[-r_\varphi,r_\varphi]$, and for a real number $u$ in this interval we have $2qu-q^2=u^2-(u-q)^2\le u^2$; since the rational numbers of the interval are dense in it, the term $(u-q)^2$ can be made arbitrarily small, and the supremum of the values of the members at $u$ is therefore $u^2$.
\end{proof}

For formulas $\varphi$ and $\psi$, the product $\varphi\psi$ abbreviates $\tfrac14(\varphi+\psi)^2-\tfrac14(\varphi-\psi)^2$, whose value is the product of the values of $\varphi$ and $\psi$.
Polynomial expressions with rational coefficients in formulas are therefore formulas.

\emph{Conditions and theories.}
The statements of the logic are not formulas but conditions, which compare the value of a sentence with a rational number.
For a sentence $\varphi$ and a rational number $r$, a \emph{condition} is an expression $\varphi>r$ or $\varphi<r$.
A structure $\fM$ satisfies $\varphi>r$ when the value $\varphi^\fM$ of $\varphi$ in $\fM$ is greater than $r$, and it satisfies $\varphi<r$ when this value is less than $r$.
A \emph{theory} is a set of conditions, and a \emph{model} of a theory is a structure that satisfies every condition in it.
Non-strict inequalities are recovered as countable theories: the \emph{closed condition} $\varphi\ge r$ denotes the theory $\{\varphi>r-\tfrac1n:n\ge1\}$, the closed condition $\varphi\le r$ denotes $\{\varphi<r+\tfrac1n:n\ge1\}$, and the closed condition $\varphi=r$ denotes the union of the two.
An \emph{axiom} is a closed condition on a formula $\varphi(\bar x)$, possibly without free variables, required of all values of $\bar x$.
The axiom $\varphi(\bar x)\le r$ is the closed condition $\sup_{\bar x}\varphi(\bar x)\le r$, and the axiom $\varphi(\bar x)\ge r$ is $\inf_{\bar x}\varphi(\bar x)\ge r$; the axiom $\varphi(\bar x)=r$ is the union of the two.
For formulas $\varphi$ and $\psi$, we write $\varphi\le\psi$ and $\varphi=\psi$ for the axioms $\varphi-\psi\le0$ and $\varphi-\psi=0$, which are closed conditions when $\varphi$ and $\psi$ are sentences.

\emph{Truth values.}
As in the work of Aubone, Dueñez and Iovino~\cite{Aubone-Duenez-Iovino:2026a}, the value $1$ stands for true and the value $-1$ for false.
A formula is \emph{two-valued} in a structure if it takes only the values $1$ and $-1$ there.
For formulas $\theta$ and $\theta'$ that are two-valued in a structure, negation, disjunction and conjunction are expressed there by $-\theta$, $\theta\vee\theta'$ and $\theta\wedge\theta'$.
The existential and the universal quantifier are expressed by $\sup_x$ and $\inf_x$, and countable disjunction and conjunction, for families whose a priori bounds have a finite supremum, by $\bigvee_n$ and $\bigwedge_n$.
A supremum of values in $\{-1,1\}$ over a nonempty set is attained, and the connectives $\sup_x$ and $\inf_x$, applied to a two-valued formula, are therefore the traditional quantifiers over the sort of $x$.
We use them in this way over the Solèr sort, on which equality and order are given by formulas that are two-valued in every model of the theory $T_0$ of \cref{sec:axioms} (see \cref{lem:lattice-truth} below).
A two-valued formula $\theta$ is asserted by the axiom $\theta\ge1$, and the implication from $\theta$ to $\theta'$ is asserted by the axiom $-\theta\vee\theta'\ge1$.
If $\theta$ is two-valued and $\psi$ is any formula, the axiom $\theta\wedge\psi\le0$ says that $\theta$ implies $\psi\le0$: where $\theta$ is false, the minimum is at most $-1$.
Norms, the measure and the real and imaginary parts of inner products keep their numerical values; every statement about these quantities is a closed condition or an axiom with a rational threshold.
A measure or an assignment that takes only the values $0$ and $1$ is called $\{0,1\}$-valued, never two-valued, since the latter term refers to the truth values $1$ and $-1$.

\emph{Truncations.}
A countable join requires a common a priori bound on its member formulas, not merely a common bound on their values.
For example, the formulas $((n\cdot\varphi)\wedge1)\vee0$ take values in $[0,1]$, but their a priori bounds are $\max(nr_\varphi,1)$.
When $r_\varphi>0$, these bounds grow without limit, so the formulas cannot be joined over $n$.

For the atomic predicates used below, we therefore name their truncations separately in the vocabulary of \cref{sec:vocabulary}, each with a priori bound $1$.
The axioms (Tr) of \cref{sec:axioms} identify each symbol with the corresponding truncated expression, so the countable joins of the symbols are legitimate formulas that give exact tests for strict inequalities in structures satisfying (Tr).\footnote{A different presentation of a normed space can use sorts for all closed balls of positive integer radii, with inclusion maps between the balls and operations between the appropriate sorts.
On a ball of radius $n$ the norm has bound $n$, and on two Hilbert balls of radii $n$ and $m$ the real and imaginary parts of the inner product have bound $nm$.
A scaling map from the unit ball to a larger ball, followed by a separately named truncated predicate of bound $1$ there, can supply a steep predicate on the unit ball with bound $1$.
Here we name the truncations directly; for the measure on the Solèr sort they are likewise new predicate symbols whose values (Tr) fixes.}

Following Aubone, Dueñez and Iovino~\cite{Aubone-Duenez-Iovino:2026a}, we call a \emph{fragment} a set of formulas that contains the first-order formulas, those in which all joins are finite, and is closed under subformulas, under substitution of terms and under the connectives $s\cdot\varphi$, $\varphi+\psi$, $\varphi\vee\psi$ and $\sup_x$.

\subsection{Henkin types and model existence}\label{sec:henkin-types}
In the model existence argument of \cref{sec:game}, we build a model from an assignment of values to sentences that respects the semantics of the logic.
We recall the definition and the Henkin conditions from the topological treatment of model existence by Aubone, Dueñez and Iovino~\cite{Aubone-Duenez-Iovino:2026a}: an assignment satisfying these conditions is a Henkin type and comes from a structure whose elements are named by closed terms.
We record these conditions for the later game and its comparison with the topological construction in \cref{app:henkin}; the geometric representation arguments that follow do not use them.

For a countable fragment $\cF$ of real-valued $L_{\omega_1\omega}$ for a countable vocabulary, such as the vocabulary $\cL$ of \cref{sec:vocabulary}, let $C$ be a set of new constants, countably many of each sort.
Let $\cF_C$ be the fragment that $\cF$ generates in the vocabulary enlarged by $C$, and let $S^{\cF}_C$ be the set of sentences of $\cF_C$.
An \emph{$S^{\cF}_C$-Henkin type} is an element $\mathfrak m$ of the product of the intervals $[-r_\varphi,r_\varphi]$, over the sentences $\varphi$ of $S^{\cF}_C$, that satisfies the following \emph{Henkin conditions}:
\begin{enumerate}[label=(H\arabic*$'$)]
\item $\mathbf 1^{\mathfrak m}=1$;
\item $(s\cdot\varphi)^{\mathfrak m}=s\,\varphi^{\mathfrak m}$ for every rational number $s$;
\item $(\varphi+\psi)^{\mathfrak m}=\varphi^{\mathfrak m}+\psi^{\mathfrak m}$;
\item $(\varphi\vee\psi)^{\mathfrak m}=\max(\varphi^{\mathfrak m},\psi^{\mathfrak m})$;
\item $(\sup_x\varphi)^{\mathfrak m}=\sup_t\varphi(t)^{\mathfrak m}$, the supremum being taken over the closed terms $t$ of the sort of $x$;
\item $(\bigvee_n\varphi_n)^{\mathfrak m}=\sup_n\varphi_n^{\mathfrak m}$.
\end{enumerate}
Here $\varphi^{\mathfrak m}$ denotes the coordinate of $\mathfrak m$ at $\varphi$, and the sentences involved range over $S^{\cF}_C$.
With the topology induced by the product topology, the $S^{\cF}_C$-Henkin types form the \emph{space of Henkin types} $\cH^{\cF}_C$, a $G_\delta$ subset of the product and therefore a \v{C}ech-complete space, hence a Baire space by Čech's form of the Baire category theorem~\cite[Proposition~6 and Corollary~7]{Aubone-Duenez-Iovino:2026a}.
Every Henkin type is the $S^{\cF}_C$-type of a structure for the enlarged vocabulary in which every element is the value of a closed term~\cite[Lemma~3]{Aubone-Duenez-Iovino:2026a}.

\historypara{Historical note.}
Carol Karp, a student of Henkin's, wrote her dissertation of 1959 on languages with expressions of infinite length and published a book on them in 1964~\cite{Karp:1964}.
In 1969 Makkai applied Smullyan's method of consistency properties to $L_{\omega_1\omega}$, building models in the manner of Henkin's completeness proof, and extended it to the infinitary logics on admissible sets~\cite{Makkai:1969,Makkai:1969b}.
In 1971 Keisler made the model existence theorem for consistency properties a central method of his book on $L_{\omega_1\omega}$~\cite[Chapter~3]{Keisler:1971}, which also contains the omitting types theorem for countable fragments~\cite[Chapter~11]{Keisler:1971}.
In 1973 Keisler built countable models by forcing with finite sets of sentences satisfiable in a given class of models~\cite[p.~100, Example~1.2]{Keisler:1973}.
In the same work he replaced principal types by $\forall\exists$ formulas over a set of formulas of a countable fragment~\cite[p.~107]{Keisler:1973}.
His general omitting types theorem gives conditions under which a $\forall\exists$ class contains a countable model in which countably many given $\forall\exists$ sentences hold~\cite[p.~108]{Keisler:1973}.
Keisler's finite conditions and $\forall\exists$ formulas have counterparts in the model existence game of \cref{sec:game} and in the contextual presentations of \cref{sec:obstructions} (see \cref{app:henkin}).
In 1974 Morley derived the omitting types theorem for countable fragments from the Baire category theorem, applied to a space of countable structures, which he showed to be a Baire space~\cite{Morley:1974}.
In 2001 Keisler and Miller carried the theorem into first-order logic without equality, where its usual formulation does not hold.
There the model is built from closed terms, and principal types give way to locally realized types~\cite[Section~4]{Keisler-Miller:2001}.

Model theory with truth values in the unit interval goes back to Chang's work of 1958 and 1961 on the infinite-valued logic of Łukasiewicz~\cite{Chang:1958,Chang:1961}.
The continuous model theory of Chang and Keisler of 1966 took its truth values in a compact Hausdorff space~\cite{Chang-Keisler:1966}.
Henson introduced the logic of positive bounded formulas in 1976 to study a Banach space through its nonstandard hull~\cite{Henson:1976}.
In metric approaches to real-valued model theory, a distance predicate takes the place of equality, and continuity requirements ensure that the other symbols respect the resulting uniform structure.\footnote{In \L{}ukasiewicz predicate logic, the congruence axioms for a similarity, or graded equality, make every predicate $1$-Lipschitz for the distance given by one minus the similarity (see \cite[Section~5.6]{Hajek:1998} and the earlier work collected there).}
Uniform continuity with moduli prescribed by the signature was the setting of Iovino's thesis of 1994 and of his papers of 1999 on stable Banach spaces~\cite{Iovino:1994a,Iovino:1999a,Iovino:1999b}; in 2002 Henson and Iovino presented the model theory of normed space structures in it~\cite[Sections~2 and~3]{Henson-Iovino:2002a}.
In 2008 Ben Yaacov, Berenstein, Henson and Usvyatsov set out the model theory of metric structures in continuous logic~\cite{BenYaacov-Berenstein-Henson-Usvyatsov:2008}.
In 2009 Ben Yaacov and Iovino introduced an infinitary continuous logic for model-theoretic forcing~\cite{BenYaacov-Iovino:2009}.
In 2010 Ben Yaacov and Usvyatsov built moduli of uniform continuity into the signatures of continuous logic~\cite[Section~2]{BenYaacov-Usvyatsov:2010}.
In 2014 Eagle proved the omitting types theorem for infinitary $[0,1]$-valued logic, whose metric structures have a master metric with respect to which the predicates are uniformly continuous~\cite[Definitions~1.4 and~1.5]{Eagle:2014}.
In the model theory for real-valued structures that Keisler set out in 2023, a general structure requires no metric signature: its function and predicate symbols may be interpreted by arbitrary mappings, subject only to the boundedness of the predicates~\cite[Sections~1 and~2.1]{Keisler:2023}.
General structures correspond to first-order structures without equality, as metric structures with a distance predicate correspond to first-order structures with equality~\cite[Section~1]{Keisler:2023}.
For vocabularies with countably many function and predicate symbols, every general structure expands to one with a distance predicate that is a uniform limit of formulas~\cite[Theorem~3.3.4]{Keisler:2023}.
The logic of this paper, whose conditions are strict inequalities, is the infinitary extension of Keisler's model theory for real-valued structures; Aubone, Dueñez and Iovino introduced it, built its models from closed terms and proved its omitting types theorem by Morley's method in the \v{C}ech-complete space of Henkin types~\cite{Aubone-Duenez-Iovino:2026a}.

\section{The model existence game and finite obstructions}\label{sec:game}

In first-order logic, compactness turns a theory without models into a finite set of axioms without models, much as finitely many vectors suffice in the Kochen--Specker theorem~\cite{Kochen-Specker:1967}.
The countable suprema and infima of real-valued $L_{\omega_1\omega}$, through which $\sigma$-additivity is expressed, carry the logic beyond the reach of compactness, and in this section we prove, for every countable theory of closed conditions, the theorem that takes its place, \cref{thm:d} of the Introduction.
To obtain finite obstructions in this setting, we organize the model existence argument of Aubone, Dueñez and Iovino~\cite[Lemma~3]{Aubone-Duenez-Iovino:2026a} as a game in which each play of a winning strategy reaches a finite contradiction.
The argument realizes Henkin types, the value assignments of \cref{sec:henkin-types}, through the Henkin construction of canonical structures whose elements are named by closed terms.
The same construction underlies the omitting types theorem of Aubone, Dueñez and Iovino, whose Baire category proof we compare with the game in \cref{app:henkin}.

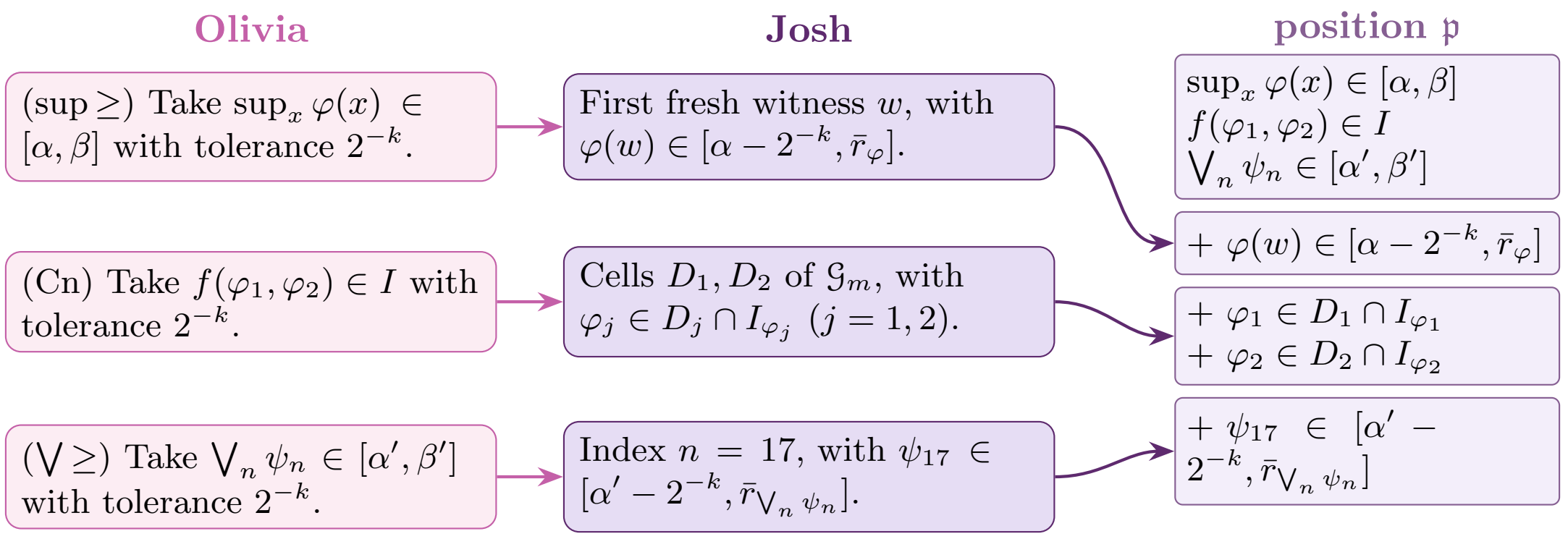
\begin{figure}[ht]
\centering
\begin{tikzpicture}[font=\footnotesize,
  ol/.style={draw=rose,fill=blush!45,rounded corners=4pt,text width=4.2cm,align=left,inner sep=4pt},
  jo/.style={draw=plum,fill=lavender!45,rounded corners=4pt,text width=4.2cm,align=left,inner sep=4pt},
  card/.style={draw=keyframe,fill=keyfill,rounded corners=2pt,text width=3.3cm,align=left,inner sep=3pt}]
\node[font=\small\bfseries,text=rose] at (0,0.9) {Olivia};
\node[font=\small\bfseries,text=plum] at (5.1,0.9) {Josh};
\node[font=\small\bfseries,text=keyframe] at (10.2,0.9) {position $\mathfrak p$};
\node[ol] (o1) at (0,0) {($\sup{\ge}$) Take $\sup_x\varphi(x)\in[\alpha,\beta]$ with tolerance $2^{-k}$.};
\node[jo] (j1) at (5.1,0) {First fresh witness $w$, with $\varphi(w)\in[\alpha-2^{-k},\bar r_\varphi]$.};
\node[ol] (o2) at (0,-1.6) {(Cn) Take $f(\varphi_1,\varphi_2)\in I$ with tolerance $2^{-k}$.};
\node[jo] (j2) at (5.1,-1.6) {Cells $D_1,D_2$ of $\cG_m$, with $\varphi_j\in D_j\cap I_{\varphi_j}$ ($j=1,2$).};
\node[ol] (o3) at (0,-3.2) {($\bigvee{\ge}$) Take $\bigvee_n\psi_n\in[\alpha',\beta']$ with tolerance $2^{-k}$.};
\node[jo] (j3) at (5.1,-3.2) {Index $n=17$, with $\psi_{17}\in[\alpha'-2^{-k},\bar r_{\bigvee_n\psi_n}]$.};
\foreach \i in {1,2,3}{\draw[-{Stealth},rose,thick] (o\i.east) -- (j\i.west);}
\node[card] (c0) at (10.2,0) {$\sup_x\varphi(x)\in[\alpha,\beta]$\\ $f(\varphi_1,\varphi_2)\in I$\\ $\bigvee_n\psi_n\in[\alpha',\beta']$};
\node[card,anchor=north] (c1) at ([yshift=-3pt]c0.south) {$+\ \varphi(w)\in[\alpha-2^{-k},\bar r_\varphi]$};
\node[card,anchor=north] (c2) at ([yshift=-3pt]c1.south) {$+\ \varphi_1\in D_1\cap I_{\varphi_1}$\\ $+\ \varphi_2\in D_2\cap I_{\varphi_2}$};
\node[card,anchor=north] (c3) at ([yshift=-3pt]c2.south) {$+\ \psi_{17}\in[\alpha'-2^{-k},\bar r_{\bigvee_n\psi_n}]$};
\foreach \i in {1,2,3}{\draw[-{Stealth},plum,thick] (j\i.east) to[out=0,in=180] (c\i.west);}
\end{tikzpicture}
\caption{Three rounds of the game.}
\label{fig:round}
\end{figure}

The two players are Olivia, who claims that a given countable theory has no models, and Josh, who attempts to describe one by answering her challenges.
Each answer adds information that all later answers must retain (see \cref{fig:round}).
Olivia wins when Josh has no legal answer or the accumulated intervals for a basic predicate value have no common point in its prescribed range; Josh wins if he can continue indefinitely.
The game follows Väänänen's model existence game for two-valued logic~\cite{Vaananen:2011}.
We compare its positions and rules with Keisler's forcing in \cref{app:henkin}.

In \cref{thm:game-existence} below we show that a nonempty countable theory of closed conditions has no models exactly when Olivia has a winning strategy.
The plays that follow such a strategy form a well-founded tree, with a finite obstruction at each leaf: statements about sentence values that no structure can satisfy together.
Only challenges to countable suprema or infima can give Josh infinitely many answers; for a fixed winning strategy, the plays whose countable answers form any one prescribed finite sequence have bounded length.
In \cref{sec:bks-game} we return to the quantum theory: for the dispersion-free extension of the theory $T$ of \cref{thm:a}, every finite set of axioms has a model (\cref{prop:bks-game}), and in the strategy that follows the proof of \cref{thm:c}, Olivia first locates an atom of large measure and then records enough of the finite geometry to run the parity argument.
Plays that reach this argument end at positions whose commitments concern a Kochen--Specker configuration; the other plays end at contradictions exposed along the way.

\subsection{The setting}\label{sec:game-setting}
In this section we work with a countable many-sorted vocabulary $\cL$ for the logic of \cref{sec:logic} and with a nonempty countable theory $T$ in $\cL$; within the section, these letters denote an arbitrary vocabulary and theory, not the vocabulary of \cref{sec:vocabulary} and the theory of \cref{thm:a}.
As usual, every sort of a structure is nonempty.

In the game, Josh commits the values of sentences to intervals, and the rules require for each sentence an interval fixed in advance that contains every possible value of the sentence.
These are the intervals $I_\varphi=[-\bar r_\varphi,\bar r_\varphi]$ of \cref{sec:logic}, whose integer endpoints are determined by the a priori bound $r_\varphi$ and do not change when terms are substituted for the free variables of $\varphi$.
In the game, the infima $\inf_x$ over a sort and the countable infima $\bigwedge_n$, which abbreviate $-\sup_x-\varphi$ and $-\bigvee_n-\varphi_n$, are treated as primitive, with rules of their own, and the interval of $\inf_x\varphi$ or $\bigwedge_n\varphi_n$ is that of $\sup_x\varphi$ or $\bigvee_n\varphi_n$.
The formulas $s\cdot\varphi$, $\varphi+\psi$ and $\varphi\vee\psi$, by contrast, are handled by the rule (Cn) of \cref{def:game} through the functions $t\mapsto st$, $(t,u)\mapsto t+u$ and $(t,u)\mapsto\max(t,u)$ that interpret them.
We call these functions the \emph{connective functions}; each of them is continuous.
Finally, the logical unity constant $\mathbf 1$, a nullary predicate symbol, takes only the value $1$.
For the definition of infeasible positions below, we write $\mathcal V_{\mathbf 1}=\{1\}$ and $\mathcal V_R=[-r_R,r_R]$ for every other predicate symbol $R$.

The axioms enter the game through one of its rules, which is simplest to state when all axioms have the same form.
As everywhere in the paper, an axiom of $T$ is a closed condition required of all values of its free variables.
If $\bar x$ lists the free variables of a formula $\varphi$, then, by the conventions of \cref{sec:logic}, an axiom $\varphi(\bar x)\le r$ is the closed condition $\sup_{\bar x}\varphi\le r$ and an axiom $\varphi(\bar x)\ge r$ is $\inf_{\bar x}\varphi\ge r$.
An axiom $\varphi(\bar x)=r$ is the union of these two closed conditions.
Every axiom of $T$ may therefore be regarded as a closed condition $\psi\le r$ or $\psi\ge r$ for a sentence $\psi$ and a rational number $r$, an axiom $\varphi(\bar x)=r$ counting as two.
The one or two closed conditions on sentences that thus correspond to an axiom constitute its \emph{normal form}.
Theories that contain strict conditions are discussed in \cref{rem:strict}.

Josh's description needs names for elements before their predicate values and function images can be tested.
For each sort $\mathbf{S}$, we therefore add to $\cL$ a set $\cW_{\mathbf S}=\{w_{\mathbf S,i}:i\in\NN\}$ of new constants of sort $\mathbf S$, the \emph{witnesses}.
We write $\cW$ for their union, $\cL_{\cW}$ for the expanded vocabulary, and $w,w'$ for witnesses.
The game treats a constant symbol of $\cL$ as a function symbol with no arguments, so that its value is named by the same rule as any other function output.
A sentence of $\cL_{\cW}$ is \emph{basic} when it has the form $R(\bar w)$, where $R$ is a predicate symbol and $\bar w$ is a tuple of witnesses of the appropriate sorts.
The tuple may be empty; in particular, the sentence $\mathbf 1$ is basic.

The description retains both numerical bounds and the function values needed to interpret them.
A commitment places a sentence value in an interval, while a record names the output of a function at given witnesses.
The axiom rule introduces commitments, and the totality rule supplies function records.
The other rules test commitments by following the definition of sentence values: a connective leads to its arguments, a countable supremum or infimum to its members, and a supremum or infimum over a sort to instances at witnesses.
For an atomic sentence, the records replace compound terms by witnesses until a basic sentence remains.
Thus a conflict among composite assertions can eventually be exposed as incompatible intervals for a basic predicate value.

\begin{definition}\label{def:position}
A \emph{commitment} is an expression $\varphi\in I$, where $\varphi$ is a sentence of $\cL_{\cW}$ in which only finitely many witnesses occur and $I=[\alpha,\beta]$ is an interval with rational endpoints $\alpha\le\beta$; the interval may be a single point.
A \emph{record} is an expression $F(\bar w)\mapsto w'$, where $F$ is a function symbol and $\bar w,w'$ are witnesses of the appropriate sorts.
A \emph{position} is a finite set $\mathfrak p$ of commitments and records that contains at most one record $F(\bar w)\mapsto w'$ for each function symbol $F$ and each tuple $\bar w$.
A witness \emph{occurs} in $\mathfrak p$ when it occurs in a commitment or a record of $\mathfrak p$.
In a structure for $\cL$, a position is \emph{satisfied} under an interpretation of the witnesses that occur in it when each of its commitments holds and, for each of its records $F(\bar w)\mapsto w'$, the term $F(\bar w)$ and the witness $w'$ have the same value.
A position $\mathfrak p$ is \emph{infeasible} when there is a basic sentence $R(\bar w)$ such that the intervals $I$ for which $\mathfrak p$ contains the commitment $R(\bar w)\in I$, together with $\mathcal V_R$, have empty intersection.
\end{definition}

An infeasible position is contradictory in the most elementary way: it assigns to a basic sentence $R(\bar w)$ finitely many rational intervals that have no common point in $\mathcal V_R$.
Infeasibility is the elementary conflict detected by the rules; satisfiability asks whether all the commitments and records can hold together in a structure.
A position may already be unsatisfiable without being infeasible, and Olivia's subsequent challenges may expose the contradiction at the level of basic sentences.

To test a commitment on a connective, Olivia needs bounds on its arguments that control the resulting value.
Josh therefore answers with small grid cells rather than exact real numbers.
The continuity of the connective on the product of the fixed intervals assigned to its argument sentences determines how fine the cells must be; a prescribed widening of the committed interval allows for their width.
For $m\in\NN$, the \emph{dyadic grid} $\cG_m$ is the set of intervals $[i2^{-m},(i+1)2^{-m}]$, for $i\in\ZZ$.
Let $\ell$ be $1$ or $2$, let $f\colon\RR^\ell\to\RR$ be a connective function, let $\varphi_1,\dots,\varphi_\ell$ be sentences, and let $k\in\NN$ be given.
The function $f$ is uniformly continuous on the compact box $I_{\varphi_1}\times\dots\times I_{\varphi_\ell}$; we write $m(f,\bar\varphi,k)$ for the least $m\in\NN$ such that $|f(\xi)-f(\xi')|\le2^{-k}$ for all points $\xi,\xi'$ of the box with $|\xi_j-\xi'_j|\le2^{-m}$ for $j=1,\dots,\ell$.
For an interval $I=[\alpha,\beta]$ and $k\in\NN$, we write $I^{(k)}=[\alpha-2^{-k},\beta+2^{-k}]$.

When Josh chooses a witness, the rules allow witnesses already in the position and one designated fresh witness.
This restriction leaves finitely many choices at each such move, while permitting a previously unused witness whenever one is needed.
The choices of an index for a countable supremum or infimum retain their countable branching; this distinction is the basis of the finite-obstruction analysis in the third subsection.

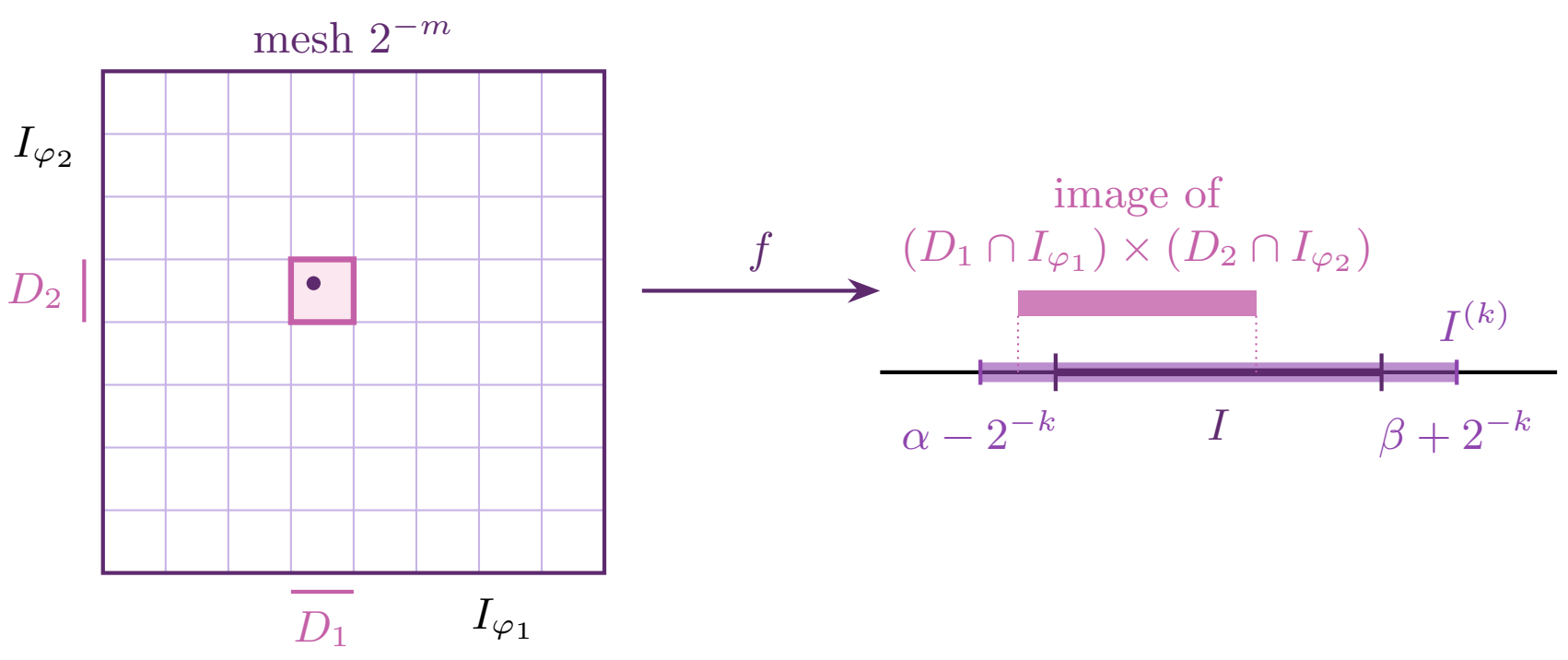
\begin{figure}[ht]
\centering
\begin{tikzpicture}[font=\footnotesize,scale=0.9]
\fill[blush!60] (1.5,2.0) rectangle (2.0,2.5);
\draw[lavender] (0,0) grid[step=0.5] (4,4);
\draw[plum,thick] (0,0) rectangle (4,4);
\draw[rose,very thick] (1.5,2.0) rectangle (2.0,2.5);
\fill[plum] (1.68,2.31) circle (1.6pt);
\draw[rose,thick] (1.5,-0.15) -- (2.0,-0.15) node[midway,below] {$D_1$};
\draw[rose,thick] (-0.15,2.0) -- (-0.15,2.5) node[midway,left] {$D_2$};
\node[below] at (3.2,-0.05) {$I_{\varphi_1}$};
\node[left] at (-0.05,3.4) {$I_{\varphi_2}$};
\node[above,plum] at (2,4) {mesh $2^{-m}$};
\draw[-{Stealth},thick,plum] (4.3,2.25) to[out=0,in=180] node[midway,above] {$f$} (6.2,2.25);
\draw[thick] (6.2,1.6) -- (11.6,1.6);
\draw[orchid,line width=4pt,opacity=0.6] (7.0,1.6) -- (10.8,1.6);
\draw[plum,line width=1.6pt] (7.6,1.6) -- (10.2,1.6);
\foreach \x in {7.6,10.2}{\draw[plum,thick] (\x,1.45) -- (\x,1.75);}
\foreach \x in {7.0,10.8}{\draw[orchid,thick] (\x,1.5) -- (\x,1.7);}
\node[plum,below] at (8.9,1.45) {$I$};
\node[orchid,below] at (7.0,1.45) {$\alpha-2^{-k}$};
\node[orchid,below] at (10.8,1.45) {$\beta+2^{-k}$};
\node[orchid,above] at (10.95,1.7) {$I^{(k)}$};
\fill[rose,opacity=0.8] (7.3,2.05) rectangle (9.2,2.25);
\node[rose,above,align=center] at (8.25,2.25) {image of\\ $(D_1\cap I_{\varphi_1})\times(D_2\cap I_{\varphi_2})$};
\draw[rose,densely dotted] (7.3,2.05) -- (7.3,1.6) (9.2,2.05) -- (9.2,1.6);
\end{tikzpicture}
\caption{The rule (Cn) for a connective of two variables.}
\label{fig:grid}
\end{figure}

\begin{definition}[The model existence game]\label{def:game}
The game $\mathrm{G}(T)$ is played by Olivia and Josh.
It starts at the empty position; at each round, Olivia makes a \emph{challenge}, and Josh \emph{answers} it by adding to the current position $\mathfrak p$ commitments and records, if any, in the manner that the rule of the challenge prescribes.
When Josh answers a challenge, the \emph{first fresh witness} of a sort $\mathbf S$ is the witness $w_{\mathbf S,i}$ with the least index $i$ that occurs neither in $\mathfrak p$ nor in the challenge.
In each rule other than (Ax) and (Tot), the commitment named by Olivia belongs to $\mathfrak p$, and its interval is written $I=[\alpha,\beta]$.
\begin{description}
\item[(Ax)] Olivia names an axiom $\psi\le r$ or $\psi\ge r$ of $T$.
Josh adds $\psi\in[\min(-\bar r_\psi,r),r]$ or $\psi\in[r,\max(r,\bar r_\psi)]$, respectively.
\item[(Tot)] Olivia names a function symbol $F$ and witnesses $\bar w$ of the appropriate sorts.
If $\mathfrak p$ contains no record $F(\bar w)\mapsto w'$, Josh adds such a record, choosing for $w'$ a witness of the sort of the values of $F$ that occurs in $\mathfrak p$ or is the first fresh witness of this sort; otherwise he adds nothing.
\item[(Tm)] Olivia names a commitment $R(\bar t)\in I$, where $R$ is a predicate symbol and $\bar t$ is a tuple of terms, together with a term $F(\bar w)$ that occurs in $\bar t$ and in which $\bar w$ is a tuple of witnesses (possibly empty).
Josh first answers as in (Tot) for $F$ and $\bar w$; then, with $w'$ the witness such that the enlarged position contains the record $F(\bar w)\mapsto w'$, he adds $R(\bar t')\in I$, where $\bar t'$ is obtained from $\bar t$ by replacing every occurrence of $F(\bar w)$ by $w'$.
\item[(Cn)] Olivia names a commitment $f(\varphi_1,\dots,\varphi_\ell)\in I$, where $f$ is a connective function, and a number $k\in\NN$.
With $m=m(f,\bar\varphi,k)$, Josh chooses intervals $D_1,\dots,D_\ell$ of $\cG_m$ such that each $D_j$ meets $I_{\varphi_j}$ and $f$ maps $(D_1\cap I_{\varphi_1})\times\dots\times(D_\ell\cap I_{\varphi_\ell})$ into $I^{(k)}$, and adds the commitments $\varphi_j\in D_j\cap I_{\varphi_j}$ (see \cref{fig:grid}).
\item[($\sup{\le}$)] Olivia names a commitment $\sup_x\varphi(x)\in I$, where $x$ is a variable of sort $\mathbf S$, and a witness $w\in\cW_{\mathbf S}$.
Josh adds $\varphi(w)\in[-\bar r_\varphi,\beta]$.
\item[($\sup{\ge}$)] Olivia names a commitment $\sup_x\varphi(x)\in I$, where $x$ is a variable of sort $\mathbf S$, and a number $k\in\NN$.
Josh chooses a witness $w\in\cW_{\mathbf S}$ that occurs in $\mathfrak p$ or is the first fresh witness of sort $\mathbf S$, and adds $\varphi(w)\in[\alpha-2^{-k},\bar r_\varphi]$.
\item[($\inf{\ge}$)] Olivia names a commitment $\inf_x\varphi(x)\in I$, where $x$ is a variable of sort $\mathbf S$, and a witness $w\in\cW_{\mathbf S}$.
Josh adds $\varphi(w)\in[\alpha,\bar r_\varphi]$.
\item[($\inf{\le}$)] Olivia names a commitment $\inf_x\varphi(x)\in I$, where $x$ is a variable of sort $\mathbf S$, and a number $k\in\NN$.
Josh chooses a witness $w$ as in ($\sup{\ge}$), and adds $\varphi(w)\in[-\bar r_\varphi,\beta+2^{-k}]$.
\item[($\bigvee{\le}$)] Olivia names a commitment $\bigvee_n\varphi_n\in I$, where $(\varphi_n)_{n\in\NN}$ is a family of sentences, and a number $n\in\NN$.
Josh adds $\varphi_n\in[-\bar r_{\bigvee_n\varphi_n},\beta]$.
\item[($\bigvee{\ge}$)] Olivia names a commitment $\bigvee_n\varphi_n\in I$ as in the preceding rule, and a number $k\in\NN$.
Josh chooses a number $n\in\NN$ and adds $\varphi_n\in[\alpha-2^{-k},\bar r_{\bigvee_n\varphi_n}]$.
\item[($\bigwedge{\ge}$)] Olivia names a commitment $\bigwedge_n\varphi_n\in I$, where $(\varphi_n)_{n\in\NN}$ is a family of sentences, and a number $n\in\NN$.
Josh adds $\varphi_n\in[\alpha,\bar r_{\bigwedge_n\varphi_n}]$.
\item[($\bigwedge{\le}$)] Olivia names a commitment $\bigwedge_n\varphi_n\in I$ as in the preceding rule, and a number $k\in\NN$.
Josh chooses a number $n\in\NN$ and adds $\varphi_n\in[-\bar r_{\bigwedge_n\varphi_n},\beta+2^{-k}]$.
\end{description}
An answer is \emph{legal} when it obeys the rule and every interval that it adds is nonempty, that is, has its left endpoint at most its right endpoint.
Josh \emph{loses} at a round if he has no legal answer to Olivia's challenge, or if his answer makes the position infeasible.
A \emph{play} is a finite or infinite sequence of rounds played according to these rules; it ends at the round, if any, at which Josh loses, a round that consists of Olivia's challenge alone when he has no legal answer to it.
Olivia wins a play if Josh loses at some round; Josh wins a play if he never loses, in which case the play is infinite.
A \emph{strategy} for a player chooses the player's next move from the moves made so far; it is \emph{winning} when the player wins every play in which the strategy is followed, whatever the opponent does.
\end{definition}

Josh's answers are of three kinds, according to the number of choices that the rules leave him.
The rules (Ax), ($\sup{\le}$), ($\inf{\ge}$), ($\bigvee{\le}$) and ($\bigwedge{\ge}$) leave him no choice: the answer is determined by the challenge.
Under (Tot), (Tm), (Cn), ($\sup{\ge}$) and ($\inf{\le}$) he has finitely many choices because only finitely many witnesses occur in a position and only finitely many intervals of $\cG_m$ meet a bounded interval.
Only the rules ($\bigvee{\ge}$) and ($\bigwedge{\le}$), which concern countable suprema and infima, leave him a choice among infinitely many answers.
We call the challenges made under these two rules \emph{countable challenges}, and Josh's answers to them \emph{countable answers}.

\subsection{The model existence theorem as a game}
Josh's survival yields a model only if every requirement is eventually tested: an axiom left unnamed or a commitment left unchallenged need not hold in the structure built from the play.
A challenge is \emph{available} at a round when the rules allow Olivia to make it there, equivalently when the commitment it names, if any, is already in the position.
Since positions only grow, an available challenge remains available.
The vocabulary $\cL_{\cW}$ and the theory $T$ are countable, and every sentence occurring in a play is atomic or is obtained from a subformula of an axiom by substituting closed terms for its free variables.
Thus only countably many challenges can occur, and at least one is always available because $T$ is nonempty.
A play is \emph{fair} if every challenge that becomes available is eventually made.
A strategy of Olivia is \emph{fair} if every infinite play in which she follows it is fair, whatever Josh does.
The following \emph{bookkeeping} gives a fair strategy.
She fixes a bijection $i\mapsto(j(i),l(i))$ from $\NN$ onto $\NN\times\NN$ with $j(i)\le i$ for every $i\in\NN$, and at each round she enumerates, with repetitions if necessary, the challenges then available.
At round $i$ she makes the $l(i)$-th challenge of her enumeration for round $j(i)$, which is still available because $j(i)\le i$.

\begin{keytheorem}[The model existence theorem as a game]\label{thm:game-existence}
Let $T$ be a nonempty countable theory of real-valued $L_{\omega_1\omega}$, in a countable vocabulary, whose axioms are closed conditions $\psi\le r$ or $\psi\ge r$ for sentences $\psi$ and rational numbers $r$, and let $\mathrm{G}(T)$ be the game of \cref{def:game}.
The following conditions are equivalent:
\begin{enumerate}
\item the theory $T$ has a model;
\item Josh has a winning strategy in $\mathrm{G}(T)$;
\item Olivia has no winning strategy in $\mathrm{G}(T)$.
\end{enumerate}
\end{keytheorem}

\begin{proof}
\emph{(1) implies (2).}
Let $\fN$ be a model of $T$; Josh answers every challenge truthfully in $\fN$.
To do so, he keeps along the play a map $h$ from the witnesses occurring in the current position to $\fN$.
The map $h$ makes every commitment of the position hold in $\fN$ when each witness $w$ is interpreted as $h(w)$.
Moreover, it satisfies the equation $F^\fN(h(\bar w))=h(w')$ for every record $F(\bar w)\mapsto w'$ of the position.
When a challenge names a witness that does not yet occur, Josh first extends $h$ to it arbitrarily, which is possible because the sorts of $\fN$ are nonempty.
The commitment that he adds under (Ax) holds because $\fN$ is a model of $T$.
Under (Tot) and (Tm), if the position contains no record $F(\bar w)\mapsto w'$, Josh adds one with $w'$ the first fresh witness and extends $h$ to $w'$ by letting $h(w')=F^\fN(h(\bar w))$.
Under (Tm), whether the record is new or not, the value of $R(\bar t')$ in $\fN$ equals that of $R(\bar t)$.
To a challenge (Cn), Josh answers with cells $D_j\in\cG_m$ that contain the values $\xi_j$ of the sentences $\varphi_j$ in $\fN$.
Every point of $(D_1\cap I_{\varphi_1})\times\dots\times(D_\ell\cap I_{\varphi_\ell})$ then has coordinates within $2^{-m}$ of those of $(\xi_1,\dots,\xi_\ell)$.
By the choice of $m$, the connective $f$ maps every such point to within $2^{-k}$ of $f(\xi_1,\dots,\xi_\ell)$, which lies in $I$.
The commitment added under ($\sup{\le}$) holds because the value of $\varphi(w)$ is at most the supremum, hence at most $\beta$.
Under ($\sup{\ge}$), Josh takes the first fresh witness $w$.
Since the supremum is at least $\alpha$, he lets $h(w)$ be an element of $\fN$ at which $\varphi$ takes a value at least $\alpha-2^{-k}$.
The remaining rules are handled in the same way; under ($\bigvee{\ge}$), for instance, Josh names a number $n$ for which the value of $\varphi_n$ is at least $\alpha-2^{-k}$.
Every answer is legal, since the true values lie in the added intervals.
Moreover, no position becomes infeasible, since the value in $\fN$ of a basic sentence $R(\bar w)$ lies in every interval committed for the sentence and in $\mathcal V_R$.
Since both of these facts hold at every round, Josh never loses.

\emph{(3) implies (1).}
We first show that every infinite fair play, that is, every fair play in which Josh never loses, yields a model of $T$.
Let $\mathfrak p_0\subseteq\mathfrak p_1\subseteq\cdots$ be the positions of such a play, indexed by $i\in\NN$, and let $\mathfrak p_\infty$ be their union.
We define a structure $\fM$ for $\cL_{\cW}$ whose elements are the witnesses: the universe of sort $\mathbf S$ is $\cW_{\mathbf S}$, and each witness is interpreted as itself.
For a function symbol $F$ and witnesses $\bar w$, Olivia makes the challenge (Tot) for $F$ and $\bar w$ at some round, because the play is fair.
The union $\mathfrak p_\infty$ therefore contains a record $F(\bar w)\mapsto w'$, which is unique because every $\mathfrak p_i$ is a position, and we define $F^\fM(\bar w)$ to be $w'$.
For a basic sentence $R(\bar w)$, the intervals $I$ for which $\mathfrak p_\infty$ contains $R(\bar w)\in I$ are compact, and so is $\mathcal V_R$.
Every finite subfamily of the family formed by these intervals and $\mathcal V_R$ has a common point, because no position $\mathfrak p_i$ is infeasible.
By the finite intersection property, the whole family has a common point, and we let $R^\fM(\bar w)$ be one.
As a point of $\mathcal V_R$, this value respects the a priori bound of $R$, and it equals $1$ when $R$ is $\mathbf 1$.
The logic requires of the interpretations $F^\fM$ and $R^\fM$ only that they have the appropriate sorts and, for $R^\fM$, that it respect the a priori bound of $R$ (see \cref{sec:logic}).
The interpretations chosen above therefore complete the definition of $\fM$.

We claim that every commitment $\varphi\in I$ of $\mathfrak p_\infty$ holds in $\fM$.
The claim is proved by induction on the construction of $\varphi$ and, for atomic sentences, on the number of occurrences of function symbols of the original vocabulary $\cL$.
For a basic sentence, the claim holds by the choice of $R^\fM$.
A non-basic atomic sentence $R(\bar t)$ contains a term $F(\bar w)$ whose arguments $\bar w$ are witnesses (possibly none, when $F$ is a constant symbol).
By fairness, Olivia made the challenge (Tm) for the commitment $R(\bar t)\in I$ and this term, and Josh added $R(\bar t')\in I$, which holds by the induction hypothesis.
The sentences $R(\bar t)$ and $R(\bar t')$ have the same value in $\fM$, because the witness $w'$ that replaces $F(\bar w)$ is $F^\fM(\bar w)$.
For $f(\varphi_1,\dots,\varphi_\ell)\in I$, fairness gives, for every $k\in\NN$, intervals $D_1,\dots,D_\ell$ such that $\mathfrak p_\infty$ contains the commitments $\varphi_j\in D_j\cap I_{\varphi_j}$ and the connective $f$ maps $(D_1\cap I_{\varphi_1})\times\dots\times(D_\ell\cap I_{\varphi_\ell})$ into $I^{(k)}$.
By the induction hypothesis, the value of $f(\varphi_1,\dots,\varphi_\ell)$ in $\fM$ lies in $I^{(k)}$ for every $k$, hence in $I$.
For $\sup_x\varphi(x)\in[\alpha,\beta]$, fairness adds to $\mathfrak p_\infty$ a commitment $\varphi(w)\in[-\bar r_\varphi,\beta]$ for every witness $w$ of the sort of $x$, that is, for every element of $\fM$ of this sort.
Fairness also adds to $\mathfrak p_\infty$, for every $k\in\NN$, a commitment $\varphi(w)\in[\alpha-2^{-k},\bar r_\varphi]$ for some witness $w$.
By the induction hypothesis, the supremum lies in $[\alpha,\beta]$.
The remaining cases of the induction are similar.
Finally, every axiom $\psi\le r$ of $T$ was named in a challenge (Ax), to which Josh answered with the commitment $\psi\in[\min(-\bar r_\psi,r),r]$.
By the claim, the structure $\fM$ satisfies $\psi\le r$; the axioms $\psi\ge r$ are treated in the same way.

Hence the reduct of $\fM$ to $\cL$ is a model of $T$.

Now suppose that $T$ has no model, and let Olivia follow a fair strategy, such as the bookkeeping described before the theorem.
By the preceding paragraphs, no play in which she follows it is infinite, since such a play would be fair and would yield a model of $T$.
Every such play therefore ends with Josh losing, and the fair strategy is a winning strategy for Olivia.
Thus, if Olivia has no winning strategy, the theory $T$ has a model.

\emph{(2) implies (3).}
The two players cannot both have winning strategies, since a play in which both follow their winning strategies would be won by both.

The game is thus determined.
Determinacy also follows from the theorem of Gale and Stewart~\cite{Gale-Stewart:1953}, since each of Olivia's wins is decided at a finite round and the set of her wins is therefore open in the space of plays; the argument above does not need it.
\end{proof}

\begin{remark}[Baire category and the Banach--Mazur game]\label{rem:banach-mazur}
Olivia's bookkeeping can be read as a Baire category argument in the space of plays, and as a winning strategy in a Banach--Mazur game on the same space.
We set out both readings in \cref{app:henkin}, beside the Baire category argument of Aubone, Dueñez and Iovino~\cite{Aubone-Duenez-Iovino:2026a} in the space of Henkin types.
\end{remark}

\subsection{Finite obstructions}
When a theory has no models, the model existence theorem supplies a winning strategy whose tree of plays organizes the resulting finite contradictions.
Let $\sigma$ be a strategy for Olivia.
The \emph{tree of $\sigma$}, denoted $\cP_\sigma$, is the set of finite plays in which Olivia follows $\sigma$ and Josh has not lost before the last round, ordered by extension.
A finite play of $\cP_\sigma$ in which Josh has not lost has as children its extensions by Olivia's next challenge and a legal answer of Josh.
When Josh has no legal answer to this challenge, the only child of the play is its extension by the challenge alone.
A \emph{leaf} of $\cP_\sigma$ is a finite play with no children, that is, one at which Josh has lost.
The strategy $\sigma$ is winning exactly when $\cP_\sigma$ has no infinite branch, that is, when $\cP_\sigma$ is \emph{well-founded}.
The \emph{rank} of a well-founded tree is defined by recursion: a leaf has rank $0$, and every other node has as rank the least ordinal exceeding the ranks of all its children.
The rank of the tree is the rank of its root.
The tree $\cP_\sigma$ is countably branching, hence countable; when it is well-founded, its rank is therefore a countable ordinal.
This rank is finite exactly when the plays in which Olivia follows $\sigma$ have bounded length.

To distinguish the countable choices from the intervening finite choices, we retain only Josh's answers to countable challenges.
In \cref{cor:finite-obstruction} below we show that only countable challenges can give a node of $\cP_\sigma$ infinitely many children.
The \emph{countable record} of a node of $\cP_\sigma$ is the finite sequence of the numbers $n$ that Josh names in his countable answers along the branch from the root to the node.
The numbers appear in this sequence in the order in which Josh names them.
The \emph{skeleton} of $\sigma$ is the set of countable records of nodes of $\cP_\sigma$, ordered by extension; it is a tree of finite sequences of natural numbers.

\begin{corollary}[Finite obstructions]\label{cor:finite-obstruction}
Let $T$ be a nonempty countable theory of real-valued $L_{\omega_1\omega}$, in a countable vocabulary, whose axioms are closed conditions $\psi\le r$ or $\psi\ge r$ for sentences $\psi$ and rational numbers $r$, and suppose that $T$ has no models.
Then Olivia has a winning strategy in $\mathrm{G}(T)$, and every winning strategy $\sigma$ of Olivia has the following properties:
\begin{enumerate}
\item The tree $\cP_\sigma$ is well-founded: every play in which Olivia follows $\sigma$ ends after finitely many rounds, at a leaf.
\item At every leaf, the last position $\mathfrak p$ is a finite set of commitments and records, involving finitely many axioms of $T$ and finitely many witnesses, that is satisfied in no structure under any interpretation of the witnesses.
Either $\mathfrak p$ is infeasible, in which case a basic sentence $R(\bar w)$ receives finitely many rational intervals with no common point in $\mathcal V_R$, or Josh has no legal answer to Olivia's last challenge, which names a commitment of $\mathfrak p$ that holds in no structure.
\item A node of $\cP_\sigma$ has infinitely many children only if $\sigma$ makes a countable challenge at it; every other node has finitely many children.
\item If every node of $\cP_\sigma$ has finitely many children, then the tree $\cP_\sigma$ is finite.
In that case there is a number $N_\sigma$ such that every play in which Olivia follows $\sigma$ ends within $N_\sigma$ rounds, and some finite set of axioms of $T$ has no models.
Every node of $\cP_\sigma$ has finitely many children when no axiom of $T$ contains a countable supremum or infimum.
\item The skeleton of $\sigma$ is well-founded, and for every element $\vec n$ of the skeleton there is a number $N_\sigma(\vec n)$ such that every play in which Olivia follows $\sigma$ and whose countable record is $\vec n$ ends within $N_\sigma(\vec n)$ rounds.
\end{enumerate}
\end{corollary}

\begin{proof}
By \cref{thm:game-existence}, Olivia has a winning strategy; let $\sigma$ be any winning strategy of Olivia.
Part (1) holds because $\sigma$ is winning.
For part (2), let $\mathfrak p$ be the last position of a leaf, a finite set of commitments and records like every position.
If Josh lost because his last answer made $\mathfrak p$ infeasible, then the position $\mathfrak p$ is satisfied in no structure, since the value of a basic sentence $R(\bar w)$ in a structure lies in $\mathcal V_R$.
Otherwise Josh had no legal answer to Olivia's last challenge, which therefore is not a challenge (Ax), (Tot) or (Tm), because under these rules a legal answer always exists.
The challenge names a commitment of $\mathfrak p$, and the argument for ``(1) implies (2)'' in the proof of \cref{thm:game-existence}, which uses the axioms of $T$ only under (Ax), produces a legal answer whenever this commitment holds in some structure under some interpretation of its witnesses.
The commitment therefore holds in no structure, and neither does $\mathfrak p$.
Part (3) follows from the count of Josh's possible answers in the paragraph after \cref{def:game}.

Parts (4) and (5) rest on K\H{o}nig's lemma: an infinite tree in which every node has finitely many children has an infinite branch.
For part (4), the tree $\cP_\sigma$ is well-founded and by hypothesis finitely branching; hence it is finite, and the length of its longest branch gives $N_\sigma$.
By part (3), the hypothesis holds when $\sigma$ makes no countable challenge.
In particular, the hypothesis holds when no axiom of $T$ contains a countable supremum or infimum.
Indeed, every sentence that occurs in a play is atomic or is an instance of a subformula of an axiom.
It remains to be shown that some finite set of axioms of $T$ has no models.
Let $T_\sigma$ be the set of axioms named in challenges (Ax) at nodes of $\cP_\sigma$, together with one axiom of $T$ if there are none.
In $\mathrm{G}(T_\sigma)$, every challenge that $\sigma$ makes at a node of $\cP_\sigma$ is allowed by the rules, and Josh has the same answers as in $\mathrm{G}(T)$.
The strategy $\sigma$ therefore wins $\mathrm{G}(T_\sigma)$, and the finite theory $T_\sigma$ has no models by \cref{thm:game-existence}.

For part (5), let $\vec n$ be a finite sequence of natural numbers, and let $\cP_\sigma^{\vec n}$ be the set of nodes of $\cP_\sigma$ whose countable record is an initial segment of $\vec n$.
The set $\cP_\sigma^{\vec n}$ is a subtree of $\cP_\sigma$, and every play in which Olivia follows $\sigma$ and whose countable record is $\vec n$ runs inside this subtree.
In $\cP_\sigma^{\vec n}$, a node at which $\sigma$ makes a countable challenge has at most one child.
If Josh has no legal answer to the challenge, its only child is the extension by the challenge alone.
Otherwise, if the countable record of the node is $\vec n$ itself, the node has no child in $\cP_\sigma^{\vec n}$.
In the remaining case, the only child of the node is the extension by the countable answer whose number is the next entry of $\vec n$.
By part (3), every other node of $\cP_\sigma^{\vec n}$ has finitely many children.
Since $\cP_\sigma^{\vec n}$ is well-founded, it is finite by K\H{o}nig's lemma, and the length of its longest branch gives $N_\sigma(\vec n)$.
Finally, suppose, towards a contradiction, that the skeleton has an infinite branch $\vec n_1\subsetneq \vec n_2\subsetneq\cdots$.
The set of nodes of $\cP_\sigma$ whose countable record is an initial segment of some $\vec n_i$ then forms a finitely branching subtree of $\cP_\sigma$.
This subtree contains nodes of every countable record $\vec n_i$ and is therefore infinite.
By K\H{o}nig's lemma, it has an infinite branch, contradicting part (1); the skeleton is therefore well-founded.
\end{proof}

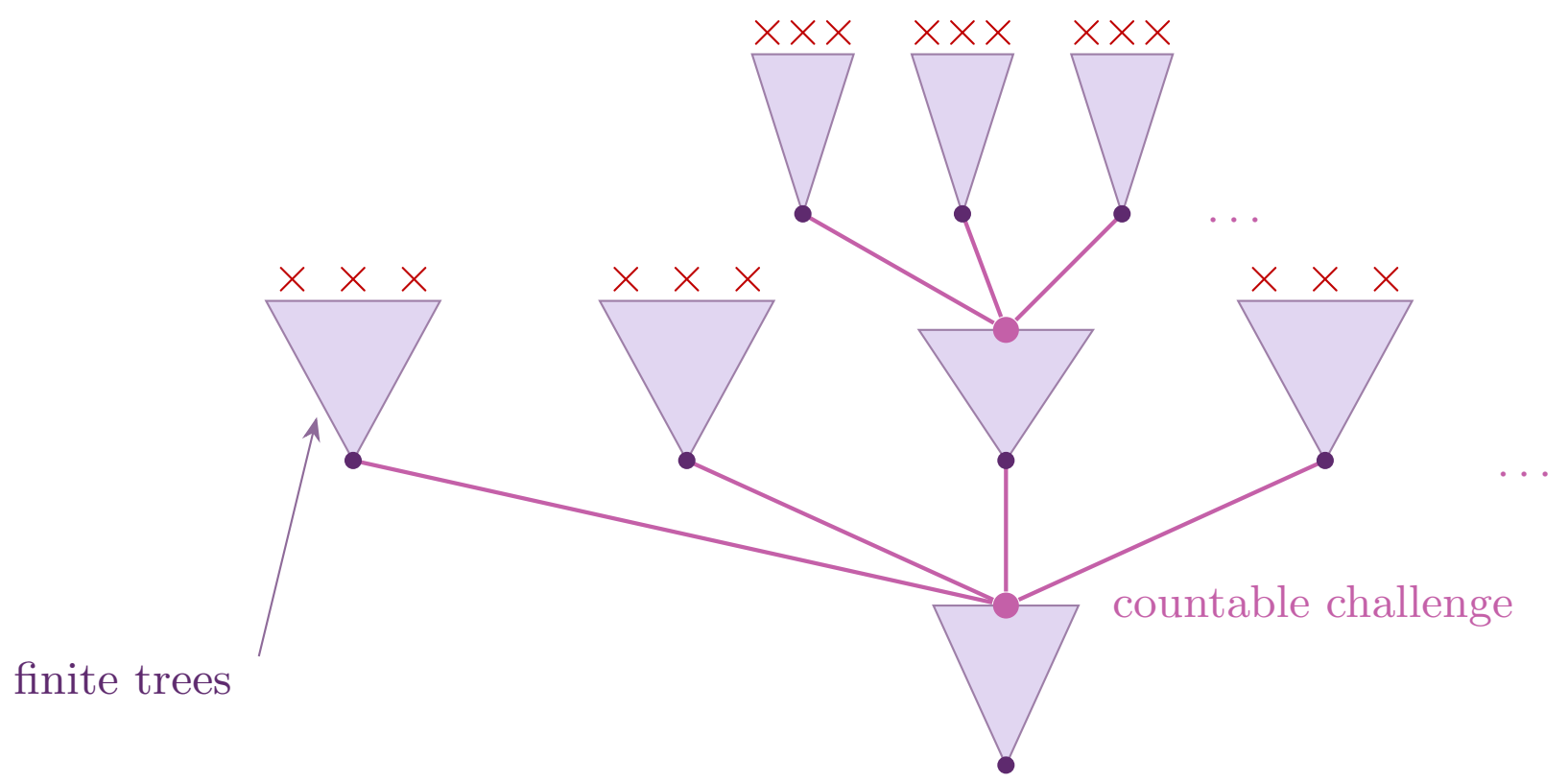
\begin{figure}[ht]
\centering
\begin{tikzpicture}[font=\footnotesize,yscale=-1,
  nd/.style={circle,fill=plum,inner sep=1.2pt},
  cn/.style={circle,fill=rose,inner sep=1.8pt},
  lf/.style={red!75!black,font=\small\bfseries}]
\fill[lavender!55] (0,0) -- (-0.5,-1.1) -- (0.5,-1.1) -- cycle;
\draw[plum!60] (0,0) -- (-0.5,-1.1) -- (0.5,-1.1) -- cycle;
\node[nd] (r) at (0,0) {};
\node[cn] (c) at (0,-1.1) {};
\foreach \x in {-4.5,-2.2,0,2.2}{\draw[rose,thick] (c) -- (\x,-2.1);}
\node[rose] at (3.6,-2.0) {$\cdots$};
\node[rose,right] at (0.6,-1.1) {countable challenge};
\finitesubtree{-4.5}{-2.1}{0.6}
\finitesubtree{-2.2}{-2.1}{0.6}
\finitesubtree{2.2}{-2.1}{0.6}
\node[nd] at (-4.5,-2.1) {};\node[nd] at (-2.2,-2.1) {};\node[nd] at (2.2,-2.1) {};
\fill[lavender!55] (0,-2.1) -- (-0.6,-3.0) -- (0.6,-3.0) -- cycle;
\draw[plum!60] (0,-2.1) -- (-0.6,-3.0) -- (0.6,-3.0) -- cycle;
\node[nd] at (0,-2.1) {};
\node[cn] (c2) at (0,-3.0) {};
\foreach \x in {-1.4,-0.3,0.8}{\draw[rose,thick] (c2) -- (\x,-3.8);}
\node[rose] at (1.6,-3.75) {$\cdots$};
\foreach \x in {-1.4,-0.3,0.8}{\finitesubtree{\x}{-3.8}{0.35}\node[nd] at (\x,-3.8) {};}
\node[plum,left] at (-5.2,-0.6) {finite trees};
\draw[-{Stealth},plum!70] (-5.15,-0.75) -- (-4.75,-2.4);
\end{tikzpicture}
\caption{The tree of a winning strategy for Olivia.}
\label{fig:skeleton}
\end{figure}

Part (4) of \cref{cor:finite-obstruction} recovers compactness for countable theories of closed conditions without countable suprema and infima: if such a theory has no models, some finite set of its axioms already has no models.
For the infinitary logic, part (3) identifies the additional choices that may prevent this conclusion, and part (5) bounds the lengths of plays with a fixed countable record (see \cref{fig:skeleton}).
Parts (1) and (2) supply a well-founded tree of finite obstructions, one at each leaf, even when no finite set of axioms is inconsistent.
The extension of the theory of \cref{thm:a} by the requirement that $\mu$ take only the values $0$ and $1$ gives an example.
This is the theory $T^{01}$ of the Introduction; in \cref{prop:bks-game} below we show that every finite set of its axioms has a model, although the whole theory has no models by \cref{thm:c}.

\begin{proof}[Proof of \cref{thm:d}]
In the model existence game, Olivia is the player who claims that the theory has no models, and the tree of plays of a strategy $\sigma$ of hers is the tree $\cP_\sigma$.
Let $T'$ be as in \cref{thm:d}, with its axioms written in the normal form of Subsection~\ref{sec:game-setting}.
By \cref{thm:game-existence}, the theory $T'$ has no models if and only if Olivia has a winning strategy in $\mathrm{G}(T')$.
If $T'$ has no models, the remaining assertions of \cref{thm:d} are those of \cref{cor:finite-obstruction}, with $T'$ in place of $T$.
At a leaf, the commitments of the last position are satisfied in no structure even when its records are disregarded, as \cref{thm:d} asserts: by part~(2) of the corollary, either they assign to a basic sentence $R(\bar w)$ intervals with no common point in $\mathcal V_R$, or one of them holds in no structure.
\end{proof}

\begin{remark}[Strict conditions]\label{rem:strict}
The conditions of the logic are strict inequalities, and a closed condition is a countable theory of strict ones (see \cref{sec:logic}).
For a theory that contains strict conditions, the game has one more rule.
\begin{description}
\item[(Str)] Olivia names a condition $\varphi<r$ of $T$, for a sentence $\varphi$.
Josh chooses a rational number $\varrho<r$ and adds $\varphi\in[-\bar r_\varphi,\varrho]$.
For a condition $\varphi>r$ of $T$, Josh chooses a rational number $\varrho>r$ and adds $\varphi\in[\varrho,\bar r_\varphi]$.
\end{description}
With one more case for this rule, the proof of \cref{thm:game-existence} applies to theories that consist of closed conditions as above and of strict conditions on sentences.
In the additional case, a truthful Josh chooses $\varrho$ between the value of $\varphi$ and $r$, and in the structure built from a fair play the commitment $\varphi\in[-\bar r_\varphi,\varrho]$ gives $\varphi\le \varrho<r$.
Under (Str), however, Josh chooses among infinitely many answers, and the challenges (Str) are accordingly counted as countable challenges in \cref{cor:finite-obstruction}.
In the countable record, Josh's rational answers are coded by natural numbers; parts (1), (3) and (5) of the corollary then hold.
Part (2) holds once the condition named in the last challenge is added to $\mathfrak p$.
Indeed, Josh has no legal answer to (Str) for $\varphi<r$ exactly when $r\le-\bar r_\varphi$, and none for $\varphi>r$ exactly when $r\ge \bar r_\varphi$.
In either case, the condition named is itself unsatisfiable.
With the same counting, part (4) holds, except for its clause on theories without countable suprema and infima.
The exception arises because a challenge (Str) with a legal answer has infinitely many children and because such theories may have no models although each finite set of their axioms has one.
For a nullary predicate symbol $R$ with $r_R>0$, the theory $\{R>0\}\cup\{R<\tfrac1n:n\ge1\}$ is such a theory.
Olivia wins by naming first $R>0$, to which Josh answers with some rational $\varrho>0$ and the commitment $R\in[\varrho,\bar r_R]$, and then a condition $R<\tfrac1n$ with $\tfrac1n<\varrho$, to which every legal answer of Josh adds an interval disjoint from $[\varrho,\bar r_R]$ and thereby makes the position infeasible.
The tree of this strategy branches infinitely at its root.
The contrast between the two kinds of conditions has a topological counterpart in the space of Henkin types of Aubone, Dueñez and Iovino~\cite{Aubone-Duenez-Iovino:2026a}, which we describe in \cref{app:henkin}.
\end{remark}

\historypara{Historical note.}
The Banach--Mazur game began as Problem~43 of the Scottish Book, posed by Mazur.
On 4 August 1935 Banach wrote in the book that Mazur's conjecture is true.
His proof was never published; Oxtoby published a proof in 1957~\cite{Mauldin:1981,Oxtoby:1957}.
Gale and Stewart proved in 1953 that open and closed games are determined~\cite{Gale-Stewart:1953}.
The analysis of finite obstructions rests, in addition, on the lemma on finitely branching trees that Dénes K\H{o}nig had published in 1927~\cite{Konig:1927}.
The witnesses of the game descend from the new constants of Henkin's completeness proof of 1949~\cite{Henkin:1949}.
For Väänänen there are essentially three kinds of games in logic, the semantic game, the model existence game and the Ehrenfeucht--Fraïssé game; his book of 2011 is built on them~\cite{Vaananen:2011}.

\section{The vocabulary}\label{sec:vocabulary}

We now return to the theory $T$ of quantum probability of \cref{thm:a} and present its vocabulary, in which the three sorts of \cref{sec:setting} are linked by projection, coordinate and span maps.

\enlargethispage{2pt}

The vocabulary $\cL$ has exactly three sorts, the Hilbert sort, the Solèr sort and the sequence sort.
The Hilbert and sequence sorts are vector sorts, intended to be the closed unit balls of complex Hilbert spaces; this restriction keeps their norm and inner-product predicates bounded.
In the intended interpretation, the elements of the Solèr sort name closed subspaces and the span function assigns a lattice element to a sequence of vectors.
The list below specifies the symbols and their individual a priori bounds, while the axioms govern their interaction and standardness remains a separate hypothesis.
\begin{itemize}
\item The \emph{Hilbert sort} $\sH$ has a binary function $(x,y)\mapsto\lambda x+\lambda'y$ for each pair $\lambda,\lambda'$ of elements of the field $\QI$ of complex numbers with rational real and imaginary parts such that $|\lambda|+|\lambda'|\le1$.
The sort $\sH$ also has the predicates $\operatorname{Re}\ip{x}{y}$, $\operatorname{Im}\ip{x}{y}$ and $\|x\|$, with a priori bound $1$, and constants $e_n$ for $n\in\NN$.
We write $0$ for the closed term $0e_0+0e_0$, and $\ip{x}{y}$ for the pair of predicates $\operatorname{Re}\ip{x}{y}$ and $\operatorname{Im}\ip{x}{y}$.
We also write $|\ip{x}{y}|^2$ for the formula $(\operatorname{Re}\ip{x}{y})^2+(\operatorname{Im}\ip{x}{y})^2$, built with the squares of \cref{lem:squares}.
\item The \emph{sequence sort} $\sK$ has the same functions and predicates as $\sH$, one constant $\varepsilon$, and unary functions $\mathrm{res}_N$ for $N\in\NN$; in this sort, we write $0$ for the closed term $0\varepsilon+0\varepsilon$.
\item The \emph{Solèr sort} $\sL$ has the constants $0$ and $1$, the binary functions $\wedge$ and $\vee$ and the unary function ${}^\perp$.
\item Function symbols connecting the sorts:
\[
\begin{gathered}
\pi\colon\sL\times\sH\to\sH,\quad \at\colon\sH\to\sL,\quad \rep\colon\sL\to\sH,\\
\iota_n\colon\sH\to\sK,\quad \kappa_n\colon\sK\to\sH\ (n\in\NN),\quad \spa\colon\sK\to\sL.
\end{gathered}
\]
\item A predicate $\mu$ on $\sL$, the \emph{measure}, with a priori bound $2$.
The bound is wider than $1$ so that the inequality $\mu(a)\le1$, like $\mu(a)\ge0$, is imposed by the axiom (M1) of \cref{sec:axioms} below, and not by the vocabulary.
\item \emph{Truncations.} For every predicate symbol $R$ listed above, every positive integer $n$ and every rational number $q$, two predicate symbols $R^{+}_{n,q}$ and $R^{-}_{n,q}$ with the arguments of $R$ and with a priori bound $1$.
\end{itemize}
The intended values of the truncations are
\[
R^{+}_{n,q}(\bar x)=\min\bigl(1,\max(0,n(R(\bar x)-q))\bigr),\qquad R^{-}_{n,q}(\bar x)=\min\bigl(1,\max(0,n(q-R(\bar x)))\bigr),
\]
which the axioms (Tr) of \cref{sec:axioms} impose.
For an atomic formula $R(\bar t)$ whose predicate symbol is one of the inner-product, norm or measure predicates listed above, and a rational number $q$, let
\[
\lceil R(\bar t)>q\rceil=\bigvee_{n\ge1}R^{+}_{n,q}(\bar t),\qquad\lceil R(\bar t)<q\rceil=\bigvee_{n\ge1}R^{-}_{n,q}(\bar t),
\]
countable joins of members with a priori bound $1$.
In every structure that satisfies (Tr), the formula $\lceil R(\bar t)>q\rceil$ takes the value $1$ where $R(\bar t)>q$ and the value $0$ elsewhere, and likewise for $\lceil R(\bar t)<q\rceil$.
Several symbols now have more than one meaning: the symbol $0$ names a sentence, a closed term of each vector sort and a constant of the Solèr sort, and the symbol $1$ names a sentence and a constant of the Solèr sort.
The symbols $\wedge$ and $\vee$, for their part, name both connectives and functions of the Solèr sort, and the arguments to which they are applied always show which is meant.
For vector terms $t_1,t_2$ of the same sort, the equation $t_1\doteqq t_2$ abbreviates the axiom $\|\tfrac12t_1-\tfrac12t_2\|\le0$, a closed condition when $t_1$ and $t_2$ are closed.
A vector $x$ \emph{belongs} to an element $a$ of the Solèr sort, or is a \emph{member} of $a$, when $\pi(a,x)\doteqq x$.

In the intended interpretation, the connecting symbols have the following roles.
Here $a$ is a lattice element, $x$ a vector, $s$ a sequence, and $n,N\in\NN$.
\begin{center}
\small
\begin{tabular}{@{}p{2.4cm}p{9.1cm}@{}}
\toprule
Symbol & Intended meaning\\
\midrule
$\pi(a,x)$ & Orthogonal projection of $x$ onto the subspace named by $a$.\\[3pt]
$\at(x)$ & The ray spanned by $x$ if $x\ne0$; the zero element otherwise.\\[3pt]
$\rep(a)$ & A chosen unit vector in a nonzero subspace $a$; the zero vector for $a=0$.\\[3pt]
$\iota_nx$ & The sequence with $x$ in position $n$ and zeros elsewhere.\\[3pt]
$\kappa_ns$ & The entry of $s$ in position $n$.\\[3pt]
$\mathrm{res}_Ns$ & The sequence obtained by keeping the entries with indices $n<N$ and replacing the others by zero.\\[3pt]
$\spa(s)$ & The lattice element naming the closed span of the entries of $s$.\\
\bottomrule
\end{tabular}
\end{center}
The constants $(e_n)_{n\in\NN}$ name an orthonormal basis of the Hilbert space.
The sequence sort consists of the square-summable sequences of vectors of the Hilbert space with $\ell^2$ norm at most $1$, and $\varepsilon$ names $(2^{-n-1}e_n)_{n\in\NN}$.
The named sequence $\varepsilon$, together with the axioms (K4) and (M5) of \cref{sec:axioms}, is used in the proof of \cref{thm:c} to locate a basis atom of measure $1$ when $\mu$ takes only the values $0$ and $1$.
The parity argument of that proof then takes place in the span of a block of four consecutive basis vectors, chosen to contain this atom.

The language has no primitive equality symbol.
For vectors, we already use zero norm of a scaled difference.
For lattice elements, the representative vector supplies a test for zero: in the intended interpretation, $\|\rep(a)\|$ is $0$ for the zero subspace and $1$ otherwise.
To test whether $a\le b$, apply this test to $a\wedge(a\wedge b)^\perp$, the orthogonal complement of $a\wedge b$ in $a$.
This element is zero exactly when $a\le b$.
We therefore define
\[
[a\ne0]=2\|\rep(a)\|-1,\qquad[a\le b]=1-2\|\rep(a\wedge(a\wedge b)^\perp)\|,\qquad[a=b]=[a\le b]\wedge[b\le a].
\]
In the intended interpretation, these are two-valued tests for nonzero subspace, inclusion and equality of subspaces.
In \cref{lem:envelope} below we associate a closed subspace of a Hilbert space with each element of the Solèr sort in every model of $T_0$.
In \cref{lem:lattice-truth} we show that the three formulas are two-valued in every such model: their value $1$ means, respectively, that the associated subspace is nonzero, that the subspace associated with $a$ is contained in the one associated with $b$, and that the two subspaces coincide.
The equality test does not assert literal identity of the elements of the Solèr sort.
For terms $t_1,t_2$ of the Solèr sort, the equation $t_1=t_2$ abbreviates the axiom $[t_1=t_2]\ge1$.

The axiom of \cref{thm:b} excludes the points of a Borel set from the range of the measure.
This axiom needs a two-valued formula that says whether a value of the measure lies in the set.
For a Borel set $B\subseteq\RR$, we write $\beta_B(a)$ for any formula of real-valued $L_{\omega_1\omega}$, with the single free variable $a$ of the Solèr sort, whose values are prescribed as follows.
At every element $a$ of every structure for $\cL$ that satisfies the axioms (Tr), the formula takes the value $1$ if $\mu(a)\in B$ and $-1$ otherwise.

\begin{lemma}\label{lem:borel}
For every Borel set $B\subseteq\RR$, there is a formula $\beta_B(a)$ as just described.
\end{lemma}

\begin{proof}
Let $\mathcal{S}$ be the class of sets $B\subseteq\RR$ for which such a formula exists with a priori bound at most $3$.
For rational numbers $q\le q'$, the formula
\[
1-2\bigl(\lceil\mu(a)>q'\rceil\vee\lceil\mu(a)<q\rceil\bigr)
\]
takes the value $1$ when $\mu(a)\in[q,q']$ and $-1$ otherwise, in every structure that satisfies (Tr).
Its a priori bound is $3$, and the class $\mathcal{S}$ therefore contains the closed intervals with rational endpoints.
If $B\in\mathcal{S}$, then the negation $-\beta_B(a)$, which has the same a priori bound, is a formula for $\RR\setminus B$.
If $B_n\in\mathcal{S}$ for every $n\in\NN$, then the countable join of the formulas $\beta_{B_n}(a)$, whose a priori bounds are at most $3$, is a formula for $\bigcup_nB_n$ with a priori bound at most $3$.
The class $\mathcal{S}$ is therefore a $\sigma$-algebra, and since every open subset of $\RR$ is a countable union of closed intervals with rational endpoints, it contains every Borel set.
\end{proof}

When a formula $\beta_B(a)$ is obtained by following, as in the proof of \cref{lem:borel}, a construction of $B$ from closed intervals with rational endpoints by complements and countable unions, its negations and countable joins are nested as in the construction.
The depth of this nesting is a countable ordinal $\alpha$, and an induction on $\alpha$ shows that the set $B$ then belongs to the class $\boldsymbol\Delta^0_{\alpha+2}$ of the Borel hierarchy.
Since this hierarchy is proper on $\RR$~\cite[Section~22]{Kechris:1995}, the least depth of such a construction of $B$ is unbounded below $\omega_1$ as $B$ varies.
The axioms of \cref{thm:b}, with formulas $\beta_B(a)$ built in this way, therefore involve formulas whose joins are nested to arbitrarily large countable depth, although the formulas $\beta_B(a)$ contain no quantifier.

\begin{definition}\label{def:standard}
Let $\fM$ be a structure for $\cL$, or for a vocabulary obtained from $\cL$ by removing the predicate $\mu$, the sequence sort with the symbols that involve it, or both, the truncations of the removed predicates being removed with them.
The structure $\fM$ is \emph{standard} when its Hilbert sort is the closed unit ball of a complex Hilbert space $H$, when its sequence sort, if it has one, is the closed unit ball of a complex Hilbert space $K$, and when the functions $\lambda x+\lambda'y$ and the predicates $\operatorname{Re}\ip{x}{y}$, $\operatorname{Im}\ip{x}{y}$ and $\|x\|$ are the restrictions of the operations, the inner products and the norms of these spaces.
\end{definition}

The theories introduced below satisfy $T_0\subseteq T_1\subseteq T$: the theory $T_0$ uses the Hilbert and Solèr sorts, $T_1$ adds sequences and their spans, and $T$ adds the measure.
Standardness concerns the two vector sorts and their Hilbert-space operations.
It does not by itself identify the Solèr sort with the full lattice of closed subspaces, or the sequence sort with all square-summable sequences of vectors of the Hilbert sort; these identifications are the results of \cref{sec:t0,sec:t1}.
Standardness is a hypothesis on the model, not a consequence of the axioms: every standard model of $T_0$ has an uncountable Hilbert sort (see \cref{prop:sep} below), whereas the downward Löwenheim--Skolem theorem gives $T_0$, $T_1$ and $T$ models whose sorts are countable (see \cref{rem:nonstandard}).
Accordingly, \cref{thm:c} concerns every model of $T$, whereas part (1) of \cref{thm:a} and \cref{thm:b} concern the standard models.

\section{The axioms}\label{sec:axioms}

The axioms build the intended structure in three stages: $T_0$ links the Hilbert and Solèr sorts by projection symbols, $T_1$ adds represented spans of sequences, and $T$ adds the measure.
In standard models, these stages yield a pre-Solèr lattice, the full lattice of closed subspaces, and a measure on the full lattice, in each case taking the quotient by the formula expressing equality.

The axiom groups follow these tasks.
Groups (H) and (OL) specify the vector and lattice operations; (P) and (Sub) identify the projections and their lattice operations; (At) links vectors with rays and chooses representatives of nonzero subspaces.
Groups (K) and (Sp) supply sequences and their spans, while (M) governs the measure.
Group (Tr) fixes the truncations used for exact threshold tests.
The formal definitions of $T_0$, $T_1$ and $T$ follow the list.

The variables $x,y,z$ range over $\sH$, the variables $s,t$ over $\sK$, and the variables $a,b,c$ over $\sL$; every axiom is a closed condition required of all values of its free variables, as in \cref{sec:logic}.
By contrast, the letters $\lambda,\lambda'$ and $n,m,k,N$ are metavariables, not variables of the logic: the pair $\lambda,\lambda'$ ranges over $\QI$ with $|\lambda|+|\lambda'|\le1$, and the indices $n,m,k,N$ range over $\NN$.
A free occurrence of one of these metavariables gives an axiom schema, one axiom for each permitted choice, whereas bound occurrences index the indicated sums, joins or infima.

Vector equations written with $\doteqq$ and lattice equations written with $=$ retain the meanings defined in \cref{sec:vocabulary}.
For two-valued formulas $\theta,\theta'$, the axiom $-\theta\vee\theta'\ge1$ expresses $\theta\Rightarrow\theta'$.
For a two-valued formula $\theta$ and any formula $\psi$, the axiom $\theta\wedge\psi\le0$ expresses that $\theta$ implies $\psi\le0$.
Italic readings describe the geometry of the standard models.
For $N\in\NN$, write $b_N=\at(e_0)\vee\dots\vee\at(e_{N-1})$, with $b_0=0$.

\begin{axiomgroup}{(H) The Hilbert sort}
\begin{enumerate}[label=(H\arabic*),leftmargin=*]
\item \emph{Inner product space.} For any two terms built from variables by the functions $\lambda x+\lambda'y$ that denote the same $\QI$-linear combination of the variables, the equation between them; the equations $\ip{\lambda x+\lambda'y}{z}=\lambda\ip{x}{z}+\lambda'\ip{y}{z}$ and $\ip{y}{x}=\overline{\ip{x}{y}}$, each written as two real equations; and $\|x\|\ge0$ and $\|x\|^2=\operatorname{Re}\ip{x}{x}$.
\item \emph{Separability.} $\ip{e_n}{e_m}=\delta_{nm}$, where $\delta_{nm}$ is the Kronecker delta, and
\[
\bigwedge_N\bigl\|\tfrac12x-\tfrac12\pi(b_N,x)\bigr\|\le0.
\]
\emph{The constants $e_n$ are orthonormal, and every vector is the limit of its projections onto the spans of $e_0,\dots,e_{N-1}$.}
\end{enumerate}
\end{axiomgroup}

\begin{axiomgroup}{(OL) The Solèr sort}
\begin{enumerate}[label=(OL),leftmargin=*]
\item The equations of ortholattices:
\begin{gather*}
a\wedge b=b\wedge a,\quad a\wedge(b\wedge c)=(a\wedge b)\wedge c,\quad a\wedge(a\vee b)=a,\\
a\vee b=b\vee a,\quad a\vee(b\vee c)=(a\vee b)\vee c,\quad a\vee(a\wedge b)=a,\\
a\wedge0=0,\quad a\vee1=1,\quad a^{\perp\perp}=a,\quad (a\wedge b)^\perp=a^\perp\vee b^\perp,\quad a\wedge a^\perp=0.
\end{gather*}
\end{enumerate}
\end{axiomgroup}

\begin{axiomgroup}{(P) Each element of the Solèr sort determines a projection}
\begin{enumerate}[label=(P\arabic*),leftmargin=*]
\item $\pi(a,\lambda x+\lambda'y)\doteqq\lambda\pi(a,x)+\lambda'\pi(a,y)$.
\item $\pi(a,\pi(a,x))\doteqq\pi(a,x)$.
\item $\ip{\pi(a,x)}{y}=\ip{x}{\pi(a,y)}$.
\end{enumerate}
\end{axiomgroup}

The projection axioms (P) concern each lattice element separately, while group (Sub) relates the projection maps to the lattice operations.
For a closed subspace $W$ of a Hilbert space, let $P_W$ denote the orthogonal projection onto $W$.
No noncommutative polynomial with fixed coefficients yields $P_{U\cap V}$ from $P_U$ and $P_V$ for every pair of closed subspaces $U$ and $V$, since such polynomials are continuous in the norm topology and the map $(P_U,P_V)\mapsto P_{U\cap V}$ is not.
The projection $P_{U\cap V}$ is the strong limit of the powers $(P_UP_V)^n$, by a theorem of von Neumann~\cite[Lemma~22]{vonNeumann:1949} (see also~\cite{vonNeumann:1950}).

\begin{theorem}[von Neumann's alternating projection theorem]\label{thm:von-neumann}
Let $U$ and $V$ be closed subspaces of a Hilbert space, and let $P_U$, $P_V$ and $P_{U\cap V}$ be the orthogonal projections onto $U$, $V$ and $U\cap V$.
Then, for every vector $x$ of the space, the sequence $\bigl((P_UP_V)^nx\bigr)_{n\in\NN}$ converges in norm to $P_{U\cap V}x$.
\end{theorem}

Axiom (Sub1) below requires the projection of $a\wedge b$ to be the limit in \cref{thm:von-neumann}.
Axiom (Sub2) requires the projection of $a\vee b$ to be the identity minus the same limit for $a^\perp$ and $b^\perp$, since the closure of $U+V$ is $(U^\perp\cap V^\perp)^\perp$.
In \cref{lem:lattice-truth} the theorem then identifies the subspace represented by $a\wedge b$ with the intersection of the subspaces represented by $a$ and $b$, and the subspace represented by $a\vee b$ with the closure of their sum.
The convergence admits no rate that is uniform in $U$ and $V$.
For $n\ge1$ the norm of $(P_UP_V)^n-P_{U\cap V}$ is $c^{2n-1}$, where $c$ is the cosine of the Friedrichs angle between $U$ and $V$~\cite{Kayalar-Weinert:1988}.
When $U+V$ is not closed the convergence is arbitrarily slow~\cite{Bauschke-Deutsch-Hundal:2009}.
Accordingly, axioms (Sub1) and (Sub2) assert convergence through the limit superior $\bigwedge_N\bigvee_{n\ge N}$, which vanishes exactly when the nonnegative errors tend to zero.

\begin{axiomgroup}{(Sub) The operations of the Solèr sort are those of subspaces}
\begin{enumerate}[label=(Sub\arabic*),leftmargin=*,start=0]
\item $\pi(0,x)\doteqq0$ and $\pi(1,x)\doteqq x$. \\
\emph{Only the zero vector belongs to $0$; every vector belongs to $1$.}
\item $\displaystyle\bigwedge_N\bigvee_{n\ge N}\bigl\|\tfrac12\pi(a\wedge b,x)-\tfrac12(P_aP_b)^nx\bigr\|\le0$, where $(P_aP_b)^nx$ is the term obtained by applying $\pi(b,\cdot)$ and then $\pi(a,\cdot)$, $n$ times. \\
\emph{The members of $a\wedge b$ are the common members of $a$ and $b$.}
\item $\displaystyle\bigwedge_N\bigvee_{n\ge N}\bigl\|\tfrac14\pi(a\vee b,x)-\tfrac14x+\tfrac14(P_{a^\perp}P_{b^\perp})^nx\bigr\|\le0$, where $(P_{a^\perp}P_{b^\perp})^nx$ is the term of (Sub1) for $a^\perp$ and $b^\perp$, and the vector inside the norm is denoted by a term built from the functions $\lambda x+\lambda'y$. \\
\emph{The closed subspace spanned by the members of $a\vee b$ is the closure of the sum of the closed subspaces spanned by the members of $a$ and by those of $b$.}
\item $\tfrac12\pi(a^\perp,x)\doteqq\tfrac12x-\tfrac12\pi(a,x)$. \\
\emph{The members of $a^\perp$ are the vectors orthogonal to all members of $a$.}
\end{enumerate}
\end{axiomgroup}

\begin{axiomgroup}{(At) Atoms and unit vectors}
\begin{enumerate}[label=(At\arabic*),leftmargin=*]
\item $\pi(\at(x),x)\doteqq x$, and $\|\pi(\at(x),y)\|^2\,\|x\|^2=|\ip{y}{x}|^2$. \\
\emph{Every vector $x$ belongs to $\at(x)$, and for $x\ne0$ the projection onto $\at(x)$ is the projection onto the line spanned by $x$.}
\item $\pi(a,\rep(a))\doteqq\rep(a)$; $\|\rep(a)\|\wedge\bigl(1-\|\rep(a)\|\bigr)\le0$; and $\|\pi(a,x)\|\le\|\rep(a)\|$. \\
\emph{The vector $\rep(a)$ belongs to $a$; its norm is $0$ or $1$; and it is $1$ unless $a$ has no nonzero member.}
\item $\|\rep(\at(x))\|\le\lceil\|x\|>0\rceil$. \\
\emph{If $\|x\|=0$, then $\at(x)$ has no nonzero member.}
\end{enumerate}
\end{axiomgroup}

\begin{axiomgroup}{(K) The sequence sort}
\begin{enumerate}[label=(K\arabic*),leftmargin=*]
\item The axioms (H1), written for $\sK$.
\item $\iota_n(\lambda x+\lambda'y)\doteqq\lambda\iota_nx+\lambda'\iota_ny$ and $\kappa_n(\lambda s+\lambda't)\doteqq\lambda\kappa_ns+\lambda'\kappa_nt$; $\ip{\iota_nx}{\iota_my}=\delta_{nm}\ip{x}{y}$; and $\ip{\kappa_ns}{x}=\ip{s}{\iota_nx}$.
\item $\displaystyle\ip{\mathrm{res}_Ns}{t}=\sum_{n<N}\ip{\kappa_ns}{\kappa_nt}$, written as two real equations, and \mbox{$\displaystyle\bigwedge_N\bigl\|\tfrac12s-\tfrac12\mathrm{res}_Ns\bigr\|\le0$}. \\
\emph{The sequence $\mathrm{res}_Ns$ keeps the first $N$ entries of $s$, and these sequences converge to $s$ as $N$ grows.}
\item $\kappa_n\varepsilon\doteqq2^{-n-1}e_n$.
\end{enumerate}
\end{axiomgroup}

\newpage
\begin{axiomgroup}{(Sp) Spans}
\begin{enumerate}[label=(Sp\arabic*),leftmargin=*]
\item $\pi(\spa(s),\kappa_ns)\doteqq\kappa_ns$. \\
\emph{Every entry of $s$ belongs to $\spa(s)$.}
\item $\displaystyle\Bigl(-\bigwedge_n[\at(\kappa_ns)\le b]\Bigr)\vee[\spa(s)\le b]\ge1$. \\
\emph{If every entry of $s$ belongs to $b$, then the span $\spa(s)$ lies below $b$.}
\end{enumerate}
\end{axiomgroup}

\begin{axiomgroup}{(Tr) Truncations}
\begin{enumerate}[label=(Tr),leftmargin=*]
\item For every predicate symbol $R$ of $\cL$ other than $\mathbf 1$ and the truncation symbols, every positive integer $n$ and every rational number $q$:
\[
R^{+}_{n,q}(\bar x)=\bigl((n\cdot(R(\bar x)-q))\wedge1\bigr)\vee0,\qquad R^{-}_{n,q}(\bar x)=\bigl((n\cdot(q-R(\bar x)))\wedge1\bigr)\vee0.
\]
\end{enumerate}
\end{axiomgroup}

We require the predicate $\mu$ to take values in $[0,1]$, to have value $1$ at the greatest element of the Solèr sort, and to agree on elements related by the equality formula $[a=b]=1$.
Axioms (M1)--(M3) impose these requirements, while (M4) imposes orthogonal additivity and (M5) makes the measure continuous along the finite joins that exhaust the represented span of a sequence.
In (M4), orthogonality is built into $b\wedge a^\perp$; for any $c$ orthogonal to $a$, the equality $[c=c\wedge a^\perp]=1$, together with (M3) and compatibility of equality with the lattice operations (see \cref{cor:quotient} below), gives the usual additive identity.
A term built from the lattice operations involves only finitely many elements of the Solèr sort, and the function $\spa$ supplies the joins of infinite sequences of atoms.
In the standard models every element of the Solèr sort is equal, in the sense of $[a=b]=1$, to the span $\spa(s)$ of a sequence $s$, and axioms (M5) and (M4) then yield $\sigma$-additivity with the help of (M1) and (M3) (see the proof of \cref{thm:a}).
For $N\in\NN$, let $c_N(s)$ be the term $\at(\kappa_0s)\vee\dots\vee\at(\kappa_{N-1}s)$, with $c_0(s)=0$.

\begin{axiomgroup}{(M) The measure}
\begin{enumerate}[label=(M\arabic*),leftmargin=*]
\item $\mu(a)\ge0$ and $\mu(a)\le1$.
\item $\mu(1)=1$.
\item $[a=b]\wedge|\mu(a)-\mu(b)|\le0$. \\
\emph{Equal elements have equal measures.}
\item $\mu\bigl(a\vee(b\wedge a^\perp)\bigr)=\mu(a)+\mu(b\wedge a^\perp)$. \\
\emph{The measure is additive on orthogonal elements.}
\item $\displaystyle\mu(\spa(s))=\bigvee_N\mu(c_N(s))$. \\
\emph{The measure of the span of a sequence is the supremum, over $N$, of the measures of the joins of the atoms of its first $N$ entries.}
\end{enumerate}
\end{axiomgroup}

\begin{definition}
The theory $T_0$ consists of the axioms (H), (OL), (P), (Sub), (At) and the axioms (Tr) for the predicates of the Hilbert sort, in the vocabulary obtained from $\cL$ by removing the sequence sort with its constant $\varepsilon$, the functions $\iota_n$, $\kappa_n$, $\mathrm{res}_N$ and $\spa$, and the predicate $\mu$, together with the truncations of the removed predicates.
The theory $T_1$ consists of $T_0$ together with (K), (Sp) and the axioms (Tr) for the predicates of the sequence sort, in the vocabulary $\cL$ without $\mu$ and its truncations; the theory $T$ consists of $T_1$ together with (M) and the axioms (Tr) for $\mu$.
\end{definition}

We now recover the geometry described by these axioms in every model of $T_0$, without assuming standardness.
First we construct a Hilbert space in which the Hilbert-sort operations are represented and the projection symbols act as orthogonal projections.
Von Neumann's alternating projection theorem (\cref{thm:von-neumann}) then shows that the projection of $a\wedge b$ has as range the intersection of the ranges of the projections of $a$ and $b$, and that the projection of $a\vee b$ has as range the closure of their sum; equality of the represented subspaces will make the quotient of the Solèr sort an ortholattice.

\begin{lemma}\label{lem:envelope}
Let $\fM$ be a model of $T_0$.
There are a complex Hilbert space $\bar E$, a subset $Y$ of the closed unit ball of $\bar E$ whose $\QI$-linear span $E$ is dense in $\bar E$, and a map $x\mapsto\hat x$ from the Hilbert sort of $\fM$ onto $Y$ with the following properties.
The map sends $\lambda x+\lambda'y$ to $\lambda\hat x+\lambda'\hat y$, and the predicates $\ip{x}{y}$ and $\|x\|$ to $\ip{\hat x}{\hat y}$ and $\|\hat x\|$.
For every element $a$ of the Solèr sort, the map $\hat x\mapsto\widehat{\pi(a,x)}$ is well defined and extends uniquely to an orthogonal projection $\bar P_a$ of $\bar E$.
If $\fM$ is standard, then the space $\bar E$ is the Hilbert space whose closed unit ball is the Hilbert sort, and the map $x\mapsto\hat x$ is the identity.
\end{lemma}

\begin{proof}
The idea is to recover a vector space from the Hilbert sort by formal multiples.
By (H1), the form given by the pair of predicates $\ip{x}{y}$ is hermitian and satisfies $\operatorname{Re}\ip{x}{x}=\|x\|^2\ge0$; the Cauchy--Schwarz inequality therefore holds for it.
Each pair $(k,x)$, with $k\ge1$ an integer and $x$ an element of the Hilbert sort, stands for $kx$.
We define the sum of two pairs by $(k,x)+(l,y)=(2kl,\tfrac1{2l}x+\tfrac1{2k}y)$ and, for an element $\lambda$ of $\QI$, the multiple $\lambda(k,x)=(km,\tfrac{\lambda}{m}x)$, where $m$ is a positive integer with $m\ge|\lambda|$; the form is given by $\ip{(k,x)}{(l,y)}=kl\ip{x}{y}$.
We identify two pairs when the form vanishes on their difference.
Axiom (H1) makes the quotient a $\QI$-vector space $E$ with a positive definite hermitian form, and we let $\hat x$ be the class of $(1,x)$; since the class of $(k,x)$ is $k\hat x$, the space $E$ is the $\QI$-linear span of the set $Y$ of the vectors $\hat x$.
The completion $\bar E$ of $E$ is a complex Hilbert space, because multiplication by the elements of $\QI$, which are dense in $\CC$, extends by continuity to multiplication by complex numbers.
For an element $a$ of the Solèr sort, consider the assignment $Q$ that sends a pair $(k,x)$ to $(k,\pi(a,x))$; by (P1)--(P3), it respects sums and multiples, it is idempotent and it is symmetric with respect to the form, in each case up to pairs on which the form vanishes.
For every pair $v$, symmetry and idempotence therefore give $\ip{Qv}{Qv}=\ip{v}{Qv}$, which is at most $\ip{v}{v}^{1/2}\ip{Qv}{Qv}^{1/2}$ by the Cauchy--Schwarz inequality, whence $\ip{Qv}{Qv}\le\ip{v}{v}$.
In particular, the assignment $Q$ sends pairs on which the form vanishes to such pairs, and it therefore induces on $E$ a $\QI$-linear, idempotent, symmetric and contractive map $Q_a$, which extends by continuity to a self-adjoint idempotent $\bar P_a$ of $\bar E$.
In a standard model, the Hilbert sort is the closed unit ball of a Hilbert space $H$, and the form is positive definite on it.
In that case, the map sending the class of $(k,x)$ to $kx$ identifies $E$ with $H$, since every vector of $H$ is a positive integer multiple of a vector of the ball.
As $H$ is complete, it follows that $\bar E=H$ and that $x\mapsto\hat x$ is the identity of the ball.
\end{proof}

For each element $a$ of the Solèr sort, write $S_a$ for the range of the projection $\bar P_a$ just constructed; in a standard model, $S_a$ is the closed subspace spanned by the members of $a$, as in the Introduction.
We must show that the lattice operations act on these ranges as intended and that the formulas of \cref{sec:vocabulary} detect their nonzeroness, inclusion and equality.
The next lemma proves this in every model of $T_0$, allowing distinct elements of the sort to name the same subspace.

\begin{lemma}\label{lem:lattice-truth}
Let $\fM$ be a model of $T_0$, and let $\bar E$ and $\bar P_a$ be as in \cref{lem:envelope}.
\begin{enumerate}[label=(\arabic*)]
\item We have $\bar P_0=0$, $\bar P_1=\mathrm{Id}$ and $\bar P_{a^\perp}=\mathrm{Id}-\bar P_a$; moreover, $S_{a\wedge b}=S_a\cap S_b$ and $S_{a\vee b}=\overline{S_a+S_b}$.
\item The norm $\|\rep(a)\|$ is $1$ if $S_a\ne\{0\}$ and $0$ otherwise; hence the formula $[a\ne0]$ is two-valued, with the value $1$ exactly when $S_a\ne\{0\}$.
\item The formula $[a\le b]$ is two-valued, with the value $1$ exactly when $S_a\subseteq S_b$, and the formula $[a=b]$ is two-valued, with the value $1$ exactly when $S_a=S_b$.
\item For every $x$ with $\|x\|>0$, the subspace $S_{\at(x)}$ is $\CC\hat x$, and for every $x$ with $\|x\|=0$ it is $\{0\}$.
\end{enumerate}
\end{lemma}

\begin{proof}
(1) The first three equalities follow from (Sub0) and (Sub3) by continuity, since $E$ is dense in $\bar E$.
By (Sub1), the projection $\bar P_{a\wedge b}$ agrees on the vectors $\hat x$, and hence on $\bar E$, with the strong limit of the operators $(\bar P_a\bar P_b)^n$, which is the projection onto $S_a\cap S_b$ by \cref{thm:von-neumann}.
Likewise, axiom (Sub2) and the equality $\bar P_{c^\perp}=\mathrm{Id}-\bar P_c$ show that $\bar P_{a\vee b}=\mathrm{Id}-P_{S_a^\perp\cap S_b^\perp}$, whose range is $(S_a^\perp\cap S_b^\perp)^\perp=\overline{S_a+S_b}$.

(2) By (At2), together with the axiom $\|x\|\ge0$ of (H1), the norm $\|\rep(a)\|$ is $0$ or $1$; by (At2) again, the vector $\widehat{\rep(a)}$ lies in $S_a$, and every $x$ satisfies $\|\bar P_a\hat x\|\le\|\rep(a)\|$.
If $S_a\ne\{0\}$, then the projection $\bar P_a$ does not vanish on the dense subspace $E$, and hence not on some $\hat x$, whence $\|\rep(a)\|>0$ and so $\|\rep(a)\|=1$.
If $S_a=\{0\}$, then the membership $\widehat{\rep(a)}\in S_a$ gives $\|\rep(a)\|=0$.

(3) By (1), the subspace $S_{a\wedge(a\wedge b)^\perp}$ is $S_a\cap(S_a\cap S_b)^\perp$, the orthogonal complement of $S_a\cap S_b$ in $S_a$.
It is $\{0\}$ exactly when $S_a\cap S_b=S_a$, that is, when $S_a\subseteq S_b$.
By (2), the formula $[a\le b]$ is therefore two-valued, with the value $1$ exactly when $S_a\subseteq S_b$, and the assertion about $[a=b]$ follows.

(4) Let $P=\bar P_{\at(x)}$ and suppose that $\|x\|>0$, whence $\hat x\ne0$.
By (At1), the vector $\hat x$ lies in $S_{\at(x)}$, and every $y$ satisfies $\|P\hat y\|^2\|\hat x\|^2=|\ip{\hat y}{\hat x}|^2$.
Let $P_x$ be the projection onto $\CC\hat x$, which satisfies $\|P_x\hat y\|^2\|\hat x\|^2=|\ip{\hat y}{\hat x}|^2$ as well; since the vector $\hat x$ lies in the range of $P$, we have $P_x\le P$.
The identity $\|P\hat y\|^2=\|P_x\hat y\|^2+\|(P-P_x)\hat y\|^2$ then gives $(P-P_x)\hat y=0$ for every $y$, whence $P=P_x$ by density.
Finally, if $\|x\|=0$, then (At3) and the axioms (Tr) for the norm give $\|\rep(\at(x))\|=0$, whence $S_{\at(x)}=\{0\}$ by (2).
\end{proof}

By part (3) of the lemma, the formula $[a=b]$ detects exactly when two elements of the Solèr sort name the same subspace, and the quotient of the sort by this relation therefore inherits the geometry of the subspaces.

\begin{corollary}\label{cor:quotient}
Let $\fM$ be a model of $T_0$, and let $\bar E$ and $S_a$ be as in \cref{lem:lattice-truth}.
The relation $[a=b]=1$ is an equivalence relation on the Solèr sort that is respected by $\wedge$, $\vee$ and ${}^\perp$.
The quotient $\sL^\fM/{=}$ by this relation is an ortholattice, and the map that sends the class of $a$ to $S_a$ is an embedding of ortholattices into the lattice $L(\bar E)$ of closed subspaces of $\bar E$ that carries the order of the quotient to inclusion.
\end{corollary}

\begin{proof}
Part (3) of \cref{lem:lattice-truth} shows that the relation $[a=b]=1$ holds exactly when $S_a=S_b$, and the relation is therefore an equivalence relation.
Since, by part (1) of the lemma, the subspaces $S_{a\wedge c}$, $S_{a\vee c}$ and $S_{a^\perp}$ are determined by $S_a$ and $S_c$, the relation is respected by the operations.
The same part (1) shows that the map sending the class of $a$ to $S_a$ preserves the operations; this map is injective by part (3), and its image is therefore a sub-ortholattice of $L(\bar E)$ isomorphic to the quotient.
Finally, the order of the quotient, in which $a\le b$ means $a\wedge b=a$, corresponds to inclusion, again by part (3).
\end{proof}

In \cref{lem:lattice-truth,cor:quotient}, the axioms (OL) are not used; they hold in every structure that satisfies the other axioms of $T_0$, by the embedding into $L(\bar E)$.
We keep them, so that the ortholattice structure of the Solèr sort has axioms of its own (see \cref{rem:redundancy}).

\historypara{Historical note.}
In the continuous logic set out by Ben Yaacov, Berenstein, Henson and Usvyatsov in 2008, the unit ball of a Banach space carries the maps $(x,y)\mapsto\alpha x+\beta y$ with $|\alpha|+|\beta|\le1$ as functions.
A Hilbert space also carries its inner product as a predicate~\cite{BenYaacov-Berenstein-Henson-Usvyatsov:2008}.
We use the affine operations and inner-product predicates in the real-valued vocabulary of \cref{sec:vocabulary}; their interpretations are governed by the axioms above, without metric or continuity requirements in the logic.
Von Neumann's alternating projection theorem, on which axioms (Sub1) and (Sub2) rest, comes from lecture notes that he distributed in mimeographed form in 1933.
The notes were printed in 1950 as the second volume of his \emph{Functional Operators}~\cite{vonNeumann:1950}.

\section{Standard models of \texorpdfstring{$T_0$}{T0}: pre-Solèr lattices}\label{sec:t0}

In a standard model of $T_0$, whose vocabulary has neither the sequence sort nor the measure, the Hilbert sort is the closed unit ball of a complex Hilbert space $H$.
Since the vocabulary has no equality symbol, distinct elements of the Solèr sort may represent the same closed subspace of $H$.
By \cref{cor:quotient}, the quotient of this sort by the relation $[a=b]=1$ of \cref{sec:vocabulary} is an ortholattice, and sending the class of each element $a$ to the closed subspace spanned by the members of $a$ embeds the quotient into the lattice $L(H)$ of closed subspaces of $H$.
We show that the quotient is a pre-Solèr lattice: it satisfies the conditions for a Solèr lattice except, possibly, completeness.
The image of the embedding is join-dense and meet-dense in the sense of \cref{sec:lattices}, and the lattice $L(H)$ is therefore isomorphic to the MacNeille completion of the quotient.

The quotient need not be complete, as the standard model of \cref{ex:finite-cofinite} shows.
Nevertheless, its $\sigma$-additive states, defined below, are given by density operators (see \cref{thm:gleason-t0}).
In \cref{sec:t1}, the sequence sort supplies the remaining subspaces in standard models of $T_1$.

Fix a standard model $\fM$ of $T_0$, and let $H$ be the Hilbert space whose closed unit ball is its Hilbert sort (see \cref{def:standard}).
In this case the Hilbert space $\bar E$ constructed in \cref{lem:envelope} is $H$ itself.
Write $\sL^\fM$ for the Solèr sort.
For each $a\in\sL^\fM$, let $P_a$ be the projection $\bar P_a$ extending $\pi(a,\cdot)$, and let $S_a$ be its range.
Thus $S_a$ is the closed subspace spanned by the members of $a$, the vectors $x$ of the unit ball with $\pi(a,x)=x$.
Write $\tilde a$ for the class of $a$ in $\sL^\fM/{=}$; by \cref{lem:lattice-truth}, two such classes coincide exactly when their subspaces coincide.

\begin{proposition}\label{prop:sep}
The space $H$ is a separable infinite-dimensional complex Hilbert space, and the constants $e_n$ of the Hilbert sort, for $n\in\NN$, form an orthonormal basis of $H$.
\end{proposition}

\begin{proof}
By (H2), the vectors $e_n$ are orthonormal.
By \cref{lem:lattice-truth}, the projection $P_{b_N}$ is the projection onto the span of $e_0,\dots,e_{N-1}$, and (H2) places every vector of the unit ball, and hence every vector of $H$, in the closed span of the vectors $e_n$.
Thus the sequence $(e_n)_{n\in\NN}$ is an orthonormal basis of $H$, which is therefore separable and infinite-dimensional; it is complex by \cref{def:standard}.
\end{proof}

\begin{proposition}\label{prop:embedding}
The map $\tau\colon\tilde a\mapsto S_a$ is an embedding of ortholattices from $\sL^\fM/{=}$ into the lattice $L(H)$ of closed subspaces of $H$: it is injective and satisfies
\[
S_0=\{0\},\quad S_1=H,\quad S_{a\wedge b}=S_a\cap S_b,\quad S_{a\vee b}=\overline{S_a+S_b},\quad S_{a^\perp}=S_a^\perp,
\]
and each of the relations $S_a\subseteq S_b$ and $[a\le b]=1$ holds exactly when $\tilde a\le\tilde b$.
\end{proposition}

\begin{proof}
The proposition restates \cref{cor:quotient} for the standard model $\fM$, in which $\bar E=H$; part~(1) of \cref{lem:lattice-truth} gives the displayed equalities, and part~(3) the description of the order.
\end{proof}

\begin{proposition}\label{prop:atoms}
For every nonzero $x$ in the unit ball, the subspace $S_{\at(x)}$ is $\CC x$.
The atoms of $\sL^\fM/{=}$ are exactly the classes $\widetilde{\at(x)}$ with $x\ne0$, and the map $x\mapsto\widetilde{\at(x)}$ induces a bijection from the unit sphere of $H$, modulo multiplication by complex numbers of modulus one, onto the set of atoms of $\sL^\fM/{=}$, whose inverse sends an atom $\tilde a$ to the class of $\rep(a)$.
\end{proposition}

\begin{proof}
The first assertion is part~(4) of \cref{lem:lattice-truth}, since in a standard model every nonzero vector $x$ has positive norm.
Since the ray $\CC x$ is an atom of $L(H)$ and the map $\tau$ is an order embedding, the class $\widetilde{\at(x)}$ is an atom of $\sL^\fM/{=}$.
Let $\tilde a$ be an atom of $\sL^\fM/{=}$.
By part~(2) of \cref{lem:lattice-truth}, the vector $\rep(a)$ is a unit vector, and it belongs to $a$ by (At2).
Hence the ray $S_{\at(\rep(a))}=\CC\rep(a)$ lies in $S_a$, and the class $\widetilde{\at(\rep(a))}$ is a nonzero element below $\tilde a$; since $\tilde a$ is an atom, we have $\tilde a=\widetilde{\at(\rep(a))}$.
For unit vectors $x$ and $y$, the classes $\widetilde{\at(x)}$ and $\widetilde{\at(y)}$ coincide exactly when $\CC x=\CC y$, that is, when $y$ is a multiple of $x$ by a complex number of modulus one.
In particular, for a unit vector $x$, the unit vector $\rep(\at(x))$ is such a multiple of $x$, because it lies in $S_{\at(x)}=\CC x$.
The inverse of the bijection is therefore as stated.
\end{proof}

\begin{theorem}\label{thm:t0}
Let $\fM$ be a standard model of $T_0$, let $H$ be the Hilbert space whose closed unit ball is its Hilbert sort, and, for each element $a$ of its Solèr sort $\sL^\fM$, let $S_a$ be the closed subspace of $H$ spanned by the members of $a$.
Then the quotient $\sL^\fM/{=}$ is a pre-Solèr lattice, and the image of the map $\tau\colon\tilde a\mapsto S_a$ is join-dense and meet-dense in the lattice $L(H)$ of closed subspaces of $H$, which is therefore isomorphic to the MacNeille completion of $\sL^\fM/{=}$.
The image of $\tau$ contains every finite-dimensional subspace of $H$ and the orthogonal complement of each, and the map $\tau$ sends every join that exists in $\sL^\fM/{=}$ to the corresponding join in $L(H)$.
\end{theorem}

\begin{proof}
The atoms already represent all rays, so their finite joins give the finite-dimensional subspaces.
We use these subspaces to verify that $\sL^\fM/{=}$ is a pre-Solèr lattice.
The zero subspace is $S_0$; now let $F$ be a nonzero finite-dimensional subspace of $H$ with a basis $x_1,\dots,x_k$ of unit vectors.
By \cref{prop:embedding,prop:atoms}, the map $\tau$ sends $(\at(x_1)\vee\dots\vee\at(x_k))^\sim$ to the closure of $F$, which is $F$ itself, and $((\at(x_1)\vee\dots\vee\at(x_k))^\perp)^\sim$ to $F^\perp$.

Condition (OM) holds because $\tau(\sL^\fM/{=})$ is a sub-ortholattice of the orthomodular lattice $L(H)$.
For (A), let $\tilde a$ be a nonzero element of $\sL^\fM/{=}$; the atoms below $\tilde a$ are the classes $\widetilde{\at(x)}$ with $x\in S_a$ and $0<\|x\|\le1$.
Any upper bound $\tilde b$ of these atoms satisfies $S_b\supseteq S_a$, hence $\tilde b\ge\tilde a$, which makes $\tilde a$ their join.

For (CP), let $\tilde a$ be an element of $\sL^\fM/{=}$, and let $\widetilde{\at(x)}$ be an atom not below $\tilde a$; then the subspace $S_a$ does not contain $x$.
By \cref{prop:embedding,lem:plus-one}, the subspace $S_{a\vee\at(x)}$ is $S_a+\CC x$, in which $S_a$ has codimension one.
Since $\tau$ is an order embedding, no element lies strictly between $\tilde a$ and $(a\vee\at(x))^\sim$.
Condition (I) holds because the join of two distinct atoms corresponds to a plane, which contains infinitely many rays, all in $\tau(\sL^\fM/{=})$.
The atoms $\widetilde{\at(e_n)}$, for $n\in\NN$, give (D).

For (S), take $p_n=\widetilde{\at(e_n)}$ and $c_n=\widetilde{\at(\tfrac12(e_n+e_{n+1}))}$; since $\tau$ is injective and preserves the lattice operations and the orthocomplementation, it suffices to verify the condition for the images of these atoms.
These images are the rays of $e_n$ and of $e_n+e_{n+1}$; by~\eqref{eq:harmonic}, the harmonic conjugate of $\CC(e_n+e_{n+1})$ with respect to $\CC e_n$ and $\CC e_{n+1}$ is $\CC(e_n-e_{n+1})$, the orthogonal complement of $\CC(e_n+e_{n+1})$ in their plane, because $\|e_n\|=\|e_{n+1}\|$.

Each closed subspace $S$ of $H$ is the join in $L(H)$ of its rays, and it equals $\bigcap\{y^\perp:y\in S^\perp,\ 0<\|y\|\le1\}$, where each $y^\perp=S_{\at(y)^\perp}$ lies in the image of $\tau$.
These two facts make the image of $\tau$ join-dense and meet-dense in $L(H)$.
By \cref{prop:dm-char}, the lattice $L(H)$ is therefore isomorphic to the MacNeille completion of $\sL^\fM/{=}$.

Finally, let $A\subseteq\sL^\fM/{=}$ have a join $\tilde c$ in $\sL^\fM/{=}$, and let $W$ be the closed span of the subspaces $S_a$, for $\tilde a\in A$.
Since every element of $A$ lies below $\tilde c$, the subspace $W$ is contained in $S_c$.
If $y$ is a vector of $S_c\cap W^\perp$ with $0<\|y\|\le1$, then every element of $A$ lies below $(\at(y)^\perp)^\sim$, and hence so does $\tilde c$.
The vector $y$, which lies in $S_c$, is then orthogonal to itself, which is absurd.
Since $S_c=W+(S_c\cap W^\perp)$ by the projection theorem, it follows that $S_c=W$.
\end{proof}

\begin{example}\label{ex:finite-cofinite}
Let $L_{\mathrm{fc}}$ (for ``finite dimension or finite codimension'') be the set of finite-dimensional subspaces of the space $\ell^2$ of square-summable complex sequences and their orthogonal complements.
For finite-dimensional subspaces $F$ and $G$ of $\ell^2$, the joins $F\vee G=F+G$, $F^\perp\vee G^\perp=(F\cap G)^\perp$ and $F\vee G^\perp=(F^\perp\cap G)^\perp$ in $L(\ell^2)$ belong to $L_{\mathrm{fc}}$, and the De~Morgan law places the meets in $L_{\mathrm{fc}}$ as well; the set $L_{\mathrm{fc}}$ is therefore a sub-ortholattice of $L(\ell^2)$.
Take the closed unit ball of $\ell^2$ as Hilbert sort and the standard unit vectors as the constants $e_n$, and take $L_{\mathrm{fc}}$, with its lattice operations and orthocomplementation, as Solèr sort.
Define $\pi(S,x)=P_Sx$ and $\at(x)=\CC x$, let $\rep(S)$ be a unit vector chosen in $S$ for each $S\ne\{0\}$, and let $\rep(\{0\})$ be $0$.
Each truncation $R^{\pm}_{n,q}$ of a predicate of the Hilbert sort is interpreted by its intended value $\min\bigl(1,\max(0,\pm n(R-q))\bigr)$.
Axioms (Sub1) and (Sub2) hold by \cref{thm:von-neumann}, and the other axioms of $T_0$ by construction; the result is therefore a standard model of $T_0$, in which $[S=S']=1$ exactly when $S=S'$.
Its Solèr sort $L_{\mathrm{fc}}$, which is its own quotient, is not complete.
Indeed, an upper bound in $L_{\mathrm{fc}}$ of the rays $\CC e_{2n}$, for $n\in\NN$, is infinite-dimensional and hence of the form $G^\perp$ for a finite-dimensional subspace $G$ orthogonal to every $e_{2n}$.
The closed span of the vectors $e_{2n+1}$ is infinite-dimensional, whereas $G$ is finite-dimensional; this span therefore contains a unit vector $y$ orthogonal to $G$.
The subspace $(G+\CC y)^\perp$, which contains every $e_{2n}$ but not $y$, is an upper bound of the same rays that lies strictly below $G^\perp$.
The rays $\CC e_{2n}$ therefore have no join in $L_{\mathrm{fc}}$.
\end{example}

We now show that Gleason's representation extends to the quotient of the Solèr sort in every standard model of $T_0$.
Its $\sigma$-additive states correspond to density operators on the model's Hilbert space $H$, which represent quantum states.
Because the quotient may be incomplete, countable additivity is required only for orthogonal families whose join exists there.

A \emph{density operator} on $H$ is a positive operator $\rho$ whose trace $\Tr\rho=\sum_k\ip{\rho u_k}{u_k}$, computed in any orthonormal basis $(u_k)_{k\in\NN}$ of $H$, equals $1$.
A \emph{$\sigma$-additive state} on an ortholattice $\Lambda$ is a map $\nu\colon\Lambda\to[0,1]$ with $\nu(1)=1$ such that $\nu(\bigvee_{n\in J}a_n)=\sum_{n\in J}\nu(a_n)$ for every family $(a_n)_{n\in J}$ of mutually orthogonal elements of $\Lambda$, indexed by a set $J\subseteq\NN$, whose join exists in $\Lambda$.

\begin{theorem}[Gleason's theorem for $T_0$]\label{thm:gleason-t0}
Let $\fM$ be a standard model of $T_0$, let $H$ be the Hilbert space whose closed unit ball is its Hilbert sort, and, for each element $a$ of its Solèr sort $\sL^\fM$, let $P_a$ be the orthogonal projection onto the closed subspace spanned by the members of $a$.
The $\sigma$-additive states on $\sL^\fM/{=}$ are exactly the maps $\tilde a\mapsto\Tr(\rho P_a)$, where $\rho$ ranges over the density operators on $H$, and each such state determines its density operator.
\end{theorem}

\begin{proof}
A density operator gives a $\sigma$-additive state by taking traces of the projections; conversely, the values of a $\sigma$-additive state on rays form a frame function, to which Gleason's theorem applies.
Let $\rho$ be a density operator on $H$.
The map $\tilde a\mapsto\Tr(\rho P_a)$ is well defined, since $P_a$ depends only on $\tilde a$.
Suppose that the elements $\tilde a_n$, for $n$ in a set $J\subseteq\NN$, are mutually orthogonal and have a join $\tilde c$ in $\sL^\fM/{=}$.
By \cref{prop:embedding,thm:t0}, the subspace $S_c$ is then the closure of the sum of the mutually orthogonal subspaces $S_{a_n}$, whence $P_cx=\sum_{n\in J}P_{a_n}x$ for every $x\in H$.
We compute the traces in an orthonormal basis of $H$ assembled from orthonormal bases of the subspaces $S_{a_n}$ and of $S_c^\perp$; since the terms $\ip{\rho u}{u}$ are nonnegative, they may be grouped by $n$, and we obtain $\Tr(\rho P_c)=\sum_{n\in J}\Tr(\rho P_{a_n})$.
Moreover, the map $\tilde a\mapsto\Tr(\rho P_a)$ takes the value $\Tr\rho=1$ at $\tilde1$, and its values lie in $[0,1]$.

Conversely, let $\nu$ be a $\sigma$-additive state on $\sL^\fM/{=}$.
Additivity on the orthogonal pair $\tilde1,\tilde0$, whose join is $\tilde1$, gives $\nu(\tilde0)=0$.
For a unit vector $x$, define $f(x)=\nu(\widetilde{\at(x)})$.
For $S_a\ne\{0\}$, let $\mathcal U_a$ be an orthonormal basis of $S_a$, which is countable because $H$ is separable.
The atoms $\widetilde{\at(u)}$, for $u\in\mathcal U_a$, are mutually orthogonal, and their join in $\sL^\fM/{=}$ is $\tilde a$.
Indeed, the element $\tilde a$ is an upper bound of these atoms, and every upper bound $\tilde b$ of them satisfies $S_b\supseteq S_a$.
Hence the value $\nu(\tilde a)$ is the sum $\sum_{u\in\mathcal U_a}f(u)$.

Taking $a=1$ shows that $f$ is a nonnegative \emph{frame function} of weight $1$: its sum over every orthonormal basis of $H$ is $1$.
Since $H$ is separable of dimension at least three, Gleason's theorem~\cite{Gleason:1957} gives a unique density operator $\rho$ with $f(x)=\ip{\rho x}{x}$ for every unit vector $x$.
It follows that $\nu(\tilde a)=\sum_{u\in\mathcal U_a}\ip{\rho u}{u}=\Tr(\rho P_a)$ whenever $S_a\ne\{0\}$; the same equality holds when $S_a=\{0\}$, since then $\tilde a=\tilde0$ and $\nu(\tilde0)=0$.
\end{proof}

\historypara{Historical note.}
Gleason proved his theorem in 1957 to answer a question of Mackey's.
Mackey had asked whether every countably additive assignment of probabilities to the closed subspaces of a separable Hilbert space is given by the Born rule~\cite{Mackey:1957,Gleason:1957}.
The Born rule gives a closed subspace, in the state described by a density operator $\rho$, the probability $\Tr(\rho P)$, where $P$ is the orthogonal projection onto the subspace.
Frame functions are the heart of Gleason's proof: he shows first, by means of spherical harmonics, that the continuous frame functions on the unit sphere of $\RR^3$ come from quadratic forms.
Gleason then shows, in the hardest step, that every nonnegative frame function on this sphere is continuous.
Finally, he passes from $\RR^3$ to real and complex separable Hilbert spaces of dimension at least three.
Gleason noted that the restriction on the dimension is essential: in dimension two the theorem fails.
Cooke, Keane and Moran gave in 1985 an elementary proof, longer than Gleason's but with fewer prerequisites~\cite{Cooke-Keane-Moran:1985}.

\section{Standard models of \texorpdfstring{$T_1$}{T1} and \texorpdfstring{$T$}{T}: Gleason's theorem, the range of the measure and definable sets of values}\label{sec:t1}

In the standard model of $T_0$ in \cref{ex:finite-cofinite}, the Solèr sort represents only the subspaces of finite dimension or finite codimension, so $T_0$ does not require every closed subspace to be represented.
Every closed subspace of a separable Hilbert space is, however, spanned by a square-summable sequence of vectors, and the theory $T_1$ adds the sequence sort and the span function $\spa$, with their axioms, to supply these sequences and their spans.
In a standard model of $T_1$ we first identify the sequence sort with the closed unit ball of the space of square-summable sequences of vectors of the Hilbert space, and then show that $\spa$ represents the closed spans of these sequences; the quotient of the Solèr sort is consequently isomorphic to the whole lattice of closed subspaces of the Hilbert space, the representation needed for \cref{thm:a,thm:b}.
In \cref{sec:definable} we use the range of the measure to show that the sets of probabilities defined uniformly in standard models are exactly the $\sigma$-projective subsets of $[0,1]$.

Until \cref{sec:definable}, fix a standard model $\fM$ of $T_1$.
Let $H$ and $K$ be the Hilbert spaces whose closed unit balls are its Hilbert and sequence sorts, and write $\sL^\fM$ for its Solèr sort.
The reduct of $\fM$ to the vocabulary of $T_0$ is a standard model of $T_0$, so the geometric representation of \cref{sec:t0} remains available.
We retain its notation: $P_a$ extends $\pi(a,\cdot)$, $S_a$ is its range, $\tilde a$ is the class of $a$ in $\sL^\fM/{=}$, and $\tau$ sends $\tilde a$ to $S_a$.
In the intended interpretation, the map $\kappa_n$ reads the entry in position $n$ of a sequence, and the map $\iota_n$ places a vector in that position (see \cref{sec:vocabulary}); the first lemma identifies the sequence sort accordingly.

\begin{lemma}\label{lem:sequence-sort}
For each $n\in\NN$, the map $\kappa_n$ from the sequence sort to the Hilbert sort extends to a bounded linear map from $K$ to $H$, and the map $s\mapsto(\kappa_ns)_{n\in\NN}$ is a linear isometry of $K$ onto the space $\ell^2(\NN,H)$ of sequences $(x_n)_{n\in\NN}$ in $H$ with $\sum_n\|x_n\|^2<\infty$.
\end{lemma}

\begin{proof}
By (K2), each map $\kappa_n$ is linear on the unit ball of $K$, and the maps $\iota_n$ extend to linear isometries from $H$ into $K$ with mutually orthogonal ranges.
By the same axiom, we have $\ip{\kappa_ms}{x}=\ip{s}{\iota_mx}$.
It follows that $\|\kappa_ms\|\le\|s\|$, and $\kappa_m\iota_nx=\delta_{mn}x$.
Each map $\kappa_m$ therefore extends to a bounded linear map from $K$ to $H$.

For $s$ and $t$ in the unit ball of $K$, axiom (K3) gives $\ip{\mathrm{res}_Ns}{t}=\sum_{n<N}\ip{\kappa_ns}{\kappa_nt}=\sum_{n<N}\ip{\iota_n\kappa_ns}{t}$, whence $\mathrm{res}_Ns=\sum_{n<N}\iota_n\kappa_ns$.
Since the identity $\ip{\kappa_ms}{x}=\ip{s}{\iota_mx}$ makes $\kappa_m$ the adjoint of the isometry $\iota_m$, the operator $\sum_{n<N}\iota_n\kappa_n$ is the orthogonal projection onto the sum of the mutually orthogonal ranges of $\iota_0,\dots,\iota_{N-1}$.
These projections increase with $N$, whence the norms $\|s-\mathrm{res}_Ns\|$ decrease, and by the second part of (K3) their infimum, and hence their limit, is $0$.
Consequently, the numbers $\ip{\mathrm{res}_Ns}{s}=\sum_{n<N}\|\kappa_ns\|^2$ converge to $\|s\|^2$, which is therefore equal to $\sum_n\|\kappa_ns\|^2$, and the map $s\mapsto(\kappa_ns)_{n\in\NN}$ is isometric on the unit ball and, being linear, on $K$.

Given $(x_n)_{n\in\NN}\in\ell^2(\NN,H)$, the partial sums $\sum_{n<N}\iota_nx_n$ form a Cauchy sequence in $K$, since the maps $\iota_n$ are isometries with mutually orthogonal ranges.
Since $K$ is complete, these partial sums converge to some $s$, which satisfies $\kappa_ns=x_n$ by the continuity of $\kappa_n$.
\end{proof}

\begin{lemma}\label{lem:span}
For every $s$ in the unit ball of $K$, the subspace $S_{\spa(s)}$ is the closed span of the entries $\kappa_ns$.
\end{lemma}

\begin{proof}
Let $W$ be this closed span.
By (Sp1), the subspace $S_{\spa(s)}$ contains every entry of $s$, and hence it contains $W$.
Conversely, let $y$ be a vector of $W^\perp$ with $0<\|y\|\le1$, and let $b=\at(y)^\perp$; by \cref{lem:lattice-truth}, the subspace $S_b$ is $y^\perp$.
Every entry $\kappa_ns$ lies in $S_b$, and so does the subspace $S_{\at(\kappa_ns)}$, which is $\CC\kappa_ns$ or $\{0\}$; by \cref{lem:lattice-truth}, we therefore have $[\at(\kappa_ns)\le b]=1$ for every $n$.
Axiom (Sp2) gives $[\spa(s)\le b]=1$, that is, $S_{\spa(s)}\subseteq y^\perp$.
Since $y$ was arbitrary, the subspace $S_{\spa(s)}$ is orthogonal to $W^\perp$, whence $S_{\spa(s)}\subseteq W^{\perp\perp}=W$.
\end{proof}

\begin{keytheorem}\label{thm:t1-iso}
Let $\fM$ be a standard model of $T_1$, let $H$ be the Hilbert space whose closed unit ball is its Hilbert sort, and, for each element $a$ of its Solèr sort $\sL^\fM$, let $S_a$ be the closed subspace of $H$ spanned by the members of $a$ and let $\tilde a$ be the class of $a$ in the quotient $\sL^\fM/{=}$ of $\sL^\fM$ by its equality, the relation $[a=b]=1$ of \cref{sec:vocabulary}.
Then $H$ is a separable infinite-dimensional complex Hilbert space, and the map $\tau\colon\tilde a\mapsto S_a$ is an isomorphism of ortholattices from $\sL^\fM/{=}$ onto the lattice $L(H)$ of closed subspaces of $H$; in particular, this quotient is a Solèr lattice.
Conversely, for every separable infinite-dimensional complex Hilbert space $H'$, there is a standard model of $T_1$ whose Hilbert sort is the closed unit ball of $H'$ and whose Solèr sort is $L(H')$.
\end{keytheorem}

\begin{proof}
By \cref{prop:sep}, the space $H$ is separable and infinite-dimensional.
By \cref{prop:embedding}, the map $\tau$ is an embedding; it remains to be shown that it is onto, for then the quotient $\sL^\fM/{=}$, being isomorphic to $L(H)$, is a Solèr lattice by \cref{thm:soler-lattices}.
Let $S$ be a closed subspace of $H$, and let $(u_n)_{n\in\NN}$ enumerate an orthonormal basis of $S$, which is countable by \cref{prop:sep}, followed by zeros if the basis is finite.
By \cref{lem:sequence-sort} there is $s$ in the unit ball of $K$ with $\kappa_ns=2^{-n-1}u_n$ for all $n$, since $\sum_n4^{-n-1}<1$, and the subspace $S_{\spa(s)}$ is $S$ by \cref{lem:span}.

For the converse, we build the model sort by sort.
The Hilbert sort is the closed unit ball of $H'$, and the constants $e_n$ are interpreted as an orthonormal basis $(e_n)_{n\in\NN}$ of $H'$.
The Solèr sort is $L(H')$, with $\pi(S,x)=P_Sx$ for the orthogonal projection $P_S$ onto $S$ and with $\at(x)=\CC x$.
For each $S\ne\{0\}$ we let $\rep(S)$ be a unit vector chosen in $S$, and we let $\rep(\{0\})$ be $0$.
The sequence sort is the closed unit ball of $\ell^2(\NN,H')$, on which the function $\iota_n$ places a vector $x$ in position $n$ of a sequence whose other entries are $0$, and the function $\kappa_n$ reads the entry of a sequence $s$ in position $n$.
The function $\mathrm{res}_N$ keeps the entries of $s$ in the positions $n<N$ and replaces the others by $0$, while the constant $\varepsilon$ is the sequence $(2^{-n-1}e_n)_{n\in\NN}$, whose norm is less than $1$.
The function $\spa$ sends each element $s$ of the sequence sort to the closed span of its entries.
Finally, each truncation $R^{\pm}_{n,q}$ is interpreted by its intended value $\min\bigl(1,\max(0,\pm n(R-q))\bigr)$, so that the axioms (Tr) hold.
Axioms (Sub1) and (Sub2) hold by \cref{thm:von-neumann}, and the others by construction; in this model, the formula $[S=S']$ takes the value $1$ exactly when $S=S'$.
\end{proof}

\begin{proof}[Proof of \cref{thm:a}]
The quotient of the Solèr sort has now been identified with $L(H)$; it remains to identify the measure.
Axiom (M4) gives finite additivity, and (M5) will turn the measure of a span into a sum over an orthonormal basis.
Together these yield the $\sigma$-additive state needed for \cref{thm:gleason-t0}.
Let $\fM$ be a standard model of $T$.
Its reduct to the vocabulary of $T_1$ is a standard model of $T_1$, to which the notation fixed at the beginning of this section applies.
By \cref{thm:t1-iso}, the space $H$ is a separable infinite-dimensional complex Hilbert space, and the map $\tau$ is an isomorphism of $\sL^\fM/{=}$ onto $L(H)$.
By (M3), the value $\mu(a)$ depends only on the class $\tilde a$, and we write $\mu(\tilde a)$ for it.

\emph{Step 1: finite additivity.}
Let $a$ and $c$ be elements of $\sL^\fM$ with $S_c\subseteq S_a^\perp$.
By \cref{lem:lattice-truth,cor:quotient}, the elements $c\wedge a^\perp$ and $c$ have the same class, and so have $a\vee(c\wedge a^\perp)$ and $a\vee c$.
Axioms (M3) and (M4) therefore give $\mu(a\vee c)=\mu(a)+\mu(c)$.
Taking $a=c=0$ gives $\mu(0\vee0)=2\mu(0)$, and since $0\vee0$ and $0$ have the same class, (M3) turns this equation into $\mu(0)=2\mu(0)$, whence $\mu(0)=0$.

\emph{Step 2: sums along orthonormal bases.}
Let $a$ be an element of $\sL^\fM$, and let $(u_k)_{k\in\NN}$ enumerate an orthonormal basis of $S_a$, followed by zeros if it is finite.
As in the proof of \cref{thm:t1-iso}, \cref{lem:sequence-sort} gives an element $s$ of the unit ball of $K$ with $\kappa_ks=2^{-k-1}u_k$ for every $k$, and by \cref{lem:span} the element $\spa(s)$ has the class $\tilde a$.
For $u_k\ne0$, the subspaces $S_{\at(\kappa_ks)}$ and $S_{\at(u_k)}$ are both $\CC u_k$ by \cref{prop:atoms}, and for $u_k=0$ both are $\{0\}$ by part~(4) of \cref{lem:lattice-truth}.
By Step~1 and (M3), the measure of $c_N(s)$ is therefore $\sum_{k<N}\mu(\at(u_k))$.
These partial sums increase with $N$, since their terms are nonnegative by (M1), and axioms (M5) and (M3) therefore give
\[
\mu(a)=\sum_k\mu(\at(u_k)). \tag{$*$}
\]

\emph{Step 3: application of Gleason's theorem.}
We show that the map $\tilde a\mapsto\mu(\tilde a)$ is a $\sigma$-additive state on $\sL^\fM/{=}$, in order to apply \cref{thm:gleason-t0}.
Its values lie in $[0,1]$ by (M1), and by (M2) it takes the value $1$ at $\tilde1$.
Let $(\tilde a_n)_{n\in J}$, for a set $J\subseteq\NN$, be a family of mutually orthogonal elements of $\sL^\fM/{=}$.
This family has a join $\tilde c$, because the map $\tau$ is onto $L(H)$.
By \cref{prop:embedding,thm:t0}, the subspaces $S_{a_n}$ are mutually orthogonal, and the closure of their sum is $S_c$; orthonormal bases of the subspaces $S_{a_n}$ therefore form together an orthonormal basis of $S_c$.
Applied to $c$, the identity $(*)$ expresses $\mu(c)$ as the sum of the values $\mu(\at(u))$ over the vectors $u$ of this basis.
Applied to each $a_n$, the identity expresses $\mu(a_n)$ as the sum of the same values over the vectors that come from $S_{a_n}$.
Since these values are nonnegative by (M1), the first sum may be grouped by $n$, whence $\mu(\tilde c)=\sum_{n\in J}\mu(\tilde a_n)$.
By \cref{thm:gleason-t0}, applied to the reduct of $\fM$ to the vocabulary of $T_0$, there is a unique density operator $\rho$ on $H$ with $\mu(a)=\Tr(\rho P_a)$ for every $a\in\sL^\fM$.

\emph{The converse.}
Let a separable infinite-dimensional complex Hilbert space $H'$ and a density operator $\rho$ on $H'$ be given.
We expand a standard model of $T_1$ built for $H'$ as in the proof of \cref{thm:t1-iso} by $\mu(a)=\Tr(\rho P_a)$, where $P_a$ is the orthogonal projection onto the closed subspace $a$ of $H'$, and we interpret the truncations of $\mu$ by their intended values, so that the axioms (Tr) for $\mu$ hold.
Axiom (M3) holds because the formula $[a=b]$ takes the value $1$ only when $a=b$.
Axioms (M1) and (M2) hold because every orthogonal projection $P$ satisfies $0\le\Tr(\rho P)\le1$ and because the trace of $\rho$ is $1$.
Axiom (M4) holds because the subspaces $a$ and $b\wedge a^\perp$ are orthogonal, whence the projection onto their join is the sum of their projections.
For (M5), let $s$ be an element of the sequence sort.
The subspaces $c_N(s)$ increase with $N$, and the closure of their union is the closed span $\spa(s)$ of the entries of $s$.
The projections onto them therefore converge strongly to $P_{\spa(s)}$, and since $\rho$ is of trace class, the numbers $\Tr(\rho P_{c_N(s)})$ increase to $\Tr(\rho P_{\spa(s)})$.
\end{proof}

The representation in \cref{thm:a} expresses the measure through a density operator, but \cref{thm:b} also needs its full range.
For a density operator $\rho$ that is not a rank-one projection, every ray probability is at most $\|\rho\|<1$, so the following lemma considers all closed subspaces, not only rays.

\begin{lemma}\label{lem:range}
For every separable infinite-dimensional complex Hilbert space $H$ and every density operator $\rho$ on $H$, the set $\{\Tr(\rho P_S):S\in L(H)\}$ is the interval $[0,1]$.
\end{lemma}

\begin{proof}
The values $0$ and $1$ are attained at $\{0\}$ and $H$, respectively.
Fix $\alpha$ with $0<\alpha<1$.
We start with a subspace spanned by eigenvectors whose probability reaches or exceeds $\alpha$, and rotate its last basis vector until the probability falls below $\alpha$.
Continuity will then give the desired value.
Since $\rho$ is a positive compact operator of trace $1$, there are orthonormal eigenvectors $v_1,v_2,\dots$ (finitely or infinitely many) whose eigenvalues $\eta_1\ge\eta_2\ge\dots>0$ are the nonzero eigenvalues of $\rho$, repeated according to multiplicity, and these eigenvalues satisfy $\sum_i\eta_i=1$.
Write $\eta_{\le k}=\eta_1+\dots+\eta_k$, with $\eta_{\le0}=0$, and choose an integer $k\ge0$ with $\eta_{\le k}<\alpha\le\eta_{\le k+1}$.

On the unit vectors of $W=\{v_1,\dots,v_{k+1}\}^\perp$, the quantity $\ip{\rho x}{x}$ has infimum $0$.
Indeed, if $\rho$ has only finitely many nonzero eigenvalues, then its kernel has finite codimension in the infinite-dimensional space $H$, and it therefore meets $W$, which also has finite codimension, in a nonzero vector; otherwise the eigenvectors $v_i$ with $i>k+1$ lie in $W$, and their eigenvalues $\eta_i$ tend to $0$.
Choose a unit vector $u\in W$ with $\ip{\rho u}{u}<\alpha-\eta_{\le k}$.

For $0\le\vartheta\le1$, let $w_\vartheta$ be the unit vector $\sqrt{1-\vartheta}\,v_{k+1}+\sqrt{\vartheta}\,u$, and let $S_\vartheta$ be the subspace spanned by $v_1,\dots,v_k$ and $w_\vartheta$.
Because $u$ is orthogonal to $v_1,\dots,v_{k+1}$, the vectors $v_1,\dots,v_k,w_\vartheta$ form an orthonormal basis of $S_\vartheta$, and the trace of $\rho P_{S_\vartheta}$ is the sum of the values of $\ip{\rho x}{x}$ at these vectors.
In the expansion of $\ip{\rho w_\vartheta}{w_\vartheta}$, the cross terms are multiples of $\ip{\rho v_{k+1}}{u}=\eta_{k+1}\ip{v_{k+1}}{u}$ and of its conjugate, and they therefore vanish.
Consequently, we have
\[
\Tr(\rho P_{S_\vartheta})=\eta_{\le k}+(1-\vartheta)\eta_{k+1}+\vartheta\ip{\rho u}{u},
\]
a continuous function of $\vartheta$ that equals $\eta_{\le k+1}\ge\alpha$ at $\vartheta=0$ and is smaller than $\alpha$ at $\vartheta=1$.
By the intermediate value theorem, it takes the value $\alpha$.
\end{proof}

\begin{proposition}[The range of the measure]\label{thm:b}
In every standard model of $T$, the range of $\mu$ is the whole interval $[0,1]$.
Consequently, for every Borel set $B\subseteq\RR$ that meets $[0,1]$, the theory $T$, together with the axiom $\beta_B(a)\le-1$, which says that no value of $\mu$ lies in $B$ (with $\beta_B$ as in \cref{lem:borel}), has no standard models.
\end{proposition}

\begin{proof}
Let $\fM$ be a standard model of $T$, and let $H$ be the Hilbert space whose closed unit ball is its Hilbert sort.
By \cref{thm:a}, there is a density operator $\rho$ on $H$ with $\mu(a)=\Tr(\rho P_a)$ for every element $a$ of the Solèr sort, and the map $\tilde a\mapsto S_a$ is onto $L(H)$.
By \cref{lem:range}, the range of $\mu$ is $[0,1]$.
Let $B\subseteq\RR$ be a Borel set that meets $[0,1]$.
Some element $a$ of the Solèr sort then satisfies $\mu(a)\in B$, and the formula $\beta_B(a)$ takes the value $1$ at this element.
Hence $\fM$ does not satisfy the axiom $\beta_B(a)\le-1$, which says that no value of $\mu$ lies in $B$.
\end{proof}

\historypara{Historical note.}
Sierpi\'nski proved in 1922 that a measure without atoms takes, on the measurable subsets of a set, every value between $0$ and the measure of the set~\cite{Sierpinski:1922}.
A $\{0,1\}$-valued measure, which is what an assignment of truth values amounts to, is the opposite extreme: the whole space is an atom.
By contrast, the measures $S\mapsto\Tr(\rho P_S)$ that density operators $\rho$ define on the closed subspaces of a separable infinite-dimensional Hilbert space take every value in $[0,1]$ (see \cref{lem:range}).
Such a measure takes all these values even when $\rho$ is a pure state, the orthogonal projection onto a ray, whose classical counterpart would be a point mass.

\subsection{Definable sets of values}\label{sec:definable}

The Borel tests of \cref{lem:borel} determine whether a probability belongs to a prescribed Borel set, using connectives alone.
For formulas depending on an additional lattice variable only through its measure, quantification over the variable performs a projection along $[0,1]$, since every value in this interval is a value of $\mu$ in every standard model of $T$ (see \cref{thm:b}).
We show in \cref{thm:definable} that the sets obtained uniformly in standard models are exactly the $\sigma$-projective subsets of the interval.
The exclusion result of \cref{thm:b} therefore extends from nonempty Borel sets of probabilities to all nonempty sets in this class.

To compare the operations on formulas with operations on sets, we use the least classes containing the Borel sets and closed under complements, countable unions and projections along $\NN^\NN$.
Aguilera, Müller and Schlicht~\cite{Aguilera-Muller-Schlicht:2021} study the resulting $\sigma$-projective sets on finite powers of $\NN^\NN$.
The proof of \cref{thm:definable} below also uses cubes of probabilities and a space of orthogonal projections, so we define the classes on Polish spaces and establish closure under Borel preimages and projections along Polish spaces.
Let $\Gamma$ be the least assignment of a class $\Gamma(X)$ of subsets of $X$ to each Polish space $X$ with the following two properties: each class $\Gamma(X)$ contains the Borel subsets of $X$ and is closed under complements in $X$ and countable unions; and, whenever $A\in\Gamma(X\times\NN^\NN)$, the set $\{x\in X:(x,z)\in A\text{ for some }z\in\NN^\NN\}$ belongs to $\Gamma(X)$.
A subset of a Polish space $X$ is \emph{$\sigma$-projective} when it belongs to $\Gamma(X)$.
In particular, every projective set~\cite[Section~37]{Kechris:1995} is $\sigma$-projective, and so is every countable union of projective sets.

\begin{definition}\label{def:definable-set}
A set $D\subseteq[0,1]$ is \emph{definable over $T$} if there is a formula $\gamma(a)$ of real-valued $L_{\omega_1\omega}$, with the single free variable $a$ of the Solèr sort, whose value, in every standard model of $T$ and at every element $a$ of its Solèr sort, is $1$ if $\mu(a)\in D$ and $-1$ otherwise; we say that $\gamma$ \emph{defines} $D$.
\end{definition}

Only subsets of $[0,1]$ matter here, because (M1) confines $\mu$ to this interval.
A formula that tests membership in $E\subseteq\RR$ in every standard model therefore defines $E\cap[0,1]$ in the sense of \cref{def:definable-set}.

\begin{theorem}\label{thm:definable}
A subset of $[0,1]$ is definable over $T$ if and only if it is $\sigma$-projective.
\end{theorem}

Consequently, if $D\subseteq[0,1]$ is a nonempty $\sigma$-projective set and $\gamma(a)$ defines $D$, then the theory $T$ together with the axiom $\gamma(a)\le-1$, which in a standard model says that no value of $\mu$ lies in $D$, has no standard models.

To define the sets, we encode membership in tuples of measure values and use the full range of $\mu$ to express projection along $[0,1]$ by quantification over the Solèr sort.
For the converse, a formula that works in every standard model must work in one whose symbols are interpreted by Borel functions on Polish spaces.
In this model we follow the descriptive complexity of the formula's values through the connectives and quantifiers, then pull back along a continuous map assigning to each probability $t\in[0,1]$ a projection of measure $t$.

\paragraph{Closure properties used in the proof.}
On the finite powers of $\NN^\NN$, the $\sigma$-projective sets are the sets of Aguilera, M\"uller and Schlicht.
Their class is the least that contains the open sets and is closed under complements, countable unions and projections along $\NN^\NN$.
This class is therefore contained in $\Gamma$, which has these properties on the finite powers of $\NN^\NN$.
Conversely, for a Polish space $X$, consider the sets $A\subseteq X$ such that $g^{-1}[A]$ belongs to the class of Aguilera, M\"uller and Schlicht for every $k\ge1$ and every Borel map $g\colon(\NN^\NN)^k\to X$.
The class of Aguilera, M\"uller and Schlicht contains the Borel sets and is closed under the operations that define $\Gamma$.
Hence these sets $A$ include the Borel sets and are closed under complements, countable unions and projections along $\NN^\NN$.
By the minimality of $\Gamma$, they include $\Gamma(X)$.
For $X=(\NN^\NN)^k$, taking for $g$ the identity map then shows that every set of $\Gamma(X)$ belongs to the class of Aguilera, M\"uller and Schlicht.

The $\sigma$-projective sets are closed under preimages by Borel maps between Polish spaces.
Indeed, the sets $A\subseteq X$ such that $g^{-1}[A]$ is $\sigma$-projective for every Borel map $g$ from a Polish space into $X$ include the Borel sets and are closed under complements and countable unions.
They are also closed under projections along $\NN^\NN$, since, for a set $A\subseteq X\times\NN^\NN$ and a Borel map $g$ from a Polish space $Z$ into $X$, the preimage under $g$ of the projection of $A$ is the projection of the preimage of $A$ under the Borel map $g\times\mathrm{id}$ from $Z\times\NN^\NN$ into $X\times\NN^\NN$.
These sets therefore include every $\sigma$-projective set.

The $\sigma$-projective sets are also closed under projections along every nonempty Polish space $Y$.
Given a continuous surjection $p\colon\NN^\NN\to Y$~\cite[Theorem~7.9]{Kechris:1995} and a $\sigma$-projective set $A\subseteq X\times Y$, the projection of $A$ to $X$ is the projection along $\NN^\NN$ of the preimage of $A$ under $\mathrm{id}\times p$, and it is therefore $\sigma$-projective.

A function $f$ from a Polish space to $\RR$ is \emph{$\Gamma$-measurable} when the set $\{f>r\}$ is $\sigma$-projective for every real number $r$.
The set $\{f\ge r\}$, which is the intersection of the sets $\{f>r-\tfrac1n\}$, is then $\sigma$-projective as well.

\begin{proof}[Proof of \cref{thm:definable}]
\emph{Every $\sigma$-projective subset of $[0,1]$ is definable.}
For $m\ge1$, consider a formula $\gamma(a_1,\dots,a_m)$ with free variables $a_1,\dots,a_m$ of the Solèr sort.
The formula $\gamma$ \emph{strongly defines} a set $A\subseteq[0,1]^m$ if it has a priori bound at most $3$, is two-valued in every model of $T$ and, in every standard model of $T$, equals $1$ at $(a_1,\dots,a_m)$ exactly when $(\mu(a_1),\dots,\mu(a_m))\in A$.
Let $\mathcal A$ be the class of sets $A\subseteq[0,1]^m$, for $m\ge1$, that some formula strongly defines.
We show that $\mathcal A$ contains every $\sigma$-projective subset of every cube $[0,1]^m$; a formula that strongly defines a set $D\subseteq[0,1]$ defines it in the sense of \cref{def:definable-set}.

First, $\mathcal A$ contains the Borel sets.
We repeat the proof of \cref{lem:borel} with the closed boxes $[q_1,q'_1]\times\dots\times[q_m,q'_m]$, for rational numbers $q_i\le q'_i$, in place of the closed intervals with rational endpoints.
The formula for such a box is the minimum
\[
\beta_{[q_1,q'_1]}(a_1)\wedge\dots\wedge\beta_{[q_m,q'_m]}(a_m)
\]
of the formulas of the original proof for the intervals $[q_i,q'_i]$.
For every Borel set $B\subseteq\RR^m$, the modified proof gives a formula $\beta_B(a_1,\dots,a_m)$.
In every structure for $\cL$ that satisfies (Tr), its value at every tuple of elements of the Solèr sort is $1$ if $(\mu(a_1),\dots,\mu(a_m))\in B$ and $-1$ otherwise.
The formula $\beta_B(a_1,\dots,a_m)$ strongly defines $B\cap[0,1]^m$.
Next, if $\gamma$ strongly defines $A\subseteq[0,1]^m$, then the formula $-\gamma$ strongly defines $[0,1]^m\setminus A$.
If $\gamma_n$ strongly defines $A_n\subseteq[0,1]^m$ for every $n\in\NN$, then the join $\bigvee_n\gamma_n$, whose members have a priori bounds at most $3$, strongly defines $\bigcup_nA_n$.
Finally, if $\gamma(a_1,\dots,a_m,b)$ strongly defines $A\subseteq[0,1]^{m+1}$, then the formula $\sup_b\gamma(a_1,\dots,a_m,b)$ is two-valued in every model of $T$.
Indeed, the Solèr sort contains the constant $0$, and a supremum of values in $\{-1,1\}$ over a nonempty set is attained.
In a standard model of $T$, the formula $\sup_b\gamma(a_1,\dots,a_m,b)$ takes the value $1$ at $(a_1,\dots,a_m)$ exactly when $(\mu(a_1),\dots,\mu(a_m),\mu(b))\in A$ for some element $b$.
By \cref{thm:b}, such an element $b$ exists exactly when $(\mu(a_1),\dots,\mu(a_m),t)\in A$ for some $t\in[0,1]$.
Thus $\mathcal A$ contains the Borel subsets of the cubes $[0,1]^m$ and is closed under complements, countable unions and projections along $[0,1]$.

The formulas above implement projections along $[0,1]$, whereas $\Gamma$ was defined using projections along $\NN^\NN$.
A Borel isomorphism between these spaces lets us compare the two constructions.
Consider the least class $\Gamma_0$ of subsets of the cubes $[0,1]^m$, for $m\ge1$, containing their Borel subsets and closed under complements in the cube, countable unions and projections from $[0,1]^{m+1}$ to $[0,1]^m$.
We have already shown that the class $\Gamma_0$ is contained in $\mathcal A$.

Fix a Borel isomorphism $\chi$ of $[0,1]$ onto $\NN^\NN$~\cite[Theorem~15.6]{Kechris:1995}.
For each Polish space $X$, let $\Gamma'(X)$ be the class of sets $A\subseteq X$ such that $g^{-1}[A]\in\Gamma_0$ for every $m\ge1$ and every Borel map $g\colon[0,1]^m\to X$.
Each class $\Gamma'(X)$ contains the Borel subsets of $X$ and is closed under complements in $X$ and countable unions.
Let $A\in\Gamma'(X\times\NN^\NN)$ be given, and let $g\colon[0,1]^m\to X$ be a Borel map.
Since $\chi$ is onto $\NN^\NN$, we have
\[
g^{-1}\bigl[\{x:\exists z\in\NN^\NN\ (x,z)\in A\}\bigr]=\{u\in[0,1]^m:\exists t\in[0,1]\ (u,t)\in(g\times\chi)^{-1}[A]\},
\]
the projection along $[0,1]$ of a set of $\Gamma_0$, because $g\times\chi$ is a Borel map on $[0,1]^{m+1}$; this projection belongs to $\Gamma_0$.
Hence the assignment $\Gamma'$ is closed under projections along $\NN^\NN$, and we have $\Gamma(X)\subseteq\Gamma'(X)$ for every Polish space $X$.
For a $\sigma$-projective set $A\subseteq[0,1]^m$, taking for $g$ the identity shows that $A\in\Gamma_0$ and therefore that $A\in\mathcal A$.

\emph{Every definable subset of $[0,1]$ is $\sigma$-projective.}
Let $\gamma(a)$ define a set $D\subseteq[0,1]$.
Let $H=\ell^2(\NN)$, let $v$ and $w$ be orthogonal unit vectors of $H$, and let $\rho$ be the projection onto $\CC v$, a density operator.
By \cref{thm:a}(2) and its proof, there is a standard model $\fM$ of $T$ that can be described as follows.
Its Hilbert sort is the closed unit ball $B_H$ of $H$, its sequence sort is the closed unit ball $B_K$ of $K=\ell^2(\NN,H)$, and its Solèr sort is $L(H)$.
The symbols are interpreted as in the proof of \cref{thm:a}, with the standard unit vectors of $\ell^2(\NN)$ as the constants $e_n$ and with $\mu(S)=\Tr(\rho P_S)$.
There, the value $\rep(S)$ may be any unit vector of $S$ when $S\ne\{0\}$, and here we let it be $P_Se_n/\|P_Se_n\|$ for the least $n$ with $P_Se_n\ne0$.
\enlargethispage{2pt}
We identify $L(H)$ with the set of orthogonal projections on $H$, with the strong operator topology.
Fix a dense sequence $(x_k)_{k\in\NN}$ in $B_H$, and give the set of operators of norm at most $1$ the strong operator topology.
For such operators, convergence at the vectors $x_k$ implies convergence at every vector.
The map $Q\mapsto(Qx_k)_{k\in\NN}$ is therefore a homeomorphism of the set of operators of norm at most $1$ onto its image in the Polish space $B_H^{\NN}$.
The image is closed: if $(y_k)_{k\in\NN}$ is the limit of sequences $(Q_jx_k)_{k\in\NN}$ in the image, then we have $\|\sum_kc_ky_k\|\le\|\sum_kc_kx_k\|$ for every finitely supported family of complex numbers $c_k$.
The map $x_k\mapsto y_k$ therefore extends to an operator of norm at most $1$.
The projections form a closed subset, since a strong limit $Q$ of projections satisfies $\ip{Qx}{x}=\|Qx\|^2=\ip{Q^*Qx}{x}$ for every $x$, whence $Q=Q^*Q$ by polarization; the operator $Q$ is therefore self-adjoint, and it is idempotent, since $Q^2=Q^*Q=Q$.
Thus the space $L(H)$ is Polish, and so are $B_H$ and $B_K$ with their norm topologies.

Every symbol of $\cL$ is interpreted in $\fM$ by a Borel function on a product of these spaces.
The vector operations, the inner products and the norms are continuous, as are the maps $\iota_n$, $\kappa_n$ and $\mathrm{res}_N$, the constants and the orthocomplement $Q\mapsto I-Q$.
The map $(Q,x)\mapsto Qx$, which interprets $\pi$, is continuous because $\|Qx-Q_0x_0\|\le\|x-x_0\|+\|(Q-Q_0)x_0\|$.
The measure, given by $\mu(Q)=\Tr(\rho Q)=\ip{Qv}{v}$, is continuous as well.

The meet $Q\wedge Q'$ is the strong limit of $(QQ'Q)^n$ as $n\to\infty$, by \cref{thm:von-neumann} and the identity $(QQ'Q)^n=(QQ')^nQ$.
As a pointwise limit of continuous maps, the meet is Borel, and the join is given by $Q\vee Q'=I-((I-Q)\wedge(I-Q'))$.
The truncations are continuous functions of the predicates that they truncate.

The map $\at$ sends $x\ne0$ to the projection $y\mapsto\ip{y}{x}x/\|x\|^2$, continuously, and $0$ to $0$.
For $\rep$, the set of projections $Q$ with $Qe_k=0$ for $k<n$ and $Qe_n\ne0$ is Borel, and on it $\rep(Q)=Qe_n/\|Qe_n\|$ is continuous.

For $\spa$, the term $c_N(s)=\at(\kappa_0s)\vee\dots\vee\at(\kappa_{N-1}s)$ of axiom (M5) is interpreted by the projection onto the span of the entries $\kappa_0s,\dots,\kappa_{N-1}s$ of $s$, which is a Borel function of $s$, being a composition of the maps $\kappa_n$, $\at$ and $\vee$.
These projections increase with $N$ and converge strongly to the projection $\spa(s)$ onto the closed span of the entries of $s$, and the map $\spa$, as a pointwise limit of Borel maps into the Polish space $L(H)$, is therefore Borel.
It follows by induction on terms that every term is interpreted by a Borel function.

We show by induction on formulas that, for every formula $\varphi(\bar x)$, the function $\bar x\mapsto\varphi^\fM(\bar x)$ is $\Gamma$-measurable on the product of the Polish spaces that interpret the sorts of $\bar x$.
Atomic formulas, those built from the truncations included, are continuous functions of Borel functions of terms, hence Borel.
For the connectives $s\cdot\varphi$, $\varphi+\psi$ and $\varphi\vee\psi$, the set where the value exceeds $r$ is the preimage of an open subset of $\RR$ or $\RR^2$ under the function given by the argument formula or pair of argument formulas.
Write the open set as a countable union of rational open intervals in the unary case, or of rational open boxes in the binary cases.
For rational endpoints $q<q'$, the preimage of $(q,q')$ under $\varphi^\fM$ is $\{\varphi^\fM>q\}\setminus\{\varphi^\fM\ge q'\}$, and box preimages are intersections of the corresponding sets for $\varphi^\fM$ and $\psi^\fM$.

For a countable join, the set $\{\bigvee_n\varphi_n>r\}$ is the union $\bigcup_n\{\varphi_n>r\}$.
For a supremum over a variable $y$, the set $\{\bar x:\sup_y\varphi(\bar x,y)>r\}$ is the projection of $\{(\bar x,y):\varphi(\bar x,y)>r\}$ along the Polish space that interprets the sort of $y$.
In each case the set is $\sigma$-projective.

Now the function $\gamma^\fM$ is $\Gamma$-measurable, and the set $\{Q\in L(H):\gamma^\fM(Q)>0\}$ is $\sigma$-projective.
For $t\in[0,1]$, let $\xi(t)$ be the projection $\at(\sqrt t\,v+\sqrt{1-t}\,w)$ onto the line spanned by the unit vector $\sqrt t\,v+\sqrt{1-t}\,w$.
The map $\xi$ is continuous and satisfies $\mu(\xi(t))=|\ip{v}{\sqrt t\,v+\sqrt{1-t}\,w}|^2=t$.
Since $\gamma^\fM$ takes the value $1$ exactly at the elements $Q$ with $\mu(Q)\in D$,
\[
D=\{t\in[0,1]:\gamma^\fM(\xi(t))>0\}=\xi^{-1}\bigl[\{\gamma^\fM>0\}\bigr],
\]
a continuous preimage of a $\sigma$-projective set, hence $\sigma$-projective.

\emph{The consequence stated after the theorem} follows as in the proof of \cref{thm:b}, with $\gamma$ in place of $\beta_B$: in a standard model, some element $a$ has $\mu(a)\in D$, and $\gamma(a)$ takes the value $1$ there.
\end{proof}

The class of $\sigma$-projective subsets of $[0,1]$ is strictly larger than that of projective ones.
For each $n\ge1$, the hierarchy theorem for projective sets~\cite[Section~37]{Kechris:1995} provides a $\boldsymbol\Sigma^1_{n+1}$ set $C_n\subseteq[0,1]$ that is not $\boldsymbol\Sigma^1_n$, and we let $h_n$ be a homeomorphism of $[0,1]$ onto $[2^{-n},2^{-n+1}]$.
The set $C=\bigcup_nh_n[C_n\cap(0,1)]$ is a countable union of projective sets, and it contains none of the points $2^{-m}$.
If $C$ were $\boldsymbol\Sigma^1_k$, then for $n>k$ the set $h_n^{-1}[C]$, which equals $C_n\cap(0,1)$, would be $\boldsymbol\Sigma^1_k$.
The set $C_n$, which differs from $C_n\cap(0,1)$ by at most two points, would then be $\boldsymbol\Sigma^1_k$ and hence $\boldsymbol\Sigma^1_n$, which is false.
If $C$ were $\boldsymbol\Pi^1_k$, it would be $\boldsymbol\Sigma^1_{k+1}$, which is excluded in the same way.
Hence $C$ is not projective, and by \cref{thm:definable} it is definable over $T$.

The upper-bound argument also relates the complexity of a defined set to the nesting of quantifiers in its defining formula, and shows why defining the nonprojective set $C$ above requires a countable join whose members have unbounded finite quantifier depths.
We define the \emph{quantifier depth} $\operatorname{qd}(\varphi)$ of a formula $\varphi$, a countable ordinal, by induction.
Atomic formulas have depth $0$, and the depth of $s\cdot\varphi$, of $\varphi+\psi$ and of $\varphi\vee\psi$ is the largest depth of the formulas to which the connective is applied.
The depth of a countable join is given by $\operatorname{qd}(\bigvee_n\varphi_n)=\sup_n\operatorname{qd}(\varphi_n)$, and that of a supremum by $\operatorname{qd}(\sup_x\varphi)=\operatorname{qd}(\varphi)+1$.

\begin{proposition}[Quantifier depth and projective complexity]\label{prop:qd}
Let $\fM$ be the standard model of $T$ constructed in the proof of \cref{thm:definable}, whose symbols are interpreted by Borel functions, and let $k<\omega$.
\begin{enumerate}[label=(\arabic*)]
\item For every formula $\varphi$ with $\operatorname{qd}(\varphi)\le k$ and every real number $r$, the set $\{\varphi^\fM>r\}$ is $\boldsymbol\Delta^1_{k+1}$.
\item Every set $D\subseteq[0,1]$ definable over $T$ by a formula of quantifier depth at most $k$ is $\boldsymbol\Delta^1_{k+1}$.
\item Every formula that defines the set $C$ above over $T$ has quantifier depth at least $\omega$, and contains a countable join whose members have unbounded finite quantifier depths.
\end{enumerate}
\end{proposition}

\begin{proof}
For part (1), the induction on formulas given in the proof of \cref{thm:definable} applies.
For $k=0$, the formula $\varphi$ contains no supremum over a variable, and these sets are Borel.
For $k\ge1$, if $\operatorname{qd}(\sup_y\varphi)\le k$, then the formula $\varphi$ has depth at most $k-1$, and the set $\{(\sup_y\varphi)^\fM>r\}$ is the projection of a $\boldsymbol\Delta^1_k$ set.
This projection is therefore $\boldsymbol\Sigma^1_k$, and hence $\boldsymbol\Delta^1_{k+1}$.
The other cases use only countable unions, countable intersections and complements, under which $\boldsymbol\Delta^1_{k+1}$ is closed.
For part (2), if a formula $\gamma$ with $\operatorname{qd}(\gamma)\le k$ defines $D$, then the set $D=\xi^{-1}[\{\gamma^\fM>0\}]$ of the proof of \cref{thm:definable} is $\boldsymbol\Delta^1_{k+1}$ by part (1).
For part (3), the set $C$ is not projective, so by part (2) every formula that defines it has quantifier depth at least $\omega$.
As the subformula relation is well founded, such a formula has a subformula of depth at least $\omega$ whose proper subformulas all have finite depth.
This subformula is a countable join whose members have unbounded finite depths.
\end{proof}

\paragraph{Uniformity across standard models.}
Definability over $T$ requires one formula to work independently of the basis, the density operator and the unit representatives chosen by $\rep$ (see \cref{def:definable-set}).
In the model $\fM$ of the proof of \cref{thm:definable}, whose symbols are interpreted by Borel functions, definability by one formula in this single model already gives definability over $T$.
Indeed, let $D\subseteq[0,1]$ be a set such that some formula takes, at every element $Q$ of the Solèr sort of $\fM$, the value $1$ if $\mu(Q)\in D$ and $-1$ otherwise.
The second half of the proof applies to this formula and shows that the set $D$ is $\sigma$-projective, and hence definable over $T$.

In an arbitrary standard model, definability by one formula cannot replace definability over $T$.
For $0<t<1$, the ray of the unit vector $\sqrt t\,e_0+\sqrt{1-t}\,e_1$ has probability $t$ in the pure state of $e_0$, and the choice of $\rep$ between this vector and its negative can encode an arbitrary subset of $(0,1)$ in one standard model.
We detect the choice by requiring both coordinates of the chosen representative to be positive and real.

To express this test, for $i\in\NN$ we write $z_i(b)$ for the pair of predicates $\ip{\rep(b)}{e_i}$, and $\nonatom(b)$ for the two-valued formula $\sup_c\bigl([c\ne0]\wedge-[c=b]\wedge[c\le b]\bigr)$, which takes the value $1$ exactly when some element lies strictly between $0$ and $b$.
For an atomic formula $R(\bar t)$ whose predicate symbol is one of the inner-product, norm or measure predicates, let $[R(\bar t)>0]=2\lceil R(\bar t)>0\rceil-1$ and $[R(\bar t)=0]=1-2\bigl(\lceil R(\bar t)>0\rceil\vee\lceil R(\bar t)<0\rceil\bigr)$, two-valued formulas in every structure that satisfies (Tr).
Let $\mathrm{eq}_\mu(a,b)$ be the formula $\beta_B(a,b)$ of the proof of \cref{thm:definable} for the diagonal $B=\{(u,u):u\in\RR\}$.
Let $\gamma_0(a)$ be the formula
\[
\begin{aligned}
\sup_b\Bigl(&[b\ne0]\wedge-\nonatom(b)\wedge\bigl[b\le\at(e_0)\vee\at(e_1)\bigr]\wedge\mathrm{eq}_\mu(a,b)\\
&\wedge[\operatorname{Im}z_0(b)=0]\wedge[\operatorname{Im}z_1(b)=0]\wedge[\operatorname{Re}z_0(b)>0]\wedge[\operatorname{Re}z_1(b)>0]\Bigr).
\end{aligned}
\]

\begin{proposition}[Definability in one standard model is not uniform]\label{prop:nonuniform}
For every set $\Omega\subseteq(0,1)$, there is a standard model of $T$ in which the formula $\gamma_0(a)$ takes the value $1$ exactly when $\mu(a)\in\Omega$, and $-1$ otherwise.
Consequently, some sets that the single formula $\gamma_0$ defines in a single standard model are not definable over $T$.
\end{proposition}

\begin{proof}
Given $\Omega\subseteq(0,1)$, take the model of the proof of \cref{thm:a}(2) on $\ell^2(\NN)$, with the standard unit vectors as the constants $e_n$ and with $\rho$ the projection onto $\CC e_0$.
In this model, let $\rep$ send the atom spanned by $\sqrt t\,e_0+\sqrt{1-t}\,e_1$, for $0<t<1$, to the vector $\sqrt t\,e_0+\sqrt{1-t}\,e_1$ itself if $t\in\Omega$ and to its negative otherwise.
For $0<t<1$, the atoms $b$ below $\at(e_0)\vee\at(e_1)$ with $\mu(b)=t$ are the lines of the vectors $\sqrt t\,e_0+\zeta\sqrt{1-t}\,e_1$ with $|\zeta|=1$.
A unit vector of such a line has two real coordinates only when $\zeta=\pm1$, and two positive real coordinates only when $\zeta=1$ and the vector is $\sqrt t\,e_0+\sqrt{1-t}\,e_1$ itself; for $t\in\{0,1\}$, one coordinate vanishes.
Hence the formula $\gamma_0(a)$ takes the value $1$ in this model exactly when $\mu(a)\in\Omega$, and $-1$ otherwise.

Since the vocabulary is countable, there are at most $2^{\aleph_0}$ formulas, and hence at most $2^{\aleph_0}$ sets definable over $T$.
The interval $(0,1)$, by contrast, has $2^{2^{\aleph_0}}$ subsets, and some sets $\Omega$ are therefore not definable over $T$.
\end{proof}

A quantified formula defining a set $D$ over $T$ need not take the value $1$ exactly when $\mu(a)\in D$ in a nonstandard model.
By \cref{rem:nonstandard}, the theory $T$ has a countable model $\fN$, so some $\alpha\in[0,1]$ is absent from its countable range of $\mu$.
The formula $\sup_b\beta_{\RR\times\{\alpha\}}(a,b)$, with the Borel test constructed in the proof of \cref{thm:definable}, searches for an element of measure $\alpha$.
It defines $[0,1]$ over $T$ by \cref{thm:b}, yet takes the value $-1$ everywhere in $\fN$.

The quantifier-free tests $\beta_B(a)$ of \cref{lem:borel} retain their intended meaning in every structure satisfying (Tr): they take the value $1$ exactly when $\mu(a)\in B$.
\Cref{rem:standard-models} below concerns all models of $T$ and therefore uses these Borel tests.
The classification in \cref{thm:definable} concerns standard models, where the full interval of probabilities is available to the quantifiers.

\section{The proof of Theorem~C}\label{sec:bks}

The preceding representation theorems concern standard models, but the exclusion of a $\{0,1\}$-valued measure holds in every model of $T$.
We prove \cref{thm:c} directly from the axioms, using the Hilbert space $\bar E$ constructed in \cref{lem:envelope} to supply the geometry.
Assuming that $\mu$ takes only the values $0$ and $1$, we apply (M5) to the named sequence $\varepsilon$ to locate a basis atom of measure $1$, then obtain a finite parity contradiction within a block of four basis vectors.

Within the model, the argument uses only the vectors denoted by the closed terms described below.
By (H2), the constants $e_n$ are orthonormal, and by (K4) the entries of the sequence constant $\varepsilon$ satisfy $\kappa_n\varepsilon\doteqq2^{-n-1}e_n$.
These equations use the relation $\doteqq$ of \cref{sec:vocabulary}, which asserts that a scaled difference has norm zero.
By \cref{lem:envelope}, the related elements have the same image in $\bar E$.
They need not be identical as elements of the model, since the equality-free vocabulary imposes no requirement that elements agreeing in all predicates be identical.
Closed terms formed from four of the constants $e_n$ using the functions $\lambda x+\lambda'y$ represent nonzero scalar multiples of the eighteen vectors of \cref{tab:ceg}, and therefore their rays.
Their images in $\bar E$ suffice for the parity argument, which requires neither the whole unit ball nor the whole lattice of closed subspaces inside the model.

The localization step needs the named orthonormal basis and the measure of its span, while the parity step needs the corresponding rays and additivity on orthogonal subspaces.
The following lemma gathers these facts in a form valid for every model of $T$.
For a model of $T_1$, we let $\bar E$, the map $x\mapsto\hat x$ and the projections $\bar P_a$ be as in \cref{lem:envelope}, applied to the reduct of the model to the vocabulary of $T_0$, and we let $S_a$ be the range of $\bar P_a$, as in \cref{lem:lattice-truth}.

\begin{lemma}\label{lem:bks-tools}
Let $\fM$ be a model of $T$, standard or not, and let $\bar E$, $\hat x$ and $S_a$ be as above.
Then the following hold:
\begin{enumerate}[label=(\arabic*)]
\item The vectors $\hat e_n$, for $n\in\NN$, form an orthonormal basis of $\bar E$.
\item For every $x$ in the Hilbert sort with $\|x\|>0$, the subspace $S_{\at(x)}$ is $\CC\hat x$, and the formula $[\at(x)=\at(y)]$ takes the value $1$ whenever $\|x\|>0$ and $\hat y=\lambda\hat x$ for a nonzero $\lambda\in\CC$.
\item The measure satisfies $\mu(0)=0$, and $\mu(a\vee c)=\mu(a)+\mu(c)$ for all elements $a,c$ of the Solèr sort with $S_c\subseteq S_a^\perp$.
\item We have $[\spa(\varepsilon)=1]=1$ and $\mu(\spa(\varepsilon))=1$.
\end{enumerate}
\end{lemma}

\begin{proof}
(1) By (H2), the vectors $\hat e_n$ are orthonormal.
By \cref{lem:lattice-truth}, the projection $\bar P_{b_N}$ is the projection onto the span of $\hat e_0,\dots,\hat e_{N-1}$, and (H2) places every vector $\hat x$ in the closed span of the vectors $\hat e_n$; since the vectors $\hat x$ span a dense subspace, the vectors $\hat e_n$ form an orthonormal basis of $\bar E$.

(2) The first assertion is part~(4) of \cref{lem:lattice-truth}.
If $\hat y=\lambda\hat x$ with $\lambda\ne0$, then $\|y\|>0$ and $S_{\at(y)}=\CC\hat x=S_{\at(x)}$, and part~(3) of the same lemma gives $[\at(x)=\at(y)]=1$.

(3) The argument of Step~1 of the proof of \cref{thm:a} applies to $\fM$, with $S_a$ the range of $\bar P_a$, because it uses, besides axioms (M3) and (M4), only \cref{lem:lattice-truth}, which holds in every model of $T_0$.

(4) By (Sp1) and (K4), the subspace $S_{\spa(\varepsilon)}$ contains every vector $2^{-n-1}\hat e_n$, and hence, by part~(1), it is $\bar E=S_1$.
Part~(3) of \cref{lem:lattice-truth} gives $[\spa(\varepsilon)=1]=1$, and axioms (M3) and (M2) give $\mu(\spa(\varepsilon))=1$.
\end{proof}

The proofs of parts~(1) and~(2), and of the equation $[\spa(\varepsilon)=1]=1$ in part~(4), use no axiom that mentions $\mu$; these assertions therefore hold in every model of $T_1$.

\begin{table}[ht]
\centering
\renewcommand{\arraystretch}{1.15}
\begin{tabular}{@{}cllll@{}}
\toprule
basis & \multicolumn{4}{c}{vectors (coordinates with respect to an orthonormal basis $g_0,g_1,g_2,g_3$)} \\
\midrule
$1$ & $(0,0,0,1)$ & $(0,0,1,0)$ & $(1,1,0,0)$ & $(1,-1,0,0)$ \\
$2$ & $(0,0,0,1)$ & $(0,1,0,0)$ & $(1,0,1,0)$ & $(1,0,-1,0)$ \\
$3$ & $(1,-1,1,-1)$ & $(1,-1,-1,1)$ & $(1,1,0,0)$ & $(0,0,1,1)$ \\
$4$ & $(1,-1,1,-1)$ & $(1,1,1,1)$ & $(1,0,-1,0)$ & $(0,1,0,-1)$ \\
$5$ & $(0,0,1,0)$ & $(0,1,0,0)$ & $(1,0,0,1)$ & $(1,0,0,-1)$ \\
$6$ & $(1,-1,-1,1)$ & $(1,1,1,1)$ & $(1,0,0,-1)$ & $(0,1,-1,0)$ \\
$7$ & $(1,1,-1,1)$ & $(1,1,1,-1)$ & $(1,-1,0,0)$ & $(0,0,1,1)$ \\
$8$ & $(1,1,-1,1)$ & $(-1,1,1,1)$ & $(1,0,1,0)$ & $(0,1,0,-1)$ \\
$9$ & $(1,1,1,-1)$ & $(-1,1,1,1)$ & $(1,0,0,1)$ & $(0,1,-1,0)$ \\
\bottomrule
\end{tabular}
\caption{The configuration of Cabello, Estebaranz and García-Alcaine~\cite{Cabello-Estebaranz-GarciaAlcaine:1996}.}
\label{tab:ceg}
\end{table}

\begin{proof}[Proof of \cref{thm:c}]
Suppose, towards a contradiction, that $\fM$ is a model of $T$ in which $\mu$ takes only the values $0$ and $1$, and let $\bar E$, $\hat x$ and $S_a$ be as in \cref{lem:bks-tools} for this model.
We first show that exactly one of the constants $e_n$ has an atom $\at(e_n)$ of measure $1$.
We then apply a parity argument to the bases of \cref{tab:ceg} within the block of four consecutive constants that contains this one.
The entries $\kappa_n\varepsilon$ satisfy $\widehat{\kappa_n\varepsilon}=2^{-n-1}\hat e_n$ by (K4), and part~(2) of \cref{lem:bks-tools} gives $[\at(\kappa_n\varepsilon)=\at(e_n)]=1$.
Since the subspaces $S_{\at(e_n)}=\CC\hat e_n$ are mutually orthogonal, repeated applications of part~(3) of \cref{lem:bks-tools}, combined with (M3), which replaces each $\mu(\at(\kappa_n\varepsilon))$ by $\mu(\at(e_n))$, give $\mu(c_N(\varepsilon))=\sum_{n<N}\mu(\at(e_n))$ for every $N\in\NN$.
Since the terms of these partial sums are nonnegative by (M1), axiom (M5) and \cref{lem:bks-tools}(4) give
\[
1=\mu(\spa(\varepsilon))=\sum_n\mu(\at(e_n)).
\]
As $\mu$ takes only the values $0$ and $1$, there is exactly one $n_0\in\NN$ with $\mu(\at(e_{n_0}))=1$.
Let $j\in\NN$ satisfy $4j\le n_0<4j+4$, write $g_i=e_{4j+i}$ for $i=0,1,2,3$, and let $F$ be the span of $\hat g_0,\dots,\hat g_3$.

Among the four constants $g_0,\dots,g_3$, whose images $\hat g_0,\dots,\hat g_3$ form an orthonormal basis of $F$, the constant $e_{n_0}$ is the only one whose atom has measure $1$.
We show that every orthogonal basis of $F$ whose vectors have coordinates in $\{0,1,-1\}$ with respect to $\hat g_0,\dots,\hat g_3$ likewise contains exactly one vector $v$ whose atom $\at(v)$, defined below, has measure $1$.
Let $\cB$ be such a basis, consisting of vectors $v_0,\dots,v_3$ of $F$.
The coordinates of each vector $\tfrac14v_i$ have absolute values with sum at most $1$; by \cref{lem:envelope}, some closed term $t_i$, built from the constants $g_0,\dots,g_3$ by the functions $\lambda x+\lambda'y$, therefore satisfies $\hat t_i=\tfrac14v_i$.
Hence the term $\tfrac12\bigl(\tfrac12\iota_0t_0+\tfrac12\iota_1t_1\bigr)+\tfrac12\bigl(\tfrac12\iota_2t_2+\tfrac12\iota_3t_3\bigr)$ denotes an element $s_\cB$ of the sequence sort.
Axiom (K2) gives $\ip{\kappa_m\iota_nx}{y}=\ip{\iota_nx}{\iota_my}=\delta_{mn}\ip{x}{y}$ for all $x$ and $y$, whence $\widehat{\kappa_m\iota_nx}=\delta_{mn}\hat x$, because the vectors $\hat y$ span a dense subspace of $\bar E$; together with the linearity of $\kappa_n$ asserted in (K2), this gives $\widehat{\kappa_ns_\cB}=\tfrac1{16}v_n$ for $n<4$.
Let $c_\cB$ be the element $\at(\kappa_0s_\cB)\vee\dots\vee\at(\kappa_3s_\cB)$, that is, $c_4(s_\cB)$.

For a nonzero vector $v$ of $F$ with coordinates in $\{0,1,-1\}$, we fix a closed term $t_v$, built in the same way, with $\hat t_v=\tfrac1{16}v$, and we write $\at(v)$ for $\at(t_v)$.
By \cref{lem:bks-tools}(2), another choice of $t_v$ changes $\at(v)$ only within the relation $[\cdot=\cdot]=1$, and therefore, by (M3), leaves $\mu(\at(v))$ unchanged.

By the same part of \cref{lem:bks-tools}, we have $[\at(\kappa_is_\cB)=\at(v_i)]=1$ for $i<4$.
By \cref{lem:lattice-truth}, the subspace $S_{c_\cB}$ is the span of $v_0,\dots,v_3$, which is $F$.
Hence $[c_\cB=c_{\cB'}]=1$ for any two such bases $\cB$ and $\cB'$, and we let $a_F$ be the element $c_\cB$ for the basis $\hat g_0,\dots,\hat g_3$.
By \cref{lem:bks-tools}(3), applied three times, and by (M3),
\[
\mu(a_F)=\mu(c_\cB)=\sum_{v\in\cB}\mu(\at(v)).
\]
For the basis $\hat g_0,\dots,\hat g_3$, part~(2) of \cref{lem:bks-tools} and (M3) give $\mu(\at(\hat g_i))=\mu(\at(e_{4j+i}))$ for $i<4$, and this sum therefore contains $\mu(\at(e_{n_0}))=1$ and three zeros; hence $\mu(a_F)=1$, and every such basis $\cB$ contains exactly one vector $v$ with $\mu(\at(v))=1$.

The nine bases of \cref{tab:ceg} are orthogonal bases of $F$ of this kind, and their eighteen vectors span pairwise distinct rays.
We conclude by counting in two ways the pairs $(\cB,v)$ that consist of one of the nine bases $\cB$ and a vector $v$ of $\cB$ with $\mu(\at(v))=1$.
By the preceding paragraph, each basis contains exactly one such vector, and there are therefore exactly nine pairs.
On the other hand, each of the eighteen vectors $v$ occurs in exactly two of the bases, and the atom $\at(v)$ fixed above is the same element in both.
Each vector $v$ with $\mu(\at(v))=1$ thus belongs to exactly two pairs, and the number of pairs is even.
This contradiction proves the theorem.
\end{proof}

In physical terms, \cref{thm:c} excludes noncontextual hidden-variable assignments with probabilities in $\{0,1\}$ that satisfy the measure axioms of $T$.
The experimentally verifiable propositions are represented by closed subspaces~\cite{Birkhoff-vonNeumann:1936}, whose role in a model of $T$ is played by the elements of the Solèr sort.
For each value of the hidden variable, the assignment gives a definite truth value to each proposition, independently of which other propositions are tested together with it.
Such an assignment would be recorded by a dispersion-free measure, taking the numerical values $0$ and $1$ for false and true respectively.
By \cref{thm:c}, the axioms of $T$, including the continuity requirement (M5) along represented spans, exclude this assignment as the interpretation of $\mu$.
Contextual theories, such as Bohm's~\cite{Bohm:1952}, allow the value revealed by a measurement to depend on the whole experimental arrangement and lie outside this conclusion.
Together with (M1), the theorem gives the following corollary.

\begin{corollary}[No dispersion-free measures]\label{cor:no-hidden}
In every model of $T$, standard or not, some element $a$ of the Solèr sort satisfies $0<\mu(a)<1$.
\end{corollary}

\begin{proof}
By (M1), the measure takes values in $[0,1]$, and by \cref{thm:c} it does not take only the values $0$ and $1$.
\end{proof}

\historypara{Historical note.}
Specker announced in 1960 that the closed subspaces of a three-dimensional Hilbert space cannot be mapped into a Boolean algebra in a way that respects the operations on orthogonal subspaces.
Of the proof he said only that it rests on an elementary geometric argument~\cite{Specker:1960}.
Bell derived from Gleason's theorem that no hidden-variable theory of the kind described above is possible, in a paper written in 1964 and published in 1966~\cite{Bell:1966}.\footnote{Bell's geometric comparisons of vectors and rays~\cite[Section~V]{Bell:1966} admit a modern topological interpretation.
For unit vectors $f_0,f_1$ of a complex Hilbert space, let $Q_i$ be the orthogonal projection onto $\CC f_i$, for $i=0,1$.
The identities $\|Q_0-Q_1\|=\|(Q_0-Q_1)f_0\|=\sqrt{1-|\ip{f_0}{f_1}|^2}$ and $\inf_{|\lambda|=1}\|f_0-\lambda f_1\|=\sqrt{2-2|\ip{f_0}{f_1}|}$, with $\lambda\in\CC$, identify both the strong operator topology and the operator norm topology on rank-one projections with the quotient of the norm topology on the unit sphere by scalar phases.}
In 1967 Kochen and Specker gave a finite proof of this impossibility with $117$ vectors in $\RR^3$~\cite{Kochen-Specker:1967}.
The configurations then shrank: Peres found configurations of $33$ vectors in $\RR^3$ and of $24$ in $\RR^4$ in 1991~\cite{Peres:1991}, and Kernaghan one of $20$ vectors in $\RR^4$ in 1994~\cite{Kernaghan:1994}.
In 1996 Cabello, Estebaranz and García-Alcaine found the eighteen vectors of \cref{tab:ceg}~\cite{Cabello-Estebaranz-GarciaAlcaine:1996}.
Peres conjectured that no smaller configuration exists, and Xu, Chen and Gühne proved in 2020 that no Kochen--Specker set in any dimension has fewer than eighteen vectors~\cite{Xu-Chen-Guhne:2020,Budroni-Cabello-Guhne-Kleinmann-Larsson:2022}.
In dimension three the minimum remains unknown at the time of writing.
Exhaustive computer searches raised the lower bound for this minimum to $22$ vectors in 2016~\cite{Uijlen-Westerbaan:2016} and to $24$ in 2023.
The bound of $24$ was reached independently by Kirchweger, Peitl and Szeider and by Li, Bright and Ganesh~\cite{Kirchweger-Peitl-Szeider:2023,Li-Bright-Ganesh:2024}.
The titles of the papers of Kernaghan and of Cabello, Estebaranz and García-Alcaine call the result the Bell--Kochen--Specker theorem, a name that records Bell's share and that \cref{thm:c} keeps.

\subsection{The Bell--Kochen--Specker obstruction in the game}\label{sec:bks-game}
The contradiction of \cref{thm:c} gives a concrete case in which the finite obstructions of \cref{sec:game} arise from an infinitary theory.
In this subsection the vocabulary and the theory are again those of \cref{sec:vocabulary} and of \cref{thm:a}.
We extend the theory $T$ of \cref{thm:a} by the assertion that its measure $\mu$ takes only the values $0$ and $1$.
The normal form of this assertion is $\sup_a\bigl(\mu(a)\wedge(\mathbf 1-\mu(a))\bigr)\le0$, which we call the \emph{dispersion-free axiom}.
By \cref{thm:c}, the extended theory has no models, so \cref{cor:finite-obstruction} applies.
In \cref{prop:bks-game} below we show that every finite set of its axioms has a model, and in \cref{cor:bks-branching} we deduce that the tree of every winning strategy of Olivia in the game for this theory has an infinitely branching node.
In \cref{prop:parity-strategy} below we construct a strategy following the proof of \cref{thm:c}: every play reaching the parity argument ends with commitments concerning the eighteen vectors and nine bases of \cref{tab:ceg} in a block of four consecutive constants $e_n$.

\begin{proposition}\label{prop:bks-game}
Let $T^{01}$ be the theory $T$ of \cref{thm:a} together with the axiom
\[
\sup_a\bigl(\mu(a)\wedge(\mathbf 1-\mu(a))\bigr)\le0.
\]
\begin{enumerate}
\item Olivia has a winning strategy in $\mathrm{G}(T^{01})$; every play in which she follows it ends after finitely many rounds, at a position that is satisfied in no structure.
\item Every finite set of axioms of $T^{01}$ has a model.
\item The skeleton of every winning strategy of Olivia in $\mathrm{G}(T^{01})$ has more than one element: in some play in which Olivia follows such a strategy, Josh answers a countable challenge.
\end{enumerate}
\end{proposition}

\begin{proof}
Part (1) follows from parts (1) and (2) of \cref{cor:finite-obstruction}, since $T^{01}$ has no models by \cref{thm:c}.

For part (2), recall that in the proof of \cref{thm:c} it is the instances of (K4), taken together, that make the entries of $\varepsilon$ span the whole space, and that only then does the $\sigma$-additivity of $\mu$ along $\varepsilon$ give some atom the measure $1$.
A finite set of axioms contains only finitely many of these instances, and we therefore look for a model of it among the structures whose sequences are all finitely supported.
In such a structure, axiom (M5) involves only finite-dimensional subspaces, and a measure that vanishes on every finite-dimensional subspace satisfies it.
Let a finite set of axioms of $T^{01}$ be given, and let $Z\subseteq\NN$ be the finite set of the numbers $n$ for which this set contains the instance $\kappa_n\varepsilon\doteqq2^{-n-1}e_n$ of axiom (K4).
A model of the given set is the following structure.
Its Hilbert sort is the closed unit ball of the space $\ell^2$ of square-summable complex sequences, and its Solèr sort is the lattice $L_{\mathrm{fc}}$ of \cref{ex:finite-cofinite}.
The symbols $\pi$, $\at$, $\rep$ and the constants $e_n$ are interpreted as in this example.
The sequence sort consists of the finitely supported sequences in the closed unit ball of the space $\ell^2(\NN,\ell^2)$ of square-summable sequences in $\ell^2$.
The interpretations of $\iota_n$, $\kappa_n$ and $\mathrm{res}_N$ are those of the proof of \cref{thm:t1-iso}, and the function $\spa$ sends each element $s$ of the sequence sort to the finite-dimensional span of its entries.
The measure $\mu$ vanishes on the finite-dimensional subspaces and equals $1$ on their orthogonal complements, and each truncation $R^{\pm}_{n,q}$ is interpreted by its intended value $\min\bigl(1,\max(0,\pm n(R-q))\bigr)$, so that the axioms (Tr) hold.
The constant $\varepsilon$ is the finitely supported sequence whose entry of index $n$ is $2^{-n-1}e_n$ for $n\in Z$ and $0$ otherwise; its norm is less than $1$.
The axioms of $T_0$ hold as in \cref{ex:finite-cofinite}, and the axioms (K1)--(K3) and (Sp) hold as in the proof of \cref{thm:t1-iso}, since the span of finitely many vectors is closed; in this structure, the formula $[a=b]$ takes the value $1$ exactly when $a=b$.
Of the instances of (K4), the only axiom in which $\varepsilon$ occurs, exactly those with $n\in Z$ hold.
Axiom (M3) follows from the description of the formula $[a=b]$ given above, and axioms (M1), (M2) and the dispersion-free axiom hold because $\mu$ takes only the values $0$ and $1$, with $\mu(1)=1$.
For (M4), the subspaces $a$ and $b\wedge a^\perp$ are orthogonal, and at most one of them is infinite-dimensional; the measure of their join is $1$ exactly when one of them has measure $1$.
For (M5), let $s$ be an element of the sequence sort; the element $\spa(s)$ and the elements $c_N(s)$ are finite-dimensional subspaces, and both sides of (M5) vanish.
The structure is therefore a model of the given finite set of axioms.

For part (3), suppose that the skeleton of a winning strategy $\sigma$ of Olivia had only one element, the empty sequence.
Then Josh would give no countable answer at any node of $\cP_\sigma$: a node at which $\sigma$ makes a countable challenge would have only one child, the extension by the challenge alone.
Every other node would have finitely many children by \cref{cor:finite-obstruction}(3).
By \cref{cor:finite-obstruction}(4), the theory $T^{01}$ would then have a finite set of axioms without models, contradicting part (2).
\end{proof}

\begin{corollary}\label{cor:bks-branching}
Every winning strategy of Olivia in the model existence game for the theory of \cref{thm:c} has a node at which its tree branches infinitely.
\end{corollary}

\begin{proof}
The theory of \cref{thm:c}, written in the normal form of Subsection~\ref{sec:game-setting}, is the theory $T^{01}$ of \cref{prop:bks-game}.
By part~(2) of the proposition, every finite set of axioms of $T^{01}$ has a model.
It follows from part (3) that Josh gives a countable answer at some node of the tree of every winning strategy of Olivia.
The interval that Josh adds under the rules ($\bigvee{\ge}$) and ($\bigwedge{\le}$) does not depend on the number he names.
Every number is therefore a legal answer at this node, and the tree branches infinitely there.
\end{proof}

The proof of \cref{thm:c} first uses countable additivity to locate an atom of measure $1$, then derives a finite parity contradiction in the block containing it.
To follow the proof in the game, Olivia needs recorded bounds with enough numerical margin for the next step.
For example, the equation $\mu(\spa(\varepsilon))=1$ uses every instance of (K4), whereas any position records only finitely much information.
We therefore first explain how Olivia can observe values to arbitrary accuracy.
In \cref{lem:finite-continuations} below we show how fair bookkeeping either exposes a contradiction or obtains the required margins after finitely many rounds in each play, without asserting a uniform bound on their number.

\paragraph{Observing values.}
The margins in the localization and parity argument concern predicate values at closed terms.
Olivia first names the required closed terms, from their innermost subterms outward, by finitely many challenges (Tot); the resulting records make the names denote the terms in every structure satisfying the position.
She then obtains intervals of arbitrarily small width for their predicate values, using the truncation axioms.
For a predicate symbol $R$ other than $\mathbf 1$ and the truncation symbols, the procedure at given witnesses is as follows.
Challenges (Ax) activate the axiom (Tr) for $R$ with $n=1$ and $q=0$; challenges ($\sup{\le}$) and ($\inf{\ge}$) then give commitments on its instance at the given witnesses.
Olivia follows the finite chain of connective functions down to the occurrence of $R$ by challenges (Cn), ending with the addition $R(\bar w)-0$.
The interval $I_{R(\bar w)}$ has length at least $2$, since the a priori bound of the predicate symbol is a positive integer, whereas the interval of the sentence $-0$ is the single point $0$.
The number $m$ prescribed by (Cn) for this addition is therefore $k$, so a tolerance $2^{-k}$ yields an interval of width at most $2^{-k}$ for $R(\bar w)$.
These challenges may expose a contradiction before the observation is complete, in which case Josh loses.
By an observation of $R(\bar t)$ we mean the naming records together with an interval commitment on the resulting basic sentence $R(\bar w)$; bounds on the value of $R(\bar t)$ follow from the records.

At the start of the play, Olivia makes the challenge (Tot) for the constant $\varepsilon$, and we write $w_\varepsilon$ for the witness of the record $\varepsilon\mapsto w_\varepsilon$ that Josh adds.
Below, Olivia also observes the value of the sentence
\[
\chi=\bigvee_{N\in\NN}\mu(c_N(w_\varepsilon)).
\]
The instance of (M5) at $w_\varepsilon$, obtained by (Ax) followed by the rules ($\sup{\le}$) and ($\inf{\ge}$), contains $\chi$ as an argument of the negation in $\mu(\spa(w_\varepsilon))-\chi$.
Challenges (Cn) along this finite chain give intervals of arbitrarily small width for $\chi$ as well.
Here the observation concerns the value of the join, not the choice of one of its members; such a choice will be made separately by ($\bigvee{\ge}$).

\begin{lemma}[Finite continuations]\label{lem:finite-continuations}
Let $S$ be a nonempty countable subtheory of $T^{01}$, regarded as a theory in the full vocabulary, and let $\mathfrak p$ be a position reached in the game.
Let $\psi_1,\dots,\psi_M$, for a positive integer $M$, be sentences whose values can be observed as above using axioms of $S$.
Let $U\subseteq\RR^M$ be open, and suppose that every model of $S$ satisfying $\mathfrak p$, under an interpretation of its witnesses, has
\[
(\psi_1,\dots,\psi_M)\in U.
\]
Olivia has a continuation, using only axioms of $S$, under which every play either ends with Josh losing or reaches, after finitely many rounds, a position whose commitments confine this tuple of values to a nonempty closed box contained in $U$.
\end{lemma}

\begin{proof}
Fix an enumeration of all the challenges that can occur in continuations of $\mathfrak p$ in which every challenge (Ax) names an axiom of $S$; among them are the challenges (Tot), the challenges of the observing procedures and the challenges to the commitments of $\mathfrak p$.
At stage $k\ge1$, Olivia makes every challenge among the first $k$ entries that is available at the start of the stage, and then observes all the listed values at tolerance $2^{-k}$, in a fixed order.
Each stage has finitely many rounds unless Josh loses, and every challenge that becomes available is eventually made, since availability persists.
For each observed sentence she intersects its recorded intervals with its fixed interval of possible values; these intersections form nested closed intervals, and the observations make their widths tend to zero.
She ends the continuation when their product is nonempty and contained in $U$.

Suppose that a play continued indefinitely without either outcome.
It would be fair for $S$, and the construction in the proof of \cref{thm:game-existence}, applied from the initial position $\mathfrak p$, would give a model of $S$ satisfying every commitment and record made in the play.
This construction does not require the initial position to be empty: its commitments are included in the same induction on formulas, and its records in the definition of the functions.
The tuple of observed values in this model belongs to $U$ by hypothesis and lies in every recorded box.
Since the boxes have widths tending to zero and $U$ is open, some box is contained in $U$, a contradiction.
An empty intersection cannot occur in such an infinite surviving play, since the constructed model supplies a point of every interval.
Thus every play has one of the asserted finite outcomes.
\end{proof}

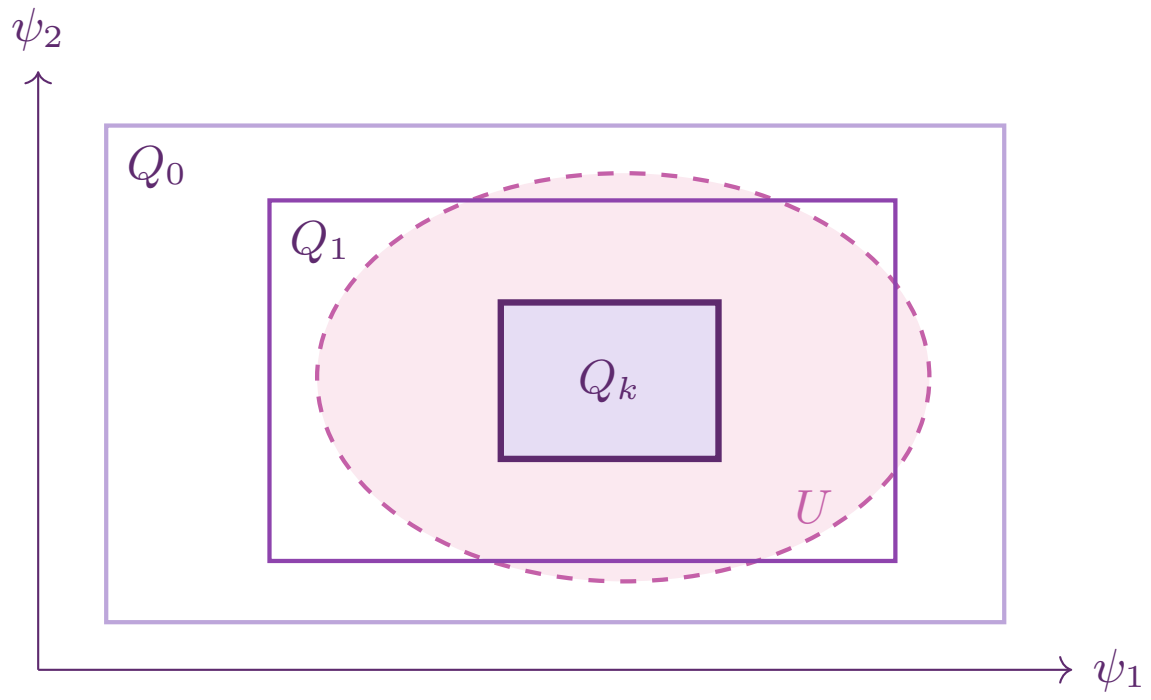
\begin{figure}[ht]
\centering
\begin{tikzpicture}[font=\footnotesize,scale=0.9]
\draw[->,plum] (0,0) -- (7.6,0) node[right] {$\psi_1$};
\draw[->,plum] (0,0) -- (0,4.4) node[above] {$\psi_2$};
\draw[rose,dashed,thick,fill=blush!55] (4.3,2.15) ellipse (2.25 and 1.5);
\node[rose] at (5.7,1.2) {$U$};
\draw[lavender!90!plum,thick] (0.5,0.35) rectangle (7.1,4.0);
\node[plum,anchor=north west] at (0.5,4.0) {$Q_0$};
\draw[orchid,thick] (1.7,0.8) rectangle (6.3,3.45);
\node[plum,anchor=north west] at (1.7,3.45) {$Q_1$};
\draw[plum,very thick,fill=lavender!45] (3.4,1.55) rectangle (5.0,2.7);
\node[plum] at (4.2,2.12) {$Q_k$};
\end{tikzpicture}
\caption{The finite continuation of \cref{lem:finite-continuations}.}
\label{fig:finite-continuation}
\end{figure}

The openness of the target lets Olivia stop when the recorded intervals themselves give the required margin, for example when a value exceeds a rational threshold or a finite sum lies within a prescribed tolerance of another value.
For the targets chosen in localization and parity, the consequences of \cref{lem:envelope,lem:lattice-truth,lem:bks-tools} verify the hypothesis of \cref{lem:finite-continuations}: every model of the chosen subtheory $S$ satisfying the current position has its tuple of values in the target open set.
The lemma converts these conclusions about all such models into finite commitments with the required margins.
At no stage must Olivia decide whether the current position has a model of $S$: she stops when its recorded box lies in $U$, or wins if Josh loses first (see \cref{fig:finite-continuation}).

\begin{proposition}[The localization and parity strategy]\label{prop:parity-strategy}
Olivia has a winning strategy in $\mathrm G(T^{01})$ such that every play in which she follows it ends after finitely many rounds and, unless Josh loses earlier, passes through the following two phases:
\begin{enumerate}
\item Josh gives a countable answer $N\ge1$, and Olivia selects an index $n_0<N$ for which Josh has committed the measure of the atom $\at(\kappa_{n_0}\varepsilon)$ to an interval whose left endpoint exceeds $15/16$;
\item the play records finite geometric commitments for the nine bases of \cref{tab:ceg} in the block containing $e_{n_0}$, obtained without any axiom that mentions $\mu$, and concludes by a finitely branching continuation that uses finitely many instances of (M1), (M3), (M4) and the dispersion-free axiom.
\end{enumerate}
\end{proposition}

\begin{proof}
\emph{Localization.}
In every model of $T$ satisfying the record for $w_\varepsilon$, the sentence $\chi$ has value $1$, by (M5) and \cref{lem:bks-tools}(4).
Apply \cref{lem:finite-continuations} with $S=T$ and the open interval $(15/16,17/16)$.
Unless Josh loses first, the intersection of the recorded intervals for $\chi$ has lower endpoint greater than $15/16$.
This endpoint is the maximum of finitely many recorded lower endpoints and the fixed lower bound $-2$; hence one actual commitment for $\chi$ has lower endpoint $\alpha>15/16$.
A challenge ($\bigvee{\ge}$) with tolerance $1/16$ then makes him name an index $N\in\NN$ and commit
\[
\mu(c_N(w_\varepsilon))\ge\alpha-1/16>7/8.
\]
If $N=0$, the selected term is $c_0(w_\varepsilon)=0$; fair bookkeeping for $T$ exposes a contradiction after finitely many rounds, since an infinite surviving play would give a model of $T$ with $\mu(0)>7/8$, contrary to \cref{lem:bks-tools}(3).
We may therefore suppose that $N\ge1$.

Write $h_n=\at(\kappa_n w_\varepsilon)$ for $n<N$.
In every model of $T$ satisfying the record $\varepsilon\mapsto w_\varepsilon$, the finite-additivity calculation in the proof of \cref{thm:c} gives
\[
\sum_{n<N}\mu(h_n)=\mu(c_N(w_\varepsilon)).
\]
{\emergencystretch=1em In every model of $T$ satisfying the current position, the commitment just made gives $\mu(c_N(w_\varepsilon))>7/8$, and the sum of the values $\mu(h_n)$ therefore exceeds $3/4$.
A further application of \cref{lem:finite-continuations}, with $S=T$ and the open set of tuples whose sum exceeds $3/4$, confines the values $\mu(h_n)$ to a closed box on which their sum exceeds $3/4$.
Let $\delta_N=\min(1/16,1/(4N))$.
We apply \cref{lem:finite-continuations} once more, now with $S$ consisting of (M1), the dispersion-free axiom and the axioms (Tr) for $\mu$, and with the open set of tuples whose coordinates all lie in $(-\delta_N,\delta_N)\cup(1-\delta_N,1+\delta_N)$.
In every model of these axioms, each value $\mu(h_n)$ is $0$ or $1$, and the hypothesis of the lemma is therefore satisfied.
The lemma gives a refinement of the recorded box in which each coordinate interval lies within one of these two intervals.
If all intervals were of the first kind, their sum would be less than $N\delta_N\le1/4$, contrary to the preceding box condition.
Choose the least index $n_0<N$ whose interval is of the second kind.
The lower endpoint of the interval for $\mu(h_{n_0})$ exceeds $1-\delta_N\ge15/16$; as for $\chi$, it is the lower endpoint of one of the commitments recorded for $\mu(h_{n_0})$, since it exceeds the fixed lower bound $-2$.
This completes the first phase.\par}

\emph{Geometric identities.}
Choose $j\in\NN$ with $4j\le n_0<4j+4$ and use the closed terms $s_\cB$ and $\at(v)$ of the proof of \cref{thm:c} for the nine bases of \cref{tab:ceg} in the block containing $e_{n_0}$.
The four contextual atom terms of a basis $\cB$ are $a_{\cB,i}=\at(\kappa_i s_\cB)$, for $i<4$.
We consider, in addition, the basis of the block formed by the four constants $e_{4j},\dots,e_{4j+3}$, ordered with $e_{n_0}$ first and the other three in increasing order of index.
For this basis we use the atom terms $\at(e_i)$ themselves, rather than contextual ones, and we now let $a_F$ be their join in this order; like the element $a_F$ of the proof of \cref{thm:c}, it names the span of the block.
To apply (M4) successively, we replace each new atom by its meet with the orthocomplement of the preceding partial join.
For any of these ten ordered bases, write its atom terms as $a_0,\dots,a_3$, and form the terms
\[
d_1=a_0,\qquad d_{i+1}=d_i\vee a_i,\qquad
a'_i=a_i\wedge d_i^\perp,\qquad d'_{i+1}=d_i\vee a'_i
\quad(1\le i\le3).
\]
Let $\mathcal D$ be the following finite list of pairs of closed terms of the Solèr sort: the two contextual terms of each shared ray; the pairs $(a'_i,a_i)$ and $(d'_{i+1},d_{i+1})$ for each ordered basis; the pairs $(d_4,a_F)$ for these bases; and the pair $(h_{n_0},\at(e_{n_0}))$.
Olivia will obtain commitments making the formula $[u=v]$ positive for each pair $(u,v)\in\mathcal D$, so that (M3) forces the two measure values to agree.
In every model of $T_1$ satisfying the records, the two terms of each pair $(u,v)$ in this list name the same closed subspace, and they therefore satisfy $[u=v]=1$ by part~(3) of \cref{lem:lattice-truth}.
For the atom terms, these subspaces are computed as in the proof of \cref{thm:c}, by \cref{lem:envelope}, (K2), (K4) and part~(2) of \cref{lem:bks-tools}; for the other terms, they follow from those of the atom terms by part~(1) of \cref{lem:lattice-truth}.
These conclusions use no measure axiom.

For such a pair, the formula $[u=v]$ is the minimum of $1-2\|\rep(u\wedge(u\wedge v)^\perp)\|$ and the corresponding expression with $u$ and $v$ interchanged.
The norms in these expressions vanish in every model of $T_1$ satisfying the records.
Apply \cref{lem:finite-continuations} with $S=T_1$ to these finitely many atomic norms and the open set of tuples in which each of these norms is less than $1/8$.
Unless Josh loses earlier, the resulting basic commitments ensure
\[
[u=v]>3/4\qquad((u,v)\in\mathcal D)
\]
in every structure satisfying the position.
Name the finitely many closed subterms involved by (Tot), retaining their records.
Thus both the contextual ray identities and the identities needed for finite additivity have been recorded to a definite margin before any probability is carried across them.

\emph{The finite parity continuation.}
Retain the finite set $\mathfrak p_0$ consisting of these function records and norm commitments, together with the interval for $\mu(h_{n_0})$ whose lower endpoint exceeds $15/16$.
Bring into play the instances of (M3) at the pairs in $\mathcal D$, the instances of (M4) at $(d_i,a_i)$ for the ten bases, and the instances of (M1) and the dispersion-free axiom at the finitely many terms occurring in these measure conditions.
Each instance is obtained by (Ax) followed by the rules ($\sup{\le}$) and ($\inf{\ge}$), applied to the witnesses that name its terms; an equation contributes its two closed inequalities.
Denote the resulting finite family of ground commitments by $\mathfrak q$.
All its sentences are built from atomic sentences by the connective functions alone, without suprema, infima or countable joins.

There is no assignment of values to its finitely many basic sentences, within their prescribed bounds and the intervals of $\mathfrak p_0$, that satisfies these ground conditions when terms are evaluated through their records.
Indeed, (M3) and the positive antecedents give $\mu(u)=\mu(v)$ for each pair in $\mathcal D$.
For each ordered basis, (M4) gives $\mu(d'_{i+1})=\mu(d_i)+\mu(a'_i)$, and the recorded identities consequently give
\[
\mu(d_{i+1})=\mu(d_i)+\mu(a_i),\qquad
\mu(a_F)=\sum_{i<4}\mu(a_i).
\]
In the named basis of the block, the first summand is $\mu(\at(e_{n_0}))=\mu(h_{n_0})>15/16$ and hence is $1$ by the dispersion-free condition and (M1).
All other summands are nonnegative and $\mu(a_F)\le1$ by (M1), whence $\mu(a_F)=1$.
Each of the nine bases of \cref{tab:ceg} therefore has exactly one contextual atom term of measure $1$.
Since the two contextual terms of each shared ray have equal measures, the double count in the proof of \cref{thm:c} yields the same contradiction, between an odd and an even number of pairs.

Olivia now fairly challenges only the finitary commitments of $\mathfrak p_0\cup\mathfrak q$ and their descendants under (Cn) and (Tm).
All their closed terms already have records; (Tm) merely replaces them by their recorded witnesses.
If this continuation survived indefinitely, the compact intervals for its finitely many basic sentences would have common values, because no position was infeasible.
The induction on finitary formulas in the proof of \cref{thm:game-existence}, using all tolerances in the fair schedule, would make these values satisfy $\mathfrak p_0\cup\mathfrak q$, contrary to the preceding paragraph.
Thus the continuation wins.
It makes no countable challenge and is finitely branching, so its tree is finite by K\H{o}nig's lemma.
Other commitments made earlier in the play may expose a contradiction sooner; they are not needed for this final argument.

Each of the preceding continuations ends after finitely many rounds in every play, and only finitely many of them follow the naming of $N$; the final parity continuation, moreover, has a finite tree.
Once the enumerations are fixed and the least index is taken at each choice, these instructions determine Olivia's move at every finite play in which Josh has not yet lost.
Every play in which she follows them is therefore finite, and the resulting strategy is a winning strategy with the two asserted phases.
\end{proof}

\begin{figure}[ht]
\centering
\begin{tikzpicture}[font=\footnotesize]
\begin{scope}[yshift=0cm]
\fill[blush!70,rounded corners=3pt] (-0.25,-0.6) rectangle (2.05,0.3);
\draw[rose,thick,rounded corners=3pt] (-0.25,-0.6) rectangle (2.05,0.3);
\draw[lavender,rounded corners=3pt] (2.55,-0.6) rectangle (4.65,0.3);
\draw[lavender,rounded corners=3pt] (5.15,-0.6) rectangle (7.45,0.3);
\fill[plum] (0.00,0) circle (1.6pt); \node[below,font=\scriptsize] at (0.00,-0.05) {$e_{0}$};
\fill[plum] (0.60,0) circle (1.6pt); \node[below,font=\scriptsize] at (0.60,-0.05) {$e_{1}$};
\fill[plum] (1.20,0) circle (1.6pt); \node[below,font=\scriptsize] at (1.20,-0.05) {$e_{2}$};
\fill[plum] (1.80,0) circle (1.6pt); \node[below,font=\scriptsize] at (1.80,-0.05) {$e_{3}$};
\fill[plum] (2.70,0) circle (1.6pt); \node[below,font=\scriptsize] at (2.70,-0.05) {$e_{4}$};
\fill[plum] (3.30,0) circle (1.6pt); \node[below,font=\scriptsize] at (3.30,-0.05) {$e_{5}$};
\fill[plum] (3.90,0) circle (1.6pt); \node[below,font=\scriptsize] at (3.90,-0.05) {$e_{6}$};
\fill[plum] (4.50,0) circle (1.6pt); \node[below,font=\scriptsize] at (4.50,-0.05) {$e_{7}$};
\fill[plum] (5.40,0) circle (1.6pt); \node[below,font=\scriptsize] at (5.40,-0.05) {$e_{8}$};
\fill[plum] (6.00,0) circle (1.6pt); \node[below,font=\scriptsize] at (6.00,-0.05) {$e_{9}$};
\fill[plum] (6.60,0) circle (1.6pt); \node[below,font=\scriptsize] at (6.60,-0.05) {$e_{10}$};
\fill[plum] (7.20,0) circle (1.6pt); \node[below,font=\scriptsize] at (7.20,-0.05) {$e_{11}$};
\node at (7.85,0) {$\cdots$};
\node[rose,above] at (0.9,0.3) {block containing $e_{n_0}$};
\end{scope}
\begin{scope}[xshift=11.0cm,yshift=0.1cm]
\draw[orchid!70] (0.000,2.100) -- (1.350,1.609);
\draw[orchid!70] (0.000,2.100) -- (0.718,-1.973);
\draw[orchid!70] (0.000,2.100) -- (2.068,0.365);
\draw[orchid!70] (0.000,2.100) -- (-1.819,-1.050);
\draw[orchid!70] (1.350,1.609) -- (0.718,-1.973);
\draw[orchid!70] (1.350,1.609) -- (-2.068,0.365);
\draw[orchid!70] (1.350,1.609) -- (1.819,-1.050);
\draw[orchid!70] (2.068,0.365) -- (1.819,-1.050);
\draw[orchid!70] (2.068,0.365) -- (-0.718,-1.973);
\draw[orchid!70] (2.068,0.365) -- (-1.819,-1.050);
\draw[orchid!70] (1.819,-1.050) -- (-0.718,-1.973);
\draw[orchid!70] (1.819,-1.050) -- (-2.068,0.365);
\draw[orchid!70] (0.718,-1.973) -- (-1.350,1.609);
\draw[orchid!70] (0.718,-1.973) -- (-0.718,-1.973);
\draw[orchid!70] (-0.718,-1.973) -- (-1.350,1.609);
\draw[orchid!70] (-1.819,-1.050) -- (-2.068,0.365);
\draw[orchid!70] (-1.819,-1.050) -- (-1.350,1.609);
\draw[orchid!70] (-2.068,0.365) -- (-1.350,1.609);
\fill[rose] (0.675,1.854) circle (2.1pt);
\fill[rose] (0.359,0.063) circle (2.1pt);
\fill[rose] (1.034,1.232) circle (2.1pt);
\fill[rose] (-0.909,0.525) circle (2.1pt);
\fill[rose] (1.034,-0.182) circle (2.1pt);
\fill[rose] (-0.359,0.987) circle (2.1pt);
\fill[rose] (1.584,0.279) circle (2.1pt);
\fill[rose] (1.943,-0.343) circle (2.1pt);
\fill[rose] (0.675,-0.804) circle (2.1pt);
\fill[rose] (0.125,-0.343) circle (2.1pt);
\fill[rose] (0.550,-1.512) circle (2.1pt);
\fill[rose] (-0.125,-0.343) circle (2.1pt);
\fill[rose] (-0.316,-0.182) circle (2.1pt);
\fill[rose] (0.000,-1.973) circle (2.1pt);
\fill[rose] (-1.034,-0.182) circle (2.1pt);
\fill[rose] (-1.943,-0.343) circle (2.1pt);
\fill[rose] (-1.584,0.279) circle (2.1pt);
\fill[rose] (-1.709,0.987) circle (2.1pt);
\node[draw=plum,fill=lavender!70,circle,inner sep=1.5pt,minimum size=0.5cm] at (0.000,2.100) {$1$};
\node[draw=plum,fill=lavender!70,circle,inner sep=1.5pt,minimum size=0.5cm] at (1.350,1.609) {$2$};
\node[draw=plum,fill=lavender!70,circle,inner sep=1.5pt,minimum size=0.5cm] at (2.068,0.365) {$3$};
\node[draw=plum,fill=lavender!70,circle,inner sep=1.5pt,minimum size=0.5cm] at (1.819,-1.050) {$4$};
\node[draw=plum,fill=lavender!70,circle,inner sep=1.5pt,minimum size=0.5cm] at (0.718,-1.973) {$5$};
\node[draw=plum,fill=lavender!70,circle,inner sep=1.5pt,minimum size=0.5cm] at (-0.718,-1.973) {$6$};
\node[draw=plum,fill=lavender!70,circle,inner sep=1.5pt,minimum size=0.5cm] at (-1.819,-1.050) {$7$};
\node[draw=plum,fill=lavender!70,circle,inner sep=1.5pt,minimum size=0.5cm] at (-2.068,0.365) {$8$};
\node[draw=plum,fill=lavender!70,circle,inner sep=1.5pt,minimum size=0.5cm] at (-1.350,1.609) {$9$};
\end{scope}
\end{tikzpicture}
\caption{The parity configuration of \cref{prop:parity-strategy}.}
\label{fig:ceg}
\end{figure}
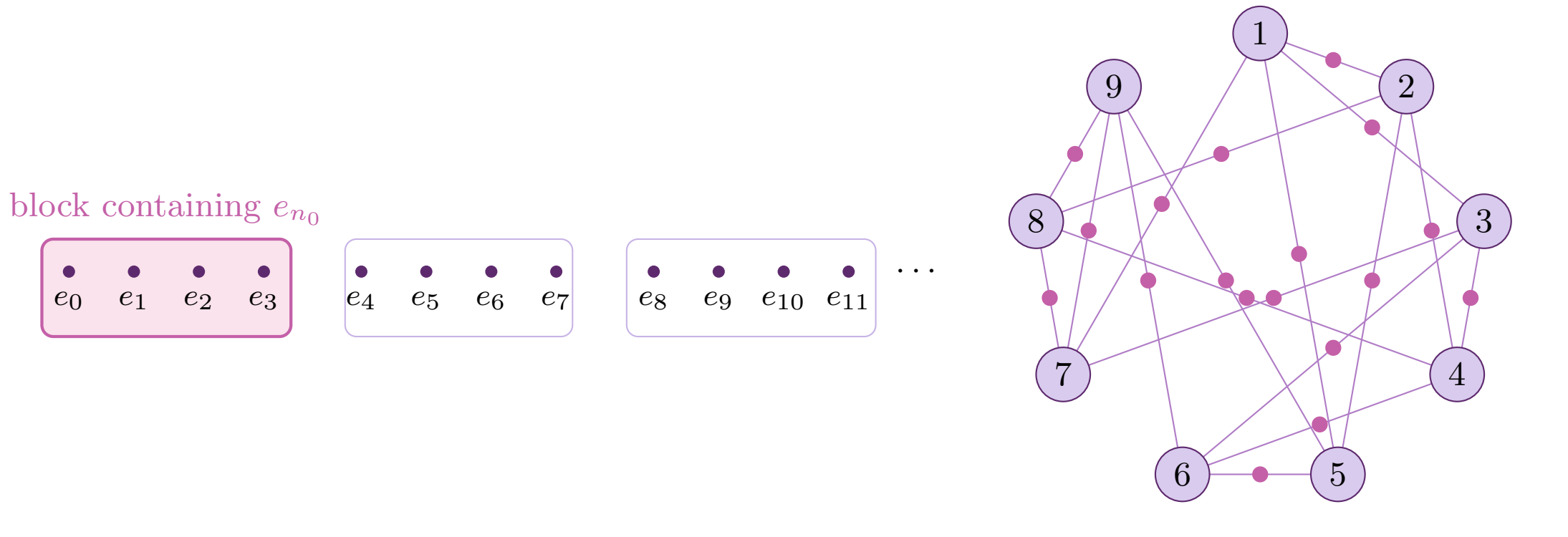

The strategy thus separates the countable search for localization and geometric margins from the final contradiction, which uses only the finite incidence pattern and finitary measure conditions.
The parity configuration is pictured in \cref{fig:ceg}.

By \cref{prop:bks-game}(3), every winning strategy of Olivia in $\mathrm{G}(T^{01})$ has some play in which Josh gives a countable answer.
For each fixed winning strategy and countable record, the lengths of plays with this record are bounded by \cref{cor:finite-obstruction}(5); neither this bound nor the finite parity continuation gives a bound independent of the record.
Whether such a bound can be achieved by a winning strategy remains open.

\begin{question}\label{q:bounded}
Does Olivia have a winning strategy in $\mathrm{G}(T^{01})$ whose tree has finite rank, that is, such that every play in which she follows it ends within a fixed number of rounds?
\end{question}

\begin{remark}[Standard models]\label{rem:standard-models}
The game $\mathrm{G}(T)$ concerns all models of $T$, standard or not, whereas \cref{thm:b} concerns the standard models, which are singled out by a hypothesis on the model rather than by axioms.
For a Borel set $B$ that meets $[0,1]$, let $T_B$ be the theory $T$ together with the axiom $\sup_a\beta_B(a)\le-1$, the normal form of the axiom of \cref{thm:b}.
Here the formula $\beta_B(a)$, given by \cref{lem:borel}, takes the value $1$ at $a$ when $\mu(a)\in B$ and $-1$ otherwise, in every structure that satisfies (Tr).
By \cref{thm:game-existence}, Olivia has a winning strategy in $\mathrm{G}(T_B)$ exactly when $T_B$ has no models at all, standard or not.
This is the case when $B$ contains $0$ or $1$, by \cref{lem:bks-tools}(3) and (M2), and when $B$ contains $(0,1)$, by \cref{cor:no-hidden}.
Some other Borel sets $B$ give theories with models.
Indeed, the range of $\mu$ in a countable model $\fN$ of $T$ (see \cref{rem:nonstandard}) is countable and contains $0$ and $1$ and, by \cref{cor:no-hidden}, a point of $(0,1)$.
For the complement $B$ of this range in $[0,1]$, the structure $\fN$ is a model of $T_B$.
We do not know which of the remaining Borel sets $B$ give theories $T_B$ with models.
\end{remark}

\section{Finite obstructions record what a Kochen--Specker configuration assumes}\label{sec:obstructions}

In the strategy of \cref{prop:parity-strategy}, countable additivity and noncontextuality enter at two distinct moments.
Countable additivity selects the four-dimensional block in which the parity argument is carried out; within the block, the geometric axioms establish that the occurrences of a ray in different bases denote the same projection, and only then does axiom (M3) carry probabilities across these identities.
We examine these two moments, after recalling what noncontextuality and determinism require of a hidden-variable model, and then give numerical tolerances under which the parity contradiction remains valid.
As in \cref{sec:bks-game}, $T^{01}$ denotes the theory $T$ of \cref{thm:a} together with the requirement that $\mu$ take only the values $0$ and $1$.

Kochen and Specker certified, by a finite configuration of rays in a finite-dimensional Hilbert space, that no quantum probability measure on the closed subspaces of the space takes only the values $0$ and $1$.
The bases of such a configuration are the sets of mutually orthogonal rays of the configuration that span the space.
An assignment of the values $0$ and $1$ to its rays is \emph{admissible} if it gives the value $1$ to exactly one ray of each basis and to no two orthogonal rays.
The configurations of Kochen and Specker have no admissible assignment.\footnote{An interactive presentation of the smallest Kochen--Specker configuration found so far in our ongoing computational search is maintained, together with the current state of this search, on the \href{https://utsa-quantum-math.github.io/assets/notes/Kochen_Specker_set.html}{UTSA Quantum Mathematics website}.}
We compare what such a configuration and the finite obstructions of the game record, how the measure enters each of them, and which topological argument gives their finiteness.

\subsection{Noncontextuality and determinism concern the responses at a single state}
One and the same projection can represent outcomes of different measurements: the projection onto a ray is an outcome of every measurement whose orthonormal basis contains a vector spanning the ray.
To interpret the measure of the game, we must therefore keep three things apart: the outcome events themselves, the projections that represent them, and the conditional probabilities that a hidden-variable model assigns to them.
Keeping them apart explains why, in the game, the geometric identities between projections and the requirements on the measure enter through separate axioms.

In the operational language of quantum foundations, an \emph{outcome event} is a specified outcome of a specified measurement procedure.
Quantum theory represents a projective measurement on a complex Hilbert space $H$ by a family of mutually orthogonal projections that sum to the identity, which we call a \emph{context}, and each of its outcome events by a member of this family.
A preparation procedure is represented by a density operator $\rho$ on $H$, a positive operator of trace $1$, and the Born rule gives an event represented by a projection $P$ the probability $\Tr(\rho P)$.
Two outcome events are \emph{operationally equivalent} when no preparation distinguishes them, that is, when their probabilities agree for every preparation; when all quantum preparations are available, this happens exactly when the two events are represented by the same projection.

A hidden-variable theory, in the form that Spekkens calls an \emph{ontological model}, supposes that each preparation selects an underlying state according to a probability distribution, and that each underlying state determines, for every outcome event, a \emph{response} in $[0,1]$, the conditional probability of the event given that state~\cite{Spekkens:2005}.
At each underlying state, the responses to the outcomes of a measurement sum to $1$, and averaging the responses over the distribution selected by a preparation must reproduce the operational probabilities.
\emph{Measurement noncontextuality} requires operationally equivalent events to receive equal responses at every underlying state, and \emph{outcome determinism} requires the responses to take only the values $0$ and $1$.
Both requirements concern the responses at individual underlying states, whereas operational equivalence concerns only their averages, which depend on the preparation.

Consider an ontological model that reproduces the quantum probabilities for all preparations and for the projective measurements on $H$ with finitely many outcomes.
At a fixed underlying state, measurement noncontextuality makes the response to an event depend only on the projection that represents it, and outcome determinism restricts the response to the values $0$ and $1$; together they yield a map from the projections of $H$ to $\{0,1\}$.
This map is normalized, since the responses within one measurement sum to $1$, and it is finitely additive: for orthogonal projections $P$ and $Q$, comparison of the measurements $\{P,Q,1-P-Q\}$ and $\{P+Q,1-P-Q\}$ shows that the response to $P+Q$ is the sum of the responses to $P$ and to $Q$.
On the whole lattice of projections, finite additivity already yields an obstruction: a fixed finite family of projections of infinite rank refutes every finitely additive $\{0,1\}$-valued measure (\cref{sec:localization}).

Countable additivity enters with the measurements that have infinitely many outcomes.
A nondegenerate measurement in an orthonormal basis $(f_n)_{n\in\NN}$ of an infinite-dimensional space $H$ has as outcomes the rays of the vectors $f_n$, and the requirement that, at each underlying state, the responses to these outcomes sum to $1$ is countable additivity along the basis.
The predicate $\mu$ of $T^{01}$ represents the map determined by the responses, and axiom (M5), which requires the measure to be continuous along the finite joins that exhaust a represented span, contains this normalization.
The strategy of \cref{sec:bks-game} uses (M5), together with the instances of (K4), along the orthonormal sequence named by the constants $e_n$ to locate a ray of measure close to $1$ and the block of four consecutive basis vectors that contains the ray.

Even on the measurements with finitely many outcomes whose events have finite rank or finite corank, no outcome-deterministic, measurement-noncontextual model reproduces the Born probabilities.
Let $F$ be any four-dimensional subspace, and complete each basis of a configuration of Kochen and Specker in $F$ by $F^\perp$.
If a dispersion-free noncontextual valuation on these events gave $F$ the value $1$, it would give $F^\perp$ the value $0$, and each of the completed bases would have its unique ray of value $1$ inside $F$, which the parity argument excludes.
Every ray lies in a four-dimensional subspace, so finite additivity gives the value $0$ to every ray, and hence to every finite-dimensional subspace.
The only such valuation is thus the measure $\mu_{\mathrm{fc}}$ of \cref{prop:fc-measure}, which gives the value $0$ even to a rank-one outcome that is certain for its eigenpreparation, and no average of such valuations reproduces the Born probabilities.
What $\mu_{\mathrm{fc}}$ shows is that normalization and finite additivity alone are consistent on this restricted lattice.
A measurement in an orthonormal basis, which has infinitely many outcomes, excludes $\mu_{\mathrm{fc}}$ by countable additivity, and \cref{sec:localization} shows how countable additivity chooses the rays that carry the obstruction.

The dispersion-free requirement thus concerns the responses at a single underlying state, whereas a density operator gives the Born probabilities of a preparation, which may take every value in $[0,1]$.

\subsection{Countable additivity chooses the configuration}\label{sec:localization}
A configuration of Kochen and Specker is fixed in advance, and one configuration serves against every assignment.
Such a configuration can be tested in the laboratory: D'Ambrosio and his coauthors implemented the eighteen yes--no tests of the configuration of \cref{tab:ceg} on four-dimensional single-photon systems in 2013~\cite{DAmbrosio-Herbauts-Amselem-Nagali-Bourennane-Sciarrino-Cabello:2013}.
Under the strategy of \cref{sec:bks-game}, the location of the ray configuration is determined by Josh's answers.
In the localization phase of \cref{prop:parity-strategy}, the constants $e_n$, for $n\in\NN$, give the orthonormal basis $(\hat e_n)_{n\in\NN}$ of the Hilbert envelope $\bar E$, and (K4) prescribes the entries of $\varepsilon$ by $\widehat{\kappa_n\varepsilon}=2^{-n-1}\hat e_n$.
During localization, a countable challenge to (M5) at the constant $\varepsilon$ makes Josh name a number $N$ and commit the measure of the join of the rays of $e_0,\dots,e_{N-1}$ to a value greater than $7/8$.
His later answers single out some $n_0<N$ for which he commits the measure of the ray of $e_{n_0}$ to a value greater than $15/16$.
The eighteen rays then lie in the span $F$ of the block $e_{4j},\dots,e_{4j+3}$ with $4j\le n_0<4j+4$.

Under this strategy, the countable-additivity scheme along the named orthonormal basis selects the four-dimensional block used for the parity argument, and the models of finite sets of axioms of $T^{01}$ explain why this localization uses the entries prescribed for $\varepsilon$.
In the models of the finite sets of axioms of $T^{01}$ constructed in the proof of \cref{prop:bks-game}(2), the measure vanishes on every finite-dimensional subspace.
These structures satisfy (M5) on the sequences present in their restricted sequence sorts.
Their measures are nevertheless not $\sigma$-additive on the lattice $L_{\mathrm{fc}}$ of \cref{ex:finite-cofinite}: the rays of the vectors $e_n$ are mutually orthogonal and have join $1$ in this lattice, yet each of them has measure $0$.
The entries that (K4) prescribes for $\varepsilon$ bring this countable family within reach of (M5).
Countable additivity along any cofinite part of the basis would also exclude these measures, since the span of such a part has a finite-dimensional orthogonal complement, and therefore measure $1$, while each of its rays has measure $0$.
In finite dimension this step is automatic, since a configuration of Kochen and Specker spans the space.

The models of the finite sets of axioms of $T^{01}$ described above rest on the following elementary fact, which isolates what finite additivity cannot exclude on finite-dimensional events.

\begin{proposition}\label{prop:fc-measure}
Let $H$ be a separable infinite-dimensional complex Hilbert space, and let $L_{\mathrm{fc}}(H)$ be the set of finite-dimensional subspaces of $H$ and their orthogonal complements.
Then $L_{\mathrm{fc}}(H)$ is a sub-ortholattice of $L(H)$, and the map $\mu_{\mathrm{fc}}$ that vanishes on the finite-dimensional subspaces and takes the value $1$ on their orthogonal complements is a finitely additive $\{0,1\}$-valued probability measure on $L_{\mathrm{fc}}(H)$: $\mu_{\mathrm{fc}}(H)=1$, and $\mu_{\mathrm{fc}}(a\vee b)=\mu_{\mathrm{fc}}(a)+\mu_{\mathrm{fc}}(b)$ whenever $a$ and $b$ are orthogonal elements of $L_{\mathrm{fc}}(H)$.
\end{proposition}

\begin{proof}
Up to unitary equivalence $H=\ell^2$, and $L_{\mathrm{fc}}(H)$ is then the sub-ortholattice $L_{\mathrm{fc}}$ of $L(\ell^2)$ of \cref{ex:finite-cofinite}.
Since $H$ is infinite-dimensional, no subspace is both finite-dimensional and the orthogonal complement of a finite-dimensional subspace, so the map $\mu_{\mathrm{fc}}$ is well defined, and $\mu_{\mathrm{fc}}(H)=1$ because $H=\{0\}^\perp$.
Let $a$ and $b$ be orthogonal elements of $L_{\mathrm{fc}}(H)$.
If $a=G^\perp$ with $G$ finite-dimensional, then $b\subseteq a^\perp=G$, so $b$ is finite-dimensional.
Hence at most one of $a$ and $b$ is infinite-dimensional, and $a\vee b$ is infinite-dimensional exactly when one of them is.
Both sides of the additivity equation therefore equal $1$ when one of $a$ and $b$ is infinite-dimensional, and $0$ otherwise.
\end{proof}

A configuration of rays inside a finite-dimensional subspace $F$ can refute a $\{0,1\}$-valued measure only when the measure gives $F$ the value $1$, and $\mu_{\mathrm{fc}}$ gives every finite-dimensional subspace the value $0$.
On every finite list of members of $L_{\mathrm{fc}}(H)$, the values of $\mu_{\mathrm{fc}}$ are the Born probabilities of a pure preparation: a unit vector orthogonal to the finite-dimensional members of the list and to the orthogonal complements of the others lies in each of these others, and so gives the prescribed values $0$ and $1$.
No density operator gives these values on all of $L_{\mathrm{fc}}(H)$, since it would give probability $0$ to every ray of an orthonormal basis $(f_n)_{n\in\NN}$.
The Born probabilities of the preparations $|f_n\rangle\langle f_n|$ nevertheless converge pointwise on $L_{\mathrm{fc}}(H)$ to $\mu_{\mathrm{fc}}$, since $\langle f_n,Pf_n\rangle\to0$ for every projection $P$ of finite rank and the probabilities of the complementary projections tend to $1$.
Each of these preparations has total probability $1$ on the basis, whereas the limit gives every ray of the basis the value $0$: normalization along a countable basis is not preserved in the limit.

On the full projection lattice, finite additivity alone already gives an obstruction using infinite-rank projections.
Write a separable infinite-dimensional Hilbert space $H$ as $\CC^4\otimes\ell^2$, and, for each vector $v$ of \cref{tab:ceg}, let $P_v$ be the projection of $\CC^4$ onto the ray of $v$.
For each basis $\cB$ of the table, the four projections $P_v\otimes1$ with $v$ in $\cB$ are pairwise orthogonal with sum $1$.
The eighteen projections $P_v\otimes1$, which have infinite rank, therefore form a configuration that is fixed in advance and refutes every finitely additive $\{0,1\}$-valued probability measure on $L(H)$ by the same parity argument.
Alda used the same amplification with a three-dimensional configuration~\cite{Alda:1981}.
More generally, Döring proved the Kochen--Specker theorem for von Neumann algebras without summands of types $\mathrm I_1$ and $\mathrm I_2$~\cite{Doring:2005}.
His argument passes through the finitely additive form of Gleason's theorem due to Christensen, Yeadon, and Bunce and Wright~\cite{Christensen:1982,Yeadon:1983,Bunce-Wright:1992} (see Hamhalter~\cite{Hamhalter:2003}).

The Solèr sort of the models of finite sets of axioms of $T^{01}$ constructed in the proof of \cref{prop:bks-game}(2) is the lattice $L_{\mathrm{fc}}$ of \cref{ex:finite-cofinite}.
This lattice contains no subspace of infinite dimension and infinite codimension, and in particular none of the ranges of the projections $P_v\otimes1$.
The $\{0,1\}$-valued measures of these models therefore do not contradict the finite-additivity obstruction on the full lattice $L(H)$.

A configuration of Kochen and Specker certifies the impossibility by the absence of an admissible assignment, whereas a finite obstruction certifies it by a finite position of commitments and function records that no structure satisfies.
The first kind of certificate rests on compactness: as recalled in the Introduction, when the Hilbert space is finite-dimensional, the condition that an orthonormal basis contain exactly one ray of value $1$ defines a closed set of assignments, and finitely many of these closed sets already have empty intersection, which turns the theorem of J.~S.~Bell into a finite configuration~\cite{Bell-JL:1996,Pitowsky:1998}.
The condition that the probabilities of the rays of the vectors $e_n$ add up to $1$ defines a set that is not closed, and conditions of this kind, expressed by the countable suprema and infima of real-valued $L_{\omega_1\omega}$, lie beyond the reach of compactness.
For this logic, the space of Henkin types is \v{C}ech-complete but in general not compact, and the Baire category theorem in that space supports the omitting types argument recalled in \cref{app:henkin}.
The finite obstructions of the game come from model construction along fair plays; the bookkeeping of this construction has its own Baire category reading, in the space of plays, which \cref{app:henkin} sets beside the argument in the space of Henkin types.

\subsection{The game records noncontextuality as the substitutivity of equality}

\begin{figure}[ht]
\centering
\begin{tikzpicture}[font=\footnotesize]
\node[draw=plum,thick,circle,minimum size=0.9cm] (u) at (0,0) {$u$};
\node[draw=plum,thick,circle,minimum size=0.9cm] (v) at (6,0) {$v$};
\node[plum,below=2pt] at (u.south) {in $\cB$};
\node[plum,below=2pt] at (v.south) {in $\cB'$};
\draw[rose,dashed,thick] (u) -- node[above,rose] {$[u=v]>3/4$} node[below,plum] {norm commitments} (v);
\node[plum] (m) at (3,-2.0) {$\mu(u)=\mu(v)$};
\draw[-{Stealth},plum] (3,-0.55) -- node[right,plum] {(M3)} (m);
\end{tikzpicture}
\caption{Probability transport across a shared ray; the label bounds the truth value of the identity of the two projections.}
\label{fig:probability-transport}
\end{figure}
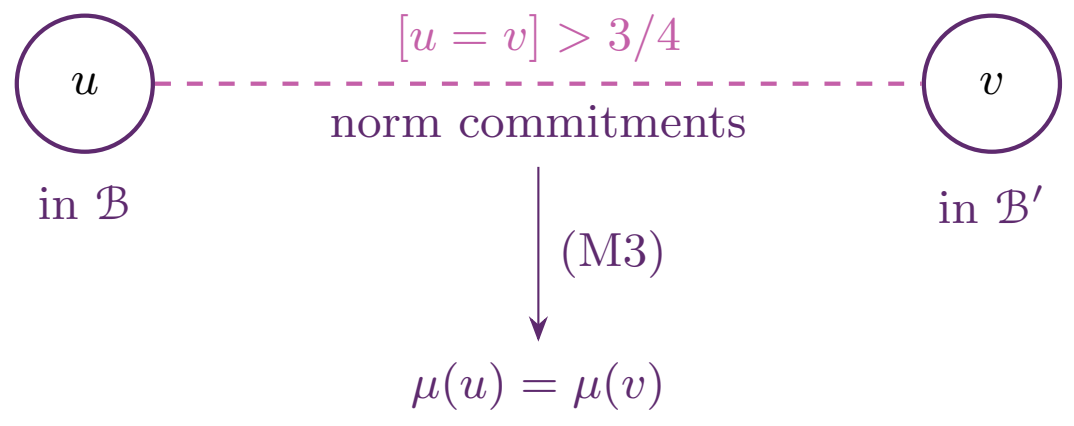

We continue with the plays of the strategy of \cref{sec:bks-game} that reach the parity argument.
The first two bases of \cref{tab:ceg} both contain the vector $(0,0,0,1)$.
If $P$ is the projection onto its ray, the two corresponding outcomes have the same Born probability $\Tr(\rho P)$ for every preparation $\rho$; measurement noncontextuality requires their responses to agree at each underlying state, and outcome determinism requires the common response to be $0$ or $1$.
In a usual Kochen--Specker configuration, a ray is one vertex in all the bases that contain it, so one value for the ray across these bases is part of the data.
A \emph{contextual vertex}, by contrast, indexes an occurrence of a term representing a ray in a specified basis; each ray of the eighteen-ray configuration has two such vertices.
In every model of $T$, the two terms denote elements $a$ and $b$ with $[a=b]=1$, meaning that the corresponding orthogonal projections $\bar P_a$ and $\bar P_b$ of the Hilbert envelope $\bar E$ are equal (see \cref{lem:envelope} and part~(3) of \cref{lem:lattice-truth}).
The elements themselves may be distinct: Keisler's general structures, on which the logic rests, have no equality symbol~\cite[Section~2.1]{Keisler:2023}, and even a standard model may contain distinct elements representing the same projection.
The finite position in \cref{prop:parity-strategy} records positive lower bounds on the formulas expressing these identities.
The corresponding ground instances of (M3), with closed terms substituted for the variables, then force equality of the measure values; a positive antecedent suffices, without requiring the value $1$.
\Cref{fig:probability-transport} illustrates this transfer between contextual vertices.
Its label $[u=v]>3/4$ is a lower bound on the truth value of the formula expressing the identity of the two projections, not a tolerance on an operational equivalence: any positive bound, through (M3), forces exact equality of the two probabilities.

Representing contextual occurrences separately has precedents in quantum foundations.
For the eighteen-ray configuration, Kunjwal and Spekkens~\cite[Fig.~2]{Kunjwal-Spekkens:2015} represent the thirty-six outcome events by distinct nodes and mark the operational equivalences between them.
Their noncontextuality inequality constrains observed statistics under preparation and measurement noncontextuality, and allows responses that are not deterministic.
In the game, the geometric axioms derive the identities of the projections that represent paired events, and (M3) carries the values of a single dispersion-free valuation across these identities.

Other frameworks also represent separately the data supplied by different contexts.
In the contextuality-by-default approach of Dzhafarov and Kujala, random variables with a common content in different contexts are distinct, and noncontextuality asks for a coupling, a jointly distributed family reproducing the distribution in each context, with specified maximal agreement on occurrences with common content~\cite{Kujala-Dzhafarov-Larsson:2015,Dzhafarov-Kujala-Larsson:2015}.
In the sheaf-theoretic framework of Abramsky and Brandenburger, noncontextuality asks for a probability distribution on assignments of outcomes to all measurements whose marginal on each context reproduces the distribution of that context~\cite{Abramsky-Brandenburger:2011}.
Both frameworks pose consistency questions about joint distributions for the empirical data, whereas our certificates record how identities between shared projections constrain the values of a single hypothetical assignment.

\Cref{tab:ks-strategy-contrast} compares the finite obstructions with the usual ray configuration, separating what holds for every winning strategy of Olivia from what belongs to the strategy of \cref{sec:bks-game}.

\begin{table}[htbp]
\centering
\footnotesize
\renewcommand{\arraystretch}{1.15}
\begin{tabular}{@{}p{2.65cm}p{\dimexpr(\linewidth-2.65cm-4\tabcolsep)*5/11\relax}p{\dimexpr(\linewidth-2.65cm-4\tabcolsep)*6/11\relax}@{}}
\toprule
\raggedright  & \raggedright Kochen--Specker configurations & \raggedright finite obstructions of the model existence game for $T^{01}$ \tabularnewline
\midrule
\raggedright framework & \raggedright partial Boolean algebras: in dimension at least three, the projections admit no homomorphism into the two-element Boolean algebra & \raggedright real-valued $L_{\omega_1\omega}$: the theory $T^{01}$ has no models, although each of its finite sets of axioms has one \tabularnewline
\midrule
 & & \raggedright \emph{for every winning strategy of Olivia, who claims that $T^{01}$ has no models} \tabularnewline
\raggedright tree of plays & \raggedright none: one configuration serves against every assignment & \raggedright well-founded, and branching infinitely at some node (\cref{cor:bks-branching}) \tabularnewline
\raggedright certificate of impossibility & \raggedright no assignment of $0$ and $1$ to the rays gives $1$ to exactly one ray of each basis and to no two orthogonal rays & \raggedright at each leaf, a finite position of commitments and function records satisfied in no structure \tabularnewline
\midrule
& & \raggedright \emph{under the strategy of \cref{sec:bks-game}, in a play that lasts until Olivia has carried out the parity argument of \cref{thm:c}} \tabularnewline
\raggedright bases & \raggedright orthonormal bases of the space & \raggedright the nine bases of \cref{tab:ceg} in the span $F$ of a block, each completed by $F^\perp$ \tabularnewline
\raggedright vertex & \raggedright a ray & \raggedright a contextual occurrence of a term representing a ray in one of the nine bases \tabularnewline
\raggedright identifications & \raggedright none needed: a ray is one vertex in all the bases that contain it (the hypothesis of noncontextuality) & \raggedright principally eighteen identifications of rays and nine of the joins of the bases with $F$, all derived from axioms of $T$ that do not mention $\mu$, along which axiom (M3) carries probabilities; each of the eighteen identifications of rays is indispensable to the inconsistency of the Boolean constraint system described below \tabularnewline
\raggedright values & \raggedright $0$ or $1$ at each ray & \raggedright $0$ or $1$ for the probability of each contextual occurrence, imposed by ground instances of (M1) and of the dispersion-free axiom \tabularnewline
\raggedright constraint on a basis & \raggedright exactly one ray of value $1$ & \raggedright sum $1$ of the probabilities of its occurrences, imposed by ground instances of (M1), (M3), (M4) and the dispersion-free axiom once an atom of the block has a recorded value above $15/16$; by \cref{thm:quantitative-parity}, the contradiction persists under tolerances $\delta$ on these values and $\eta$ on these sums with $4\delta+\eta<1$ \tabularnewline
\raggedright location of the rays & \raggedright fixed in advance & \raggedright the eighteen rays of \cref{tab:ceg}, in the block determined by Josh's answers \tabularnewline
\raggedright space & \raggedright finite dimension, spanned by the rays & \raggedright infinite dimension; the rays lie in the span of a block $e_{4j},\dots,e_{4j+3}$ of the orthonormal basis $(e_n)_{n\in\NN}$ named by constants of $T$ \tabularnewline
\raggedright role of $\sigma$-additivity & \raggedright none: every basis is finite; in infinite dimension, the eighteen projections $P_v\otimes1$ on $\CC^4\otimes\ell^2$, described in \cref{sec:localization}, need only finite additivity & \raggedright locates the unit-weight block: Josh names a number $N$ and commits the measure of the join of the rays of $e_0,\dots,e_{N-1}$ to a value greater than $7/8$, and the eighteen rays lie in the block that contains some $e_{n_0}$ with $n_0<N$ \tabularnewline
\bottomrule
\end{tabular}
\caption{The certificates of the game beside the certificate of Kochen and Specker.}
\label{tab:ks-strategy-contrast}
\end{table}

The geometric identifications in \cref{tab:ks-strategy-contrast} come from the terms used to construct each basis.
As in the proof of \cref{thm:c}, each of the nine bases $\cB$ of \cref{tab:ceg}, with vectors $v_0,\dots,v_3$, gives a closed term $s_\cB$ of the sequence sort whose entries $\kappa_is_\cB$ satisfy $\widehat{\kappa_is_\cB}=\tfrac1{16}v_i$ for $i<4$.
Since the subspace $S_{\at(x)}$ is the ray $\CC\hat x$ whenever $\|x\|>0$ (\cref{lem:lattice-truth}), the vertices of $\cB$ are indexed by the terms $\at(\kappa_is_\cB)$, for $i<4$.
Every ray of the table lies in two bases and is therefore represented by two contextual occurrences of terms.

The instances of (K2), together with those of the axioms behind \cref{lem:envelope,lem:lattice-truth}, give $[\at(\kappa_is_\cB)=\at(v_i)]=1$, where $\at(v_i)$ is as in the proof of \cref{thm:c}; hence the two contextual terms of a ray denote elements that are equal in the sense of the formula $[a=b]$.
The same instances give $[c_\cB=a_F]=1$ for the nine elements $c_\cB=\at(\kappa_0s_\cB)\vee\dots\vee\at(\kappa_3s_\cB)$, where the element $a_F$ stands for the span $F$ of the block.
All these identifications follow from axioms that do not mention $\mu$.

Axiom (M3), $[a=b]\wedge|\mu(a)-\mu(b)|\le0$, which every quantum probability measure satisfies, says that $\mu$ respects the equality of the Solèr sort.
This axiom requires the measure to take the same value at the elements denoted by the two contextual terms of a ray.
Noncontextuality, which a configuration of Kochen and Specker builds into the identification of a ray across the bases that contain the ray, therefore enters the certificate exactly as these ground instances of (M3).
Outcome determinism enters separately, through the ground instances of the dispersion-free axiom.

For the following minimality assertion, retain only the Boolean variables attached to the thirty-six contextual occurrences, the nine exactly-one basis constraints and the eighteen equalities between paired occurrences.
Each equality is indispensable to the inconsistency of this constraint system.
Indeed, the nine bases can be arranged in the cycle $1,2,4,3,6,5,9,8,7$ of their numbers in \cref{tab:ceg}, in which any two consecutive bases (the last and the first included) share a ray.
If the identification attached to a ray of a basis $\cB$ is omitted, give the value $1$ to the vertex of the ray in $\cB$.
Next, pair off the other eight bases, which follow one another along the cycle, by four of the rays that consecutive bases share.
Finally, give the value $1$ to both vertices of each of these four rays and the value $0$ to all the other vertices.
Each basis then has exactly one vertex of value $1$.
Without identifications, a vertex of value $1$ could be chosen in each basis independently.
An independent choice is possible in the dense families of orthonormal bases that Clifton and Kent constructed in each finite dimension.
In each of these families, any two bases generate Boolean algebras of projections with only $0$ and the identity in common~\cite{Clifton-Kent:2000}.

\subsection{The parity obstruction survives numerical tolerances}
The interval commitments of the game admit a quantitative stability statement for the parity argument.
For a hypothetical valuation, we attach a value to each occurrence of a ray in a basis and allow these values to be near $0$ or $1$, with each basis sum near $1$.
The following theorem gives sufficient numerical tolerances while keeping the rays fixed and the values at the two occurrences of a shared ray exactly equal.

\begin{keytheorem}[Quantitative parity obstruction]\label{thm:quantitative-parity}
Let $\mathfrak B$ be the family of the nine bases of \cref{tab:ceg}, each regarded as a set of four of the eighteen rays determined by the table.
Let $\delta$ and $\eta$ be nonnegative real numbers with $4\delta+\eta<1$.
There is no family of numbers $(p_{B,r})_{B\in\mathfrak B,\,r\in B}$ in $[0,1]$ satisfying all three conditions:
\begin{enumerate}[label=(\roman*)]
\item $p_{B,r}=p_{B',r}$ whenever $B,B'\in\mathfrak B$ and $r\in B\cap B'$;
\item $p_{B,r}\in[0,\delta]\cup[1-\delta,1]$ for every $B\in\mathfrak B$ and every $r\in B$;
\item $\bigl|\sum_{r\in B}p_{B,r}-1\bigr|\le\eta$ for every $B\in\mathfrak B$.
\end{enumerate}
\end{keytheorem}

\begin{proof}
Since $4\delta+\eta<1$, the intervals $[0,\delta]$ and $[1-\delta,1]$ are disjoint.
Call an occurrence, that is, a pair $(B,r)$ with $B\in\mathfrak B$ and $r\in B$, high when $p_{B,r}$ belongs to $[1-\delta,1]$.
A basis with no high occurrence has sum at most $4\delta<1-\eta$, contrary to (iii).
A basis with at least two high occurrences has sum at least $2-2\delta>1+\eta$, again contrary to (iii).
Thus each of the nine bases has exactly one high occurrence.
Condition (i) gives both occurrences of each ray the same value, and each ray belongs to exactly two bases.
Counting high occurrences basis by basis therefore gives nine, whereas counting them ray by ray gives an even number, a contradiction.
\end{proof}

The bound $4\delta+\eta<1$ is sharp for these three conditions.
If $4\delta+\eta\ge1$ and $\eta\le1$, the family with every $p_{B,r}$ equal to $(1-\eta)/4$ satisfies them, since $(1-\eta)/4\le\delta$ and every basis sum is $1-\eta$; if $\eta\ge1$, the family that vanishes identically satisfies them.
The sharpness concerns the numerical model of \cref{thm:quantitative-parity}, in which the rays and the equalities across bases are exact, and not the noise tolerance of an experiment.

In the finite parity continuation of \cref{prop:parity-strategy}, the three conditions of the theorem appear with $\delta=\eta=0$.
Condition~(i) comes from the ground instances of (M3) and the recorded geometric antecedents.
The dispersion-free ground conditions, together with (M1), give condition~(ii).
Localization has recorded an atom value greater than $15/16$, which these conditions raise to $1$; finite additivity, nonnegativity and the upper bound on the block measure then force each of the nine basis sums to be $1$, which is condition~(iii).

The quantitative theorem isolates the same parity obstruction when endpoint and basis-sum conditions are allowed numerical tolerances: for bases of $k$ rays, with an integer $k\ge2$, the bound $k\delta+\eta<1$ forces exactly one high value per basis.
Any configuration whose exactly-one constraints have no Boolean solution then gives the same contradiction.

The parameters $\delta$ and $\eta$ bound errors in commitments about a hypothetical dispersion-free valuation.
Born probabilities may take intermediate values, so condition (ii) is a restriction on this valuation rather than an estimate of experimental frequency errors.
An experimental noncontextuality inequality would additionally require specified preparation and measurement procedures, operational equivalences and a statistical bound for the relevant class of ontological models.
\Cref{thm:quantitative-parity} thus establishes the quantitative stability of the obstruction to a dispersion-free valuation when the rays are exact; the distinct problem of imperfect measurements has its own analyses~\cite{Simon-Brukner-Zeilinger:2001,Larsson:2002}.

\historypara{Historical note.}
In 1990 Peres gave a proof of the Kochen--Specker theorem by observables on two qubits, valid for a particular state.
Building on Peres's proof, Mermin found a finite obstruction of another kind, formed by nine such observables arranged in three rows and three columns of commuting observables.
No assignment of the values $\pm1$ to these nine observables respects the multiplicative relations along the rows and the columns~\cite{Peres:1990,Mermin:1990}.
The array thus gives a state-independent proof of the theorem.
In 1999 Meyer observed that a coloring of Godsil and Zaks~\cite{Godsil-Zaks:1988} assigns $0$ and $1$ to the unit vectors of $\RR^3$ with rational coordinates, which are dense in the unit sphere.
In this assignment, each orthogonal triple contains exactly one vector of value $1$, and no two orthogonal vectors have the value $1$~\cite{Meyer:1999}.
Kent extended the result in the same year to dense sets of resolutions of the identity in every finite dimension~\cite{Kent:1999}, and the families of Clifton and Kent date from 2000~\cite{Clifton-Kent:2000}.

\section{Remarks}\label{sec:remarks}

The vocabulary and axioms determine how the theory describes subspaces and their measures, and which representation results require standardness.
The following remarks explain the use of projections, the role of (K4) in \cref{thm:c}, and the choices of scalars and logical operations.
They also distinguish the conclusions proper to standard models from the expressive limits of the axioms and record two redundancies among the axioms.
None is needed for the proofs of the main theorems.

\begin{remark}[Projections rather than membership]
The join of two closed subspaces is their closed span, so a description of the join must account for closure.
The vocabulary does this through projections.
A vector $x$ belongs to an element $a$ of the Solèr sort exactly when $\|\tfrac12x-\tfrac12\pi(a,x)\|=0$.
One might instead try to describe membership in $a\vee b$ by asking whether $x$ is a sum $y+z$, with $y$ in $a$ and $z$ in $b$.
In the logic, this would mean that the infimum of a scaled norm of $x-y-z$ vanishes.
The vectors quantified over lie in the unit ball, and this restriction prevents the description from recovering the closed span.
Consider a standard model whose Hilbert sort is the closed unit ball of a Hilbert space $H$.
In this model, the infimum of the scaled norm of $x-y-z$ vanishes exactly at the sums $y+z$ of such members, since these sums form a closed set.
Indeed, bounded sequences in the closed subspaces $S_a$ and $S_b$ spanned by the members of $a$ and of $b$ have weakly convergent subsequences.
When $S_a+S_b$ is not closed, these sums therefore miss the vectors of the unit ball that lie in the join $\overline{S_a+S_b}$ of $S_a$ and $S_b$ in $L(H)$ but not in $S_a+S_b$.
Axiom (Sub2) reaches the join through the projections instead, by von Neumann's alternating projection theorem.
The atoms are reached in the same manner, since axiom (At1) describes $\at(x)$ through the projection onto it: in a standard model of $T_0$, the functions $\at$ and $\rep$ provide the bijection between the rays of $H$, which correspond to the unit vectors up to a phase, and the atoms of the quotient of the Solèr sort by its equality (see \cref{prop:atoms}).
Completeness of this quotient, by contrast, is imposed by no axiom: in the standard models of $T_1$ it is a theorem, which the sequence sort yields through the function $\spa$ (see \cref{thm:t1-iso}).
\end{remark}

\begin{remark}[Finitary definability]\label{q:finitary}
Every closed interval with rational endpoints is definable over $T$ by a formula with countable joins and without quantifiers (see the proof of \cref{lem:borel}), and with quantifiers the definable sets are exactly the $\sigma$-projective sets.
We do not know whether a formula without countable joins can define over $T$ any subset of $[0,1]$ other than $\emptyset$ and $[0,1]$.
An argument through ultraproducts and the continuity of formulas would have to overcome two features of the standard models of $T$.
An ultraproduct of standard models of $T$ need not satisfy axiom (H2), which places every vector in the closed span of the vectors $e_n$;
and the lattice operations are not continuous for the distance between orthogonal projections in operator norm, since the meet of two distinct rays is the zero subspace however close the rays are.
\end{remark}

\begin{remark}[Nonstandard models]\label{rem:nonstandard}
The representation theorems identify the lattice and the measure in standard models.
Their conclusions cannot extend to every model of $T$: the downward Löwenheim--Skolem theorem gives models whose sorts are countable~\cite[Lemma~4]{Aubone-Duenez-Iovino:2026a}.
In such a model, the range of $\mu$ is countable and hence differs from $[0,1]$.
The quotient of the Solèr sort by $[a=b]=1$ is countable as well, so it cannot be isomorphic to the uncountable lattice of closed subspaces of an infinite-dimensional Hilbert space.
Thus standardness is essential to \cref{thm:b} and to part (1) of \cref{thm:a}.

In a nonstandard model of $T_0$, even the conclusion of Gleason's theorem for $T_0$ (see \cref{thm:gleason-t0}) can fail, as the pre-Solèr lattice of \cref{ex:e0} shows.
Recall from this example that $E_0$ is the space of finitely supported complex sequences with the inner product of $\ell^2$, and that $\Lambda_0$ is the lattice formed by the finite-dimensional subspaces of $E_0$ and their orthogonal complements in $E_0$.
Take the closed unit ball of $E_0$ as Hilbert sort, the standard unit vectors as the constants $e_n$ for $n\in\NN$, and the lattice $\Lambda_0$ as Solèr sort.
The function $\pi$ sends $(a,x)$ to the orthogonal projection of $x$ onto $a$, and the function $\at$ sends $x$ to $\CC x$.
Finally, the function $\rep$ assigns to each $a\ne\{0\}$ a unit vector $\rep(a)$ chosen in $a$, and to $\{0\}$ the vector $0$.
With the truncations interpreted by their intended values, the resulting structure is a model of $T_0$ that is not standard, since $E_0$ is not complete.
Passage to closures in $\ell^2$ embeds the ortholattice $\Lambda_0$ in the lattice $L_{\mathrm{fc}}$ of \cref{ex:finite-cofinite}, and each map $\pi(a,\cdot)$ is the restriction to the unit ball of $E_0$ of the orthogonal projection $P_a$ of $\ell^2$ onto the closure of $a$; axioms (Sub1) and (Sub2) therefore hold by \cref{thm:von-neumann}, applied in $\ell^2$, and the other axioms of $T_0$ hold by construction.
Every maximal orthonormal family of $E_0$ is countable, and its rays have join $1$ in $\Lambda_0$.
A $\sigma$-additive state $\nu$ on $\Lambda_0$ would therefore give a function $f$ on the unit vectors of $E_0$, defined by $f(x)=\nu(\CC x)$, with values in $[0,1]$.
Its sum over every maximal orthonormal family would be $1$, and its sum over an orthonormal basis of any finite-dimensional subspace $M$ would be $\nu(M)$, independently of the basis.
Dvurečenskij calls a nonnegative function on the unit vectors of an inner product space a \emph{weak frame function} when its sums over all maximal orthonormal families equal one positive constant and its restriction to the unit vectors of every finite-dimensional subspace has the same sum over all the orthonormal bases of the subspace~\cite{Dvurecenskij:1995}.
The function $f$ would therefore be a weak frame function, and by his theorem an inner product space that carries a weak frame function is complete~\cite[Theorem~3.1]{Dvurecenskij:1995} (compare Hamhalter and Pták~\cite{Hamhalter-Ptak:1987}).
The space $E_0$, however, is not complete.
Hence the lattice $\Lambda_0$ has no $\sigma$-additive state, although every density operator $\rho$ on $\ell^2$ defines a map $a\mapsto\Tr(\rho P_a)$ on $\Lambda_0$.
The identification of $\sigma$-additive states with all such trace maps therefore fails in this nonstandard model.

In standard models, the axioms of \cref{thm:b} exclude every range other than $[0,1]$, one Borel set at a time, or more generally one $\sigma$-projective set at a time (see \cref{thm:definable}).
No sentence, however, and indeed no countable theory, says that the range of $\mu$ is disconnected.
Indeed, suppose that a countable theory $T'$ of real-valued $L_{\omega_1\omega}$ in the vocabulary of $T$ were satisfied, among the models of $T$, exactly by those in which the range of $\mu$ is disconnected; a sentence $\psi$ that took the value $1$ exactly in these models would give such a theory, namely the closed condition $\psi=1$.
Let $\fM$ be a standard model of $T$, which exists by \cref{thm:a}(2).
By \cref{thm:b}, the range of $\mu$ in $\fM$ is the connected set $[0,1]$, which excludes $\fM$ from the models of $T'$.
By the downward Löwenheim--Skolem theorem~\cite[Lemma~4 and its proof]{Aubone-Duenez-Iovino:2026a}, applied to a countable fragment that contains the sentences of the axioms of $T$ and of the conditions of $T'$, the model $\fM$ has a substructure $\fN$ with countable sorts in which every sentence of the fragment keeps its value.
The substructure $\fN$ is therefore a model of $T$ but not of $T'$.
The range of $\mu$ in $\fN$, however, is countable and contains $\mu(0)=0$ and $\mu(1)=1$; it is therefore disconnected, which contradicts the choice of $T'$.

We leave open whether the measure $\mu$ of a nonstandard model of $T$ can differ from every map $a\mapsto\Tr(\rho\bar P_a)$ with $\rho$ a density operator on a complex Hilbert space $\bar E$ as in \cref{lem:envelope}.
Here $\bar P_a$ is the orthogonal projection of $\bar E$ that extends $\pi(a,\cdot)$.

\end{remark}

\begin{remark}[Axiom (K4)]
Unlike the representation theorems, \cref{thm:c} excludes a dispersion-free measure in every model of $T$, standard or not.
Axiom (K4) prescribes the entries of the constant $\varepsilon$ by $\widehat{\kappa_n\varepsilon}=2^{-n-1}\hat e_n$, for $n\in\NN$; together with (H2) and (Sp1), these identities make the subspace represented by $\spa(\varepsilon)$ the whole Hilbert envelope (see \cref{lem:bks-tools}(4)).
Without (K4), a nonstandard model of the remaining axioms may contain no sequence whose entries span its Hilbert envelope.
The structure in the proof of \cref{prop:bks-game}(2), with the zero sequence as $\varepsilon$, satisfies every other axiom of $T$ and has a measure taking only the values $0$ and $1$.
In standard models, Gleason's theorem already excludes such a measure by \cref{thm:a}, whose proof does not use (K4).
\end{remark}
\begin{remark}[Scalars]
The vocabulary fixes the complex numbers as scalars.
For real Hilbert spaces, replace $\QI$ by $\QQ$ and omit the predicate $\operatorname{Im}\ip{x}{y}$.
For quaternionic spaces, use the quaternions with rational coordinates, acting on the left, and four predicates for the real components of the form.
In \cref{prop:atoms}, ``complex number of modulus one'' then becomes ``real number of modulus one'' or ``quaternion of modulus one''.
Each choice yields the analogue of \cref{thm:t1-iso} for the corresponding case of \cref{thm:soler-lattices}, and the proof of \cref{thm:c} applies in all three cases, since the vectors of \cref{tab:ceg} have real coordinates.
\end{remark}

\begin{remark}[The logical form of the axioms]
Every axiom of $T$ is a closed condition of real-valued infinitary logic, built with the connectives $s\cdot\varphi$, $\varphi+\psi$, $\varphi\vee\psi$, the suprema over the sorts and countable joins of families with bounded a priori bounds.
The axioms (H1), (OL), (P), (Sub0), (Sub3), (At1), (At2), (K1), (K2), (K4), (Sp1), (Tr) and (M1)--(M4) involve no countable join, except through the squares of \cref{lem:squares} in (H1), (K1) and (At1).
Axioms (H2), (Sub1), (Sub2), (K3), (Sp2) and (M5) use countable suprema or infima, and (At3) uses the countable join $\lceil\|x\|>0\rceil$.
The hypothesis of the implication in (M3) is the formula $[a=b]$ of \cref{sec:vocabulary}, and that of (Sp2) is the countable conjunction of the formulas $[\at(\kappa_ns)\le b]$.
These formulas are built from the norms of the vectors $\rep(c)$, and they are two-valued in every model of $T_1$, by \cref{lem:lattice-truth} applied to its reduct to the vocabulary of $T_0$; hence so is their conjunction.
In every such model they express the order and the equality of the Solèr sort.
Each axiom (Tr), for given $R$, $n$ and $q$, consists of closed conditions on formulas without countable joins, whose a priori bounds grow with $n$; the formulas $\lceil R(\bar t)>q\rceil$ and $\lceil R(\bar t)<q\rceil$, which gather the truncations for all $n$, are countable joins of members with a priori bound $1$.
The functions $\at$ and $\spa$ take values in the Solèr sort, on which the formula $d(a,b)=\tfrac12(1-[a=b])$ defines, in every model of $T$, a pseudometric with the values $0$ and $1$.
Its topology induces the discrete topology on the quotient by the equality.
On the subspace of rays $\CC\hat x$ with $\hat x\ne0$, represented by the classes of $\at(x)$, this topology differs from the strong operator topology that \cref{sec:definable} places on the projections.
In a standard model, the function $\at$ sends arbitrarily close nonzero, nonproportional vectors to distinct quotient classes, so it is not continuous from the norm topology to the discrete quotient topology.
The function $\spa$ has the same discontinuity along the sequences $\iota_0x$ and $\iota_0y$, since \cref{lem:span} identifies the class of $\spa(\iota_0x)$ with the class of $\at(x)$.
The logic accommodates this discontinuity, since it places no condition on the interpretation of a function symbol beyond its sorts.
This observation concerns the norm topology and the discrete quotient topology just specified.
It does not preclude a pre-metric expansion in Keisler's sense: for a countable vocabulary such an expansion exists with a distance predicate that is a uniform limit of formulas~\cite[Theorem~3.3.4]{Keisler:2023}, and its topology need not be the norm topology.
\end{remark}

\begin{remark}[Redundancy]\label{rem:redundancy}
The axioms (OL) hold in every structure that satisfies the other axioms of $T_0$, since the proof of \cref{lem:lattice-truth}, which embeds the quotient of the Solèr sort into the lattice of closed subspaces of a Hilbert space, does not use them.
In standard models, axiom (Sub2) likewise follows from (Sub1), (Sub3) and (OL).
The two arguments do not allow both groups to be omitted together, because the proof of \cref{lem:lattice-truth} uses (Sub2) to compute joins, whereas the derivation of (Sub2) uses (OL); both groups of axioms are kept nonetheless, so that each property has an axiom of its own.
\end{remark}

\historypara{Historical note.}
Löwenheim proved in 1915 that a first-order sentence with a model has a countable model~\cite{Lowenheim:1915}.
Skolem extended the result to countable sets of sentences in 1920, and in 1922 he observed that axiomatic set theory, in which uncountable sets exist, therefore has countable models~\cite{Skolem:1970}.
The countable models of \cref{rem:nonstandard} are of the same kind: they satisfy every axiom of a theory whose standard models are uncountable.
That an incomplete inner product space can be recognized by the absence of states goes back to Hamhalter and Pták.
In 1987 they proved that a separable inner product space is complete exactly when the lattice of its orthogonally closed subspaces admits a $\sigma$-additive state~\cite{Hamhalter-Ptak:1987}.
The theorem of Dvurečenskij of 1995 used in \cref{rem:nonstandard}, which detects completeness through frame functions rather than through states on the whole lattice of orthogonally closed subspaces, belongs to the same line of work~\cite{Dvurecenskij:1995}.

\appendix

\section{Lattices, hermitian spaces and projective geometry}\label{app:background}

In this appendix we give the geometric definitions and representation results used in \cref{sec:t0,sec:t1}, together with their proofs: the MacNeille completion, hermitian spaces, coordinatization, Solèr's condition and the completion of a pre-Solèr lattice.
These results describe the abstract lattices represented by the vector and event axioms; the model existence game uses only the logic of \cref{sec:logic}.

\subsection{Lattices}\label{sec:lattices}

Birkhoff and von Neumann proposed in 1936 that the propositions about a quantum system be identified with the closed subspaces of its Hilbert space, ordered by inclusion~\cite{Birkhoff-vonNeumann:1936}; the order-theoretic notions below describe this geometry without reference to vectors.
The completion that MacNeille constructed in 1937~\cite{MacNeille:1937} supplies the joins and meets that an incomplete lattice of subspaces lacks; the lattices that the axioms on vectors and events alone describe may be incomplete in this way.

\begin{definition}
For a subset $S$ of a partially ordered set, an \emph{upper bound} of $S$ is an element $b$ with $a\le b$ for all $a\in S$, and the \emph{join} $\bigvee S$ of $S$ is its least upper bound, when it exists; lower bounds and the \emph{meet} $\bigwedge S$ are defined dually.
A \emph{lattice} is a partially ordered set in which every pair of elements $a,b$ has a join $a\vee b$ and a meet $a\wedge b$.
A lattice is \emph{bounded} when it has a least element $0$ and a greatest element $1$, and \emph{complete} when every subset has a join.
\end{definition}

\begin{definition}
An \emph{orthocomplementation} on a bounded lattice $\Lambda$ is a map $a\mapsto a^\perp$ from $\Lambda$ to $\Lambda$ such that, for all $a,b\in\Lambda$, the equations $a^{\perp\perp}=a$, $a\wedge a^\perp=0$ and $a\vee a^\perp=1$ hold, and the inequality $a\le b$ implies $b^\perp\le a^\perp$.
A bounded lattice with an orthocomplementation is an \emph{ortholattice}.
Two elements $a,b$ of an ortholattice are \emph{orthogonal}, written $a\perp b$, when $a\le b^\perp$.
An ortholattice is \emph{orthomodular} when, for all elements $a$ and $b$, the inequality $a\le b$ implies $b=a\vee(b\wedge a^\perp)$.
\end{definition}

The lattice $L(H)$ of closed subspaces of a Hilbert space $H$, ordered by inclusion and equipped with the orthocomplementation $M\mapsto M^\perp$, is an ortholattice.
Its orthomodularity says that every closed subspace $N$ containing a closed subspace $M$ is the orthogonal sum of $M$ and $N\cap M^\perp$.

\begin{definition}\label{def:atoms}
An \emph{atom} of a bounded lattice is a nonzero element $p$ such that every element $b$ with $b\le p$ is equal to $0$ or to $p$.
A bounded lattice is \emph{atomistic} when every element $a$ is the join of the atoms below it (this join being required to exist).
It has the \emph{covering property} when, for every atom $p$ and every element $a$ with $p\not\le a$, no element lies strictly between $a$ and $a\vee p$.
It is \emph{irreducible} when, for any two distinct atoms $p$ and $q$, some atom $r$ distinct from $p$ and $q$ lies below $p\vee q$.
\end{definition}

\begin{definition}
A \emph{sub-ortholattice} of an ortholattice $\Lambda$ is a subset that contains $0$ and $1$ and is closed under the operations $\vee$, $\wedge$ and ${}^\perp$ of $\Lambda$.
An \emph{embedding} of ortholattices is an injective map that preserves $\vee$, $\wedge$, ${}^\perp$, $0$ and $1$.
\end{definition}

Let $\Gamma$ be a complete lattice, and let $\Lambda$ be a subset of $\Gamma$.
We say that $\Lambda$ is \emph{join-dense} in $\Gamma$ when every $g\in\Gamma$ satisfies $g=\bigvee\{a\in\Lambda:a\le g\}$, the join being computed in $\Gamma$, and that $\Lambda$ is \emph{meet-dense} in $\Gamma$ when every $g\in\Gamma$ satisfies $g=\bigwedge\{a\in\Lambda:a\ge g\}$.

For a partially ordered set $\Lambda$ and a subset $X$ of $\Lambda$, let $X^u$ and $X^l$ be the sets of upper and lower bounds of $X$ in $\Lambda$.
We call $X$ a \emph{cut} when $X=(X^u)^l$.
Ordered by inclusion, the cuts form a complete lattice $\DM(\Lambda)$, the \emph{MacNeille completion} of $\Lambda$.
The map $a\mapsto\{b\in\Lambda:b\le a\}$ embeds $\Lambda$ into this lattice, preserving all joins and meets that exist in $\Lambda$.

When $\Lambda$ is an ortholattice, its orthocomplementation induces one on the cuts.
For a cut $X$, let $X^\perp$ be the set $\{b\in\Lambda:b\le a^\perp\text{ for every }a\in X\}$ of lower bounds of $\{a^\perp:a\in X\}$.
We show after the next proposition that the map $X\mapsto X^\perp$ is an orthocomplementation of $\DM(\Lambda)$ that extends that of $\Lambda$; with it we regard $\DM(\Lambda)$ as an ortholattice.

\begin{proposition}[Banaschewski~\cite{Banaschewski:1956}, Schmidt~\cite{Schmidt:1956}]\label{prop:dm-char}
Let $\Lambda$ be a partially ordered set, and let $\Gamma$ be a complete lattice that contains $\Lambda$ as a join-dense and meet-dense subset, the order of $\Gamma$ inducing on $\Lambda$ the order of $\Lambda$.
Then the map that sends each $g\in\Gamma$ to the set $\{a\in\Lambda:a\le g\}$ is an isomorphism of $\Gamma$ onto $\DM(\Lambda)$ that restricts on $\Lambda$ to the canonical embedding $a\mapsto\{b\in\Lambda:b\le a\}$.
If $\Lambda$ is an ortholattice, then the lattice $\Gamma$ has at most one orthocomplementation that extends that of $\Lambda$; any such orthocomplementation is given by $g^\perp=\bigwedge\{a^\perp:a\in\Lambda,\ a\le g\}$, and the isomorphism carries it to the orthocomplementation $X\mapsto X^\perp$ of $\DM(\Lambda)$.
\end{proposition}

The meet determines the order through the equivalence $a\le b$ if and only if $a\wedge b=a$, which we use in \cref{lem:lattice-truth} to recover the order from the lattice operations.
If every subset has a join, its meet is the join of all its lower bounds.
Reversing the order gives the converse, so completeness may equally be stated in terms of meets.

Orthomodular or not, every ortholattice satisfies the De~Morgan law $(\bigvee S)^\perp=\bigwedge\{a^\perp:a\in S\}$ for every subset $S$ whose join exists, since the map $a\mapsto a^\perp$ is an order-reversing bijection that is its own inverse.
It is through this law that orthocomplements of joins are computed below, in the MacNeille completion and in the coordinatization theorem of \cref{sec:projective}.

Let $X$ be a cut of an ortholattice $\Lambda$.
The set $X^\perp$, as a set of lower bounds, is a cut, and the map $X\mapsto X^\perp$ reverses inclusion.
An element $b$ lies in $X^\perp$ exactly when $b^\perp$ is an upper bound of $X$, whence $X^{\perp\perp}=(X^u)^l=X$.
Moreover, the intersection $X\cap X^\perp$ is $\{0\}$, because an element $b$ with $b\le b^\perp$ satisfies $b=b\wedge b^\perp=0$.
Since $X\mapsto X^\perp$ is an order-reversing bijection that is its own inverse, it carries the meet $X\cap X^\perp$ to the join of $X^\perp$ and $X$, which is therefore $\{0\}^\perp=\Lambda$.
Thus the map $X\mapsto X^\perp$ is an orthocomplementation of $\DM(\Lambda)$; it extends that of $\Lambda$, since it sends the cut $\{b\in\Lambda:b\le c\}$ to the cut $\{b\in\Lambda:b\le c^\perp\}$.

\begin{proof}[Proof of \cref{prop:dm-char}]
For $g\in\Gamma$, write $\Lambda_{\le g}=\{a\in\Lambda:a\le g\}$; we show that each $\Lambda_{\le g}$ is a cut and that the map $g\mapsto\Lambda_{\le g}$ is an order embedding whose image contains every cut.
The assertions on orthocomplementations then follow from the De~Morgan law.

By join-density, we have $g=\bigvee\Lambda_{\le g}$; hence the upper bounds of $\Lambda_{\le g}$ in $\Lambda$ are the elements $b\in\Lambda$ with $b\ge g$.
By meet-density, the common lower bounds of these elements are the elements $a\in\Lambda$ with $a\le g$; thus the set $\Lambda_{\le g}$ is a cut.
For $g,g'\in\Gamma$, join-density also shows that $g\le g'$ if and only if $\Lambda_{\le g}\subseteq\Lambda_{\le g'}$; hence the map is an order embedding.
Because the order of $\Gamma$ induces that of $\Lambda$, the map restricts on $\Lambda$ to the canonical embedding.

If $X$ is a cut and $g=\bigvee X$ in $\Gamma$, then the cut $X$ is contained in $\Lambda_{\le g}$.
Conversely, an element $a\in\Lambda$ with $a\le g$ lies below every $b\in X^u$ because each such $b$ is an upper bound of $X$ in $\Gamma$, whence $a\in X=(X^u)^l$.
Thus $X=\Lambda_{\le g}$, and the map is onto $\DM(\Lambda)$.

The formula for $g^\perp$ follows from $g=\bigvee\Lambda_{\le g}$ and the De~Morgan law.
By this formula, an element $b\in\Lambda$ lies below $g^\perp$ exactly when $b\le a^\perp$ for every $a\in\Lambda_{\le g}$; hence $\Lambda_{\le g^\perp}=(\Lambda_{\le g})^\perp$.
\end{proof}

\historypara{Historical note.}
Birkhoff and von Neumann rejected the lattice of closed subspaces of an infinite-dimensional Hilbert space as the logic of quantum mechanics.
Letters that von Neumann wrote to Birkhoff in 1935, while the two were preparing their paper of 1936, show why: the lattice is not modular.
Von Neumann, however, wanted the probabilities of quantum propositions to behave like a dimension function, a real-valued function $d$ on the lattice that is strictly increasing and satisfies $d(a\vee b)+d(a\wedge b)=d(a)+d(b)$.
Only a modular lattice carries such a function~\cite{Redei:2007}.
In the paper itself, Birkhoff and von Neumann postulated instead an abstract projective geometry, whose lattice of subspaces is modular~\cite{Birkhoff-vonNeumann:1936}.
A year later Husimi observed that the closed subspaces of a Hilbert space satisfy the orthomodular law~\cite{Husimi:1937}, which was rediscovered several times in the 1950s and 1960s.
In his thesis of 1964, Piron took completeness, orthomodularity, atomicity and the covering property as axioms for the propositions about a quantum system.
He showed that the irreducible lattices with these properties are isomorphic to lattices of orthogonally closed subspaces of hermitian spaces when their dimension is at least four~\cite{Piron:1964}.

\subsection{Hermitian spaces}\label{sec:hermitian}

A one-dimensional subspace of a vector space is called a \emph{ray}.
By the coordinatization theorem of projective geometry, every projective geometry with four independent points is isomorphic to the geometry of the rays of a vector space over a division ring~\cite{Baer:1952}.
The theorem has an orthocomplemented form, which represents an ortholattice with the properties listed below by the orthogonally closed subspaces of a hermitian space (\cref{thm:coord}); the division ring that it yields need not be $\RR$, $\CC$ or $\HH$.

\begin{definition}
A \emph{division ring} is a ring with $1\ne0$ in which every nonzero element is invertible; multiplication need not be commutative.
An \emph{involution} of a division ring $\Delta$ is a map $\alpha\mapsto\alpha\inv$ such that $(\alpha+\beta)\inv=\alpha\inv+\beta\inv$, $(\alpha\beta)\inv=\beta\inv\alpha\inv$ and $\alpha^{**}=\alpha$ for all $\alpha,\beta\in\Delta$.
Fix an involution $\alpha\mapsto\alpha\inv$ of $\Delta$, and let $E$ be a left vector space over $\Delta$, that is, a left $\Delta$-module.\footnote{Physicists usually let quaternionic scalars act on the right; their quaternionic Hilbert spaces are therefore right modules. The choice of side is immaterial: a left vector space over $\Delta$ is a right vector space over the opposite division ring, and an involution is an isomorphism of $\Delta$ onto its opposite.}
A \emph{hermitian form} on $E$ is a map $\ip{\cdot}{\cdot}\colon E\times E\to\Delta$ such that $\ip{\alpha x+\beta y}{z}=\alpha\ip{x}{z}+\beta\ip{y}{z}$ and $\ip{y}{x}=\ip{x}{y}\inv$ for all $x,y,z\in E$ and $\alpha,\beta\in\Delta$; it is \emph{anisotropic} when $\ip{x}{x}\ne0$ for every nonzero vector $x\in E$.
A \emph{hermitian space} is a left vector space over a division ring with an involution, together with an anisotropic hermitian form.
\end{definition}

\begin{definition}
In a hermitian space $E$, two vectors $x$ and $y$ are \emph{orthogonal}, written $x\perp y$, when $\ip{x}{y}=0$.
For $S\subseteq E$, let $S^\perp$ be the set $\{x\in E:\ip{x}{y}=0\text{ for all }y\in S\}$.
A subspace $M$ is \emph{orthogonally closed} when $M^{\perp\perp}=M$.
\end{definition}

\begin{definition}[The lattice $\cC(E)$]
For a hermitian space $E$, we write $\cC(E)$ for the set of orthogonally closed subspaces of $E$, ordered by inclusion.
\end{definition}

The set $\cC(E)$, with the map $M\mapsto M^\perp$, is a complete ortholattice, as we verify below.

\begin{definition}
A hermitian space $E$ is \emph{orthomodular} when $M+M^\perp=E$ for every $M\in\cC(E)$, where $M+M^\perp$ is the set of sums $v+y$ with $v\in M$ and $y\in M^\perp$.
\end{definition}

\begin{definition}
A \emph{Hilbert space} is a hermitian space $H$ over the real numbers $\RR$, the complex numbers $\CC$ or the quaternions $\HH$, with the usual involution (the identity, complex conjugation or quaternionic conjugation), such that $\ip{x}{x}>0$ for every nonzero vector $x$ and the space is complete in the norm $\|x\|=\ip{x}{x}^{1/2}$.
\end{definition}

\begin{lemma}\label{lem:plus-one}
In a hermitian space $E$ over a division ring $\Delta$, the subspace $M+\Delta x$ belongs to $\cC(E)$ for every $M\in\cC(E)$ and every $x\in E$.
\end{lemma}

It follows from the identities that define a hermitian form that $\ip{x}{\alpha y}=\ip{x}{y}\alpha\inv$ for all $x,y\in E$ and $\alpha\in\Delta$.
The bracket is linear in its first argument, as is customary in mathematics; in the Dirac notation used by physicists it is linear in the second.

Nothing in the definition of a hermitian space forces the scalars to be real, complex or quaternionic: over $\QQ$ with the identity involution, the form $\sum_nx_ny_n$ on finitely supported rational sequences $x=(x_n)_{n\in\NN}$ and $y=(y_n)_{n\in\NN}$ is anisotropic, since a sum of rational squares vanishes only when every square does.

Orthogonality is symmetric, since $\ip{y}{x}=\ip{x}{y}\inv$ and the involution fixes $0$.
In particular, every subset $S$ of $E$ satisfies $S\subseteq S^{\perp\perp}$.

Every set of the form $S^\perp$ is an orthogonally closed subspace: applying the inclusion-reversing operation ${}^\perp$ to $S\subseteq S^{\perp\perp}$ gives $S^{\perp\perp\perp}\subseteq S^\perp$, and the reverse inclusion is the same fact for $S^\perp$.

If $(M_j)_{j\in J}$ is a family of orthogonally closed subspaces indexed by an arbitrary set $J$, then their intersection, which equals $\bigcap_jM_j^{\perp\perp}=\bigl(\bigcup_jM_j^\perp\bigr)^\perp$, is orthogonally closed.
Hence the set $\cC(E)$ is a complete lattice with least element $\{0\}$ and greatest element $E$.
In this lattice the family has meet $\bigcap_jM_j$ and join $\bigl(\sum_jM_j\bigr)^{\perp\perp}$, which may be strictly larger than $\sum_jM_j$.

For $M\in\cC(E)$, anisotropy gives $M\cap M^\perp=\{0\}$, whence $M\vee M^\perp=(M^\perp\cap M)^\perp=E$.
Since the map $M\mapsto M^\perp$ also reverses inclusion and is its own inverse on $\cC(E)$, the lattice $\cC(E)$ is, with this map, a complete ortholattice.

In a Hilbert space $H$, the projection theorem supplies the decomposition that defines orthomodular spaces above: every closed subspace $M$ satisfies $M+M^\perp=H$.
For every subspace $M$, it also gives $M^{\perp\perp}=\overline M$.
Thus orthogonal closure is ordinary norm closure, and $L(H)=\cC(H)$.
The space $H$ is orthomodular, and \cref{prop:om-equiv} below shows that its lattice is orthomodular as well.

\begin{lemma}\label{lem:finite}
In a hermitian space $E$, every finite-dimensional subspace $F$ is orthogonally closed and satisfies $F+F^\perp=E$.
\end{lemma}

\begin{proof}
The idea is to write down the orthogonal projection onto $F$.
Let $n$ be the dimension of $F$.
By anisotropy, the Gram--Schmidt process gives $F$ an orthogonal basis $w_1,\dots,w_n$ with $\ip{w_i}{w_i}\ne0$ for $i=1,\dots,n$.
For $z\in E$, the vector $Pz=\sum_i\ip{z}{w_i}\ip{w_i}{w_i}^{-1}w_i$ lies in $F$, and the difference $z-Pz$ lies in $F^\perp$; hence $F+F^\perp=E$.
If $y\in F^{\perp\perp}$, the vector $y-Py$ lies in $F^\perp\cap F^{\perp\perp}$ and is therefore orthogonal to itself; by anisotropy it is $0$, whence $y\in F$.
\end{proof}

\begin{proof}[Proof of \cref{lem:plus-one}]
If $x\in M$, there is nothing to prove; we therefore assume that $x\notin M$.
Since $M+\Delta x$ is contained in $(M+\Delta x)^{\perp\perp}$, it suffices to show that every vector of $(M+\Delta x)^{\perp\perp}$ differs from a multiple of $x$ by a vector of $M$.
We compare the functional that such a vector induces on $M^\perp$ through the form with the one induced by $x$.
Their kernels will show that they differ by multiplication on the right by a scalar.
Write $W=M^\perp$, and let $\xi_x$ be the functional $z\mapsto\ip{z}{x}$ on $W$.
This functional is nonzero, since $x\notin W^\perp=M$, and we may choose $z_0\in W$ with $\xi_x(z_0)=1$.

Let $y$ be an element of $(M+\Delta x)^{\perp\perp}$, and let $\xi_y$ be the functional $z\mapsto\ip{z}{y}$ on $W$.
Since $W\cap x^\perp=(M+\Delta x)^\perp$, the functional $\xi_y$ vanishes on the kernel of $\xi_x$; as $z-\xi_x(z)z_0$ lies in this kernel for every $z\in W$, we obtain $\xi_y(z)=\xi_x(z)\alpha$ with $\alpha=\xi_y(z_0)$.
Because $\ip{z}{\alpha\inv x}=\ip{z}{x}\alpha$, the vector $y-\alpha\inv x$ is orthogonal to $W$ and therefore lies in $W^\perp=M$, whence $y\in M+\Delta x$.
\end{proof}

\begin{proposition}\label{prop:om-equiv}
For a hermitian space $E$, the lattice $\cC(E)$ is orthomodular if and only if $E$ is an orthomodular space.
\end{proposition}

\begin{proof}
Suppose first that $E$ is orthomodular, and let $M\subseteq M'$ be elements of $\cC(E)$.
Every vector $z\in M'$ can be written as $v+y$ with $v\in M$ and $y\in M^\perp$; since $y=z-v$ lies in $M'$, we have $M'=M+(M'\cap M^\perp)\subseteq M\vee(M'\wedge M^\perp)\subseteq M'$, which is the orthomodular law.

Conversely, suppose that $\cC(E)$ is orthomodular, and let $\Delta$ be the division ring of scalars of $E$.
Let $M$ be an element of $\cC(E)$, and let $x$ be a vector of $E\setminus M$.
By \cref{lem:plus-one}, the subspace $M'=M+\Delta x$ belongs to $\cC(E)$, and orthomodularity gives $M'=M\vee(M'\cap M^\perp)$.
Since $M'\ne M$, the subspace $M'\cap M^\perp$ is nonzero.
As this subspace meets $M$ only in $0$, and as $M$ has codimension one in $M'$, the subspace $M'\cap M^\perp$ is one-dimensional, spanned by a vector $u$.
Then $M\subsetneq M+\Delta u\subseteq M'$, and the codimension count gives $M'=M+\Delta u$; hence $x\in M+M^\perp$.
\end{proof}

\historypara{Historical note.}
In 1964 Piron published a theorem asserting that every orthomodular complex inner product space is complete, but his proof contained an error~\cite{Piron:1964}.
Amemiya and Araki gave a correct proof in 1966, and their method also covers real and quaternionic scalars~\cite{Amemiya-Araki:1966}.
It was then speculated that, among infinite-dimensional hermitian spaces over arbitrary division rings, orthomodularity alone characterizes the Hilbert spaces over $\RR$, $\CC$ and $\HH$.
Keller ended the speculation in 1980 with a ``non-classical Hilbert space'', an infinite-dimensional orthomodular hermitian space over a division ring other than these three~\cite{Keller:1980}.
Keller was a student of Herbert Gross, whose papers and book on forms in infinite-dimensional vector spaces shaped the theory; the story is told in Holland's survey of 1995~\cite{Holland:1995}.

\subsection{Projective geometry, Solèr lattices and pre-Solèr lattices}\label{sec:projective}

The lattice of closed subspaces remembers more than inclusion: its atoms are rays, their joins describe incidence, and orthocomplementation records orthogonality.
The coordinatization theorem recovers a vector space and a hermitian form from this information.
It gives an embedding into the lattice of orthogonally closed subspaces; completeness is needed to make the embedding onto.
Solèr's condition then supplies the infinite orthonormal sequence that, together with orthomodularity of the representing space, forces real, complex or quaternionic Hilbert geometry.

We develop these steps to distinguish the Solèr lattices from the pre-Solèr lattices, which may lack completeness.
In the axiomatization of \cref{sec:axioms}, the Hilbert sort provides vectors on which lattice elements act as orthogonal projections, while the sequence sort supplies all closed spans in standard models (see \cref{thm:t1-iso}).
We assume familiarity with Hilbert spaces but none with projective geometry, whose incidence axioms we recall before recovering the vectors from an abstract lattice.

\subsubsection{Rays and coordinates}

In a complex Hilbert space $H$, the atoms of $L(H)$ are the rays $\CC x$ spanned by nonzero vectors $x$.
Two distinct rays have as their join the plane they span, and the rays in the plane are exactly the atoms below the join.
Thus the lattice records which rays lie in which planes without naming vectors or scalars.

To recover vectors from an abstract lattice, we first describe incidence without assuming a vector space.
Coordinatization will identify this abstract incidence structure with the geometry of rays of a vector space.

\begin{definition}
A \emph{projective geometry} is a set of \emph{points} together with a family of subsets, called \emph{lines}, such that:
\begin{enumerate}[label=(PG\arabic*)]
\item any two distinct points lie on exactly one line;
\item every line contains at least three points;
\item (\emph{the Veblen--Young axiom}) if $a,b,a',b'$ are distinct points and the line through $a$ and $b$ meets the line through $a'$ and $b'$, then the line through $a$ and $a'$ meets the line through $b$ and $b'$.
\end{enumerate}
\end{definition}

A set of points is a \emph{projective subspace} when it contains the whole line through any two of its points.
The projective subspace \emph{generated} by a set of points is the intersection of all projective subspaces containing it.
A \emph{projective plane} is the projective subspace generated by three points not on one line.
Axiom (PG3) says that two lines lying in a common projective plane always meet; a projective geometry has no parallel lines.
The word \emph{projective} comes from perspective drawing, in which parallel rails, projected onto the canvas, meet at the horizon.

The basic example comes from linear algebra.
Let $V$ be a left vector space of dimension at least two over a division ring $\Delta$.
Take as points the rays of $V$, and as lines the sets of rays contained in a two-dimensional subspace.
Axiom (PG1) holds because two distinct rays span a unique two-dimensional subspace.
Axiom (PG2) holds because the two-dimensional subspace spanned by linearly independent vectors $x$ and $y$ contains the three distinct rays $\Delta x$, $\Delta y$ and $\Delta(x+y)$.
For (PG3), suppose that the two-dimensional subspaces spanned by $a,b$ and by $a',b'$ have a ray in common.
The rays $a,b,a',b'$ therefore lie in a subspace of dimension at most three.
By counting dimensions, any two two-dimensional subspaces of this subspace have a ray in common, in particular those spanned by $a,a'$ and by $b,b'$.

Points are \emph{independent} when none of them lies in the projective subspace generated by the others.
In the basic example, rays are independent exactly when nonzero vectors spanning them, one in each ray, are linearly independent.
Two projective geometries are \emph{isomorphic} when there is a bijection between their sets of points under which the lines of the one correspond exactly to the lines of the other.
For vector spaces $E$ and $E'$ over division rings $\Delta$ and $\Delta'$, respectively, a \emph{semilinear isomorphism} from $E$ onto $E'$ is an additive bijection $f\colon E\to E'$ for which there is an isomorphism $\sigma\colon\Delta\to\Delta'$ with $f(\alpha x)=\sigma(\alpha)f(x)$ for all $\alpha\in\Delta$ and $x\in E$.
The following theorem says that the basic example is, up to isomorphism, the only one as soon as there are four independent points.

\begin{theorem}[Coordinatization of projective geometries; see Baer~{\cite{Baer:1952}}]\label{thm:fundamental}
Every projective geometry that contains four independent points is isomorphic to the geometry of rays of a left vector space $E$ over a division ring $\Delta$.
The pair $(\Delta,E)$ is unique up to a semilinear isomorphism.
\end{theorem}

The division ring $\Delta$ is thus recovered from the geometry alone, whatever it may be: every division ring occurs.
The hypothesis of four independent points, on the other hand, is indispensable: a projective plane regarded as a projective geometry in its own right need not come from a division ring.
The planes that do are exactly those in which Desargues's theorem holds (if two triangles are in perspective from a point, the three intersection points of corresponding sides lie on one line).
Planes in which Desargues's theorem fails, the \emph{non-Desarguesian} planes, exist; in a geometry with four independent points, by contrast, the theorem follows from the axioms.

\subsubsection{Recovering a hermitian space from the lattice}

The following properties of an ortholattice $\Lambda$ enter the representation theorems:
\begin{enumerate}[label=(\textup{\Alph*}),align=left,leftmargin=*]
\item[(C)] the lattice is complete;
\item[(OM)] it is orthomodular;
\item[(A)] it is atomistic;
\item[(CP)] it has the covering property;
\item[(I)] it is irreducible;
\item[(D)] it contains four mutually orthogonal atoms.
\end{enumerate}

\begin{theorem}[Coordinatization; Baer~\cite{Baer:1952}, Maeda and Maeda~{\cite[\S34]{Maeda-Maeda:1970}}, Holland~{\cite[\S3]{Holland:1995}}]\label{thm:coord}
Let $\Lambda$ be an ortholattice that satisfies (OM), (A), (CP), (I) and (D), but not necessarily (C).
Then there are a division ring $\Delta$ with an involution, a hermitian space $E$ over $\Delta$, and a bijection $p\mapsto \hat p$ from the atoms of $\Lambda$ onto the rays of $E$, such that two atoms $p$ and $q$ are orthogonal exactly when $\hat p\perp \hat q$.
The map
\[
\gamma\colon\Lambda\to\cC(E),\qquad \gamma(a)=\{0\}\cup\bigcup\{\hat p: p\text{ an atom with } p\le a\},
\]
is an embedding of ortholattices whose image contains every finite-dimensional subspace of $E$ and the orthogonal complement of each.
The triple formed by $\Delta$, $E$ and the form is unique up to a semilinear isomorphism followed by a rescaling $\ip{x}{y}'=\ip{x}{y}\zeta$, where $\zeta\inv=\zeta\ne0$, with the new involution $\alpha\mapsto\zeta^{-1}\alpha\inv\zeta$.
\end{theorem}

To apply \cref{thm:fundamental} to an orthomodular lattice that is atomistic and has the covering property, we take the atoms of the lattice as the points of a geometry.
The lines of this geometry are the sets $\{r: r\text{ an atom with } r\le p\vee q\}$ for distinct atoms $p$ and $q$; in $L(H)$, these sets consist of the rays in the plane spanned by $p$ and $q$.
The covering property gives this geometry (PG1) and, with orthomodularity, (PG3); irreducibility in the sense of \cref{def:atoms} is exactly (PG2) (see Maeda and Maeda~\cite[\S\S15--16]{Maeda-Maeda:1970} and Holland~\cite[\S4]{Holland:1995}).

Orthocomplementation adds a relation between points: by the definition of orthogonality, atoms $p$ and $q$ are orthogonal when $p\le q^\perp$.
In $L(H)$ this is orthogonality of rays, and it determines the inner product up to the ambiguity described at the end of \cref{thm:coord}.
Condition (D) supplies the four independent points of \cref{thm:fundamental}, since mutually orthogonal atoms are independent.
The coordinatization theorem, \cref{thm:coord}, goes back, in finite dimensions, to Birkhoff and von Neumann~\cite{Birkhoff-vonNeumann:1936}.
The theorem carries the atoms onto all rays of a hermitian space and preserves their orthogonality, while embedding the lattice into the orthogonally closed subspaces of the space.
Condition (C) is absent from its hypotheses, so the embedding applies to the incomplete lattices of \cref{sec:t0}; surjectivity requires a separate argument, given in the proof of \cref{thm:soler-lattices}.

\begin{proof}[Proof of the statements of \cref{thm:coord} about $\gamma$]
By \cref{thm:fundamental}, incidence among the atoms recovers a vector space $E$ whose rays correspond to all the atoms.
Orthocomplementation induces a \emph{polarity} in the sense of Holland~\cite[Definition~3.2]{Holland:1995}, sending a projective subspace to the points orthogonal to all its points.
It reverses inclusion and sends the whole geometry to the empty set.
The polar of a point $p$ consists of the atoms below $p^\perp$.
It is a hyperplane, or maximal proper projective subspace, whose polar is $\{p\}$~\cite[\S4]{Holland:1995}.

The projective-dual construction for Baer's representation~\cite[Chapter~IV]{Baer:1952} gives a form $B$ inducing the polarity, linear in its first variable and $g$-semilinear in its second for an anti-automorphism $g$ of $\Delta$ (see Holland~\cite[proof of Theorem~3.6]{Holland:1995}).
Thus $B(x,y)=0$ for nonzero vectors $x,y\in E$ exactly when the corresponding atoms are orthogonal.
The form is reflexive, meaning that $B(x,y)=0$ if and only if $B(y,x)=0$, and it is anisotropic because $p\wedge p^\perp=0$ for every atom $p$.
To obtain a hermitian form, choose $e\ne0$, let $d=B(e,e)$, and define
\[
 B_0(x,y)=B(x,y)d^{-1},\qquad \tau(\alpha)=d\,g(\alpha)d^{-1}.
\]
The forms $B_0$ and $B$ have the same zero set, so $B_0$ induces the same polarity.
Moreover, $B_0$ is hermitian for the involution $\tau$, in every characteristic.\footnote{For $y\ne0$, the left-linear functionals $L_y(x)=B_0(x,y)$ and $H_y(x)=\tau^{-1}(B_0(y,x))$ have the same kernel, so $H_y=L_y\lambda_y$ for a right scalar $\lambda_y$.
For independent $y,z$, the functionals $L_y$ and $L_z$ are independent over the scalars acting on the right, since a relation between them would, by the nondegeneracy of $B_0$, give a relation between $y$ and $z$; comparing $H_{y+z}=H_y+H_z=L_y\lambda_y+L_z\lambda_z$ with $H_{y+z}=L_{y+z}\lambda_{y+z}=L_y\lambda_{y+z}+L_z\lambda_{y+z}$ therefore gives $\lambda_y=\lambda_z=\lambda_{y+z}$.
When $y$ and $z$ are dependent, an auxiliary vector independent of both gives the same equality, so $\lambda_y$ is constant for $y\ne0$.
Since $B_0(e,e)=1$, this constant is $1$.
Hence $B_0(y,x)=\tau(B_0(x,y))$; applying this twice and using $B_0(e,\alpha e)=\tau(\alpha)$ gives $\tau^2=\mathrm{id}$.}
We verify that the resulting map $\gamma$ represents the whole lattice and preserves its operations.
If two atoms lie below an element $a$ of $\Lambda$, every atom on the line they determine lies below $a$; hence the set $\gamma(a)$ is a subspace.
For an atom $r$ the following are equivalent: $r\le a^\perp$; $a\le r^\perp$; $p\le r^\perp$ for every atom $p\le a$ (because (A) makes $a$ the join of these atoms); $\hat r\perp \hat p$ for every atom $p\le a$.
Thus the map $\gamma$ carries $a^\perp$ to $\gamma(a)^\perp$, and applying this twice gives $\gamma(a)=\gamma(a)^{\perp\perp}$.
By (A), elements $a,b\in\Lambda$ satisfy $a\le b$ if and only if $\gamma(a)\subseteq\gamma(b)$, which makes $\gamma$ injective.
Meets are preserved because an atom lies below $a\wedge b$ exactly when it lies below $a$ and $b$, and joins are preserved by the De~Morgan law.
Finally, let $F$ be a finite-dimensional subspace of $E$ with basis $x_1,\dots,x_k$, and let $p_i$ be the atom with $\hat p_i=\Delta x_i$, for $i=1,\dots,k$.
The join $p_1\vee\dots\vee p_k$ exists, as every finite join does in a lattice.
Because $\gamma$ preserves joins, the image of this join is $(\Delta x_1+\dots+\Delta x_k)^{\perp\perp}=F^{\perp\perp}$, which is $F$ by \cref{lem:finite}.
The orthogonal complement $F^\perp$ is therefore $\gamma\bigl((p_1\vee\dots\vee p_k)^\perp\bigr)$.
\end{proof}

Piron~\cite{Piron:1964} proved this representation for complete lattices; completeness is not needed for the embedding, and in the proof of \cref{thm:soler-lattices} below it serves only to make $\gamma$ surjective.

\subsubsection{Harmonic conjugates and Solèr's condition}\label{sec:soler-pre-soler}

Solèr proved in 1995 that an orthomodular space containing an infinite orthonormal sequence is a Hilbert space over $\RR$, $\CC$ or $\HH$~\cite{Soler:1995}, and Holland restated the existence of such a sequence as a condition on the lattice, through harmonic conjugates, a construction of classical projective geometry that uses joins and meets alone~\cite{Holland:1995}.
This condition, together with completeness and the conditions of the coordinatization theorem, singles out the lattices of closed subspaces of infinite-dimensional Hilbert spaces (\cref{thm:soler-lattices}).

Let $E$ be a left vector space of dimension at least three over a division ring $\Delta$.
For distinct rays $p$ and $q$ of $E$, regarded as projective points, write $p+q$ for their projective line.
Fix a point $c$ on this line other than $p$ and $q$.
Choose a point $r$ not on $p+q$.
All intersections below lie in the ray geometry of the three-dimensional subspace spanned by $p,q,r$, where two distinct projective lines meet in one point; we denote the point by their intersection.
Choose a point $r'$ on $r+p$ other than $r$ and $p$; let $c'=(c+r')\cap(q+r)$ and $c''=(p+c')\cap(q+r')$.
The point $h=(r+c'')\cap(p+q)$ is the \emph{harmonic conjugate} of $c$ with respect to $p$ and $q$.
The point $h$ does not depend on the choices of $r$ and $r'$, and, for vectors $u$ and $v$ spanning $p$ and $q$, respectively,
\begin{equation}\label{eq:harmonic}
\text{the harmonic conjugate of }\Delta(u+v)\text{ with respect to }\Delta u\text{ and }\Delta v\text{ is }\Delta(u-v).
\end{equation}
The bijection of \cref{thm:coord} between atoms and rays transfers the harmonic conjugate to the atoms of an ortholattice that satisfies the hypotheses of that theorem.

\emph{Solèr's condition}, in the form we use, asks for an infinite orthogonal sequence of atoms along which the orthocomplement and the harmonic conjugate agree at every step:
\begin{enumerate}[label=(S),leftmargin=*]
\item There are mutually orthogonal atoms $p_n$, for $n\in\NN$, and atoms $c_n\le p_n\vee p_{n+1}$ with $c_n\ne p_n,p_{n+1}$, such that, for every $n\in\NN$, the harmonic conjugate of $c_n$ with respect to $p_n$ and $p_{n+1}$ is $c_n^\perp\wedge(p_n\vee p_{n+1})$.
\end{enumerate}

\begin{definition}
A \emph{Solèr lattice} is an ortholattice that satisfies (C), (OM), (A), (CP), (I), (D) and (S).
A \emph{pre-Solèr lattice} is an ortholattice that satisfies (OM), (A), (CP), (I), (D) and (S), but not necessarily (C).
\end{definition}

\begin{theorem}[Piron, Amemiya and Araki, Solèr; in this form Holland~{\cite[Theorem~4.1]{Holland:1995}}]\label{thm:soler-lattices}
An ortholattice is a Solèr lattice if and only if it is isomorphic to $L(H)$ for some infinite-dimensional Hilbert space $H$ over $\RR$, $\CC$ or $\HH$.
\end{theorem}

Coordinatization leaves the division ring unrestricted.
For example, the lattice of all subspaces of $\QQ^4$, orthocomplemented by the form $\sum_{i<4}x_iy_i$, satisfies the hypotheses of \cref{thm:coord} with $\Delta=\QQ$.
To reach Hilbert spaces, we need a further condition on the lattice that supplies an infinite orthonormal sequence in its representing space.

Solèr's condition does this through harmonic conjugates.
For two orthogonal rays, the harmonic conjugate of a third ray in their plane agrees with its orthocomplement in the plane exactly when the third ray is spanned by the sum of two vectors of equal ``length'', one from each orthogonal ray.
Applying this requirement along an infinite orthogonal sequence permits a single rescaling of the form that makes the sequence orthonormal (see \cref{lem:s-orthonormal}).
Completeness of the lattice will supply orthomodularity of the whole representing space, the other hypothesis needed for Solèr's theorem.

In the vector space $E$, the construction of the harmonic conjugate has a simple coordinate description, which shows that the harmonic conjugate does not depend on the choices of the auxiliary points $r$ and $r'$.
Let $u$ and $v$ be vectors spanning $p$ and $q$, respectively, scaled so that $c=\Delta(u+v)$; this is possible because $c$ lies on $p+q$ and differs from $p$ and $q$.
Likewise, let $w$ be a vector spanning $r$, scaled so that $r'=\Delta(u+w)$, which is possible because $r'$ lies on $r+p$ and differs from $r$ and $p$.
A direct computation gives $c'=\Delta(v-w)$, $c''=\Delta(u-v+w)$ and $h=\Delta(u-v)$.
All coefficients are $\pm1$, so the computation is valid over every division ring, which proves~\eqref{eq:harmonic}.

Now let $\Lambda$ be an ortholattice satisfying the hypotheses of \cref{thm:coord}, and apply the construction in its coordinatizing space $E$.
The bijection between atoms and rays transfers the harmonic conjugate to the atoms, with $\vee$ in place of $+$ and $\wedge$ in place of $\cap$.
The coordinate calculation proves independence of the auxiliary choices, so the relation ``$h$ is the harmonic conjugate of $c$ with respect to $p$ and $q$'' is expressed by a first-order formula in the language of lattices.

\begin{figure}[ht]
\centering
\begin{tikzpicture}[scale=0.9]
\begin{scope}
  \draw[black!40,->] (-2,0) -- (2.3,0); \draw[black!40,->] (0,-2) -- (0,2.3);
  \draw[plum, very thick, -{Stealth}] (0,0) -- (1,0) node[below right] {$u$};
  \draw[plum, very thick, -{Stealth}] (0,0) -- (0,1) node[above left] {$v$};
  \draw[rose, thick] (-1.6,-1.6) -- (1.6,1.6) node[above right] {$c$};
  \draw[orchid, thick] (-1.6,1.6) -- (1.6,-1.6) node[below right] {$h=c^\perp$};
  \node at (0,-2.5) {$\ip{u}{u}=\ip{v}{v}$};
\end{scope}
\begin{scope}[xshift=6.3cm]
  \draw[black!40,->] (-2,0) -- (2.3,0); \draw[black!40,->] (0,-2.3) -- (0,2.3);
  \draw[plum, very thick, -{Stealth}] (0,0) -- (1,0) node[below right] {$u$};
  \draw[plum, very thick, -{Stealth}] (0,0) -- (0,2) node[above left] {$v$};
  \draw[rose, thick] (-0.9,-1.8) -- (1.05,2.1) node[above right] {$c$};
  \draw[orchid, thick] (1.05,-2.1) -- (-0.9,1.8) node[left] {$h$};
  \draw[plum, thick, dashed] (-1.9,0.95) -- (2.1,-1.05) node[right] {$c^\perp$};
  \node at (0,-2.7) {$\ip{u}{u}\ne\ip{v}{v}$};
\end{scope}
\end{tikzpicture}
\caption{The harmonic conjugate and the orthocomplement of $\Delta(u+v)$.}
\label{fig:bisect}
\end{figure}
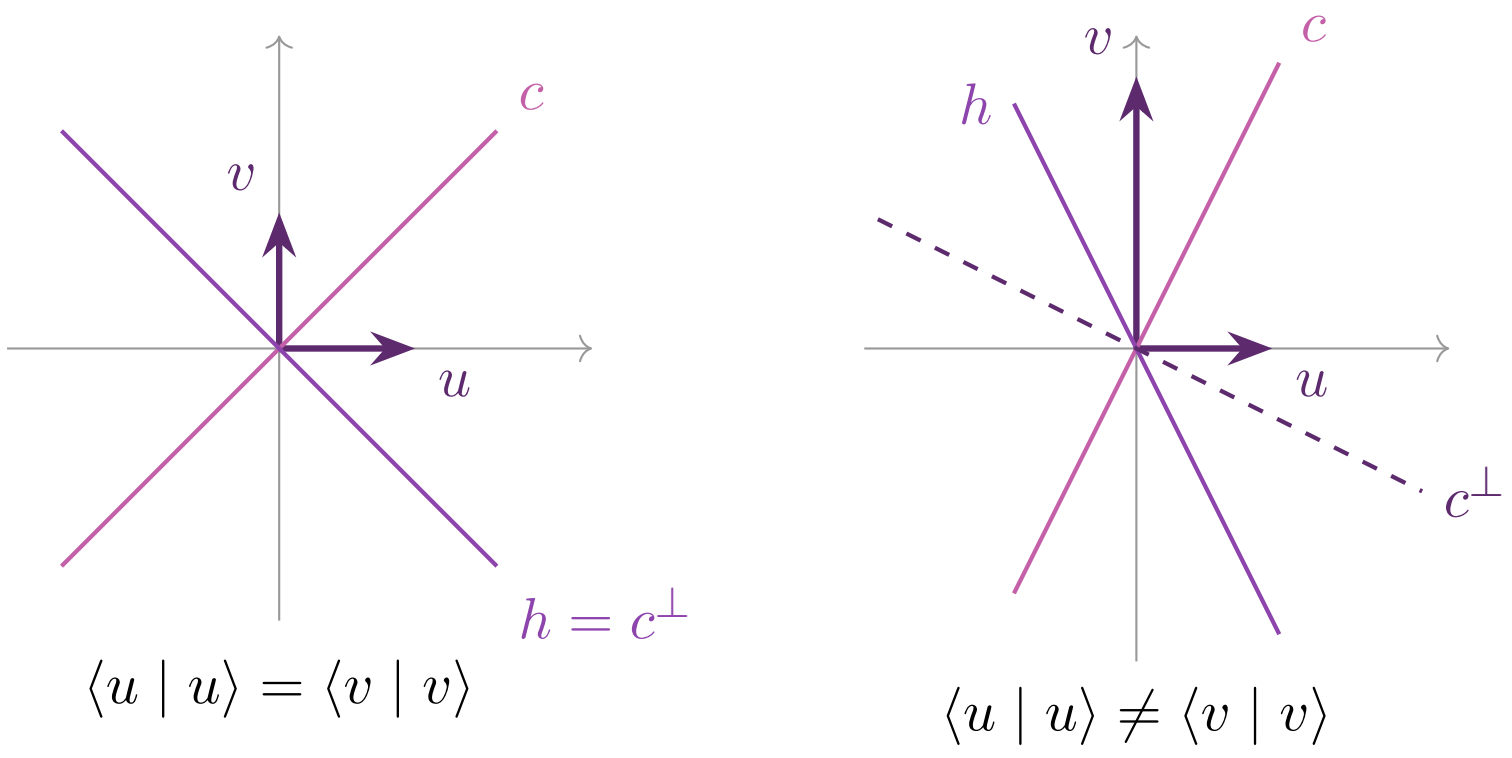

Suppose now that the vector space carries an anisotropic hermitian form, as the space $E$ of \cref{thm:coord} does, and that $p\perp q$; then the vectors $u$ and $v$ are orthogonal.
The orthocomplement of $c$ within the plane $p\vee q$, the ray $c^\perp\wedge(p\vee q)$, consists of the vectors $\alpha u+\beta v$ with $\ip{\alpha u+\beta v}{u+v}=\alpha\ip{u}{u}+\beta\ip{v}{v}=0$.
It coincides with the harmonic conjugate $\Delta(u-v)$ exactly when $\ip{u}{u}=\ip{v}{v}$, that is, when $u$ and $v$ have the same ``length'' (see \cref{fig:bisect}).

This formulation of Solèr's condition through harmonic conjugates is due to Holland~\cite[Theorems~3.7 and~4.1]{Holland:1995}.
Morash's angle-bisecting axiom~\cite{Morash:1973}, which asks for bisectors as in the left half of \cref{fig:bisect}, may be used instead (see Holland~\cite[\S4]{Holland:1995}).
The next lemma verifies the role of (S): agreement of the harmonic conjugate with the orthocomplement provides a common normalization of the infinite orthogonal sequence.

\begin{lemma}\label{lem:s-orthonormal}
Let $\Lambda$ be an ortholattice that satisfies (OM), (A), (CP), (I), (D) and (S), and let $E$ be a hermitian space over a division ring $\Delta$ and $p\mapsto \hat p$ a bijection from the atoms of $\Lambda$ onto the rays of $E$ that \cref{thm:coord} provides for $\Lambda$.
After a rescaling of the form, the space $E$ contains an infinite orthonormal sequence.
\end{lemma}

\begin{proof}
Let $p_n$ and $c_n$, for $n\in\NN$, be atoms as in (S).
Choose $e_0$ with $\hat p_0=\Delta e_0$, and inductively choose $e_{n+1}$ with $\hat p_{n+1}=\Delta e_{n+1}$ and $\hat c_n=\Delta(e_n+e_{n+1})$.
The choice of $e_{n+1}$ is possible because, by \cref{thm:coord}, the ray $\hat c_n$ lies in the plane $\hat p_n+\hat p_{n+1}$ and differs from $\hat p_n$ and $\hat p_{n+1}$.
By~\eqref{eq:harmonic} and the computation of the orthocomplement of $c$ in the plane $p\vee q$ above, condition (S) gives $\ip{e_n}{e_n}=\ip{e_{n+1}}{e_{n+1}}$ for every $n\in\NN$; call this common value $\zeta$, which satisfies $\zeta\inv=\zeta\ne0$.
The form $\ip{x}{y}'=\ip{x}{y}\zeta^{-1}$ is hermitian for the involution $\alpha\mapsto\zeta\alpha\inv\zeta^{-1}$, has the same orthogonality relation, and satisfies $\ip{e_n}{e_n}'=1$ for every $n$; as the atoms $p_n$ are mutually orthogonal, the sequence $(e_n)_{n\in\NN}$ is orthonormal for the new form.
\end{proof}

\begin{proof}[Sketch of the proof of \cref{thm:soler-lattices}]
Let $\Lambda$ be a Solèr lattice, with $\gamma$ and $E$ as in \cref{thm:coord}.
Given $M\in\cC(E)$, completeness provides the join $a$ of the atoms $p$ with $\hat p\subseteq M$, and this join satisfies $\gamma(a)^\perp=\gamma(a^\perp)=M^\perp$, whence $\gamma(a)=M$.
Thus the embedding $\gamma$ is an isomorphism onto $\cC(E)$; the lattice $\cC(E)$ is therefore orthomodular, and by \cref{prop:om-equiv} so is the space $E$.
By \cref{lem:s-orthonormal}, the space $E$ contains, after a rescaling of the form that changes neither $\cC(E)$ nor its orthocomplementation, an infinite orthonormal sequence.
Solèr's theorem~\cite{Soler:1995} (see also Prestel~\cite{Prestel:1995}) asserts that an orthomodular space containing an infinite orthonormal sequence is a Hilbert space over $\RR$, $\CC$ or $\HH$.
The space $E$, with its form so rescaled, is therefore such a Hilbert space.
Conversely, the lattice $L(H)$ is complete and satisfies the other conditions.
For (S), take an orthonormal sequence $(e_n)_{n\in\NN}$ in $H$, and let $p_n$ be the ray of $e_n$ and $c_n$ the ray of $e_n+e_{n+1}$.
\end{proof}

Condition (S) cannot be omitted from \cref{thm:soler-lattices}.
Keller~\cite{Keller:1980} constructed infinite-dimensional orthomodular hermitian spaces over division rings other than $\RR$, $\CC$ and $\HH$.
For such a space $E$ the lattice $\cC(E)$ satisfies (C), (OM), (A), (CP), (I) and (D) (see \cref{prop:om-equiv,lem:plus-one,lem:finite}).
Nevertheless, by the uniqueness clause of \cref{thm:coord}, the lattice $\cC(E)$ is isomorphic to no lattice $L(H)$ of a Hilbert space $H$ over $\RR$, $\CC$ or $\HH$.

\historypara{Historical note.}
Von Staudt's \emph{Geometrie der Lage} of 1847 was the first treatment of projective geometry free of any metric.
His later calculus of ``throws'' outlined the algebraization of an axiomatic geometry.
Hilbert showed in his \emph{Grundlagen der Geometrie} of 1899 that Desargues's theorem follows from the axioms of space but cannot be proved in the plane without the axioms of congruence.
On this theorem he founded an algebra of segments that does not use the axioms of congruence~\cite{Hilbert:1899}.

Veblen and Young gave axioms for projective geometry in 1908~\cite{Veblen-Young:1908}.
In their book of 1910 they introduced coordinates through an algebra of points that descends from von Staudt's algebra of throws~\cite{Veblen-Young:1910}.
Solèr, a student of Herbert Gross, had just begun work on Gross's orthomodular problem when Gross died suddenly in 1989.
She completed the work under the informal guidance of Alexander Prestel of Konstanz, a friend of Gross's, and her theorem appeared in 1995~\cite{Soler:1995}.
In a survey of the same year, Holland restated her condition through harmonic conjugates, constructed as in the book of Veblen and Young, in the form used in this subsection~\cite{Holland:1995}.

\subsection{What a pre-Solèr lattice is missing}\label{sec:missing}

The preceding subsection used completeness to pass from a pre-Solèr lattice to the full lattice of closed subspaces of a Hilbert space.
We now identify its MacNeille completion with the ortholattice of orthogonally closed subspaces of the hermitian space obtained by coordinatization, and show that this completion is a Solèr lattice exactly when it is orthomodular.
A condition on sets of atoms characterizes orthomodularity of the completion without referring to the representing space.
Neither this condition nor completeness is expressible among pre-Solèr lattices by a sentence of two-valued $L_{\omega_1\omega}$ about lattice elements (see \cref{prop:geometric-nonexpressibility}).

A pre-Solèr lattice need not be complete, as \cref{ex:e0} below shows; the following criterion and that example concern the MacNeille completion of such a lattice.

\begin{theorem}\label{thm:missing}
Let $\Lambda$ be a pre-Solèr lattice, and let $E$ be a hermitian space that \cref{thm:coord} provides for $\Lambda$.
Then the following are equivalent:
\begin{enumerate}[label=(\arabic*)]
\item the MacNeille completion of $\Lambda$ is orthomodular;
\item the space $E$ is orthomodular;
\item the lattice $\Lambda$ is isomorphic to a join-dense and meet-dense sub-ortholattice of $L(H)$ for some infinite-dimensional Hilbert space $H$ over $\RR$, $\CC$ or $\HH$.
\end{enumerate}
If these conditions hold, then the space $E$, with its form suitably rescaled, is a Hilbert space as in (3).
The lattice $\Lambda$ is isomorphic to $L(H)$ for some infinite-dimensional Hilbert space $H$ over $\RR$, $\CC$ or $\HH$ if and only if $\Lambda$ is complete.
\end{theorem}

\begin{example}\label{ex:e0}
Let $E_0$ be the space of finitely supported complex sequences $x=(x_n)_{n\in\NN}$, with the inner product $\ip{x}{y}=\sum_nx_n\overline{y_n}$ of the space $\ell^2$ of square-summable sequences, and let $\Lambda_0$ be the set of finite-dimensional subspaces of $E_0$ and their orthogonal complements in $E_0$.
Then $\Lambda_0$ is a pre-Solèr lattice whose MacNeille completion is isomorphic to $\cC(E_0)$ and is not orthomodular.
\end{example}

Orthomodularity of the completion of a pre-Solèr lattice can be stated through the atoms of the lattice alone, without reference to the completion or to a coordinatizing space.
For a set $\cA$ of atoms, let $\cA^{\perp}$ be the set of atoms orthogonal to every element of $\cA$, and call $\cA$ \emph{orthoclosed} when $\cA=\cA^{\perp\perp}$.
Consider the following condition, named (MOM) for MacNeille orthomodularity:
\begin{enumerate}[label=(MOM),leftmargin=*]
\item for every orthoclosed set $\cA$ of atoms and every atom $p\notin\cA$, some atom $q\in\cA^\perp$ lies in $(\cA\cup\{p\})^{\perp\perp}$.
\end{enumerate}

\begin{theorem}\label{thm:mom}
The MacNeille completion of a pre-Solèr lattice is orthomodular if and only if the lattice satisfies (MOM).
\end{theorem}

Throughout this subsection we fix a pre-Solèr lattice $\Lambda$, together with a hermitian space $E$ and an embedding $\gamma\colon\Lambda\to\cC(E)$ that \cref{thm:coord} provides for $\Lambda$.
We identify $\Lambda$ with its image $\gamma(\Lambda)\subseteq\cC(E)$.

\begin{theorem}\label{thm:dense}
Let $\Lambda$ be a pre-Solèr lattice, let $E$ be a hermitian space and $\gamma\colon\Lambda\to\cC(E)$ an embedding that \cref{thm:coord} provides for $\Lambda$, and identify $\Lambda$ with its image $\gamma(\Lambda)$.
Then $\Lambda$ is join-dense and meet-dense in $\cC(E)$, and $\cC(E)$, with $M\mapsto M^\perp$, is isomorphic as an ortholattice to the MacNeille completion of $\Lambda$.
\end{theorem}

\begin{proof}
Every $M\in\cC(E)$ is the join in $\cC(E)$ of the rays it contains, which belong to $\Lambda$ by \cref{thm:coord}.
If $M^\perp=\{0\}$, then the subspace $M$ is $E$, which belongs to $\Lambda$.
Otherwise $M$ equals $M^{\perp\perp}=\bigcap\{y^\perp:0\ne y\in M^\perp\}$, a meet of hyperplanes, each of which belongs to $\Lambda$.
Since $\gamma$ is an embedding of ortholattices, the second assertion now follows from \cref{prop:dm-char}.
\end{proof}

The completion supplies the orthogonally closed subspaces absent from $\Lambda$.
Any such subspace $M$ must have both $M$ and $M^\perp$ infinite-dimensional, since $\Lambda$ already contains every finite-dimensional subspace and its orthogonal complement (see \cref{thm:coord}).
\Cref{thm:missing} gives the criterion for the completion to be a Solèr lattice, and \cref{ex:e0} exhibits a pre-Solèr lattice for which the criterion fails.

\begin{proof}[Proof of \cref{thm:missing}]
Conditions (1) and (2) are equivalent by \cref{thm:dense} and \cref{prop:om-equiv}.
If (2) holds, then \cref{lem:s-orthonormal} and Solèr's theorem (see the proof of \cref{thm:soler-lattices}) show that the space $E$, with its form rescaled, is a Hilbert space.
Rescaling changes neither $\cC(E)$ nor its orthocomplementation, and condition (3) follows from \cref{thm:dense}.
If (3) holds, then by \cref{prop:dm-char} the orthomodular lattice $L(H)$ is isomorphic as an ortholattice to the MacNeille completion of $\Lambda$; this is (1).
For the last assertion, a complete pre-Solèr lattice is a Solèr lattice and by \cref{thm:soler-lattices} is therefore isomorphic to $L(H)$ for an infinite-dimensional Hilbert space $H$ over $\RR$, $\CC$ or $\HH$.
Conversely, every lattice $L(H)$ is complete.
\end{proof}

In \cref{ex:e0}, the finite-dimensional subspaces of an incomplete complex inner product space and their orthogonal complements form a pre-Solèr lattice whose MacNeille completion is not orthomodular.

\begin{proof}[Verification of \cref{ex:e0}]
By \cref{lem:finite}, the finite-dimensional elements of $\Lambda_0$ are orthogonally closed; the others are orthogonal complements, hence orthogonally closed as well, and $\Lambda_0$ is closed under ${}^\perp$.
Moreover, for finite-dimensional $F$ and $G$ the joins $F\vee G=F+G$, $F^\perp\vee G^\perp=(F\cap G)^\perp$ and $F\vee G^\perp=(F^\perp\cap G)^\perp$ belong to $\Lambda_0$, whence the De~Morgan law places the meets in $\Lambda_0$ as well.
Thus $\Lambda_0$ is a sub-ortholattice of $\cC(E_0)$.
The lattice $\Lambda_0$ is orthomodular, because every $M\in\Lambda_0$ satisfies $M+M^\perp=E_0$ by \cref{lem:finite}, and the argument of \cref{prop:om-equiv} applies.
As a sub-ortholattice of $\cC(E_0)$ that contains every ray of $E_0$, the lattice $\Lambda_0$ inherits (A), (CP) and (I) from $\cC(E_0)$, in which (CP) holds by \cref{lem:plus-one}.
Conditions (D) and (S) hold with $p_n=\CC e_n$ and $c_n=\CC(e_n+e_{n+1})$, where $(e_n)_{n\in\NN}$ is the sequence of standard unit vectors.
Consequently $\Lambda_0$ is a pre-Solèr lattice, and the space $E_0$, with $\hat p=p$ for every atom $p$ of $\Lambda_0$, is a hermitian space that \cref{thm:coord} provides for $\Lambda_0$.
The corresponding embedding $\gamma$ is then the inclusion of $\Lambda_0$ in $\cC(E_0)$.
The MacNeille completion of $\Lambda_0$, which is $\cC(E_0)$ by \cref{thm:dense}, is not orthomodular.
Indeed, by the Amemiya--Araki theorem~\cite{Amemiya-Araki:1966}, a complex inner product space whose lattice of orthogonally closed subspaces is orthomodular is complete, which $E_0$ is not.
In sum, the lattice $\Lambda_0$ is orthomodular while its MacNeille completion is not; Harding~\cite{Harding:1991} studies such lattices.
\end{proof}

\begin{proof}[Proof of \cref{thm:mom}]
Let $\hat p$ denote, for each atom $p$, the ray of $E$ that corresponds to $p$ under \cref{thm:coord}.
The map $\cA\mapsto\{0\}\cup\bigcup\{\hat p:p\in\cA\}$ is a bijection from the orthoclosed sets of atoms onto $\cC(E)$, carrying $\cA^\perp$ to the orthogonal complement of the image of $\cA$.
If $M$ is the image of $\cA$ and $\hat p=\Delta x$, where $\Delta$ is the division ring of scalars of $E$, then by \cref{lem:plus-one} the set $(\cA\cup\{p\})^{\perp\perp}$ corresponds to $M+\Delta x$.
Condition (MOM) therefore says that $(M+\Delta x)\cap M^\perp\ne\{0\}$ whenever $M\in\cC(E)$ and $x\notin M$.
If $y=v+\alpha x$ is nonzero, with $v\in M$ and $y\in M^\perp$, then the scalar $\alpha$ is nonzero because $y\notin M$, whence the vector $x=\alpha^{-1}(y-v)$ lies in $M+M^\perp$.
Conversely, if $x=v+y$ with $v\in M$ and $y\in M^\perp$, then the vector $y=x-v$ is a nonzero element of $(M+\Delta x)\cap M^\perp$.
Thus (MOM) is equivalent to the orthomodularity of $E$, and the theorem follows from \cref{thm:missing}.
\end{proof}

By \cref{lem:finite}, the requirement in (MOM) holds automatically unless the subspaces of $E$ that correspond to $\cA$ and to $\cA^\perp$ in the proof of \cref{thm:mom} are both infinite-dimensional.

\historypara{Historical note.}
MacNeille introduced the completion by cuts in 1937~\cite{MacNeille:1937}, generalizing the cuts by which Dedekind had constructed the real numbers from the rationals in 1872~\cite{Dedekind:1872}.
Banaschewski and Schmidt characterized it in 1956 by join-density and meet-density~\cite{Banaschewski:1956,Schmidt:1956}.
In 1969 Adams solved negatively Problem~36 of Birkhoff's \emph{Lattice Theory}, in its third edition of 1967, by giving an orthocomplemented modular lattice whose completion by cuts is not orthomodular~\cite{Adams:1969}.
His lattice consists of the finite-dimensional subspaces of an incomplete complex inner product space and their orthogonal complements, as in \cref{ex:e0}, and his argument also rests on the Amemiya--Araki theorem.
Harding gave in 1991, without recourse to the Amemiya--Araki theorem, a construction of orthomodular lattices whose completions by cuts are not orthomodular.
Following a method that Kalmbach introduced in 1977, he glued together the Boolean algebras generated by the chains of a lattice~\cite{Harding:1991}.
Harding noted that whether every orthomodular lattice can be embedded in a complete one was still open.

\subsection{Nonexpressibility over pre-Solèr lattices}

Completeness is the one condition of a Solèr lattice that a pre-Solèr lattice may fail, and among pre-Solèr lattices no sentence of $L_{\omega_1\omega}$ about the lattice alone expresses it.
In classical two-valued $L_{\omega_1\omega}$, the same holds for condition (MOM), which by \cref{thm:mom} characterizes the pre-Solèr lattices whose MacNeille completion is orthomodular; (MOM) quantifies over sets of atoms, whereas the logic quantifies only over individual elements of the lattice.
For both conditions the obstruction is countability, as the proof of the following proposition shows.

\begin{proposition}\label{prop:geometric-nonexpressibility}
Neither (C) nor (MOM) is equivalent, among pre-Solèr lattices, to a sentence of two-valued $L_{\omega_1\omega}$ in the vocabulary of ortholattices expanded by constants for atoms $p_n$ and $c_n$, for $n\in\NN$, that witness (S).
\end{proposition}

\begin{proof}
Suppose that a sentence $\varphi$ of $L_{\omega_1\omega}$ were equivalent to (MOM) among pre-Solèr lattices.
By \cref{thm:mom}, the same sentence would then be equivalent to orthomodularity of their MacNeille completions.
The axioms of ortholattices and the conditions (OM), (A), (CP), (I), (D) and (S) are expressed by countably many sentences of $L_{\omega_1\omega}$ in this vocabulary.
These sentences hold in $L(\ell^2)$, with $p_n=\CC e_n$ and $c_n=\CC(e_n+e_{n+1})$ for the standard unit vectors $e_n$, $n\in\NN$.
Indeed, by~\eqref{eq:harmonic}, the harmonic conjugate of $c_n$ with respect to $p_n$ and $p_{n+1}$ is $\CC(e_n-e_{n+1})$, which is $c_n^\perp\wedge(p_n\vee p_{n+1})$.
The sentence $\varphi$ also holds in this lattice, since the lattice is complete and orthomodular and is therefore its own MacNeille completion.
By the downward Löwenheim--Skolem theorem for a countable fragment, these sentences hold in a countable structure $\Lambda'$.
Thus $\Lambda'$ is a pre-Solèr lattice whose MacNeille completion is orthomodular.
By \cref{thm:missing}, a hermitian space $E$ provided for $\Lambda'$ by \cref{thm:coord} is, after rescaling its form, an infinite-dimensional Hilbert space.
By \cref{thm:coord}, the atoms of $\Lambda'$ correspond to all the rays of $E$.
There are uncountably many such rays, contradicting the countability of $\Lambda'$.

If instead $\varphi$ expressed completeness among pre-Solèr lattices, the same argument would give a countable complete pre-Solèr lattice $\Lambda'$.
By \cref{thm:soler-lattices}, it would be isomorphic to $L(H)$ for an infinite-dimensional real, complex or quaternionic Hilbert space $H$, again contradicting countability.
\end{proof}

The argument for completeness needs only the downward Löwenheim--Skolem theorem and \cref{thm:soler-lattices}, and it applies equally in real-valued $L_{\omega_1\omega}$, where equality of lattice elements is represented by a two-valued predicate.
Take a countable fragment containing sentences that express the ortholattice axioms and the conditions defining pre-Solèr lattices, together with a sentence proposed to express completeness.
The expansion of $L(\ell^2)$ used in the proof has a countable substructure in which every sentence of this fragment keeps its value.
The quotient of this substructure by the equality predicate is a countable pre-Solèr lattice, so by \cref{thm:soler-lattices} it cannot be complete.
This sentence takes the same value in this quotient as in the complete lattice $L(\ell^2)$, and so cannot express completeness.
Every countable set of sentences about the lattice alone thus keeps its values in some countable structure, whereas the representation of a Solèr lattice requires all the rays of an infinite-dimensional Hilbert space.

\section{The space of Henkin types and the model existence game}\label{app:henkin}

In a position of the game of \cref{sec:game}, finitely many commitments and records describe a structure that is still being built.
In this appendix we attach to each position the set of Henkin types compatible with it, in the space of Henkin types of \cref{sec:henkin-types}.
Josh's truthful answers and Olivia's terminal obstructions then acquire a topological meaning, and the positions and rules of the game can be compared with Keisler's forcing conditions and $\forall\exists$ formulas.
We close by setting the Baire category argument behind Olivia's bookkeeping, which takes place in the space of plays, beside the argument by which Aubone, Dueñez and Iovino omit types in the space of Henkin types~\cite{Aubone-Duenez-Iovino:2026a}; the two arguments take place in different spaces, and their relation is an analogy, not an equivalence.

\subsection*{The space of Henkin types}
We retain the setting of \cref{sec:game}: the vocabulary $\cL$ is countable, and $T$ is a nonempty countable theory of closed conditions $\psi\le r$ or $\psi\ge r$ for a sentence $\psi$ and a rational number $r$.
The fragment $\cF$ generated by the sentences $\psi$ of the axioms of $T$ is therefore countable (see~\cite[Section~6]{Aubone-Duenez-Iovino:2026a}).
We take as the new constants $C$ of \cref{sec:henkin-types} the witnesses $\cW$ of the game, countably many of each sort.
The fragment $\cF_C$ that $\cF$ generates in $\cL_{\cW}$ then has a countable set $S^{\cF}_C$ of sentences, and this set contains every sentence that occurs in a play of the game $\mathrm G(T)$ of \cref{def:game}, since such a sentence is an atomic sentence of $\cL_{\cW}$ or arises from a subformula of an axiom by the substitution of closed terms.
We write $\prod S^{\cF}_C$ for the product of the closed intervals $[-r_\varphi,r_\varphi]$, for $\varphi\in S^{\cF}_C$, and $\varphi^{\mathfrak m}$ for the coordinate at $\varphi$ of a point $\mathfrak m$ of the product.
This product, the space of all candidate assignments of values to sentences, is compact Hausdorff by Tychonov's theorem (see Engelking~\cite[Theorem~3.2.4]{Engelking:1989a}), and the Henkin types form its subspace $\cH^{\cF}_C$, the set of assignments that satisfy the Henkin conditions.
For a set $T'$ of closed conditions on sentences of $S^{\cF}_C$, we write $\langle\!\langle T'\rangle\!\rangle$ for the set of points that satisfy every condition of $T'$; in particular, $\langle\!\langle\emptyset\rangle\!\rangle$ is the whole product.

We use three results of Aubone, Dueñez and Iovino~\cite{Aubone-Duenez-Iovino:2026a} for the logic of this paper:
\begin{enumerate}[label=(\roman*)]
\item A structure for $\cL_{\cW}$ is \emph{canonical} when every element is the value of a closed term.
The family of the values of the sentences of $S^{\cF}_C$ in a canonical structure, its \emph{$S^{\cF}_C$-type}, belongs to $\cH^{\cF}_C$.
Conversely, every member $\mathfrak m$ of $\cH^{\cF}_C$ is the $S^{\cF}_C$-type of a canonical structure, the \emph{term structure} of $\mathfrak m$~\cite[Lemma~3]{Aubone-Duenez-Iovino:2026a}.
The elements of the term structure are the closed terms of $\cL_{\cW}$, on which the function symbols act by forming terms, and the value of each atomic sentence in it is the value that $\mathfrak m$ assigns to the sentence.
\item For a countable vocabulary and a countable fragment, every structure has a substructure whose sorts are countable and in which every sentence of the fragment keeps its value~\cite[Lemma~4 and its proof]{Aubone-Duenez-Iovino:2026a}.
A structure for $\cL$ whose sorts are countable has an expansion to $\cL_{\cW}$ that names every element by a witness, and this expansion is canonical~\cite[proof of Proposition~5]{Aubone-Duenez-Iovino:2026a}.
\item The space $\cH^{\cF}_C$ is a $G_\delta$ subset of $\prod S^{\cF}_C$~\cite[Proposition~6]{Aubone-Duenez-Iovino:2026a}.
Every Henkin condition of \cref{sec:henkin-types} defines a closed set, except the upper bounds on suprema that (H5$'$) and (H6$'$) impose; each of these upper bounds defines a countable intersection of open sets.
For a supremum $\sup_x\varphi(x)$, the upper bound says that its value does not exceed the supremum of the values of the sentences $\varphi(t)$, for the closed terms $t$ of $\cL_{\cW}$ of the sort of $x$.
For a countable join, the upper bound says that the value of the join does not exceed the supremum of the values of its members.
\end{enumerate}

From~(iii), Aubone, Dueñez and Iovino deduce that the space of Henkin types is \v{C}ech-complete~\cite[Corollary~7]{Aubone-Duenez-Iovino:2026a}.
The category arguments below require the same of the Henkin types that satisfy the theory.
Each axiom $\psi\le r$ or $\psi\ge r$ of $T$ defines the intersection of the open sets $\{\mathfrak m:\psi^{\mathfrak m}<r+\tfrac1n\}$ or $\{\mathfrak m:\psi^{\mathfrak m}>r-\tfrac1n\}$, for $n\ge1$, and since $T$ is countable, the set $\langle\!\langle T\rangle\!\rangle\cap\cH^{\cF}_C$, like $\cH^{\cF}_C$ itself, is a $G_\delta$ subset of the compact product, and hence \v{C}ech-complete~\cite[Theorems~3.1.2 and~3.9.1]{Engelking:1989a}.
In a \v{C}ech-complete space, by the Baire category theorem, the union of countably many nowhere dense subsets has dense complement~\cite[Theorem~3.9.3]{Engelking:1989a}.

\subsection*{Positions and plays}
Every position reached in a play has its sentences in $S^{\cF}_C$.
With such a position $\mathfrak p$ we associate the set $\langle\!\langle\mathfrak p\rangle\!\rangle$ of the points $\mathfrak m$ of the product that obey its commitments and its records: $\varphi^{\mathfrak m}\in I$ for every commitment $\varphi\in I$ of $\mathfrak p$, and $\theta(F(\bar w))^{\mathfrak m}=\theta(w')^{\mathfrak m}$ for every record $F(\bar w)\mapsto w'$ of $\mathfrak p$ and every atomic formula $\theta(x)$ of $\cL_{\cW}$ whose only free variable is a variable $x$ of the sort of $w'$.
These equations make $F(\bar w)$ and $w'$ agree in every atomic formula with parameters.
Since the coordinate projections are continuous, the sets $\langle\!\langle\mathfrak p\rangle\!\rangle$ and $\langle\!\langle T'\rangle\!\rangle$ are closed: the closed intervals of a position give a closed set, as finitely many strict conditions would give a basic open set.
Satisfaction of the position in a structure requires more than agreement in atomic formulas: under an interpretation of the witnesses, its commitments must hold and every record $F(\bar w)\mapsto w'$ must make the term $F(\bar w)$ and the witness $w'$ denote the same element (see \cref{def:position}).
Over the Henkin types, as the following proposition shows, the difference disappears.

\begin{proposition}\label{prop:app-meets}
Let $\mathfrak p$ be a position whose sentences belong to $S^{\cF}_C$, and let $T'$ be a subset of $T$.
The set $\langle\!\langle\mathfrak p\rangle\!\rangle\cap\langle\!\langle T'\rangle\!\rangle\cap\cH^{\cF}_C$ is nonempty if and only if $\mathfrak p$ is satisfied in some model of $T'$ under some interpretation of the witnesses.
\end{proposition}

\begin{proof}
Let $\mathfrak m$ be a point of the set, and let $\fN$ be, by~(i), a canonical structure whose $S^{\cF}_C$-type is $\mathfrak m$.
The commitments of $\mathfrak p$ hold in $\fN$, and for each record $F(\bar w)\mapsto w'$ of $\mathfrak p$, the equations defining $\langle\!\langle\mathfrak p\rangle\!\rangle$ make $F(\bar w)$ and $w'$ agree in every atomic formula with parameters.
A quotient turns this agreement into the literal equality that the record requires.
Call two elements of the same sort of $\fN$ equivalent when they agree in every atomic formula with parameters.
This relation is a congruence for every function symbol, since substituting a term for a variable of an atomic formula yields an atomic formula; the quotient map therefore preserves the value of every atomic formula and, by induction on formulas, of every formula.
Since $\mathfrak m$ lies in $\langle\!\langle T'\rangle\!\rangle$, the reduct of the quotient to $\cL$ is a model of $T'$.
When each witness is interpreted as the class of its value, this model satisfies $\mathfrak p$: its commitments hold because the quotient map preserves values, and the equations defining $\langle\!\langle\mathfrak p\rangle\!\rangle$ place the values of $F(\bar w)$ and $w'$ in one class.

Conversely, let a model $\fM$ of $T'$ satisfy $\mathfrak p$ under an interpretation $h$.
By~(ii), applied to the expansion of $\fM$ by constants for the elements $h(w)$ and to the fragment generated by $\cF$ and the sentences of $\mathfrak p$, the structure $\fM$ has a substructure $\fM_0$ whose sorts are countable, which contains these elements and satisfies $T'$ and $\mathfrak p$ under $h$, the records included.
The unused witnesses of each sort, infinitely many since $\mathfrak p$ is finite, can name every element of $\fM_0$.
The resulting expansion is canonical, and by~(i) its $S^{\cF}_C$-type belongs to $\langle\!\langle\mathfrak p\rangle\!\rangle\cap\langle\!\langle T'\rangle\!\rangle\cap\cH^{\cF}_C$.
\end{proof}

The proposition is the counterpart, for positions, of the theorem that Aubone, Dueñez and Iovino call Henkin Realizability~\cite[Proposition~5]{Aubone-Duenez-Iovino:2026a}.
Taking for $\mathfrak p$ the empty position and $T'=T$, we see that $T$ has a model exactly when $\langle\!\langle T\rangle\!\rangle\cap\cH^{\cF}_C$ is nonempty.

The proposition gives a topological meaning to both outcomes of the game.
In the proof that (1) implies (2) in \cref{thm:game-existence}, Josh answers truthfully in a model of $T$, and every position that he reaches is satisfied in this model; by \cref{prop:app-meets}, his answers keep the set $\langle\!\langle\mathfrak p\rangle\!\rangle\cap\langle\!\langle T\rangle\!\rangle\cap\cH^{\cF}_C$ nonempty at every position $\mathfrak p$ of the play.
At a leaf of the tree of a winning strategy of Olivia, the last position $\mathfrak p$ is satisfied in no structure (see \cref{cor:finite-obstruction}(2)), and \cref{prop:app-meets}, with $T'$ empty, shows that the closed set $\langle\!\langle\mathfrak p\rangle\!\rangle$ does not meet $\cH^{\cF}_C$.

Survival has the complementary meaning: a fair play that Josh survives determines a nonempty subset of $\langle\!\langle T\rangle\!\rangle\cap\cH^{\cF}_C$ that is closed in $\cH^{\cF}_C$.
Let $\mathfrak p_0\subseteq\mathfrak p_1\subseteq\cdots$ be its positions, indexed by $i\in\NN$, and let $\mathfrak p_\infty$ be their union.
The subset in question is $\bigcap_i\langle\!\langle\mathfrak p_i\rangle\!\rangle\cap\cH^{\cF}_C$, and it lies in $\langle\!\langle T\rangle\!\rangle$: by fairness, every axiom of $T$ is named in a challenge (Ax), and Josh's answer commits the sentence of the axiom to an interval on which the axiom holds.
The structure $\fM$ built from this play as in the proof that (2) implies (1) in \cref{thm:game-existence} is canonical, since its elements are the witnesses, each named by itself.
It satisfies every commitment and every record of $\mathfrak p_\infty$, and by~(i) its $S^{\cF}_C$-type lies in the set, which is therefore nonempty.
The play need not determine a single Henkin type.
Suppose, for instance, that a predicate symbol $R$ other than $\mathbf 1$, with $r_R>0$, occurs in no axiom of $T$; then no sentence that occurs in the play contains $R$.
Each value of $R^{\fM}$ remains free in $[-r_R,r_R]$, and different choices give $S^{\cF}_C$-types in the set that differ at $R(\bar w)$.

Positions play the part of the conditions in Keisler's model-theoretic forcing.
Keisler adjoins countably many new constants and takes as forcing conditions the finite sets of sentences of the enlarged vocabulary that are satisfiable in a given class of models, ordered by inclusion~\cite[pp.~98 and~100, Example~1.2]{Keisler:1973}.
The witnesses serve as the new constants, and a position grows by commitments and records as a condition grows by sentences.
By \cref{prop:app-meets}, a position $\mathfrak p$ is satisfied in some model of $T$ under some interpretation of the witnesses exactly when $\langle\!\langle\mathfrak p\rangle\!\rangle\cap\langle\!\langle T\rangle\!\rangle\cap\cH^{\cF}_C$ is nonempty, and these satisfiable positions correspond to the conditions for the class of models of $T$.
A record $F(\bar w)\mapsto w'$ is the counterpart, in a logic without equality, of an atomic sentence $c=\tau$, where $c$ is a new constant and $\tau$ is a closed term.
Keisler requires that, for every closed term $\tau$, each condition extend to one that contains such a sentence for some added constant $c$~\cite[p.~99, clause~(iv)]{Keisler:1973}, and the rule (Tot) meets the corresponding requirement in the game.
A fair play that Josh survives yields an increasing sequence of satisfiable positions and a canonical model, as Keisler's construction yields a generic model~\cite[p.~101, Lemma~1.4]{Keisler:1973}; fairness alone, however, does not impose his requirement that every sentence of the fragment be decided.

The rules also separate universal from existential choices, as Keisler's $\forall\exists$ formulas do~\cite[p.~107]{Keisler:1973}.
These formulas have the form $\forall\bar x\bigvee_n\exists\bar y(\varphi_{n1}\wedge\dots\wedge\varphi_{nj_n})$, where $n$ ranges over $\NN$, each $j_n$ is a positive integer, and the formulas $\varphi_{ni}$, for $1\le i\le j_n$, belong to a fixed set of formulas of a countable fragment.
Consider an axiom $\inf_{\bar x}\bigvee_n\sup_{\bar y}\theta_n(\bar x,\bar y)\ge r$ of the same shape, in which each formula $\theta_n$ is a minimum of finitely many formulas and the a priori bounds of the formulas $\theta_n$ have a finite supremum, as the formation of the countable join requires.
After Josh commits the value of the sentence of this axiom under (Ax), the rule ($\inf{\ge}$) lets Olivia choose the witnesses for $\bar x$, whereas ($\bigvee{\ge}$) and ($\sup{\ge}$) require Josh to choose the index $n$ and the witnesses for $\bar y$, up to the tolerances that Olivia names: the universal choices are Olivia's, and the existential choices are Josh's.
Keisler's general omitting types theorem asks that every condition remain satisfiable together with every instance of the given $\forall\exists$ formulas at tuples of the added constants~\cite[p.~108]{Keisler:1973}.
Its counterpart in the game is the requirement that Josh answer each challenge at a position satisfied in a model of $T$ by a position with the same property, and the proof that (1) implies (2) in \cref{thm:game-existence} shows that he can do so whenever $T$ has a model.
When $T$ has no model, a winning strategy of Olivia ends each play at a finite obstruction instead.
The last position at each leaf holds in no structure under any interpretation of its witnesses, so that no extension of it is satisfied in a model of $T$ or corresponds to a forcing condition for the class of models of $T$.
The correspondence is not literal in three respects: commitments confine the values of sentences to closed intervals instead of asserting sentences, the axioms of $T$ need not have the $\forall\exists$ shape, and Keisler's logic has equality~\cite[p.~97]{Keisler:1973}, whose role in the game is taken by the records.
For the strategy of \cref{sec:bks-game}, the contextual presentations of \cref{sec:obstructions} picture the terminal commitments and records of the plays that reach the parity argument as vertices, identifications and intervals (see \cref{tab:ks-strategy-contrast}).

\subsection*{Closed conditions, strict conditions and compactness}
The contrast between closed and strict conditions drawn in \cref{rem:strict} has a topological counterpart.
The Henkin types that satisfy a closed condition form a closed subset of $\cH^{\cF}_C$, and those that satisfy a strict condition an open subset.
In the compact space $\prod S^{\cF}_C$, a family of closed sets with empty intersection has a finite subfamily with empty intersection; a family of open sets need not.
For the theory $\{R>0\}\cup\{R<\tfrac1n:n\ge1\}$ of \cref{rem:strict}, each finite part has a model, so that the open set that it defines meets the space of Henkin types, by~(i) and~(ii); the whole theory defines the empty set, and the tree of the winning strategy of Olivia described in \cref{rem:strict} for this theory branches infinitely at its root.

Closed conditions alone do not restore compactness, because the space $\cH^{\cF}_C$ need not be closed in the product.
Were it closed and had $T$ no model, compactness would give a finite subset $T'$ of $T$ with $\langle\!\langle T'\rangle\!\rangle\cap\cH^{\cF}_C=\emptyset$, since the axioms of $T$ define closed sets, and by \cref{prop:app-meets} the set $T'$ would be a finite set of axioms without models.
The theory $T^{01}$ has no model, although every finite set of its axioms has one (see \cref{prop:bks-game}); for its fragment, therefore, the space $\cH^{\cF}_C$ is not closed.
What the category arguments use is not closedness but \v{C}ech-completeness, which $\cH^{\cF}_C$ owes to being a $G_\delta$ subset of the compact product.
By~(iii), the parts of the Henkin conditions that need not define closed sets are the upper bounds on suprema, and in the game Josh certifies such a bound by naming an instance.
Under ($\sup{\ge}$) and ($\inf{\le}$) the instance is a witness, chosen among the finitely many that the rules allow him: those that occur in the position and the first fresh witness.
Under ($\bigvee{\ge}$) and ($\bigwedge{\le}$) it is an index, chosen among infinitely many.
Parts (3) and (4) of \cref{cor:finite-obstruction} record these two cases; they follow from the rules of the game alone, without the topology of $\cH^{\cF}_C$.

\subsection*{Baire category}
The category argument behind Olivia's bookkeeping takes place in another space, whose points are plays.
Suppose that $T$ has a model, so that, by \cref{thm:game-existence}, Josh has a winning strategy $\sigma'$.
Let $X_{\sigma'}$ be the set of infinite plays in which Josh follows $\sigma'$, with the topology whose basic open sets consist of the plays that extend a given finite play.
Since the moves range over countable sets, the space $X_{\sigma'}$ is homeomorphic to a closed subspace of the Baire space $\NN^\NN$; it is therefore Polish and satisfies the Baire category theorem.
For each challenge and each round, the plays in which the challenge, if it first becomes available at that round, is made then or later form a dense open subset of $X_{\sigma'}$.
Since only countably many challenges occur in plays (see \cref{sec:game}), the fair plays form a countable intersection of dense open subsets of $X_{\sigma'}$, which is dense by the Baire category theorem and in particular nonempty.

The bookkeeping can also be read as a strategy in a \emph{Banach--Mazur game}.
In such a game, two players alternately choose nonempty open subsets of a space, each contained in the one before, and the second player wins if the intersection of the chosen sets is contained in a given set.
Consider the Banach--Mazur game on $X_{\sigma'}$ in which the players choose basic open sets, that is, finite plays in which Josh follows $\sigma'$, each extending the one before, and the second player wins when the infinite play that these finite plays determine is fair.
The two players of this game are not Olivia and Josh: a move of either player may extend the play by many rounds, each consisting of a challenge of Olivia and an answer of Josh.
Olivia's bookkeeping of \cref{sec:game}, indexed by the moves of the second player rather than by the rounds of the play, supplies a winning strategy for the second player (compare Oxtoby~\cite{Oxtoby:1957}).
Recall that Olivia fixes a bijection $i\mapsto(j(i),l(i))$ from $\NN$ onto $\NN\times\NN$ with $j(i)\le i$ for every $i\in\NN$, and that at each round she enumerates the challenges then available.
At the move of index $k$ of the second player, for $k\in\NN$, the strategy extends the current finite play by one round, in which Olivia makes the $l(k)$-th challenge of her enumeration for round $j(k)$ and Josh answers by $\sigma'$.
Since each move of the second player adds a round, at least $k$ rounds have been played, so round $j(k)$ has been reached and the challenge is still available.
The finite plays therefore grow without bound and determine a unique infinite play, and the play is fair, since every pair $(j,l)$ is of the form $(j(k),l(k))$.
In the proof of \cref{thm:game-existence}, the category argument reduces to the bookkeeping, because only countably many challenges are available at each round.

Aubone, Dueñez and Iovino prove the omitting types theorem for their logic by Morley's Baire category method, applied in the \v{C}ech-complete space of Henkin types, whose points are the types of structures built from closed terms~\cite[Theorem~8]{Aubone-Duenez-Iovino:2026a}.
To state it in the terms of this appendix, we call a \emph{type} in a finite tuple $\bar x$ of variables a set $\Xi(\bar x)$ of expressions $\varphi(\bar x)\ge r$, with $r$ rational and with $\varphi$ a formula of $\cF$ whose free variables are among $\bar x$~\cite[Section~11]{Aubone-Duenez-Iovino:2026a}.
Substituting for $\bar x$ a tuple $\bar t$ of closed terms of $\cL_{\cW}$ of the appropriate sorts turns these expressions into closed conditions on sentences of $S^{\cF}_C$, which form the \emph{instance} $\Xi(\bar t)$ of the type.
For a sequence $(\Xi_n(\bar x_n))_{n\in\NN}$ of types, the theorem reads as follows.
Suppose that $T$ has a model and that, for every $n\in\NN$ and every tuple $\bar t$ of closed terms of $\cL_{\cW}$ of the appropriate sorts, the set of points of $\langle\!\langle T\rangle\!\rangle\cap\cH^{\cF}_C$ that satisfy the instance $\Xi_n(\bar t)$ has empty interior in this space.
Then $T$ has a model whose sorts are countable and in which no tuple satisfies any $\Xi_n(\bar x_n)$.
Indeed, the complements of these sets are countably many dense open subsets of the nonempty space $\langle\!\langle T\rangle\!\rangle\cap\cH^{\cF}_C$, and by~(i) a point of their intersection, which the Baire category theorem provides, is the $S^{\cF}_C$-type of a canonical structure whose reduct to $\cL$ is such a model.

The hypothesis of empty interior explicates, as a statement about interiors in the \v{C}ech-complete space of Henkin types that satisfy $T$, the hypothesis of Keisler's general omitting types theorem recalled above, in the case in which his $\forall\exists$ formulas are taken over the whole fragment (see \cite[Introduction]{Aubone-Duenez-Iovino:2026a}).
The interior is tested on basic open sets given by finitely many strict conditions, not on the closed sets $\langle\!\langle\mathfrak p\rangle\!\rangle$ attached to positions.

The two category arguments thus select different objects in different spaces: the first a fair play, the second a Henkin type that omits the prescribed types.
A fair play that Josh survives determines a closed set of Henkin types, which need not reduce to a single point, and we know of no continuous map from the space of plays to the space of Henkin types that carries the dense open requirements of the one argument to those of the other.
Their relation remains an analogy.

\historypara{Historical note.}
The omitting types theorem goes back to Henkin, in 1954~\cite{Henkin:1954}, and to Orey, in 1956~\cite{Orey:1956}.
Robinson introduced forcing into model theory in an article published in 1971~\cite{Robinson:1971a}, adapting the method that Cohen had created for set theory in 1963~\cite{Cohen:1963}.
Barwise extended Robinson's forcing to infinitary logic and used it to give a new proof of the omitting types theorem (see Keisler~\cite[p.~96]{Keisler:1973}).
Keisler described his forcing of 1973, recalled in the historical note of \cref{sec:henkin-types}, as a reformulation, in terms of ordinary models, of the Boolean-valued models of Rasiowa and Sikorski~\cite{Rasiowa-Sikorski:1963}, and he regarded its basic result, the Generic Model Theorem, as essentially due to them~\cite[p.~97]{Keisler:1973}.

\bibliographystyle{amsalpha}
\bibliography{GKS}

\end{document}